\documentclass[10pt,a4paper,final]{article}
\usepackage[left=2.250cm, right=2.250cm, top=3.750cm, bottom=3.75cm]{geometry}

\usepackage{titlesec}

\titleformat*{\section}{\normalsize\bfseries}
\titleformat*{\subsection}{\normalsize\itshape} 
\titleformat*{\subsubsection}{\normalsize\itshape}  
\titleformat*{\paragraph}{\normalsize\itshape}
\titleformat*{\subparagraph}{\normalsize\itshape}
\usepackage[font={footnotesize},labelfont=bf]{caption}

\usepackage{graphics}
\usepackage{graphicx}
 \usepackage{epsfig}

\usepackage{amssymb}
 \usepackage{amsthm}
  \usepackage{amsmath}
\usepackage{tikz-cd}

 \usepackage{float}
 \usepackage{newfloat}
\usepackage{lineno}

\usepackage{ccicons}
\usepackage{epstopdf}

\usepackage{multirow}
\usepackage{color}
\usepackage{hhline}
\usepackage{booktabs}

\graphicspath{{./}}

\definecolor{nverde}{RGB}{0,61,0} 
\definecolor{cr1}{RGB}{200,0,0}
\definecolor{cr2}{RGB}{0,0,200}
\definecolor{cr12}{RGB}{100,0,100}

\newcommand{\xx}{\mathbf{x}}
\newcommand{\nn}{\mathbf{\tilde{n}}}				
\newcommand{\tnn}{\mathbf{n}}  		

\newcommand{\tang}{\mathbf{t}}  				

\newcommand{\element}{T}					
\newcommand{\face}{\Gamma}					
\newcommand{\edge}{\Lambda}					
\newcommand{\vertex}{\mathtt{V}}			
\newcommand{\bary}{\mathtt{x}}		    	
\newcommand{\VT}[1]{\left|T_{#1}\right|}	
\newcommand{\pbf}[1]{\partial T_{#1}}  
\newcommand{\AF}[1]{\left|\pbf{#1}\right|}	
\newcommand{\edgeD}{{\edge}^*}
\newcommand{\faceD}{\face^*}
\newcommand{\length}[1]{\lvert #1 \rvert}

\newcommand{\mesh}{\mathcal{T}}

\newcommand{\myfunc}{\mathcal{U}}
\newcommand{\Lagrange}{\mathcal{L}}

\DeclareMathOperator{\dive}{\nabla \cdot}
\DeclareMathOperator{\gra}{\nabla}
\DeclareMathOperator{\grae}{\nabla}
\DeclareMathOperator{\curl}{\nabla \times}

\DeclareMathOperator{\NF}{\mathrm{NF}}
\DeclareMathOperator{\RS}{\mathrm{Rus}}
\DeclareMathOperator{\Duc}{\mathrm{Duc}}

\DeclareMathOperator{\dS}{\mathrm{dS}}
\DeclareMathOperator{\dV}{\mathrm{dV}}
\DeclareMathOperator{\Dt}{\Delta t\,}
\DeclareMathOperator{\dx}{\mathrm{d\xx}}
\newcommand{\halb}{\frac{1}{2}}

\newcommand{\vel}{u}						
\newcommand{\Vel}{{u}}
\newcommand{\bvel}{\mathbf{\vel}}			
\newcommand{\bvele}{\mathbf{\vel}^{e}}		

\newcommand{\bVel}{\mathbf{\Vel}}
\newcommand{\bbVel}{{\bVel}}		
\newcommand{\bVele}{\bVel^{e}}
\newcommand{\WW}{\rho\bVel}		    
\newcommand{\WWe}{\WW^{e}}				
\newcommand{\bWW}{{\WW}}		
\newcommand{\tWW}[1]{\WW^{\ast}_{#1}}
\newcommand{\ww}{\rho\bvel}			
\newcommand{\wwe}{\ww^{e}}				
\newcommand{\tww}{{\ww}^{\ast}}			

\newcommand{\mue}{\mu^{e}}	
\newcommand{\etae}{\eta^{e}}	
\newcommand{\bx}{\mathbf{x}}			

\newcommand{\Q}{\mathbf{Q}}					
\newcommand{\tQ}{\Q^{\ast}}					
\newcommand{\bQ}{{\Q}}	
\newcommand{\Qe}{\mathbf{Q}^{e}}			
\newcommand{\press}{p}						
\newcommand{\Press}{p}
\newcommand{\presse}{\press^{e}}
\newcommand{\mom}{\rho\bvel}
\newcommand{\bvar}{B}   
\newcommand{\evar}{E}   
\newcommand{\avar}{a}   
\newcommand{\bbvar}{\mathbf{\bvar}}   
\newcommand{\bevar}{\mathbf{\evar}}   

\newcommand{\Bvar}{B}   
\newcommand{\Evar}{E}   
\newcommand{\Avar}{A}
\newcommand{\bBvar}{\mathbf{\Bvar}} 
\newcommand{\bEvar}{\mathbf{\Evar}} 
\newcommand{\Bnvar}{\overline{\Bvar}} 
\newcommand{\bAvar}{\mathbf{\Avar}}
\newcommand{\Etvar}{\overline{\Evar}} 
\newcommand{\Atvar}{\overline{\Avar}}
\newcommand{\tbBvar}{\tilde{\bBvar}}
\newcommand{\magfieldnc}{\Bnvar}

\newcommand{\magfield}{\bBvar}
\newcommand{\bmagfield}{{\bBvar}}

\newcommand{\magfielde}{\magfield^{e}}
\newcommand{\elfield}{\bEvar}
 
\newcommand{\magpotential}{\mathbf{A}}
\newcommand{\etares}{\eta}
\newcommand{\polflux}{\Psi}
\newcommand{\Flux}{\boldsymbol{\mathcal{F}}}	
\newcommand{\Fluxcv}{\Flux_{\bvel}}	
\newcommand{\Fluxvv}{\Flux_{\mu}}	
\newcommand{\Fluxcb}{\Flux_{\magfield}}	
\newcommand{\Fluxvb}{\Flux_{\eta}}	
\newcommand{\bF}{\mathbf{F}}	
\newcommand{\bFcv}{\bF_{\bvel}}	
\newcommand{\bFvv}{\bF_{\mu}}	
\newcommand{\bFcb}{\bF_{\magfield}}	
\newcommand{\bFvb}{\bF_{\eta}}	

\newcommand{\alphacv}{\alpha^{c}_{ij}} 
\newcommand{\alphavv}{\alpha^{v}_{ij}} 
\newcommand{\lambdacv}{\lambda^{c}}
\newcommand{\lambdavv}{\lambda^{v}}
\newcommand{\domain}{\Omega}
\newcommand{\ddomain}{\domain}

\DeclareFloatingEnvironment[fileext=lop]{Diagram}

\usepackage{authblk}
\usepackage{fancyhdr}
\usepackage{changepage}
\newcommand{\address}[2]{\affil[#1]{#2}}
\newcommand{\corref}{\footnote{Corresponding author}}
\newcommand{\MPI}{1}
\newcommand{\DMUNITN}{2}
\newcommand{\UVigo}{3}
\newcommand{\LAM}{4}
\newcommand{\DISI}{5}
\newcommand{\TUM}{6}

\title{A well-balanced and exactly divergence-free staggered semi-implicit\\ hybrid finite volume / finite element scheme\\ for the incompressible MHD equations}  
\author[\MPI]{\underline{F. Fambri}\corref}
\author[\DMUNITN]{E. Zampa}
\author[\UVigo]{S. Busto}
\author[\LAM,\DISI]{L. R\'io-Mart\'in}
\author[\MPI]{F. Hindenlang}
\author[\MPI,\TUM]{E. Sonnendr\"ucker}
\author[\LAM]{M. Dumbser}
\address{\MPI}{Max-Planck-Institut für Plasmaphysik, Boltzmannstraße 2, D-85748 Garching, Germany}
\address{\DMUNITN}{Department of Mathematics, University of Trento, Via Sommarive 14, 38123 Trento, Italy} 

\address{\UVigo}{%
	Department of Applied Mathematics I, Universidade de Vigo, Campus As Lagoas Marcosende s/n, 36310 Vigo, Spain
	}  
    
\address{\LAM}{Laboratory of Applied Mathematics, DICAM, University of Trento, via Mesiano 77, 38123 Trento, Italy} 
\address{\DISI}{DISI, University of Trento, via Sommarive 9, Povo, 38123 Trento, Italy} 
\address{\TUM}{Technische Universität München, Zentrum Mathematik, D-85748, Garching, Germany}

\affil[ ]{\footnote{\textit{
francesco.fambri@ipp.mpg.de, 
enrico.zampa@unitn.it,
saray.bustoulloa@unitn.it,
laura.delrio@unitn.it,
florian.hindenlang@ipp.mpg.de,\newline
sonnen@ipp.mpg.de,
michael.dumbser@unitn.it}}}

\date{1 May, 2023}    
\newcommand{\keywords}[1]{\textbf{keywords}-- #1}

\allowdisplaybreaks

\begin{document} 
\newgeometry{left=1.00cm, right=1.00cm, top=2.00cm, bottom=2.00cm} 
\maketitle 
 
\begin{abstract}
We present a new exactly divergence-free and well-balanced hybrid finite volume/finite element scheme for the numerical solution of the incompressible viscous and resistive magnetohydrodynamics (MHD) equations on staggered unstructured mixed-element meshes in two and three space dimensions. The equations are split into several subsystems, each of which is then discretized with a particular scheme that allows to preserve some fundamental structural features of the underlying governing PDE system also at the discrete level. The pressure is defined on the vertices of the primary mesh, while the velocity field and the normal components of the magnetic field are defined on an edge-based/face-based dual mesh in two and three space dimensions, respectively. This allows to account for the divergence-free conditions of the velocity field and of the magnetic field in a rather natural manner. The non-linear convective and the viscous terms in the momentum equation are solved at the aid of an explicit finite volume scheme, while the magnetic field is evolved in an exactly divergence-free manner via an explicit finite volume method based on a discrete form of the Stokes law in the edges/faces of each primary element. The latter method is stabilized by the proper choice of the numerical resistivity in the computation of the electric field in the vertices/edges of the 2D/3D elements. To achieve higher order of accuracy, a piecewise linear polynomial is reconstructed for the magnetic field, which is guaranteed to be exactly divergence-free via a constrained $L^2$ projection. Finally, the pressure subsystem is solved implicitly at the aid of a classical continuous finite element method in the vertices of the primary mesh and making use of the staggered arrangement of the velocity, which is typical for incompressible Navier-Stokes solvers. In order to maintain non-trivial stationary equilibrium solutions of the governing PDE system exactly, which are assumed to be known \textit{a priori}, each step of the new algorithm takes the known equilibrium solution explicitly into account so that the method becomes exactly well-balanced. We show numerous test cases in two and three space dimensions in order to validate our new method carefully against known exact and numerical reference solutions. In particular, this paper includes a very thorough study of the lid-driven MHD cavity problem in the presence of different magnetic fields and the obtained numerical solutions are provided as free supplementary electronic material to allow other research groups to reproduce our results and to compare with our data. We finally present long-time simulations of Soloviev equilibrium solutions in several simplified 3D tokamak configurations, showing that the new well-balanced scheme introduced in this paper is able to maintain stationary equilibria exactly over very long integration times even on very coarse unstructured meshes that, in general, do not need to be aligned with the magnetic field lines. 
\end{abstract}


\restoregeometry


\section{Introduction}\label{sec:intro}

The equations of viscous and resistive magnetohydrodynamics (MHD) describe non-ideal plasma flows in the continuum limit and can be seen as a particular case of the more general Maxwell-Vlasov-Boltzmann system. The structure of viscous and resistive MHD is the one of the Navier-Stokes equations coupled with the Faraday law (induction equation), supplemented by appropriate closure relations for the stress tensor in the momentum equation and for the electric field in the induction equation. 
In the frame of classical parabolic models for the non-ideal (dissipative) effects in viscous and resistive MHD, the viscous stresses inside the fluid are modeled via the usual rheology of Newtonian fluids, while the resistivity can be described via an extra term in the electric field that is proportional to the curl of the magnetic field. If present, heat conduction is modeled via the usual Fourier law. More general first order hyperbolic formulations that describe general electromagnetic wave propagation in moving dielectric continuous media 
and which reduce to the viscous and resistive MHD equations in their stiff relaxation limit can be found in \cite{Rom1998,GPRmodelMHD}. Potential practical applications of the MHD equations range from astrophysics over  solar physics to nuclear fusion applications on Earth, such as magnetic confinement fusion (MCF) in tokamak or stellarator devices. 

It is well-known that the MHD equations are endowed with a stationary differential constraint which states that the divergence of the magnetic field must remain zero for all times, i.e., that no magnetic monopoles can exist. Typically, the full compressible MHD equations are considered, and there is an extensive literature on the topic. It is out of the scope of this paper to give an exhaustive overview. We only list some main representatives of typical discretization schemes. Numerical methods for the MHD equations can essentially be divided into two main classes: i) the first class which preserves the divergence-free condition of the magnetic field \textit{exactly} at the discrete level, using appropriately staggered meshes, see, e.g. \cite{BalsaraSpicer1999,BalsaraSpicer1999b,BalsaraAMR,Balsara2004,balsarahlle2d} and references therein, and ii) the second class which satisfies the divergence-free condition only \textit{approximately} at the discrete level. A very well-known representative of the second class is the Powell approach \cite{PowellMHD1,PowellMHD2}, which consists in adding multiples of the divergence constraint to the original governing PDE system, making use of the original ideas of Godunov on the symmetrization and thermodynamic compatibility of the MHD equations \cite{God1961,God1972MHD}. However, when using the Powell method on the discrete level, exact momentum and total energy conservation are, in general, lost. Another classical representative of the second class of methods is the generalized Lagrangian multiplier (GLM) approach of Munz \textit{et al.}, \cite{MunzCleaning,Dedneretal}, where the original MHD system is augmented via an additional PDE for an artificial cleaning scalar that is coupled to the induction equation via a grad - div pair of differential operators and which allows to transport divergence errors in the magnetic field away via acoustic-type waves so that they cannot accumulate locally. Unlike the Powell approach, the original GLM divergence cleaning for MHD conserves momentum and total energy exactly at the discrete level, but it is not thermodynamically compatible in the sense of Godunov. Very recently, in \cite{HTCMHD} it has been shown that the GLM approach can also be rigorously cast into a thermodynamically compatible and fully conservative form, while previous thermodynamically compatible versions of GLM for MHD violate the conservation of momentum and total energy at the discrete level, see, e.g. \cite{GassnerEntropyGLM}. In this paper, we will develop an \textit{exactly divergence-free} scheme, but we also show how parts of our algorithm can be used in the context of the hyperbolic GLM divergence cleaning framework.  

As already mentioned before, most of the well-established numerical schemes used for the solution of the MHD equations address the full compressible MHD system and are typically \textit{explicit}, of the upwind finite volume or finite difference type, based on Riemann solvers and use \textit{collocated meshes}, see, e.g. \cite{GardinerStone,BrioWu,RyuJones,fallemhd3,fallemhd2,HLLCMHD,HLLD,torrilhon,torrilhonbalsara,HelzelRossmanith,HiPa18}. More recently, also high order discontinuous Galerkin finite element schemes have been developed for MHD, see  \cite{divdg1,divdg2,divdg3,WarburtonVRMHD,EntropyDGMHD,GassnerEntropyGLM,Zanotti2015c}. 

A totally different numerical framework is semi-implicit schemes on \textit{staggered meshes}, which are by now a well-established standard for the numerical solution of the incompressible Navier-Stokes equations \cite{markerandcell,chorin1,chorin2,patankar,patankarspalding,BellColellaGlaz,vanKan}, and there exist also a few extensions of semi-implicit schemes to magnetohydrodynamics, see, e.g. \cite{mhdsi1,mhdsi2,mhdsi3,mhdsi4,Amari,Lerbinger,Harned,Finan,EULAG,SIMHD,Fambri20}. For implicit structure-preserving 
schemes for the incompressible MHD equations, the reader is referred to \cite{Ma2016,Hu2021,HuMaXu2017,Gawlik2021,Gawlik2022} and references therein. 

In order to preserve stationary equilibrium solutions exactly, a numerical method must be \textit{well-balanced}. The concept of well-balancing was introduced for the first time in \cite{Bermudez1994} for the shallow water equations and in \cite{BottaKlein} for the compressible Euler equations with gravity. Subsequently, many substantial contributions have been made to the literature on well-balanced methods, and an exhaustive overview is impossible and clearly out of the scope of the present paper. For an incomplete list of references related to the topic of well-balancing, the interested reader is referred to \cite{leveque,greenbergleroux,Noelle1,Noelle2,gosse2000well,pares2004well,castro2007wellb,Kapelli2014,Klingenberg2015,xing2013high,KlingenbergPuppo,GCD18,ThomannWB1,ThomannWB2,GCD21,WBGravity} and references therein. So far, there are only very few well-balanced schemes available for the MHD equations, see, e.g. \cite{WBMHD}, but the method presented there is explicit, and it is only available on Cartesian meshes. 

The objective of this paper is therefore to design a new structure-preserving scheme for the non-ideal viscous and resistive magnetohydrodynamics equations that is at the same time exactly divergence-free, well-balanced in order to maintain arbitrary, but \textit{a priori} known, equilibrium solutions exactly at the discrete level and which is able to run on general unstructured mixed-element meshes in two and three space dimensions, see also \cite{MixedWENO2D,MixedWENO3D}. The general unstructured grids are necessary in order to deal with complex geometries like those of tokamaks or stellarators. To ease the development and the presentation of the method, in this paper, we restrict ourselves to the incompressible MHD equations, while the extension to the full, compressible case is left to future work. 

The rest of this paper is organized as follows. In Section \ref{sec:goveq}, we briefly present the governing partial differential equations. In Section \ref{sec:numdisc}, a new well-balanced and exactly divergence-free semi-implicit hybrid finite volume / finite element scheme is designed for the solution of the incompressible viscous and resistive MHD equations on general unstructured mixed-element meshes. Numerical results are shown in Section~\ref{sec:numericalresults} for a large set of test problems in two and three space dimensions, including also an application of the method to simplified 3D tokamak geometries. The paper closes with some concluding remarks and a brief outlook to future work given in Section~\ref{sec:conclusions}. 

\section{Governing partial differential equations} \label{sec:goveq}
The incompressible viscous and resistive MHD equations have been shown to describe many dynamical properties of hot, strongly magnetized plasma in the low Mach number regime. This model consists of the coupling of the conservation of momentum with the Faraday (induction) equation of electrodynamics:  
\begin{align}
&\partial_t\left( \mom \right) \ + \dive \left(\Fluxcv - \Fluxvv\right) +\gra \press = 0, \label{eq:mom}\\
& \partial_t \bbvar  + \curl \bevar = 0.
\label{eq:Faraday} 
\end{align} 
Above, $\rho$ is the constant fluid density, $\press$ denotes the pressure, $\bevar$ corresponds to the electric field vector and $\bvel=(\vel_1,\vel_2,\vel_3)$ and $\bbvar=(\bvar_1,\bvar_2,\bvar_3)$ are the fluid velocity and the magnetic field vectors, respectively. The inviscid flux tensor $\Fluxcv$ and the viscous stress tensor $\Fluxvv$ will be specified later. Since we assume an incompressible fluid, the velocity field must be  divergence-free, i.e. 
\begin{equation}
\dive \bvel = 0. \label{eq:divu} 
\end{equation} 	
Eqn.~\eqref{eq:divu} is what remains from the pressure equation of the compressible Navier-Stokes equations when the Mach number tends to zero, see \cite{KlaMaj,KlaMaj82,munzMPV,Klein2001}, and the fluid pressure $\press$ in the momentum equation \eqref{eq:mom} can be seen as a Lagrange multiplier that assures that the divergence-free condition of the velocity field \eqref{eq:divu} is always satisfied. The coupling of \eqref{eq:mom} with \eqref{eq:divu} leads to the well-known elliptic pressure Poisson equation, i.e. the pressure adjusts globally and instantaneously so that \eqref{eq:divu} is guaranteed everywhere. 

Instead, a direct consequence of the Faraday law \eqref{eq:Faraday} is that the magnetic field must remain divergence-free for all times, if it was initially divergence-free, hence
\begin{align}
	\dive \bbvar = 0. \label{eq:divB}
\end{align}
The physical meaning of \eqref{eq:divB} is that magnetic monopoles cannot exist, and the mathematical justification is that the divergence of the curl is zero. 
Although apparently similar, the two divergence-free conditions above are totally different in nature. While the magnetic field is automatically guaranteed to remain divergence-free for all times as a mere consequence of \eqref{eq:Faraday} and standard vector analysis identities, the divergence-free condition of the velocity field is the PDE that is needed to determine the pressure in \eqref{eq:mom}. 
The governing equations \eqref{eq:mom}-\eqref{eq:Faraday} with \eqref{eq:divu} are a coupled system of non-linear hyperbolic-parabolic-elliptic PDE. 

To close the system, we still need to provide the definitions for the electric field and the two flux tensors. 
The electric field vector, $\bevar$, is given by 
\begin{equation}
\bevar= - \bvel \times \bbvar + \eta \curl \bbvar,
\label{eq:E} 
\end{equation}
where the first term on the right-hand side corresponds to the electric field in ideal MHD, and the second term includes the dissipative effects in non-ideal, resistive MHD, with $\eta$ the electric resistivity. 
In~\eqref{eq:mom}, $\Fluxcv=\Fluxcv(\rho,\bvel,\bbvar)$ is the tensor of inviscid fluxes, which include the non-linear advection and the Lorentz force, while the tensor $\Fluxvv=\Fluxvv(\mu, \nabla \bvel)$ contains all the purely dissipative fluxes, i.e., the viscous stresses. 
More specifically, we have 
\begin{equation}  \label{eqn.ffluxmom}
	\Fluxcv:=  
\mom \otimes \bvel + \halb \bbvar^2 \, \mathbf{I} -  \bbvar \otimes \bbvar, 
\quad\quad \Fluxvv= \Fluxvv\left(\mu, \nabla \bvel \right) :=    \mu \left( \gra \bvel + \gra \bvel^T \right),  
\end{equation}
where $\mathbf{I}$ is the  identity matrix, and $\mu$ is the dynamic viscosity. 
For simplicity, in this work, the fluid viscosity and resistivity are supposed to be some non-negative constants. Moreover, the density is assumed to be constant in space and time, so that \eqref{eq:divu} can be replaced by
\begin{equation}
	\dive \ww = 0, \label{eq:divw}
\end{equation}
with  $\ww:=\rho\bvel$ the linear momentum.

It will turn out to be useful to introduce also the corresponding conservative form of the Faraday equation, i.e. 
\begin{align}
\partial_t \bbvar + \dive \left( \Fluxcb -\Fluxvb\right) = {0}, \qquad \Fluxcb :=\bbvar \otimes \bvel - \bvel \otimes\bbvar , \qquad \Fluxvb = \Fluxvb(\eta, \nabla \magfield) := \eta \left( \gra \bbvar - \gra \bbvar^T \right). \label{eq:FaradayCons} 
\end{align}
Therefore, the PDE system \eqref{eq:mom}-\eqref{eq:FaradayCons} may also be written as
\begin{subequations} \label{eq:PDEq} 
\begin{eqnarray}
   \partial_t \ww  + \dive \left( \Fluxcv(\Q) - \Fluxvv(\mu, \nabla \bvel) \right) + \gra \press & = &  0, \\   
   \partial_t \bbvar + \dive \left( \Fluxcb(\Q)-\Fluxvb(\eta, \nabla \bbvar) \right) & = & 0,
\end{eqnarray}
\end{subequations}
where $\Q := \left( \ww, \bbvar\right)^T$ is the vector of conservative variables, and $\rho$ is a given constant parameter.

\section{Well-balanced numerical method} \label{sec:numdisc}
To develop the new well-balanced semi-implicit hybrid finite volume / finite element method for the discretization of the system \eqref{eq:mom}-\eqref{eq:divB}, we first introduce an equilibrium solution 
$\Qe:=\left(\WWe,\magfielde\right)^{T}$, with \mbox{$\WWe=\rho^{e}\bvele$,} $\rho^{e}=\rho$ and $\presse$. Since $\Qe$ is assumed to be a stationary solution of the governing partial differential equations, that must be preserved exactly also at the discrete level, \cite{BVC94,Pares2006,Castro2008,GDC17,GCD18,BCK21,GCD21,CG23}, we have $\partial_t \Qe = 0$, hence 
\begin{align}
			\dive \left( \wwe \right) = 0, \label{eq:divwe}\\
			\dive \Fluxcv\left(\Qe\right) - \dive \Fluxvv\left(\mue, \nabla \bvele \right) + \gra \presse = 0,
			\label{eqn.mome}\\
			\dive \Fluxcb\left(\Qe\right) - \dive \Fluxvb\left(\etae, \nabla \magfielde \right) = 0 .
			\label{eqn.mage}
\end{align}
Subtracting the former relations from \eqref{eq:divw} and \eqref{eq:PDEq}, we get
\begin{align}
	\dive \left( \ww - \wwe \right) = 0, \label{eq:divwe0}\\
	\partial_t \ww + \dive \Fluxcv\left(\Q\right)  - \dive \Fluxcv\left(\Qe\right) - \dive \Fluxvv\left(\mu, \nabla \bvel \right) + \dive \Fluxvv\left(\mue, \nabla \bvele \right) + \gra \left( \press-\presse\right)=0,
	\label{eqn.convvissubsv0}\\
	\partial_t \magfield + \dive \Fluxcb\left(\Q\right)- \dive \Fluxcb\left(\Qe\right) - \dive \Fluxvb\left(\eta, \nabla \magfield \right)+ \dive \Fluxvb\left(\etae, \nabla \magfielde \right)=0.
	\label{eqn.convvissubsB0}
\end{align}

Next, after performing a semi-discretization in time, see, e.g. \cite{BFSV14,BFTVC17,Hybrid1,Hybrid2}, system  \eqref{eq:divwe0}-\eqref{eqn.convvissubsB0} is split into the following subsystems: 

\paragraph{Convective and viscous subsystem; auxiliary conservative evolution system of the magnetic field}
\begin{subequations}\label{eqn.convvissubs}
	\begin{align}
		\frac{1}{\Delta t} \left( \tWW{}-\WW^{n}\right) = - \dive \Fluxcv\left(\Q^{n}\right)  + \dive \Fluxcv\left(\Qe\right) + \dive \Fluxvv\left(\mu, \nabla \bvel^{n} \right) - \dive \Fluxvv\left(\mue, \nabla \bvele \right) - \gra \left( \press^{n}-\presse\right),
		\label{eqn.convvissubsv}\\
		\frac{1}{\Delta t} \left(\magfield^{n+1} - \magfield^{n} \right) =- \dive \Fluxcb\left(\Q^{n}\right)+ \dive \Fluxcb\left(\Qe\right) + \dive \Fluxvb\left(\eta, \nabla \magfield^{n} \right)- \dive \Fluxvb\left(\etae, \nabla \magfielde \right).
		\label{eqn.convvissubsB}
	\end{align}
\end{subequations}

\paragraph{Faraday subsystem}
\begin{equation}
	\frac{1}{\Delta t} \left( \magfield^{n+1} - \magfield^{n}\right) =   - \curl \bevar(\magfield^n,\bvel^n) + \curl \bevar(\magfielde,\bvele),
	\label{eq:Faraday2} 
\end{equation}

\paragraph{Pressure subsystem}
\begin{subequations}\label{eq:press_sys}
	\begin{align}
		&\dive  \left( \WW^{n+1} - \WWe \right)  = 0, \label{eqn.pressubs1}\\
		&\frac{1}{\Delta t} \left( \WW^{n+1}-\tWW{}\right)  +\gra  \left( \press^{n+1}-\press^{n}\right)  = 0.
		\label{eqn.pressubs2}
	\end{align}
\end{subequations}
Combination of the two pressure equations above leads to the pressure Poisson equation 
\begin{equation}
	 - \nabla^2 \left( \press^{n+1}-\press^{n} \right) = - \frac{1}{\Delta t} \dive \left( \tWW{} - \WWe \right), 
	 \label{eqn.pressurePoisson}
\end{equation}
which now also takes into account the stationary equilibrium solution in order to obtain the well-balance property of our new scheme. 

The auxiliary variable $\tWW{}$ introduced above denotes an intermediate momentum accounting only for the convective and viscous terms, without the contribution of the pressure gradient at the new time. In general, the intermediate velocity field $\tWW{}$ is not yet divergence-free. Therefore, equation \eqref{eqn.pressubs2} will be employed not only to build the pressure system but also, to obtain the final momentum $\WW^{n+1}$.

To introduce the spatial discretization, we make use of unstructured staggered mixed-element grids, where the staggered dual mesh is of the edge-type in 2D and of the face-type in 3D, see also \cite{USFORCE}. One valuable improvement with respect to former hybrid FV/FE methods, presented in \cite{BFSV14,BFTVC17,Hybrid1,Hybrid2,HybridMPI,HybridALE,HybridSWE,HybridImplicit}, is that the unstructured grids considered are now of mixed elements, i.e., they may potentially contain triangles and quadrilaterals in 2D and tetrahedra, hexahedra, square pyramids and triangular prisms in 3D simultaneously. This new feature helps on the design of meshes for complex domains and when requiring the grid to be aligned with specific fields as classically done in the simulation of tokamaks, see \cite{Tamain2016,Leddy2017,Dingfelder2020}, although more recently also non field-aligned unstructured meshes have been used in the context of high order discontinuous Galerkin finite element schemes, see, e.g. \cite{Giorgiani2020}. 

Regarding the numerical method applied to the hydrodynamics part, i.e., without the magnetic field, a major difference with respect to \cite{BFSV14,BFTVC17,HybridMPI,HybridNNT} are the spaces used to discretize the different variables. Even if staggered dual grids are also employed in the method proposed in this paper, in order to avoid stability issues such as checkerboard phenomena, all variables are computed on the primal grid, and then an interpolation between the primal and the dual grid is applied to the momentum unknowns. Besides, the main contributions of the method proposed in this paper are that it is designed to work on general unstructured mixed-element meshes, it is built to be well balanced for stationary equilibria, and that it is exactly divergence-free for the magnetic field. 

In what follows, we first introduce the employed staggered grids in two and three space dimensions and the related notation (Section~\ref{sec:mesh}), we then present the overall algorithm and describe each of its stages. 


\subsection{Staggered mixed-element unstructured meshes} \label{sec:mesh} 
To discretize the spatial domain $\Omega \subset \mathbb{R}^d$, with $d$ the number of space dimensions and $\partial \Omega$ its boundary, general unstructured mixed-element staggered faced-based grids are employed in this work. The use of staggered face-type grids is a well-documented approach when employing semi-implicit schemes, both in the framework of discontinuous Galerkin and hybrid FV/FE methods, see, e.g. \cite{TD14,TD16,HybridALE,HybridImplicit}. Here, we extend this type of staggered grids to a wider variety of primal elements.

\begin{figure}
	\centering
	\includegraphics[width=0.3\linewidth]{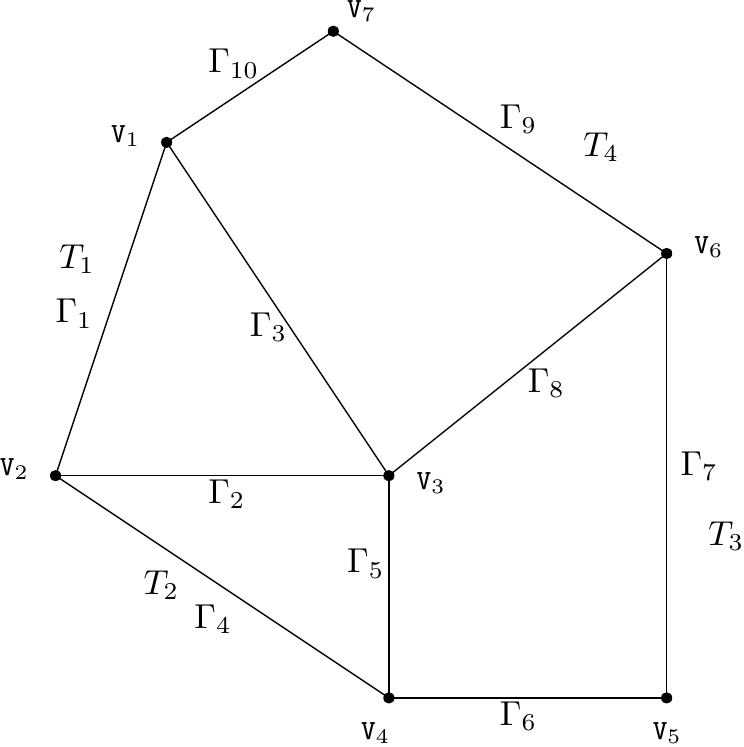}\hfill
	\includegraphics[width=0.3\linewidth]{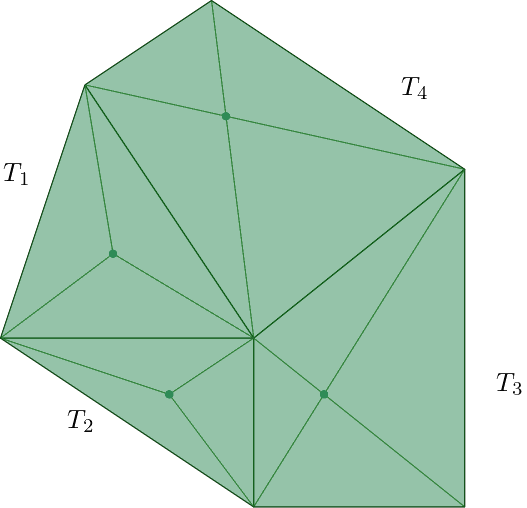}\hfill
	\includegraphics[width=0.3\linewidth]{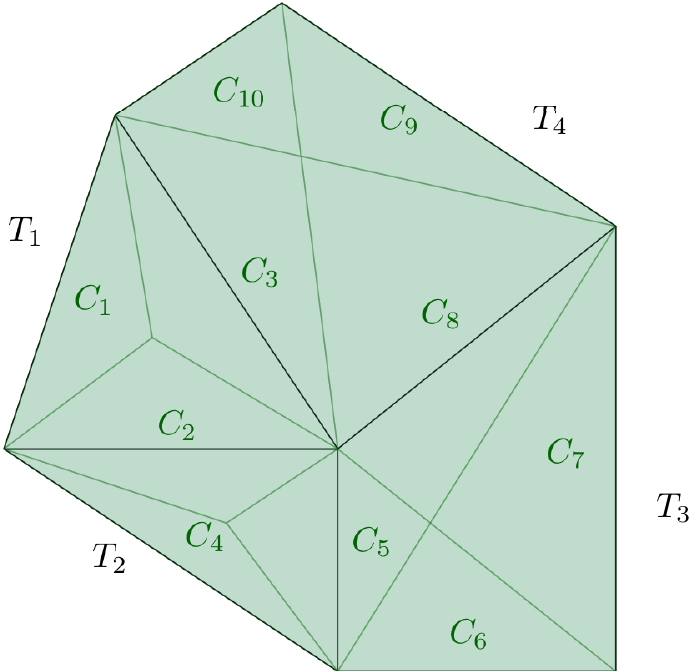}
	\caption{Construction of a dual mixed mesh. Left: primal mesh made of two triangles, $T_1$, $T_2$, and two quadrilaterals, $T_3$, $T_4$. The vertex are identified with $\vertex_k$ while the edges are denoted by $\Gamma_f$. Centre: the barycenters of primal elements are obtained and used to define the triangular subelements related to the faces. Right: merging of the two subelements related to each face results on the corresponding dual element, $C_i$. Interior cells are denoted by $C_2$, $C_3$, $C_5$, $C_8$. Boundary dual cells, $C_1$, $C_4$, $C_6$, $C_7$, $C_9$, $C_{10}$,  are constructed to be equal to the boundary subelement.}
	\label{fig:buildmesh2d}
\end{figure}

For the 2D case, we assume the primal elements to be triangles or quadrilaterals, see the left image in Figure~\ref{fig:buildmesh2d}. Meanwhile, in the 3D case, we have tetrahedra, hexahedra, quadrilateral pyramids, and triangular prisms, as shown in the top left plots in Figures \ref{fig:buildmesh3d}, \ref{fig:buildmesh3d_HexaPrism}, and \ref{fig:buildmesh3d_PiramidTetra}. We denote by $T_i$ a generic primal element, $i\in\left\lbrace 1,\dots,N_{\mathrm{el}}\right\rbrace$, with $N_{\mathrm{el}}$ the total number of primal elements, composed of vertices  $\vertex_{i_p}\equiv \vertex_k$, with $p$ the local vertex number inside element $T_i$ and $k\in\left\lbrace 1,\dots,N_{\mathrm{ver}}\right\rbrace$ the associated global vertex number, while $N_{\mathrm{ver}}$ is the total number of primal vertices. The faces in 3D are denoted by $\face_{i_m}\equiv \face_f$, with $m$ the local face number in element $T_i$, $f\in\left\lbrace 1,\dots,N_{\mathrm{face}}\right\rbrace$ the global face number, and $N_{\mathrm{face}}$ the total number of primal faces. The common face shared by two elements $T_i$ and $T_j$ is also denoted by $\partial T_{ij} = \face_{ij} = T_i \cap T_j$, while the boundary of element $T_i$ is given by the union of its faces, i.e. $\pbf{i}= \underset{m}{\cup} \face_{i_m} = \underset{j}{\cup} \face_{ij}$, with $T_j$ the neighbours of $T_i$. 

The edges in 3D are denoted by $\edge_{i_q}\equiv \edge_e$, with $q$ the local number of the edge inside element $T_i$, $e\in\left\lbrace 1,\dots,N_{\mathrm{edge}}\right\rbrace$ the global edge number, and $N_{\mathrm{edge}}$ the total number of primal edges. In two space dimensions, the faces of the elements reduce to the edges and the edges reduce to the vertices. 

The local number of vertices, edges, and faces by element depends on the element type and has been recalled in Table~\ref{tab.loclvertexedgesbyelement}. Note that, to ease the 2D/3D presentation of the method, we will employ the term faces even when referring to the edges in 2D if their role corresponds to the one played by the faces in 3D. The volume/area of a primal element is denoted by $\VT{i}$, the outward pointing normal to its boundary, $\pbf{i}$, is $\tnn_i$, and $\nn_{ij}:=\AF{ij}\tnn_{ij}$ refers to the normal at face $\pbf{ij}$, the common face between $T_i$ and $T_j$, with $\tnn_{ij}$ the unit normal and $\AF{ij}$ the area/length of the face. Moreover, $\mathcal{N}_{i}$ is the set of neighbour elements of $T_i$ on the primal mesh which share a common face.
\begin{table}[h]
	\begin{center}
		\begin{tabular}{|c||c||c|c|c|}
			\hline
			Dimension & Element type & $N_{i_v}$  &$N_{i_e}$ & $N_{i_f}$ \\
			\hline\hline
			\multirow{2}{*}{2D}& Triangle & 3 & 3 & 3\\
			\hhline{|~||-|-|-|-|}
			& Quadrilateral & 4 & 4 & 4\\
			\hline\hline
			\multirow{4}{*}{3D}& Tetrahedron & 4 & 6 & 4\\
			\hhline{|~||-|-|-|-|}
			& Hexahedron & 8 & 12 & 6\\
			\hhline{|~||-|-|-|-|}
			& Square pyramid & 5 & 8 & 5 \\
			\hhline{|~||-|-|-|-|}
			& Triangular prism & 6 & 9 & 5 \\
			\hline
		\end{tabular}
	\end{center}
	\caption{Number of local vertex, $N_{i_v}$, local edges, $N_{i_e}$, and local faces, $N_{i_f}$, for each primal element type.}\label{tab.loclvertexedgesbyelement}
\end{table}

To construct the dual grid, we proceed as follows. First, the barycenters of all primal elements, $\bary_i$, are computed. This allows the division of each element $T_i$ in as many subelements, $T_{i_f}$, as it has faces, $N_{i_f}$. Therefore, for the 2D case (center image of Figure~\ref{fig:buildmesh2d}) the subelements are triangles with basis one of the faces of the primal element and the opposite vertex the barycenter of the primal cell. In 3D, if the corresponding face is a triangle, then the resulting subcell is a tetrahedron, having the triangular face as a basis and as the opposite vertex, the barycenter of the primal cell, see Figure~\ref{fig:buildmesh3d_PiramidTetra}. Instead, for a quadrilateral face, we get a pyramid with the quadrilateral primal face as basis and the opposite vertex being the barycenter of the element, see Figures~\ref{fig:buildmesh3d}-\ref{fig:buildmesh3d_HexaPrism}. We denote by $T_{ij}$ the subelement of $T_i$ associated with the face $\face_{ij}$ that is in common with the neighbour primal element $T_j$, and by $T_{ji}$ the subelement of the primal element $T_j$ having a the same common face $\face_{ji} = \face_{ij}$ with the element $T_i$. In consequence, $|T_{ij}|$ represents the volume/area of the subelement $T_{ij}$. The dual elements, $C_{k}$, $k\in\left\lbrace 1,\dots,N_{\mathrm{dual}}\right\rbrace$, are then built by merging the two subelements generated at each side of a face, see Figures~\ref{fig:buildmesh2d},  \ref{fig:buildmesh3d}, \ref{fig:buildmesh3d_HexaPrism} and \ref{fig:buildmesh3d_PiramidTetra}. For the faces on the boundary of the computational domain, the dual element will simply be the subelement related to the face, see right panel of Figure~\ref{fig:buildmesh2d}. Hence, the vertex of each dual element will be the vertex of the face used to build it and the two barycenters of the elements containing that face. In 2D, all interior dual elements will be quadrilaterals, see Figure~\ref{fig:buildmesh2d}, while in 3D, we have two types of interior elements: one coming from merging two pyramids, see Figures~\ref{fig:buildmesh3d}-\ref{fig:buildmesh3d_HexaPrism}, and another one when gluing two tetrahedra, see Figure~\ref{fig:buildmesh3d_PiramidTetra}. In this case, the dual elements on the boundary of $\Omega$ are either pyramids or  tetrahedra, depending on the shape of the boundary face.

\begin{figure}
	\centering
	\includegraphics[width=0.33\linewidth]{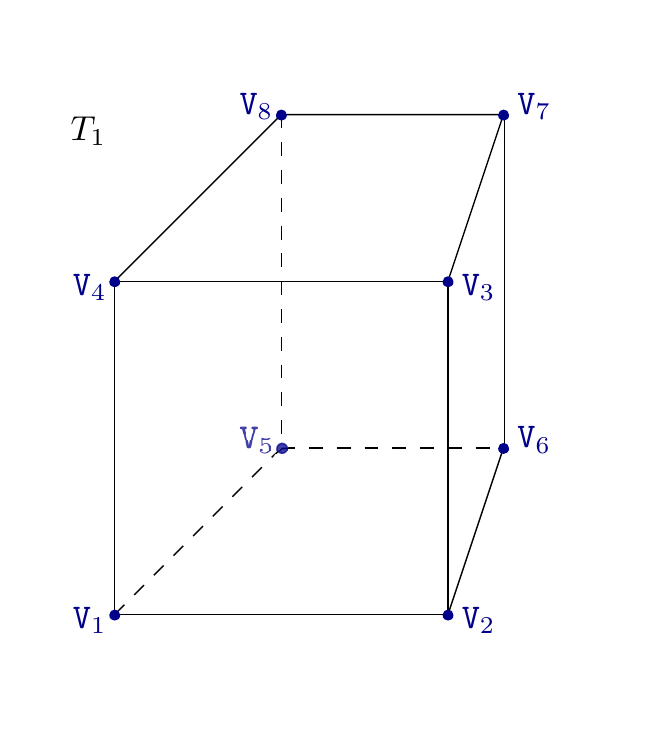}\hfill
	\includegraphics[width=0.33\linewidth]{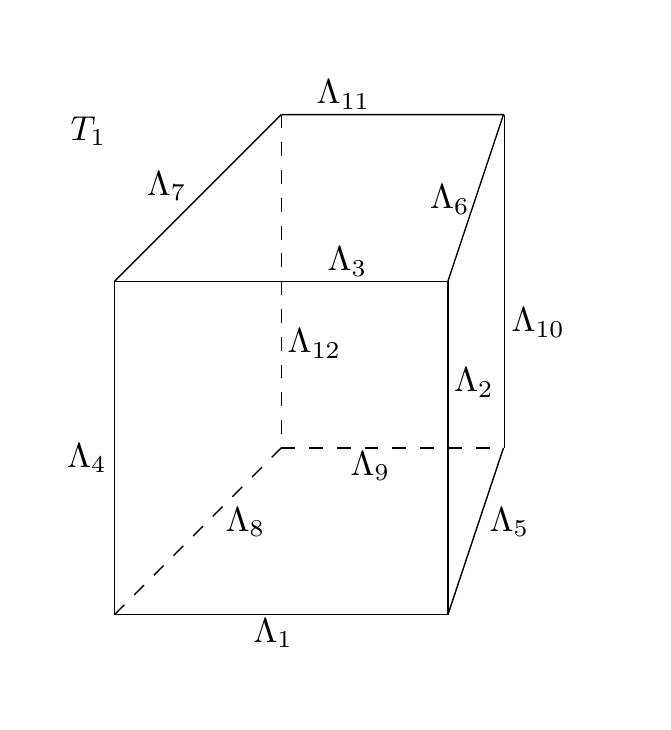}\hfill
	\includegraphics[width=0.33\linewidth]{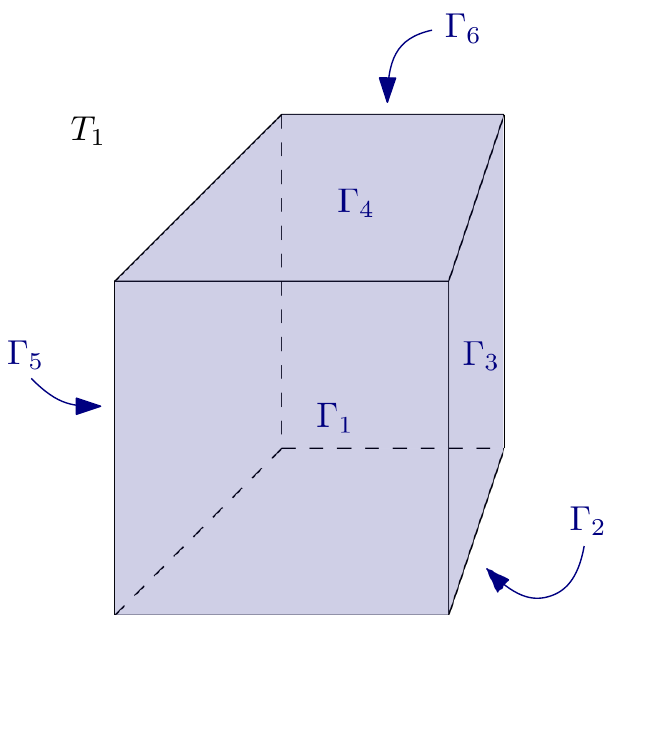}
	\caption{Notation of a tridimensional primal element $T_{1}$: the vertex are identified with $\vertex_k$ (left), the edges are denoted by $\edge_e$ (centre) and the faces are denoted by $\Gamma_f$ (right).}
	\label{fig:primalmeshnotation3d}
\end{figure}

\begin{figure}
	\centering
	\includegraphics[width=0.25\linewidth]{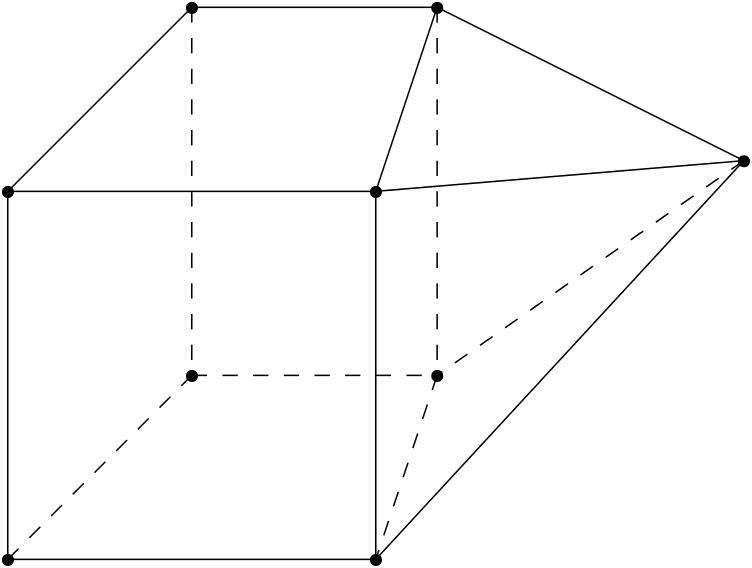}\hfill
	\includegraphics[width=0.25\linewidth]{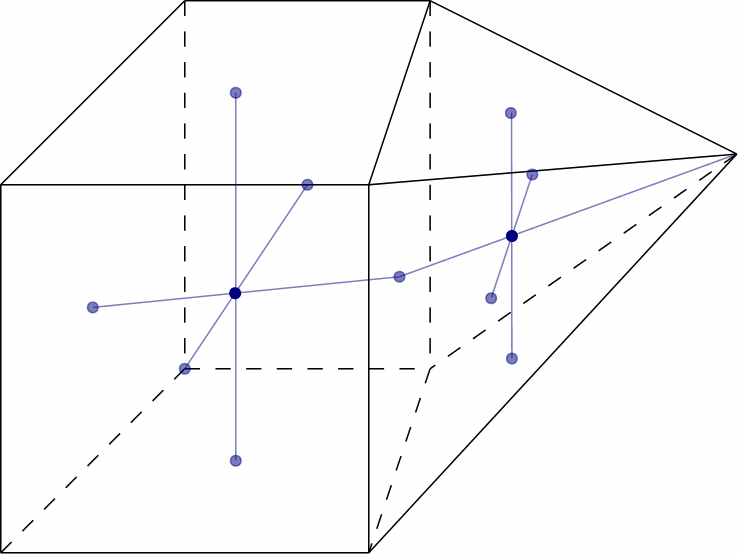}\hfill
	\includegraphics[width=0.25\linewidth]{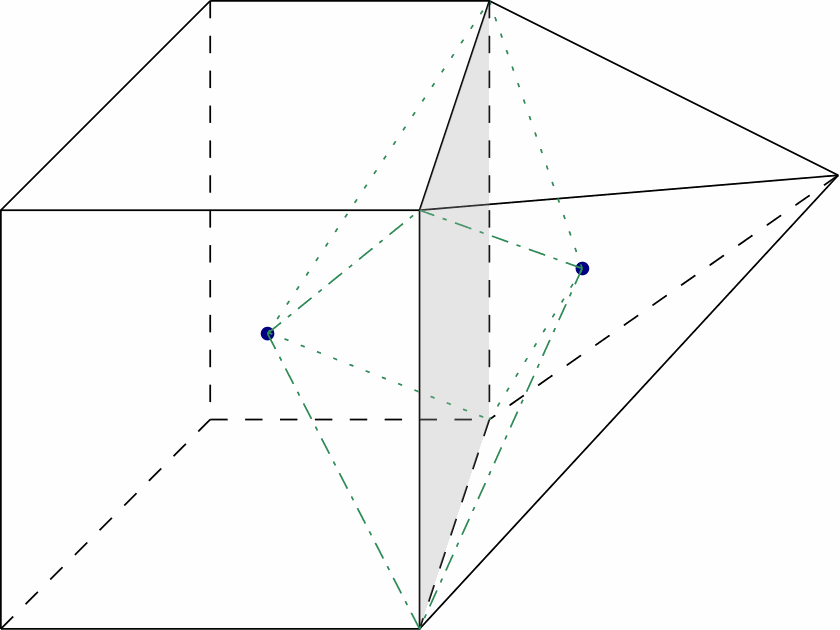} \\
	\vspace{0.05\linewidth}
	\includegraphics[width=0.25\linewidth]{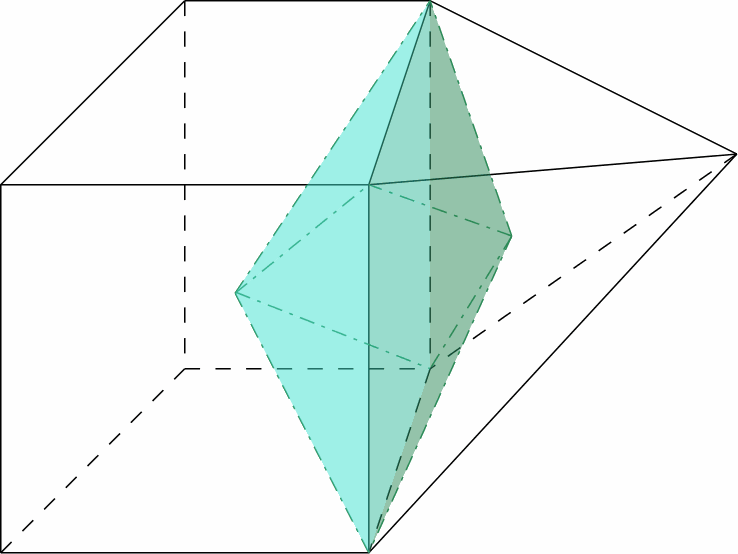}\hfill
	\includegraphics[width=0.25\linewidth]{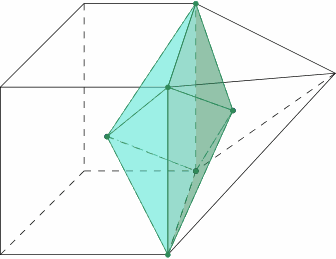}\hfill
	\includegraphics[width=0.25\linewidth]{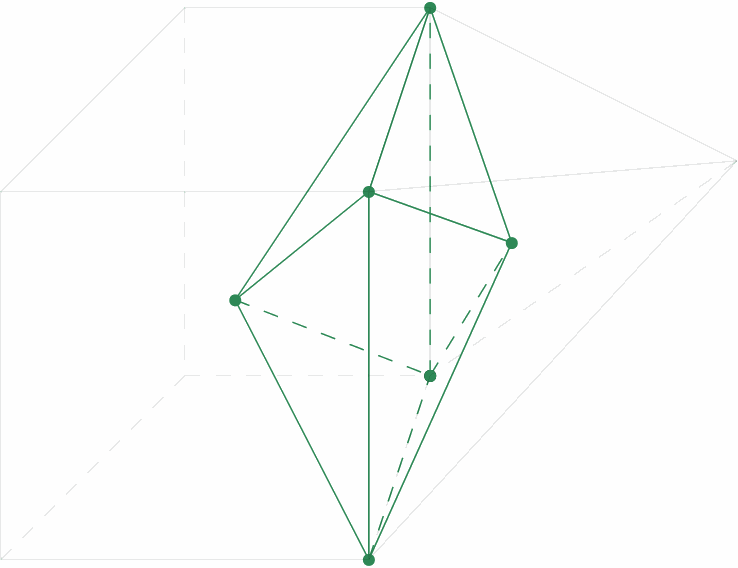}
	\caption{Construction of a dual element in 3D based on a face belonging to a hexahedron and a square pyramid. Top left: primal hexahedron (left), primal square pyramid (right), and primal vertex (black dots). Top centre: construction of the barycenters of the primal elements. Top right: each barycenter is connected to the vertex of the common face (shadowed in grey), generating a pyramid on each side of the face with basis the grey face and opposite vertex the barycenters. Bottom left: the lateral faces of the generated pyramids are shaded in light green (left, inside the hexahedron) and sea green (right, inside the primal pyramid). Bottom centre: the vertex of the two new pyramids are marked with green dots; continuous and discontinuous green lines indicate visible and shadow edges when merging both pyramids. Bottom right: the new dual element corresponds to the volume generated by merging the two pyramids constructed.}
	\label{fig:buildmesh3d}
\end{figure}

\begin{figure}
	\centering
	\includegraphics[width=0.25\linewidth]{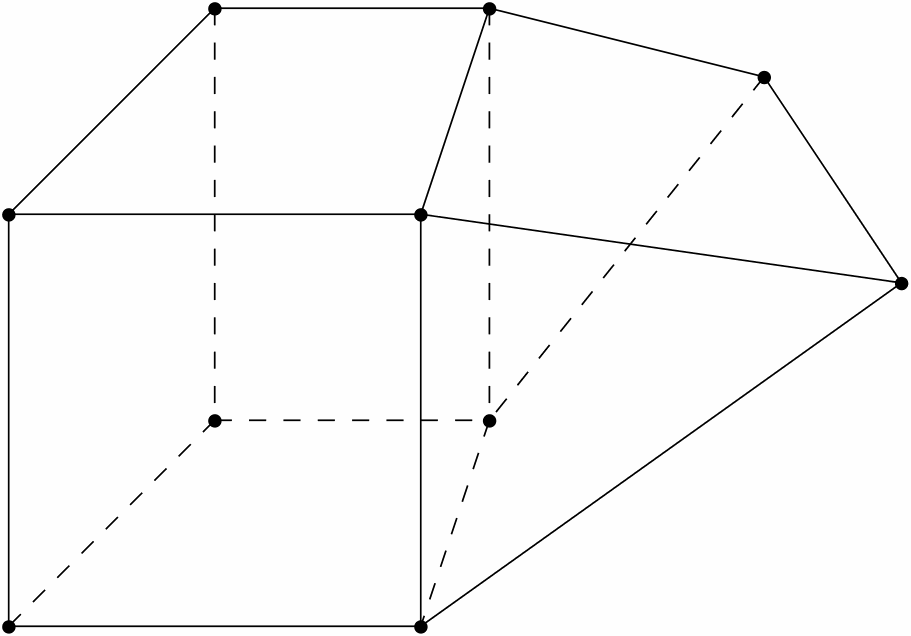}\hfill
	\includegraphics[width=0.25\linewidth]{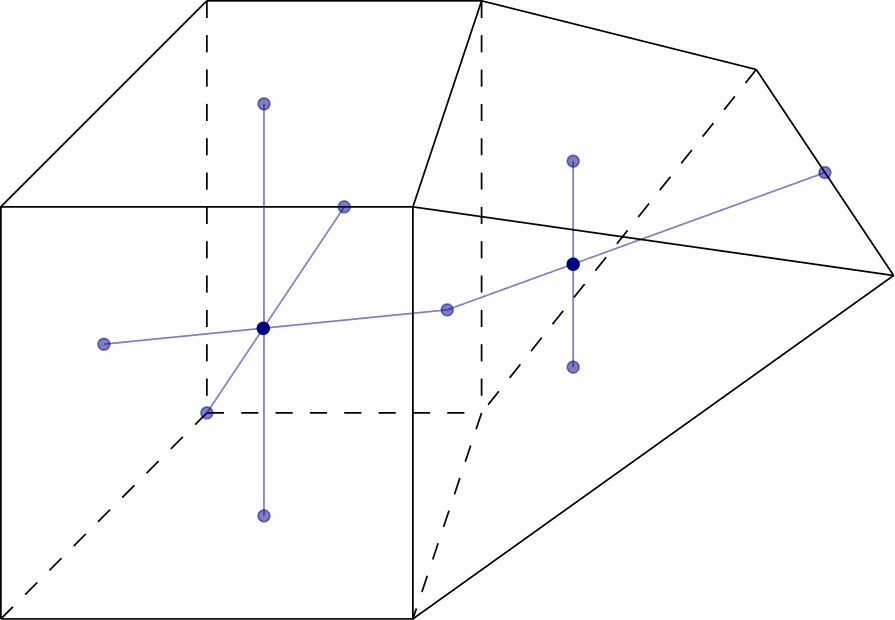}\hfill
	\includegraphics[width=0.25\linewidth]{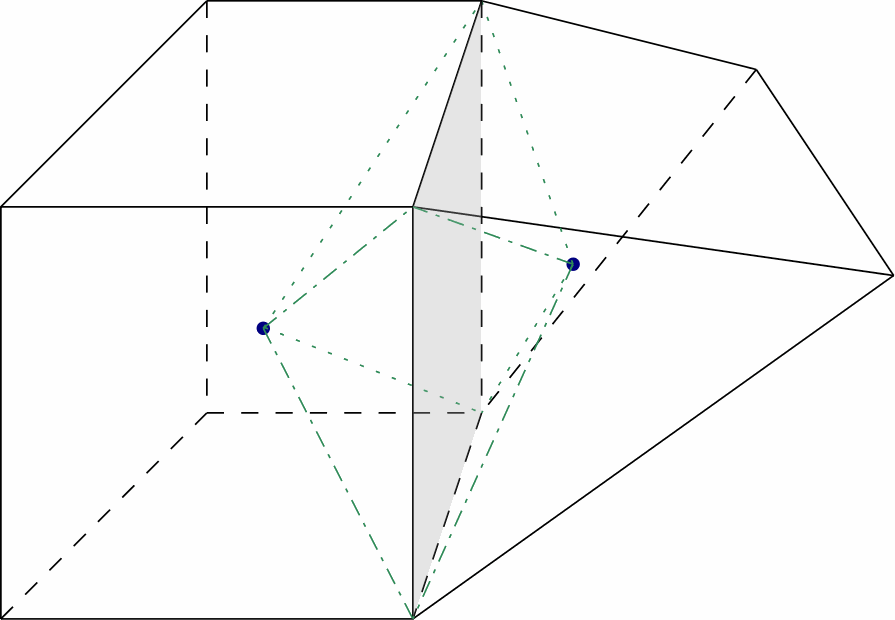} \\
	\vspace{0.05\linewidth}
	\includegraphics[width=0.25\linewidth]{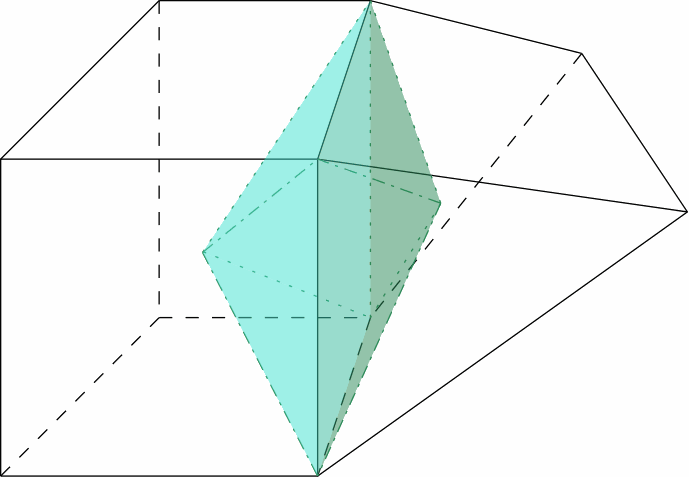}\hfill
	\includegraphics[width=0.25\linewidth]{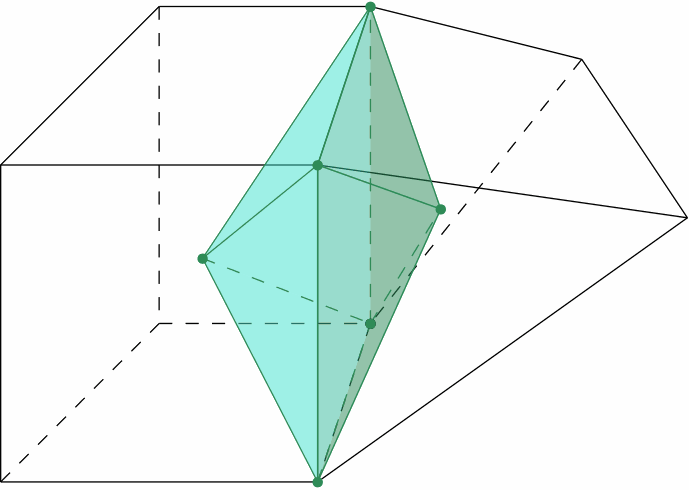}\hfill
	\includegraphics[width=0.25\linewidth]{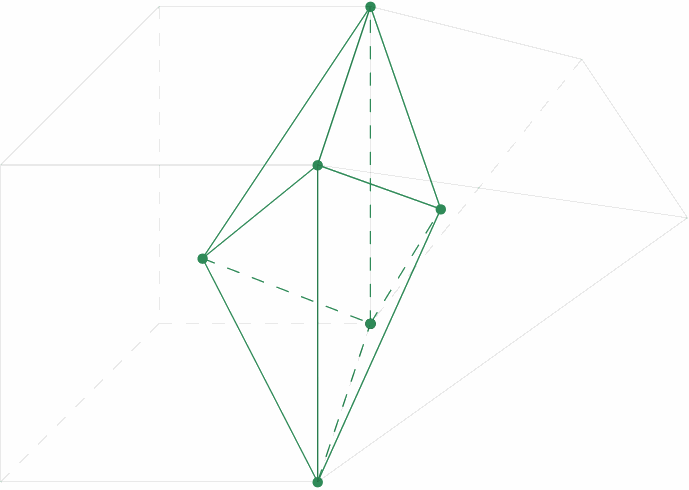}
	\caption{Construction of a dual element in 3D based on a face belonging to a hexahedron and a triangular prism. Top left: primal hexahedron (left), primal triangular prism (right), and primal vertex (black dots). Top centre: construction of the barycenters of the primal elements. Top right: each barycenter is connected to the vertex of the common face (shadowed in grey), generating a pyramid on each side of the face with basis the grey face and opposite vertex the barycenters. Bottom left: the lateral faces of the generated pyramids are shadowed in light green (left, inside the hexahedron) and sea green (right, inside the primal pyramid). Bottom centre: the vertex of the two new pyramids are marked with green dots; continuous and discontinuous green lines indicate visible and shadow edges when merging both pyramids. Bottom right: the new dual element corresponds to the volume generated by merging the two pyramids constructed.}
	\label{fig:buildmesh3d_HexaPrism}
\end{figure}

\begin{figure}
	\centering
	\includegraphics[width=0.25\linewidth]{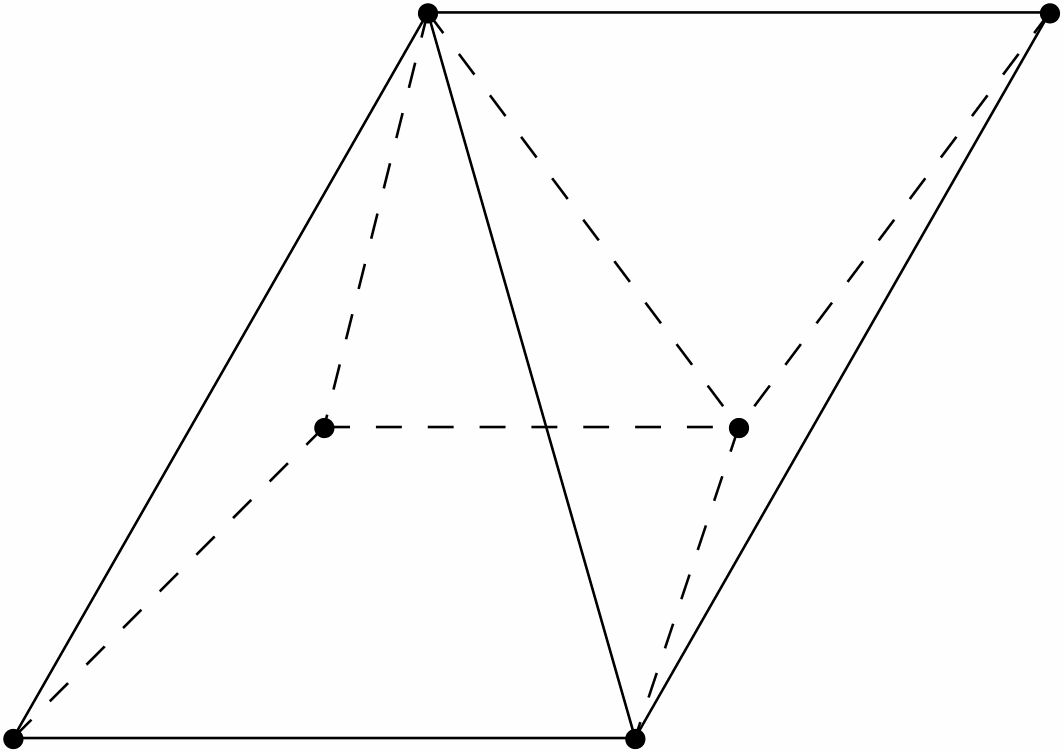}\hfill
	\includegraphics[width=0.25\linewidth]{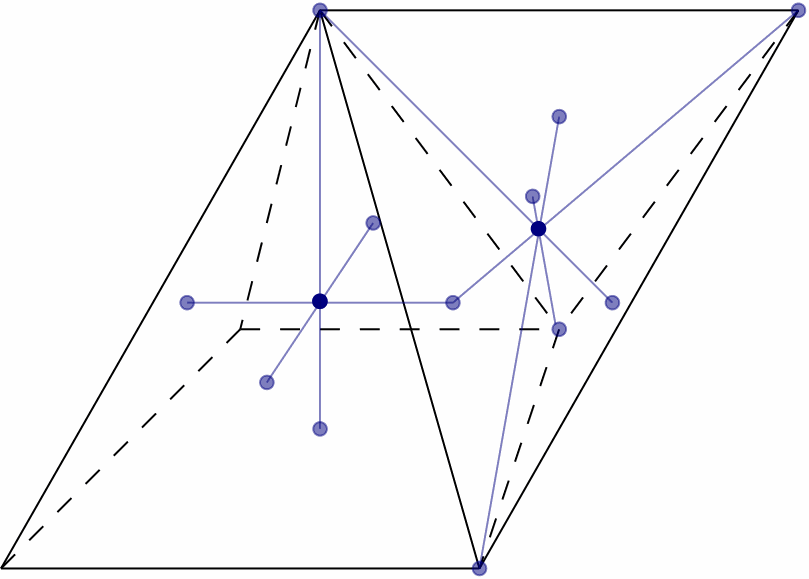}\hfill
	\includegraphics[width=0.25\linewidth]{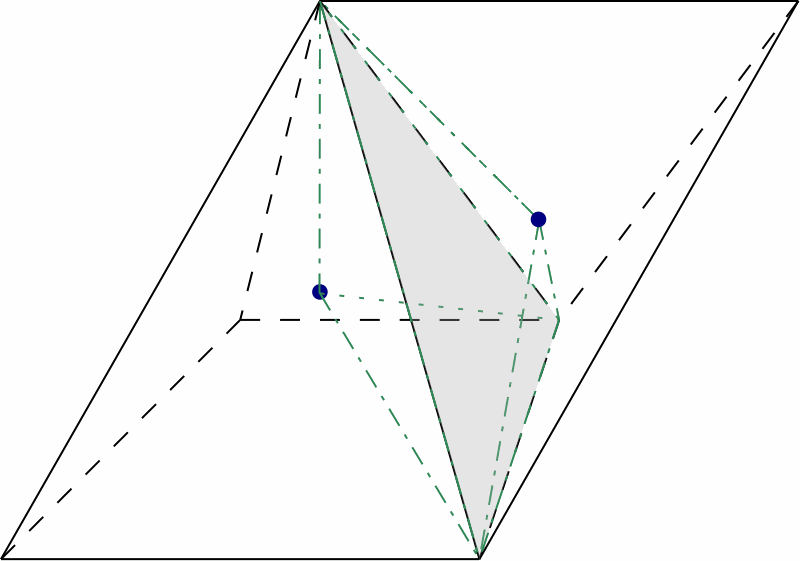} \\
	\vspace{0.05\linewidth}
	\includegraphics[width=0.25\linewidth]{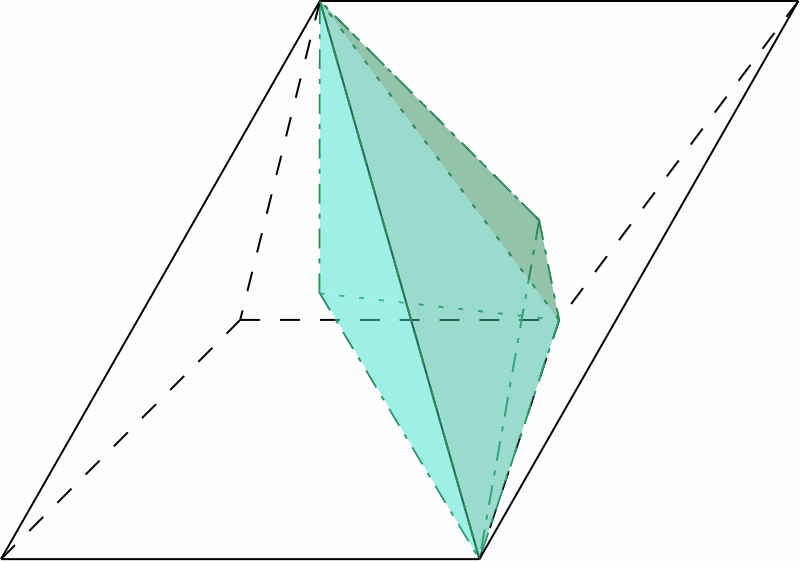}\hfill
	\includegraphics[width=0.25\linewidth]{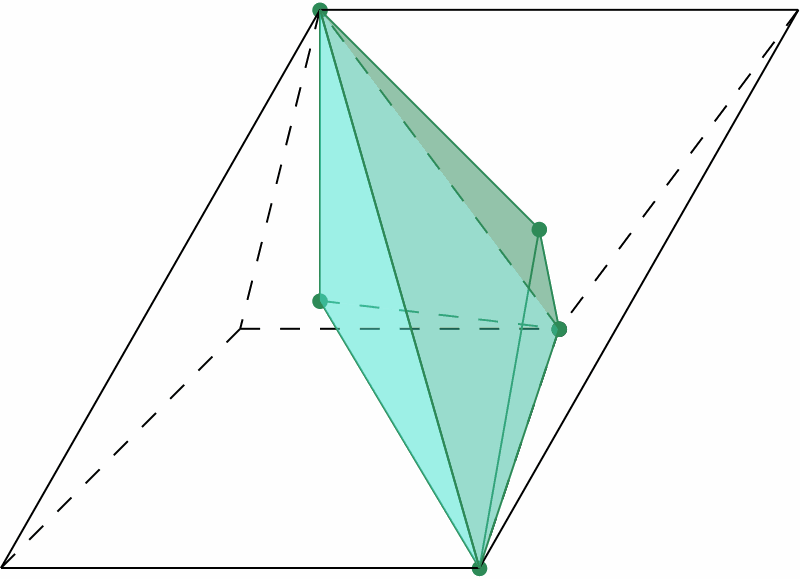}\hfill
	\includegraphics[width=0.25\linewidth]{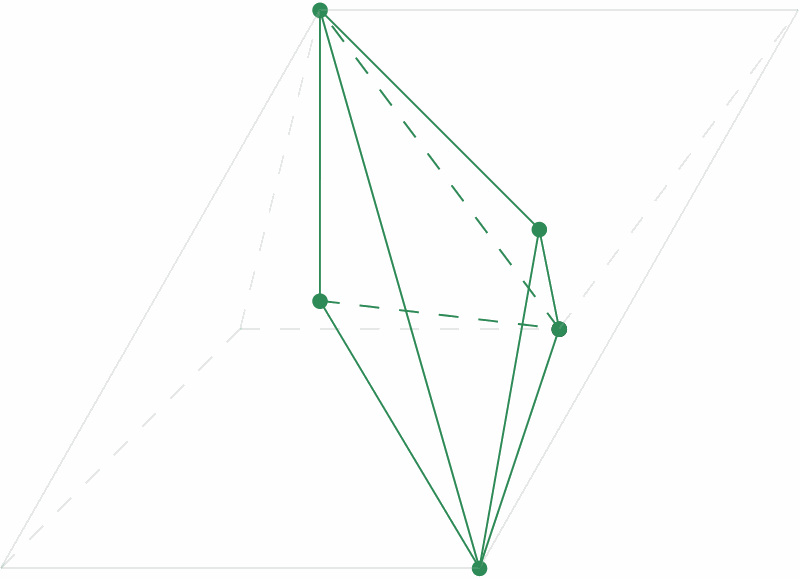}
	\caption{Construction of a dual element in 3D based on a face belonging to a square pyramid and a tetrahedron. Top left: primal square pyramid (left), primal tetrahedron (right), and primal vertex (black dots). Top centre: construction of the barycenters of the primal elements. Top right: each barycenter is connected to the vertex of the common face (shaded in grey), generating a tetrahedron on each side of the face with basis the grey face and opposite vertex the barycenters. Bottom left: the lateral faces of the generated tetrahedra are shadowed in light green (left, inside the hexahedron) and sea green (right, inside the primal pyramid). Bottom centre: the vertex of the two new tetrahedra are marked with green dots; continuous and discontinuous green lines indicate visible and shadow edges when merging both tetrahedra. Bottom right: the new dual element corresponds to the volume generated by merging the two tetrahedra constructed.}
	\label{fig:buildmesh3d_PiramidTetra}
\end{figure}

\subsection{Overview of the overall algorithm}
Before presenting all the details of the new method, we first give a brief overview of the overall algorithm. Taking into account the performed splitting of the equations, the dual mesh structure, and the combination of finite volume (FV) and finite element (FE) discretization, the proposed algorithm can be divided into four main stages:
\begin{itemize}
	
	\item Divergence-free reconstruction and transport-diffusion stage (Section~\ref{sec:transdiff}). System \eqref{eqn.convvissubs} is solved at the aid of an explicit FV method in the primal grid, obtaining an intermediate momentum, $\tww$, which in general, does not yet verify the divergence-free condition of the velocity field. The explicit transport-diffusion stage also includes the effects of the Maxwell stress tensor and of the viscous stresses in the momentum equation. Furthermore, an auxiliary cell-centered value of the magnetic field is computed, which is later corrected and overwritten. 
	To reach second order of accuracy, piecewise linear polynomials are reconstructed from the cell-centered velocity and from the auxiliary cell-centered magnetic field. In order to guarantee a locally and globally exactly divergence-free reconstruction for the magnetic field, the constrained $L^2$ projection algorithm of \cite{divfree2015} is employed, which subsequently overwrites the auxiliary magnetic field in the cell center.  
	  
	\item Divergence-free update of the magnetic field (Section \ref{sec:divfree}). In order to verify the divergence-free condition for the magnetic field exactly at the discrete level, the induction equation \eqref{eq:Faraday2} is solved inside the faces of the primal mesh, providing a natural and exactly divergence-free evolution of the normal component of the magnetic field via a discrete Stokes theorem inside each face, see \cite{BalsaraSpicer1999,divfree2015}. For the discrete Stokes law, the electric field is needed at the edges of each face and is provided by a multi-dimensional Riemann solver and a discrete curl term, including the physical resistivity. This allows computing the new magnetic field $\magfield^{n+1}$.  
	As an alternative to the exactly divergence-free evolution, also a simpler hyperbolic GLM divergence cleaning approach \cite{MunzCleaning,Dedneretal} can be applied, see Section~\ref{sec:divcleaning}. 
	
	\item Projection stage (Section~\ref{sec:projection}). The new pressure, $\press^{n+1}$, is computed at the vertices of the primal grid by solving the pressure Poisson equation \eqref{eqn.pressurePoisson} using classical continuous $\mathbb{P}_1$ Lagrange finite elements. To build the corresponding right-hand side term, the intermediate momentum $\tww$ is interpolated from the primal elements to the dual cells. 
	
	\item Post-projection stage (Section~\ref{sec:postpro}). The intermediate momentum $\tww$ is now corrected using the pressure gradient coming from the previous projection stage. As a consequence, the momentum at the new time step, $\ww^{n+1}$, is obtained. 
\end{itemize}
The workflow of the proposed well-balanced divergence-free methodology is depicted in Diagram \ref{dia:algorithm}.
\begin{Diagram}[h]
	\begin{equation*}
		\begin{tikzcd}
			\press^n \arrow[rrrrd]           &                  &           &                        &                                                                                                                         & \press^{n+1} \arrow[rd, "\textrm{Post-proj. stage}" ] &           \\
			\Q^n_i \arrow[rr, "\mathrm{Rec.}" description]       &                          & \Q^n(x) \arrow[rr, "\mathrm{FV}" description]  &  & \tQ_i \arrow[ru, "\textrm{Proj. stage}" ] \arrow[rr, "\magfield^{n+1}_i = \magfield^*_i" description] \arrow[rrd, "E = -\bVel \times \magfield + \eta\curl \magfield" description] &                                      & \Q^{n+1}_i \\
			\magfield^n_{\face} \arrow[rru, "L^2 \textrm{ proj.}" description] \arrow[rrrrrr, "\mathrm{Faraday}" description] &   &    &                   &                                                                                                                         &                                      & \magfield^{n+1}_{\face}
		\end{tikzcd}
	\end{equation*}
	\caption{Overall well-balanced and divergence-free FV/FE scheme for the magnetohydrodynamics equations.} \label{dia:algorithm}
\end{Diagram}

\subsection{Transport-diffusion stage on the primal mesh} \label{sec:transdiff}
The transport diffusion system  \eqref{eqn.convvissubs} is solved by employing explicit finite volumes in the primal grid. 
Accordingly, we integrate the corresponding PDEs on a control volume $T_i$ and apply Gauss theorem, obtaining:
\begin{align}
	\tWW{i} &= \WW^{n}_{i} 
	- \frac{\Delta t}{\VT{i}} \int\limits_{\partial T_i} \Fluxcv\left(\Q^{n}\right)\cdot \tnn_{i} \dS
	+ \frac{\Delta t}{\VT{i}} \int\limits_{\partial T_i} \Fluxcv\left(\Qe\right)\cdot \tnn_{i} \dS 
	+  \frac{\Delta t}{\VT{i}} \int\limits_{\partial T_i} \Fluxvv\left(\mu, \nabla \bVel^{n}\right) \cdot\tnn_{i} \dS \nonumber\\
	& \hspace{1cm} -  \frac{\Delta t}{\VT{i}} \int\limits_{\partial T_i} \Fluxvv\left(\mue, \nabla \bvele\right)\cdot \tnn_{i} \dS
	- \frac{\Delta t}{\VT{i}}\int\limits_{T_{i}} \gra \left( \press^{n} - \presse  \right) \dV, \label{eq:FV} \\
	\magfield^{*}_{i} &= \magfield^{n}_{i} 
	- \frac{\Delta t}{\VT{i}} \int\limits_{\partial T_i} \Fluxcb\left(\Q^{n}\right) \cdot\tnn_{i} \dS
	+ \frac{\Delta t}{\VT{i}} \int\limits_{\partial T_i} \Fluxcb\left(\Qe\right) \cdot\tnn_{i} \dS  
	+  \frac{\Delta t}{\VT{i}} \int\limits_{\partial T_i} \Fluxvb\left(\eta, \nabla \magfield^{n}\right)\cdot \tnn_{i} \dS\nonumber\\
	& \hspace{1cm} -  \frac{\Delta t}{\VT{i}} \int\limits_{\partial T_i} \Fluxvb\left(\etae, \nabla \magfielde\right)\cdot \tnn_{i} \dS,  \label{eq:FVB}
\end{align}
where $\tWW{i}$ denotes the discrete intermediate momentum and $\magfield^{*}_{i}$ is the auxiliary magnetic field in cell $T_i$, while $\WW^{n}_{i} $ and $\magfield^{n}_{i} $,   are the cell-averaged momentum and magnetic fields at time $t^{n}$,
\begin{equation*}
	\WW^{n}_{i} =  \frac{1}{\VT{i}}\int\limits_{T_i}\ww(\xx,t^{n})\dV, \quad
\magfield^{n}_{i} = \frac{1}{\VT{i}}\int\limits_{T_i}\magfield(\xx,t^{n})\dV.
\end{equation*}
To compute the convective flux terms contribution, the boundary integrals are divided into the sum of the fluxes over the primal element faces 
\begin{align*}
	\int\limits_{\partial T_i} \Fluxcv\left(\Q^{n}\right)\cdot \tnn_{i} \dS \approx \sum\limits_{T_j \in \mathcal{N}_{i}} \AF{ij} \Fluxcv^{\NF}\left(\bQ_{ij}^-,\bQ_{ij}^+,\tnn_{ij}\right),
	\\ 
	\int\limits_{\partial T_i} \Fluxcb\left(\Q^{n}\right)\cdot \tnn_{i} \dS \approx \sum\limits_{T_j \in \mathcal{N}_{i}} \AF{ij} \Fluxcb^{\NF}\left(\bQ_{ij}^-,\bQ_{ij}^+,\tnn_{ij}\right),
\end{align*}
with $\pbf{ij}=\Gamma_{ij}$ the common face of two neighbouring cells, $T_{i}$ and $T_{j}$, and a numerical flux function $\Flux^{\NF}$ that depends on the left and right boundary extrapolated values $\bQ_{ij}^-$ and $\bQ_{ij}^+$, respectively, and which can be chosen, for instance, as the classical Rusanov numerical flux function,
\begin{gather}
	\Fluxcv^{\RS}(\bQ_{ij}^-,\bQ_{ij}^+,\tnn_{ij}) = \halb \left( \bFcv(\bQ_{ij}^-) + \bFcv(\bQ_{ij}^+)\right)\tnn_{ij} - \alphacv \left(\bWW_{ij}^+ - \bWW_{ij}^- \right),  \label{eqn.RSflux}\\
	\Fluxcb^{\RS}(\bQ_{ij}^-,\bQ_{ij}^+,\tnn_{ij}) = \halb \left( \bFcv(\bQ_{ij}^-) + \bFcv(\bQ_{ij}^+)\right)\tnn_{ij} - \alphacv \left(\bmagfield_{ij}^+ - \bmagfield_{ij}^- \right),  \label{eqn.RSfluxb}
\end{gather}
or the Ducros flux function,
\begin{gather}
	\Fluxcv^{\Duc}(\bQ_{ij}^-,\bQ_{ij}^+,\tnn_{ij}) = 
	  \halb \left( \bWW_{ij}^- + \bWW_{ij}^+ \right) {\Vel}_{ij}
	- \halb \left( \bmagfield_{ij}^- + \bmagfield_{ij}^+ \right) \bmagfield_{ij}\cdot\tnn_{ij}
	+ \halb \bmagfield_{ij}^{2} \, \tnn_{ij}
	- \alphacv \left(\bWW_{ij}^+ - \bWW_{ij}^- \right), \label{eqn.Ducflux}\\
	\Fluxcb^{\Duc}(\bQ_{ij}^-,\bQ_{ij}^+,\tnn_{ij}) = 
	  \halb \left( \bmagfield_{ij}^- + \bmagfield_{ij}^+ \right) {\Vel}_{ij}
	- \halb \left( \bWW_{ij}^- + \bWW_{ij}^+ \right) \bmagfield_{ij}\cdot\tnn_{ij}
	- \alphacv \left(\bmagfield_{ij}^+ - \bmagfield_{ij}^- \right), \label{eqn.Ducfluxb}\\
	%
	 {\Vel}_{ij} = {\bvel}_{ij}\cdot\tnn_{ij}, \quad
	\bvel_{ij} = \halb \left( \bvel_{ij}^- + \bvel_{ij}^+ \right), \quad 
	\bmagfield_{ij} = \halb \left( \bmagfield_{ij}^- + \bmagfield_{ij}^+ \right) ,\quad
\end{gather}
where, instead of directly considering the conservative variables coming from the previous time step, we introduce their reconstructed values, $\bQ$, $\bWW$, $\bmagfield$, which come from a well-balanced reconstruction in space and time, a detailed description of which is provided in Section~\ref{sec:WB-recons}.

In the above expressions, the upwind coefficient, $\alphacv = \max\left\lbrace\left|\lambdacv_{i}\right|,\left|\lambdacv_{j}\right|\right\rbrace$, represents the maximum signal speed at the interface and depends on the maximum absolute value of the eigenvalues of the convective system for the left and right states to the face, $\left|\lambdacv_{i}\right|$, $\left|\lambdacv_{j}\right|$, respectively. Since we are considering the complete MHD system, those eigenvalues depend both on the velocity field and the magnetic field, being:
\begin{gather}
	\lambdacv_{1} = u-c_f,\quad \lambdacv_{2} =u-c_a,\quad 
	\lambdacv_{3,4,5} =u,\quad \lambdacv_{6} =2 u,\quad
	\lambdacv_{7} =u+c_a,\quad \lambdacv_{8} =u+c_f,
\end{gather}
where 
\begin{gather}
	u   = \left|\bvel\right|, \qquad 
	c_a  = \sqrt{\frac{1}{\rho}\left( \magfield\cdot\tnn \right)^2}, \qquad
	c_f  = \sqrt{\frac{1}{\rho} \magfield\cdot\magfield}.
	%
\end{gather}

The flux term contribution related to the equilibrium solution is then computed by directly evaluating the convective physical flux at each face barycenter as
\begin{equation}
	\int\limits_{\partial T_i} \Fluxcv\left(\Qe\right) \tnn_{i} \dS \approx \sum\limits_{T_j \in \mathcal{N}_{i}} \AF{ij} \bFcv(\Qe) \cdot \tnn_{ij},\qquad
	\int\limits_{\partial T_i} \Fluxcb\left(\Qe\right) \tnn_{i} \dS \approx \sum\limits_{T_j \in \mathcal{N}_{i}} \AF{ij} \bFcb(\Qe) \cdot \tnn_{ij},
\end{equation}
with 
\begin{equation}
	\Qe_{ij}=\left(\WWe_{ij},\magfielde_{ij}\right)^{T}, \qquad \rho^{e}_{ij} = \rho, \qquad \WWe_{ij} = \wwe\left(\xx_{ij}\right),\qquad \magfielde_{ij} = \magfielde\left(\xx_{ij}\right)
\end{equation}
and $\xx_{ij}$ are the coordinates of the barycenter of face $\face_{ij}$.

On the other hand, the viscous terms are approximated as 
\begin{gather}
	\int\limits_{\partial T_i} \Fluxvv\left(\mu, \nabla \bVel^{n}\right) \tnn_{i} \dS \approx \sum\limits_{T_j \in \mathcal{N}_{i}} \AF{ij} \Fluxvv^{\NF}(\mu,\nabla \bbVel_{ij}^-, \nabla \bbVel_{ij}^+, \bbVel_{ij}^-,\bbVel_{ij}^+,\tnn_{ij}), \\
	\int\limits_{\partial T_i} \Fluxvb\left(\eta, \nabla \magfield^{n}\right) \tnn_{i} \dS \approx \sum\limits_{T_j \in \mathcal{N}_{i}} \AF{ij} \Fluxvb^{\NF}(\eta,\nabla \bmagfield_{ij}^-, \nabla \bmagfield_{ij}^+, \bmagfield_{ij}^-,\bmagfield_{ij}^+,\tnn_{ij}),
\end{gather}
where the numerical flux functions for the viscous flux are given accordingly to \cite{Gas07} by 
\begin{eqnarray}
	\Fluxvv^{\NF}(\mu,\nabla \bbVel_{ij}^-, \nabla \bbVel_{ij}^+, \bbVel_{ij}^-,\bbVel_{ij}^+,\tnn_{ij})  &=& \halb \left( \bFvv(\mu, \nabla \bbVel_{ij}^-) + \bFvv(\mu, \nabla \bbVel_{ij}^+)\right)\tnn_{ij} - \alphavv \left(\bbVel_{ij}^+ - \bbVel_{ij}^- \right), \label{eqn.Viscflux}\\
	\Fluxvb^{\NF}(\eta,\nabla \bmagfield_{ij}^-, \nabla \bmagfield_{ij}^+, \bmagfield_{ij}^-,\bmagfield_{ij}^+,\tnn_{ij})  &=& \halb \left( \bFvb(\eta,\nabla \bmagfield_{ij}^-) + \bFvb(\eta,\nabla \bmagfield_{ij}^+)\right)\tnn_{ij} - \alphavv \left(\bmagfield_{ij}^+ - \bmagfield_{ij}^- \right),
\end{eqnarray}
with the penalty coefficient $\alphavv = \frac{2}{h_{i}+h_{j}} \max\left\lbrace\left|\lambdavv_{ i}\right|,\left|\lambdavv_{j}\right|\right\rbrace$ and the upper bound of the viscous eigenvalues $\lambdavv = \dfrac{\mu}{\rho} + \etares$. Once again, to compute the numerical fluxes, the reconstructed values are employed. 

In case the well-balanced scheme or a pressure correction approach is considered, the gradient of the pressure at the previous time step must be accounted for in the explicit stage. To approximate it, the integral of the gradient on the control volume is also transformed into the sum of integrals on the element faces. Consequently, we have
\begin{equation}
	\int\limits_{T_{i}} \gra \left( \press^{n} - \presse  \right) \dV =  \int\limits_{\partial T_i} \left( \press^{n} - \presse  \right) \tnn_{i}\dS
	= \sum\limits_{T_j \in \mathcal{N}_{i}} \AF{ij} \left(  \press^{n}_{ij}- \presse_{ij}  \right) \tnn_{ij}
\end{equation}
where the pressure at the cell faces, $\press^{n}_{ij}$, is computed as the average pressure defined at the face vertices:
\begin{equation}
	 \press^{n}_{ij}  = \frac{1}{|\mathcal{V}_{ij}|} \sum\limits_{\vertex_k \in \mathcal{V}_{ij}}^{} \press^{n}_{k},
\end{equation}
with $\mathcal{V}_{ij}$ the set of vertices of face $\face_{ij}$ and $|\mathcal{V}_{ij}|$ its cardinality.

As described above, the scheme would only be of first order in space and time. To improve the accuracy of the method, we compute an extrapolation of the data at the neighbouring of the faces and perform a half in time evolution as usual in the MUSCL-Hancock method, \cite{VL97,Toro}, and the ADER approach, \cite{Toro,TMN01}. Moreover, special attention is paid to have also a well-balanced reconstruction and to guarantee the divergence-free condition of the magnetic field.

\paragraph{Well-balanced reconstruction}\label{sec:WB-recons}
For attaining a second order scheme, we consider the half in time evolved extrapolated values of the left and right states at a face $\face_{ij}$, $\bQ_{ij}^-$ and $\bQ_{ij}^+$, respectively, when approximating the fluxes, \eqref{eqn.RSflux}, \eqref{eqn.Ducflux}, \eqref{eqn.Viscflux}, instead of the cell-averaged states coming from the previous time step, $\Q^{n}_{i}$, $\Q^{n}_{j}$. 
Accordingly, at each cell face, $\face_{ij}$, we define the two piecewise linear space-time reconstruction polynomials inside the elements $T_i$ and $T_j$ as 
\begin{gather}
 	\bQ_{i}(\xx,t) = \Qe_{ij} + \Q'_i + \gra \Q'_{i} \cdot \left(\xx-\xx_{i}\right) +\partial_t \Q'_{i} 
 	\left( t - t^n \right), \qquad \forall \xx \in T_i, \label{eqn.Qipoly} \\
 	\bQ_{j}(\xx,t) = \Qe_{ij} + \Q'_j + \gra \Q'_{j} \cdot \left(\xx-\xx_{j}\right) + \partial_t \Q'_{j} \left( t - t^n \right), \qquad \forall \xx \in T_j, 
\end{gather} 
with $\Q'_i=\Q_{i}^{n}-\Qe_{i}$,  $\Qe_{ij}$ the equilibrium solution at the face barycenter $\xx_{ij}$, and $\xx_i$ the barycenter of the primal element $T_i$. Moreover, we denote $\gra \Q'_{i}$ the spatial gradient of the fluctuation $\Q'$ in $T_i$ and  $\gra \Qe_{ij}$ the gradient of the equilibrium solution at the face barycenter.
In case $\gra \Qe_{ij}$ is not analytically known, it is approximated using finite differences. 
The boundary-reconstructed values are then given by 
\begin{equation}
	\bQ_{ij}^- = \bQ_i \left( \xx_{ij},t^{n} + \halb \Delta t \right), \qquad \textnormal{and} \qquad 
	\bQ_{ij}^+ = \bQ_j \left( \xx_{ij},t^{n} + \halb \Delta t \right).
\end{equation} 
The slopes at cell $T_{i}$, $\gra \Q'_{i}$, are approximated using a least square approach.
First, we compute the difference between the states at cell $T_{i}$ and its neighbours, taking into account the equilibrium solution and the distance between the neighbouring cell barycenters as
\begin{equation}
	\Delta\Q'_{ij} =\frac{1}{\left| T_j \right|}\int\limits_{T_j} \nabla \Q'_i\cdot(\xx-\xx_i)\dx = \nabla\Q'_{i}\cdot\Delta\xx_{ij}, \qquad \forall\; T_j\in \mathcal{S}_{i},
\end{equation}
with
\begin{equation}
		\Delta\Q'_{ij} := \left( \Q_j^{n}-\Qe_j\right) - \left( \Q_i^{n} - \Qe_i\right) ,\qquad \Delta\xx_{ij} := \xx_j-\xx_i
\end{equation} 
and $ \mathcal{S}_{i}$ the stencil for the reconstructed solution, that, due to the use of a least squares approach, may involve more cells than the minimum necessary number to reconstruct a polynomial, 
see \cite{BarthFrederickson1990,Abgrall1994,DumbserKaeser07} for further details. 
Besides, the integral has been approximated via the mid-point rule, which is enough for up to second order accurate integration. Then, applying least squares, we get the associated system of normal equations
\begin{equation}
	\mathcal{A} \nabla \Q'_{i}  = \mathcal{B},\qquad   \mathcal{A} = \left(\Delta\xx_{{ij}} \right) \cdot \left(\Delta\xx_{{ij}}\right)^{T}  , \qquad  \mathcal{B} = \left(\Delta\xx_{{ij}} \right) \cdot \left(\Delta\Q'_{{ij}}\right)^{T}
	\label{eqn.lsq} 
\end{equation}
which provides an optimal approximation of the gradient at the cell minimizing the errors committed with respect to the local-face approximated gradients. In general, the spatial reconstruction polynomial for the magnetic field 
$\magfield_i(\xx,t^n)$ obtained via the reconstruction detailed in \eqref{eqn.Qipoly}-\eqref{eqn.lsq} is not yet divergence-free, hence we denote it for the moment by the preliminary symbol $\tilde{\magfield}_i(\xx,t^n)$. Later, in Section \ref{sec:L2projection} we show how to make the reconstruction polynomial of the magnetic field truly divergence-free, which will then be denoted by $\magfield_i(\xx,t^n)$.  

Finally, following the classical MUSCL-Hancock and ADER approach \cite{Toro,titarevtoro,BTVC16}, and taking into account the equilibrium solution, the contribution of the time derivative, $\partial_t \Q'_{i}$,  is obtained by applying the Cauchy-Kovalevskaya procedure. Accordingly, it is replaced by spatial derivatives taking into account the momentum equation \eqref{eqn.convvissubsv0}, as
\begin{equation}
	\partial_t \WW'_{i} = -\frac{1}{\VT{i}}\sum\limits_{T_j \in \mathcal{N}_{i}} \AF{ij}\left(\bFcv\left(\Q_{ij} \right) +\bFvv\left(\mu, \nabla \Q_{ij} \right) -\bFcv\left(\Q^{e}_{ij} \right)-\bFvv\left(\mue, \nabla \Qe_{ij}\right) + \press_{ij}\right) \cdot \tnn_{ij} ,
\end{equation}
and the conservative form of the Faraday equation \eqref{eqn.convvissubsB0},
\begin{equation}
	\partial_t \magfield'_{i} = -\frac{1}{\VT{i}}\sum\limits_{T_j \in \mathcal{N}_{i}} \AF{ij}\left(\bFcb\left(\Q_{ij} \right) +\bFvb\left(\eta, \nabla \Q_{ij} \right) -\bFcb\left(\Q^{e}_{ij} \right)-\bFvb\left(\etae, \nabla \Qe_{ij} \right)\right) \cdot \tnn_{ij},
\end{equation}
with the boundary-reconstructed data from \textit{within} element $T_i$ at time $t^n$ given by 
\begin{equation}
	\Q_{ij} = \Qe_{ij} + \Q'_i + \gra\Q'_{i} \cdot\left( \xx_{ij} - \xx_{i}\right).
\end{equation}

To introduce the non-linearity needed to avoid spurious oscillations related to the second order approach, a limiting strategy should be employed. For instance, the classical Barth and Jespersen limiter, \cite{BJ89}, or ENO and WENO methods, \cite{Abgrall99ENO,HuShuTri,DumbserKaeser07,cweno2017} could be selected. 
Furthermore, as an alternative to the use of classical ADER schemes, the local space-time predictor ADER approach could have been used, thus avoiding the Cauchy-Kovalevskaya procedure, \cite{DET08,DBTM08}.

\subsection{Divergence-free evolution of the magnetic field and constrained $L^2$ projection}  \label{sec:divfree}

In this section, we propose a novel well-balanced algorithm for the divergence-free evolution of the Faraday equation at second order of accuracy. For every cell $T_i$, we consider a  magnetic field $\bBvar$ in the space  $\mathcal{P}_r(T_i)^d$, which is the space of $d$-vector valued polynomials of degree $r$ on $T_i$. 
It is said to be divergence-free if, for each element of the domain $T_i\in\ddomain$, we have the pointwise identity 
\begin{equation}\label{eq:divfree}
\dive \bBvar = 0, \qquad \forall \mathbf{x} \in T_i.
\end{equation}
 Then, considering the full domain $\ddomain$, the discrete solution, which consists of piecewise polynomials, is  represented in a non-conforming space of solutions.
Inspired by the work of \cite{divfree2015}, we introduce a divergence-free reconstruction operator that is conforming in the sense of the face-averaged magnetic fluxes, i.e., prescribing a discrete scalar field for the face-averaged normal component $\Bnvar_f$ on face $\face_f$, which is defined as  
\begin{equation}\label{eq:Bface}
\Bnvar_f =  \frac{1}{\length{\face_f}} \int \limits_{\face_f} \bBvar \cdot \tnn_f \, dS, \quad \forall \face_f.
\end{equation}

For this section we need some auxiliary notation, that we now introduce. Each face $\face_f$ has an intrinsic orientation induced by the face normal vector $\tnn_f$ pointing from the left element $T_{\ell}$ attached to the face to the right element $T_r$ associated with face $\face_f = \face_{\ell r}$. Moreover, if $T_i$ is an element and $\face_{f} = \face_{ij}$ is one of its faces, we define the relative orientation $\sigma_{if}\in \{-1,1\}$ of $T_i$ and $\face_f$ via the formula $\tnn_{ij} = \sigma_{ij}\tnn_{f}$, or, equivalently, $\sigma_{ij} = \tnn_{ij}\cdot \tnn_{f}$. Similarly each edge $\edge_e$ has an intrinsic orientation induced by the direction of its tangent vector $\tang_e$. If the edge $\edge_e$ belongs to the boundary of the face $\face_f$, then we define their relative orientation $\sigma_{fe}$ as equal to $1$ if $\edge_e$ is oriented counterclockwise with respect to $\face_e$, and $-1$ if it is oriented clockwise. Given an edge $\edge_e$, we define $\mathcal{T}_e$ (resp.  $\mathfrak{F}_e$) as the set of elements (resp. faces) that contain $e$ and we denote by $\lvert \mathcal{T}_e\rvert$ (resp. $\lvert \mathfrak{F}_e \rvert$) its cardinality. Moreover, if $T_i$ is an element containing $e$, we denote by $\mathfrak{F}_{ei}$ the subset of $\mathfrak{F}_e$ supported in $T_i$. Note that $\mathfrak{F}_{ei}$ contains exactly two faces, say $\face_{f^+_e}$ and $\face_{f^-_e}$, such that 
\begin{equation*}
	\sigma_{if^+_e}\sigma_{f^+_ee} + \sigma_{if^-_e}\sigma_{f^-_ee} = 0.
\end{equation*}
Finally, let $\bary_e$ denote the coordinates of the barycenter of the edge $\edge_e$.

\subsubsection{Divergence-free evolution of the normal components of the magnetic field}



First, we derive a discrete equation for the (scalar) averaged magnetic fluxes $\Bnvar_f$ in each face $\face_f$.
Given any face $\face_f$ of the primary mesh, the Faraday equation reveals the time-evolution of $\Bnvar_f$ as
\begin{equation*}
\frac{d}{dt} \int \limits_{\face_f} \bBvar \cdot \tnn_f \, dS + \int \limits_{\face_f} \curl \bEvar \cdot \tnn_f \, dS = 0.
\end{equation*} 
Then, after the application of Stokes theorem and time integration, one may derive a consistent discrete equation for the face-averaged magnetic fluxes in the form of 
\begin{equation}
	\Bnvar_{f}^{n+1} = \Bnvar_{f}^{n} - \frac{\Delta t}{\length{\face_f}}\sum_{\edge_e \in \face_f} \sigma_{fe} \,  \length{\edge_e} \, \Etvar_{e},
	\label{eqn.newBface}
\end{equation}
where the sum is taken over the edges $\edge_e$ which delimit the boundary of the face $\face_f$. 

%

Here, for any edge $\edge_e$ in the mesh, we define the averaged  tangential component of the electric field along $\edge_e$ as  
\begin{equation}\label{eq:edge}
 \Etvar_e=  \frac{1}{\length{\edge_e}} \int \limits_{\edge_e} \bEvar \cdot \tang_e \, d\lambda.
\end{equation}
 In the next section, we describe how to compute the edge-averaged electric field $\Etvar_e$ from the given velocity and magnetic fields. 
We now show that the averages obtained with equation \eqref{eqn.newBface} are divergence-free in the sense
\begin{equation}
	\sum_{\face_f\in \element_i}\sigma_{if}\lvert \face_f\rvert \Bnvar_{f}^{n+1} = 0, \quad \forall \element_i. 
	\label{eq:av.div.free}
\end{equation} 
This follows from a standard induction argument provided that the initial condition satisfies
\begin{equation}
	\sum_{\face_f \in \element_i}\sigma_{if}\lvert \face_f\rvert \Bnvar_{f}^{0} = 0, \quad \forall T_i. 
	\label{eq:av.div.free.initialcondition}
\end{equation} 
and from the straightforward computation
\begin{align*}
	\sum_{\face_f\in \element_i}\sigma_{if}\sum_{\edge_e \in \face_f}\sigma_{fe}\lvert \edge_e\rvert \Etvar_e &= \sum_{\edge_e \in \element_i}\sum_{\face_f\in\mathfrak{F}_{ei}}\sigma_{fe}\lvert \edge_e\rvert \Etvar_e\\
	& =	\sum_{\edge_e \in \element_i}(\sigma_{if}\sigma_{f^+_ee}+\sigma_{if}\sigma_{f^-_ee})\lvert \edge_e\rvert \Etvar_e = 0.
\end{align*}
To ensure that the initial face averages of the magnetic field are divergence-free in the sense \eqref{eq:av.div.free.initialcondition}, we initialize them as 
\begin{equation*}
	\Bnvar^0_f  = \frac{1}{\length{\face_f}}\sum_{\edge_e \in \face_f} \sigma_{fe} \,  \length{\edge_e} \, \Atvar_{e},
\end{equation*}
where 
\begin{equation*}
		\Atvar_e=  \frac{1}{\length{\edge_e}} \int \limits_{\edge_e} \bAvar(x,0)\cdot \tang_e \, d\lambda,
\end{equation*}
is the averaged tangential component of the magnetic potential $\bAvar$.

\subsubsection{Computation of the electric field} 
\begin{figure}
		\centering
		\includegraphics[width=0.3\linewidth]{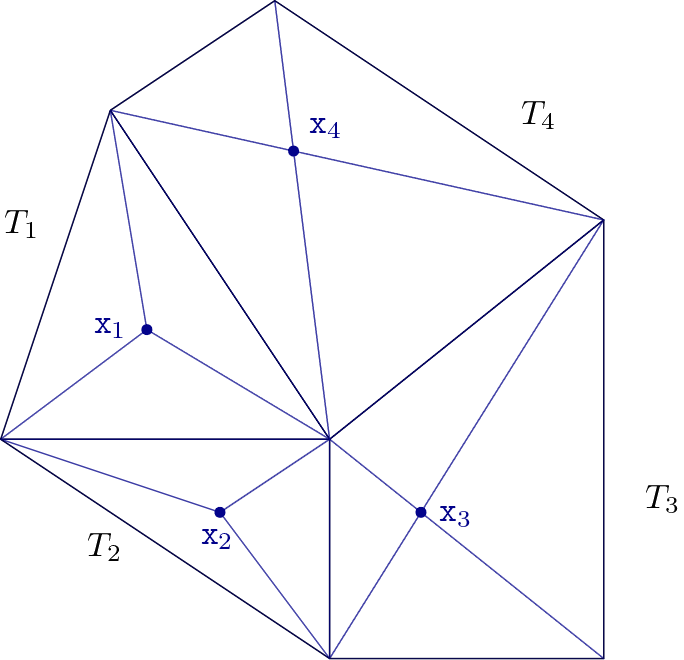}\hfill
		\includegraphics[width=0.3\linewidth]{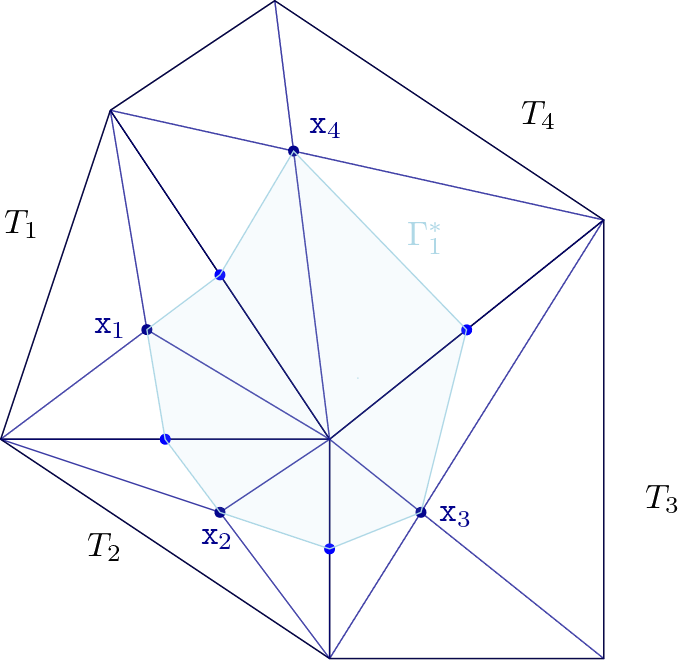}\hfill
		\includegraphics[width=0.3\linewidth]{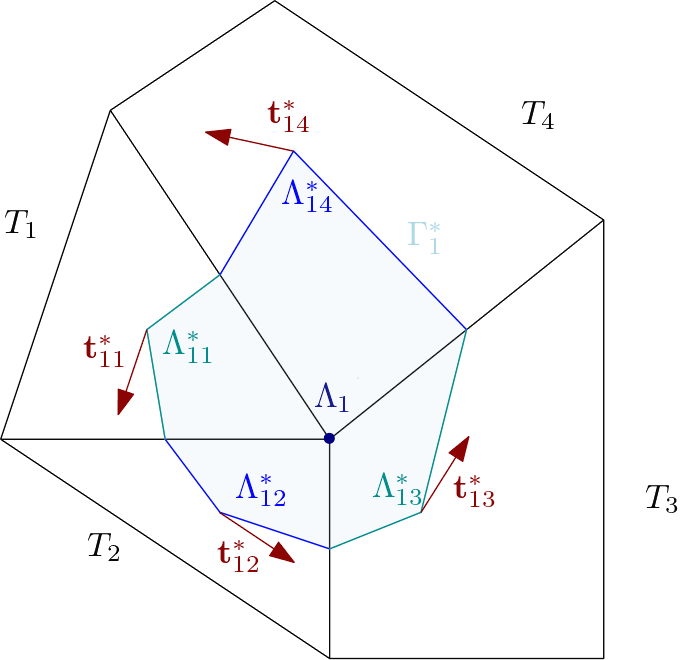}
		\caption{Dual vertex-based mesh structure employed for the computation of the electric field in 2D. Left: the primal elements are the black triangles and quadrilaterals $T_1$, $T_2$, $T_3$, and $T_4$  with $\bary_1$, $\bary_2$, $\bary_3$, and $\bary_4$ their corresponding barycenters (dark blue). Center: the blue points denote the barycenters of the primal faces related to $\edge_1$ and are connected with the barycenters of the primal elements to construct $\faceD_1$ (light blue). Right: $\edgeD_{11}$, $\edgeD_{12}$, $\edgeD_{13}$, $\edgeD_{14}$ identify the four edges, each one defined by gathering the two linear pieces inside each primal element of the boundary of $\faceD_1$, and the corresponding unitary corner tangent vectors $\tang_{11}^{\ast}$, $\tang_{12}^{\ast}$, $\tang_{13}^{\ast}$, $\tang_{14}^{\ast}$.}
		\label{fig.2dmeshvertexdual}
	\end{figure}

First, we recall how to construct the barycentric dual face $\face^*_e$ with respect to the edge $\Lambda_e$. For $\Gamma_f\in\mathfrak{F}_e$, denote by $\Gamma^*_{lfe}$ (resp. $\Gamma^*_{rfe}$) the triangle with vertices the barycenters of $T_l$ (resp. $T_r$), $\face_f$ and $\edge_e$. Then we set
\begin{equation*}
	\face^*_e = \bigcup_{f\in \mathfrak{F}_e}\face^*_{rfe}\cup \face^*_{lfe},
\end{equation*}
see Figure~\ref{fig.2dmeshvertexdual}. We denote by $\tang^*_e$ the unit tangent vector to the boundary of $\face^*_e$. 
From the definition of the electric field \eqref{eq:E}, we  evaluate the average of the tangent component along $\edge_e$ as 
\begin{eqnarray}
	\label{eqn.efield} 
	\Etvar_e & = &  -  \frac{1}{ |{\mesh_{e}}|  }\sum_{T_i \in \mesh_{e}}\frac{1}{\length{\edge_e}} 
	\int \limits_{\edge_e}\bVel_{i}\left( \xx,t^{n+\halb} \right) \times \magfield_{i}\left( \xx,t^{n+\halb} \right)\cdot\tang_{e}\,d\lambda 
	\nonumber \\ 
	&& + \frac{\eta}{\length{\faceD_e}}\sum_{T_i \in \mesh_{e}}\int \limits_{\edgeD_{ei}}\magfield_{i}\left( \xx,t^{n+\halb} \right) \cdot \tang^*_e \,d\lambda
 	+ \frac{\eta_h}{\length{\faceD_e}}\sum_{T_i \in \mesh_{e}} \length{\edgeD_{ei}} \, \magfield_i\left( \bary_{e},t^{n+\halb} \right) \cdot \tang_{ei}^*,
\end{eqnarray}
with the corner tangent vector $\tang_{ei}^*$ defined as 
\begin{equation}
  	\length{\edgeD_{ei}} \tang_{ei}^* = \int \limits_{\edgeD_{ei}} \tang^* \,d\lambda, 
  	\qquad  
  	\length{\edgeD_{ei}} = \int \limits_{\edgeD_{ei}}  \,d\lambda. 
\end{equation}
	The electric field $\Etvar_e$ is the sum of three components. The first is the averaged electric field for ideal MHD, where the average accounts for all elements $T_i$ attached to the edge. The second term is a discretization of the physical resistivity term, which is proportional to the curl of $\magfield$. It is computed as a line integral on the edges of the dual face $\faceD_e$, which is delimited by piecewise linear edges $\edgeD_{ei}$ that connect the barycenters of all ${T_i \in \mesh_{e}}$ around $\edge_e$. Each edge $\edgeD_{ei}$ of the dual face $\faceD_e$ consists of two pieces inside element $T_i$, i.e. the dual edge $\edgeD_{ei}$ itself is piecewise linear inside $T_i$: one piece starting in the barycenter $\xx_i$ of $T_i$, pointing towards the barycenter of the face in common between $T_i$ and the previous element $T_i^-$ in the ordered set $\mesh_{e}$. Likewise, the second piece of $\edgeD_{ei}$ inside $T_i$ points from $\xx_i$ to the barycenter of the common face between $T_i$ and the next element $T_i^+$ in the ordered set $\mesh_{e}$. We assume that all elements in $\mesh_{e}$ are ordered counter-clockwise around the edge $\edge_e$, i.e. according to the right hand rule with respect to the vector $\tang_e$. For a sketch of the notation related to this edges in the 2D framework see Figure~\ref{fig.2dmeshvertexdual}. The third term in \eqref{eqn.efield} is a stabilization term, written under the form of an artificial resistivity. It is inspired by the expressions obtained from approximate multi-dimensional Riemann solvers for MHD, see e.g. \cite{balsarahlle2d,balsarahllc2d,balsarahlle3d,BalsaraMultiDRS,MUSIC1,MUSIC2}. The numerical resistivity $\eta_h$ is similar to the artificial viscosity of a standard Rusanov (or local Lax-Friedrich) numerical flux, which provides an upwind-type stabilization. %
	Here $\eta_h$ is an \emph{artificial resistivity} defined as 
$	\eta_h = s_{\max} \, {\length{\faceD_e}} / {\length{\partial \faceD_e}}$,
where $s_{\max} = \max \left\lbrace\left|\lambdacv_{i}\right| \right\rbrace$ is the maximum wave speed estimated over all elements $T_i \in \mesh_{e}$. Note that, for uniform meshes, $\eta_h$ scales with the mesh size.
	
\subsubsection{Constrained $L^2$ projection}\label{sec:L2projection}

For every cell $T_i$, we define a functional as the $L^2$-distance between the final divergence-free polynomial $\bBvar_i=\bBvar_i\left(\xx,t^{n+1}\right)$ and the polynomial reconstruction $\tbBvar_i=\tbBvar_i\left(\xx,t^{n+1}\right)$, described previously in Section \ref{sec:transdiff} and based on the auxiliary magnetic field $\magfield_i^{*}$, which in general is not yet divergence-free, as follows: 
\begin{equation}
\myfunc \left( \bBvar_i \right) :=  \int \limits_{T_i} \, \left\lVert \bBvar_i - \tbBvar_i \right\rVert^2 dV.
\label{eqn.func} 
\end{equation}
We then state the following local optimization problem: look for a polynomial $\magfield_i\in \mathcal{P}_r(T_i)^d$ that minimizes  $\myfunc \left( \bBvar_i \right)$ subject to the following 
constraints:
\begin{align}
	\dive \magfield_{i} &= 0,  \qquad  \quad  \qquad \qquad \forall \xx \in T_i,  \label{eq:div-free} \\
    \int \limits_{\face_{ij}} \magfield_{i}\cdot \tnn_{ij} \, dS &= \sigma_{if} \length{\face_{ij}} \, \magfieldnc_{f}^{n+1},   \qquad \forall \face_{ij} \in \partial T_i^{\subset}, 
 \label{eq:nc}
\end{align}
with $\partial T_i^{\subset} = \partial T_i \setminus \widehat{\face}_{ij}$ the set of faces of $T_i$ from which one arbitrary face $\widehat{\face}_{ij}$ has been removed. Note that for the entire boundary $\partial T_i$ of $T_i$ the constraints are not linearly independent, as we now show. Fix an arbitrary face $\widehat{\face}_{ij}$ in $\partial T_i$, then  
\begin{eqnarray}
	 \int \limits_{\widehat{\face}_{ij}} \magfield_{i} \cdot \tnn \, dS &=& \sum_{\face_{ij} \in \partial T_i} \, \int \limits_{\face_{ij}} \magfield_{i}\cdot \tnn \, dS - \hspace{-2mm} \sum_{\face_{ij} \in \partial T_i^{\subset}} \, \int \limits_{\face_{ij}}\magfield_{i}\cdot \tnn \, dS \nonumber \\ 
	& = & \int \limits_{T_i} \dive \magfield_{i} \, dV -  \sum_{\face_{ij} \in \partial T_i^{\subset}} \, \int \limits_{\face_{ij}}\magfield_{i}\cdot \tnn \, dS \nonumber
	= 0 \quad - \sum_{\face_{ij} \in \partial T_i^{\subset}} \, \int \limits_{\face_{ij}}\magfield_{i}\cdot \tnn \, dS.  \nonumber
\end{eqnarray}
Note that the second equality follows from the Gauss theorem. Therefore, we impose the constraint \eqref{eq:nc} only on  the subset $\partial T_i^{\subset} = \partial T_i \setminus \widehat{\face}_{ij}$ and not on all faces of $T_i$, see also \cite{divfree2015}. 
The optimization problem is then easily solved by the following Lagrangian-multiplier method: 
\begin{align}\label{eq:Lagrange}\left\{
\begin{array}{l} 
\Lagrange \left( \bBvar_{i}, \lambda, \lambda_{ij} \right) := \myfunc \left( \bBvar_{i} \right)  + \lambda \dive \magfield_{i} + {\displaystyle \sum \limits_{\face_{ij} \in \partial T_i^{\subset} }} \!\!\! \lambda_{ij}  \left(
{\displaystyle \int\limits_{\face_{ij}}} \magfield_{i}\cdot \tnn_{ij} \, dS - \sigma_{if} \length{\face_{ij}} \, \magfieldnc_{f}^{n+1}   \right), \\
\cfrac{\partial \Lagrange}{\partial \bBvar_i} = 0,  \qquad 
\cfrac{\partial \Lagrange}{\partial \lambda} = 0, \qquad 
\cfrac{\partial \Lagrange}{\partial \lambda_{ij}} = 0, \quad \forall \face_{ij} \in \partial T_i^{\subset}.
\end{array}\right.
\end{align}

The magnetic field $\bBvar_i(\xx,t^{n+1})$ is assumed to be piecewise linear (i.e. $r=1$) and can be written with respect to its barycentric value and the corresponding slopes, i.e. for every element $T_i \in \ddomain$ we may write
\begin{equation} \label{eq:Brec}
\magfield_{i} = \magfield_{i}\left(\xx,t^{n+1}\right) = \magfield_{i}^{n+1} + \nabla \magfield_{i}^{n+1} \left( \xx - \xx_i \right).
\end{equation}
%
Note that, for up to second order accurate methods, the barycentric value $\magfield_{i}^{n+1}$ coincides with the corresponding cell-average.
In general $\magfield_i$ can we written as
\begin{equation} \label{eq:Brec2}
\begin{split}
\magfield_{i} = \magfield_{i}\left(\xx,t^{n+1}\right) = \bBvar^{n+1}_{i,\ell} \, \psi_\ell(\mathbf{x}) 
\end{split}
\end{equation}
where $\psi_\ell(\mathbf{x})$ form a  basis for $ \mathcal{P}_r(T)^d$ and $\bBvar^{n+1}_{i,\ell}$ is the corresponding array of degrees of freedom, which includes the barycentric value $\magfield_{i}^{n+1}$ as well as the slope $\nabla \magfield_{i}^{n+1}$. In particular, we have chosen a \emph{modal} basis, hence the degrees of freedom are the first moments of $\bBvar$, i.e. the cell average and the local spatial-slopes. 
Here, $\bx$ is the vector position in Cartesian coordinates in the physical space.
Hence, in three space dimensions, the approximation of the final divergence-free magnetic field has 12 \emph{free-parameters}, namely the three cell-averages and the nine slopes. 
Using \eqref{eq:Brec2}, the extremum conditions may be rewritten as 
 \begin{equation}\label{eq:Lagrange2}
\cfrac{\partial\Lagrange}{\partial \bBvar^{n+1}_{i,\ell}} = 0, \qquad 
\cfrac{\partial \Lagrange}{\partial \lambda} = 0, \qquad 
\cfrac{\partial \Lagrange}{\partial \lambda_{ij}} = 0 , \quad \quad 
\forall \face_{ij} \in \partial T_i^{\subset}.  
\end{equation}
 Since the selected functional $\myfunc$ is convex and we are imposing only linear constraints, the resulting matrix is symmetric and invertible and the unique solution of the linear system will be a minimizer. Moreover, since such matrix depends only on the mesh, its local inverse can be pre-computed, stored and directly used in the simulation, see \cite{divfree2015} for further details. At the end, the auxiliary cell-average of the magnetic field
 $\magfield_i^*$ is overwritten by the result of the new cell-average coming from the constrained $L^2$ projection. Then the scheme can proceed with the next time-step. 
\subsubsection{Well-balancing}
The previously-described divergence-free algorithm is not yet well-balanced. To obtain a well-balanced divergence-free scheme, we must subtract the discretization of the stationary equilibrium equation. Indeed, on the continuous level one has 
\begin{align}  
	\curl \bEvar( \eta^e, \magfielde, \bVele) = 0,   \qquad \textnormal{since} \qquad \partial_t  \magfielde = 0. 
	\label{eqn.stationary.Faraday}
\end{align}
The well-balancing procedure proposed in this paper is extremely simple, and it may drastically reduce the diffusive and dispersive numerical errors around the chosen analytical or (numerically) approximated equilibrium solution, as shown later in our numerical results. As a consequence, long-time stability is obtained with respect to any analytical or numerical solution, for any prescribed structured or unstructured grid, and almost independently on the order of accuracy of the method. It is important to notice that this is not a linearization, since the PDE are solved in their original non-linear form, from which only a discretization of the stationary version of the PDE \eqref{eqn.stationary.Faraday} has been subtracted. 

The final well-balanced discretization of \eqref{eq:Faraday2} reads 
\begin{equation}
	\Bnvar_{f}^{n+1} = \Bnvar_{f}^{n} - \frac{\Delta t}{\length{\face_f}}\sum_{\edge_e \in \face_f} \sigma_{fe} \,  \length{\edge_e} \, \left( \Etvar_{e} - \Etvar_{e}^e \right),
	\label{eqn.newBface.wb} 
\end{equation}
with the electric field in the equilibrium defined as
\begin{eqnarray}
	\label{eqn.efield.wb} 
	\Etvar_e^e & = &  -  \frac{1}{ |{\mesh_{e}}|  }\sum_{T_i \in \mesh_{e}}\frac{1}{\length{\edge_e}} 
	\int \limits_{\edge_e} \bVel^e\left( \xx \right) \times \magfield^e\left( \xx \right)\cdot\tang_{e}\,d\lambda 
	\nonumber \\ 
	&& + \frac{\eta}{\length{\faceD_e}}\sum_{T_i \in \mesh_{e}}\int \limits_{\edgeD_{ei}}\magfield^e \left( \xx \right) \cdot \tang^* \,d\lambda
	+ \frac{\eta_h}{\length{\faceD_e}}\sum_{T_i \in \mesh_{e}} \length{\edgeD_{ei}} \, \magfield^e \left( \bary_{e} \right) \cdot \tang_{ei}^*.
\end{eqnarray}
It is obvious that when the discrete solution coincides with the \textit{a priori} known equilibrium, we have $\Etvar_e = \Etvar_e^e$ and therefore the numerical method produces $\Bnvar_{f}^{n+1} = \Bnvar_{f}^{n}$, hence the method is well-balanced. To the best of our knowledge, this is the first well-balanced exactly divergence-free finite volume scheme for the discretization of the MHD equations on general mixed-element unstructured meshes. 



\subsection{GLM divergence cleaning}\label{sec:divcleaning}
A hyperbolic Generalized Lagrangian Multiplier (GLM) divergence cleaning technique can be used as a simpler alternative to the former exactly divergence-free scheme. Although GLM does not lead to an exactly divergence-free discretization, it avoids the accumulation of divergence errors in the magnetic field and also allows long-time stable simulations. Due to its simplicity, we present it here, although the main focus of this paper is on the exactly divergence-free scheme introduced above.   

Following the ideas of Munz \textit{et al.} \cite{MunzCleaning,Dedneretal}, we construct an augmented system of  \eqref{eq:mom}-\eqref{eq:divB} by introducing the cleaning variable $\psi$. It serves as a generalized Lagrangian multiplier and it is coupled to the induction equation via a grad - div pair of differential operators. This generates artificial acoustic-type waves in $\psi$ and $\magfield$ that propagate away divergence errors that will arise due to numerical discretization errors inside a general-purpose scheme. The augmented GLM-MHD system reads: 
\begin{subequations}\label{eqn.GLMsystem}
	\begin{align}
		\dive \ww = 0, \\
		\partial_t\ww \ + \dive \Fluxcv\left(\Q \right)  - \dive\Fluxvv\left(\mu, \nabla \Q \right)  +\gra \press = {0},\\
		\partial_t \bbvar + \dive \Fluxcb \left(\Q \right)  -\dive\Fluxvb\left(\eta, \nabla \Q \right)  +\gra  \psi = {0},\\
		\partial_t \psi + \dive \left(c_{h}^2 \magfield\right)   = 0,
	\end{align}
\end{subequations}
with $\Q=\left(\ww, \magfield, \psi \right)^{T}$ the extended vector of unknowns, $\Fluxcv$, $\Fluxvv$, $\Fluxcb$, and $\Fluxvb$, the fluxes given in \eqref{eqn.ffluxmom} and \eqref{eq:FaradayCons}, and $c_{h}$ the divergence cleaning speed.
For a thermodynamically compatible formulation of GLM for MHD, the reader is referred to \cite{HTCMHD}. 

The main novelty of this section is the introduction of the well-balanced feature inside the hyperbolic GLM  divergence-cleaning method. To the best of our knowledge, a well-balanced GLM method for MHD does not exist yet. As done for the original system, to develop a well-balanced scheme, we assume a divergence-free stationary solution to be known, hence $\dive \magfielde = 0$ and therefore $\psi^{e}=0$. Then, after subtracting the equilibrium solution \eqref{eq:divwe}-\eqref{eqn.mage} from the MHD-GLM system \eqref{eqn.GLMsystem}, we obtain
\begin{subequations}\label{eqn.GLMsysteme}
	\begin{align}
		\dive \left( \ww - \wwe\right) = 0, \\
		\partial_t\ww  + \dive \Fluxcv\left(\Q \right) - \dive \Fluxcv\left(\Qe \right)   - \dive\Fluxvv\left(\mu, \nabla \Q \right) + \dive\Fluxvv\left(\mue, \nabla \Qe \right) +\gra \left( \press -\presse\right) = {0},\\
		\partial_t \bbvar + \dive \Fluxcb \left(\Q \right) - \dive \Fluxcb \left(\Qe \right)  -\dive\Fluxvb\left(\eta, \nabla \Q \right)  +\dive\Fluxvb\left(\etae, \nabla \Qe \right)  +\gra  \psi = {0},\\
		\partial_t \psi + \dive \left(c_{h}^2 \magfield\right)  - \dive \left(c_{h}^2 \magfielde\right) = 0. 
	\end{align}
\end{subequations}
As for the original system, the augmented GLM-MHD system from which the equilibrium solution has been subtracted, \eqref{eqn.GLMsysteme}, is now solved with exactly the same explicit finite volume approach introduced in Section~\ref{sec:transdiff}. We only need to modify the eigenvalues, which for the augmented system also include $c_h$ and read:
\begin{gather}
	\lambdacv_{1} = u-c_f,\quad \lambdacv_{2} =u-c_a,\quad 
	\lambdacv_{3,4} =u,\quad \lambdacv_{5} =2 u, 
	\lambdacv_{6} =u+c_a,\quad \lambdacv_{7} =u+c_f,\quad \lambdacv_{8} = - c_h,\quad \lambdacv_{9} =  c_h.
\end{gather} 

As shown in Diagram~\ref{dia:GLM}, within this approach the transport diffusion stage does not make use of the divergence-free reconstruction of the magnetic field any more and there is also no exactly divergence-free evolution of the magnetic field in the faces as the one introduced in Section~\ref{sec:divfree}. Instead, the scheme directly employs the value computed at the previous time step for the well-balanced reconstruction performed within the finite volume method (Section~\ref{sec:WB-recons}). Once the solution of the transport-diffusion stage is obtained, the projection stage is performed to get the new pressure and the momentum is corrected with the contribution of the pressure gradient. As such, the implementation of the GLM-MHD system is much simpler than the exactly divergence-free scheme presented before. 
\begin{Diagram}[h]
	\begin{equation*}
		\begin{tikzcd}
			\press^n \arrow[rrd]                                      &                        &                                                                                                                         & \press^{n+1} \arrow[rd, "\textrm{Post-proj. stage}" ] &           \\
			\Q^n_i, \psi^{n}_{i} \arrow[rr, "\mathrm{FV}" description]                                 & \arrow[r, "" ] & \tQ_i,  \psi^{n+1}_{i} \arrow[ru, "\textrm{Proj. stage}" ] \arrow[rr, "\magfield^{n+1}_i = \magfield^*_i" description]  &                                      & \Q^{n+1}_i 
		\end{tikzcd}
	\end{equation*}
	\caption{Overall well-balanced FV/FE scheme for the magnetohydrodynamics equations with a GLM divergence cleaning approach.} \label{dia:GLM}
\end{Diagram}

\subsection{Projection stage}\label{sec:projection}
The projection stage is devoted to the solution of the pressure subsystem~\eqref{eq:press_sys}, i.e., the pressure Poisson equation \eqref{eqn.pressurePoisson}. To this end, we need the intermediate solution $\tWW{}$ computed in the transport-diffusion stage, which contains the convective and viscous terms. Moreover, as has been indicated before, its interpolation between the primal and the dual mesh is essential to avoid stability issues. Since in the previous stage $\tWW{}$ was computed on the primal mesh, where also the pressure will be approximated, its value is now interpolated onto the dual mesh. In general, knowing the value of a certain variable in the primal mesh $Q_i$, $\forall\, i\in\left\lbrace1,\dots,N_{\mathrm{el}}\right\rbrace$, its value on a dual cell $C_{j}$, $Q_j$, is computed as 
\begin{equation}
Q_j = \frac{1}{|C_j|}\sum_{i\in\mathcal{K}_j}Q_i\,|T_{ij}|,
\label{eq:primal2dual}
\end{equation}
where $\mathcal{K}_j$ identifies the set of subelements of the dual element $C_{j}$, $\left|T_{ij}\right|$ is the area/volume of subelement $T_{ij}$, $T_{ij}\subset T_{i}$, $T_{ij}\subset C_{j}$. Analogously, to interpolate data from the dual to the primal cells, we compute
\begin{equation}
Q_i = \frac{1}{|T_i|}\sum_{j\in\mathcal{D}_i}Q_j\,|T_{ji}|.
\label{eq:dual2primal}
\end{equation}
Here $\mathcal{D}_i$ is the set of subelements contained in the primal element $T_i$. 
Directly applying \eqref{eq:primal2dual} to $\tWW{}$ inevitably leads to an excessive numerical dissipation of the Lax-Friedrichs-type, see \cite{SIMHD} for details. To avoid it, instead of directly interpolating $\tWW{}$ we will only account for the variation with respect to the previous time step, i.e., we average only $\delta Q_i := \delta \rho \bvel_{i}=\tWW{i}-\rho\bvel^{n}_{i}$ from the primal mesh to the dual one, hence obtaining 
$\rho\bvel_{j}^* = \delta \rho \bvel_{j} + \rho\bvel_{j}^{n}$ for the intermediate momentum on the dual mesh. 
%
Once the intermediate momentum is obtained on the dual cells we can efficiently address the pressure system using continuous finite elements.

Inserting \eqref{eqn.pressubs1} into \eqref{eqn.pressubs2} leads to the pressure Poisson equation \eqref{eqn.pressurePoisson} to be solved for each time step, which can be more compactly written as 
\begin{equation}
	\nabla^2 \delta \press = \frac{1}{\Delta t } \left(\dive \tWW{}-\dive \wwe\right),   \label{eq:poisson}  
\end{equation}
where $\delta \press= \press^{n+1}-\press^{n}$ is the sought pressure correction. 

Multiplying \eqref{eq:poisson} by a test function $z\in V_{0}=\left\lbrace z\in H^{1}(\Omega): \int_{\Omega} z \dV = 0 \right\rbrace$, integrating on the computational domain $\Omega$, and applying Green's formula, leads to
\begin{align}
	&\textit{Find $ \delta \press \in V_0$ such that} \nonumber\\
	&\Delta t \int \limits_{\Omega} \grae \delta \press \cdot \grae z \dV = \int \limits_{\Omega} (\tWW{}-\wwe) \cdot \grae z \dV -  \int \limits_{\partial\Omega} (\tWW{}-\wwe) \cdot \mathbf{n}\, z \dS,\qquad\forall z\in  V_{0}.\label{eq:weak_poisson} 
\end{align}
This weak problem is discretized using classical $\mathbb{P}_1$ Lagrange finite elements in the primal mesh. Regarding the right-hand side, the integrals are computed by employing the value of the momentum in the dual control volumes. The resulting linear system is solved using a matrix-free conjugate gradient algorithm. 

\subsection{Post-projection stage}\label{sec:postpro}

A final correction step must be carried out at the end of each time step, since the contribution of the term related to $\press^{n+1}$ is not yet included in the intermediate approximations obtained within the transport-diffusion stage. In order to accomplish this, first the gradient of the pressure jump is calculated on the half dual cells, using the corresponding basis functions, and weighted averaged to get its approximation on the dual cell, $\nabla \delta\press_{j}$. These values are then interpolated to the primal cells using \eqref{eq:dual2primal}, thus setting $Q_j:=\nabla \delta\press_{j}$ in \eqref{eq:dual2primal} in order to obtain $\nabla \delta\press_{i}$ on the primal mesh. Consequently, the final momentum on the dual grid is given by
\begin{equation}
	\WW^{n+1}_{j} = \tWW{j} - \Delta t\,\nabla \delta\press_{j},
	\label{eqn.postpro_mom_dual}
\end{equation} 
while the momentum on the primal mesh is updated as 
\begin{equation}
	\WW^{n+1}_{i} = \tWW{i} - \Delta t\,\nabla \delta\press_{i}.
	\label{eqn.postpro_mom_primal}
\end{equation} 
This closes the description of the algorithm. 

\subsection{Non well-balanced scheme and algorithm without pressure correction}
The numerical method presented in the former sections is designed to be well balanced in each step. To get a non well-balanced scheme with a pressure correction approach, it suffices to neglect the terms which depend on the equilibrium solution $\Qe,\,\presse$ throughout the proposed well-balanced scheme. Moreover, to get a non well-balanced scheme without the pressure-correction approach we would also cancel the explicit pressure term $p^n$ in the transport-diffusion equation \eqref{eq:FV}, solve the pressure system in Section \ref{sec:projection} directly for the new pressure, $\Press^{n+1}$, instead than for $\delta\Press=\Press^{n+1}-\Press^{n}$, and then the momentum will be directly updated as $\WW^{n+1} = \tWW{} - \Delta t\,\nabla \Press^{n+1}$.


\section{Numerical results}\label{sec:numericalresults} 
The proposed methodology is assessed at the aid of a large set of test problems with known exact or numerical reference solutions. The numerical convergence and the well-balanced property are both analysed within a pure hydrodynamics framework (i.e., setting $\magfield = 0$) and also for the complete MHD system. Then, some classical test problems, including a current sheet test, a magnetic field loop advection, the classical lid-driven cavity for the incompressible Navier-Stokes equations as well as a new version including also viscous and resistive MHD, the double shear layer, and the viscous and resistive Orszag-Tang vortex are run. Finally, a more complex 3D MHD test case of the Grad-Shafranov equilibrium in simplified 3D tokamak-type torus geometries is studied. 

Unless stated otherwise, the density is set to $\rho=1$, and the time step for all test cases is calculated according to the CFL-type condition
\begin{equation}
	\Delta t = \mathrm{CFL} \, \min\limits_{T_{i}} \left\lbrace \frac{\; h_{i}}{\max  \left|\lambdacv_{i}\right| \; + \frac{2}{h_i}\;\max \left|\lambdavv_{i}\right|}  \right\rbrace,
\end{equation}
with $\mathrm{CFL} < \frac{1}{d}$, $d$ the number of space dimensions, $h_{i}$ the incircle diameter of the primal element $T_{i}$ and $\left|\lambdacv_{i}\right|$, $\left|\lambdavv_{i}\right|$ its corresponding absolute eigenvalues for the convective and viscous subsystems, respectively.

\subsection{Numerical convergence tests and verification of the well-balanced property} 
To check the order of convergence of the proposed hybrid FV/FE method via numerical experiments, we perform two types of simulations: the first type, where we solve the incompressible Navier-Stokes equations, and the second, where we check the experimental order of convergence by solving the incompressible ideal MHD equations. We also carry out a series of numerical experiments in order to verify the well-balanced property of our new family of hybrid FV/FE schemes.  

\subsubsection{2D Taylor-Green vortex without magnetic field}\label{sec.TGV2D}

\begin{figure}[h!]
	\begin{center}
	\includegraphics[width=0.45\linewidth]{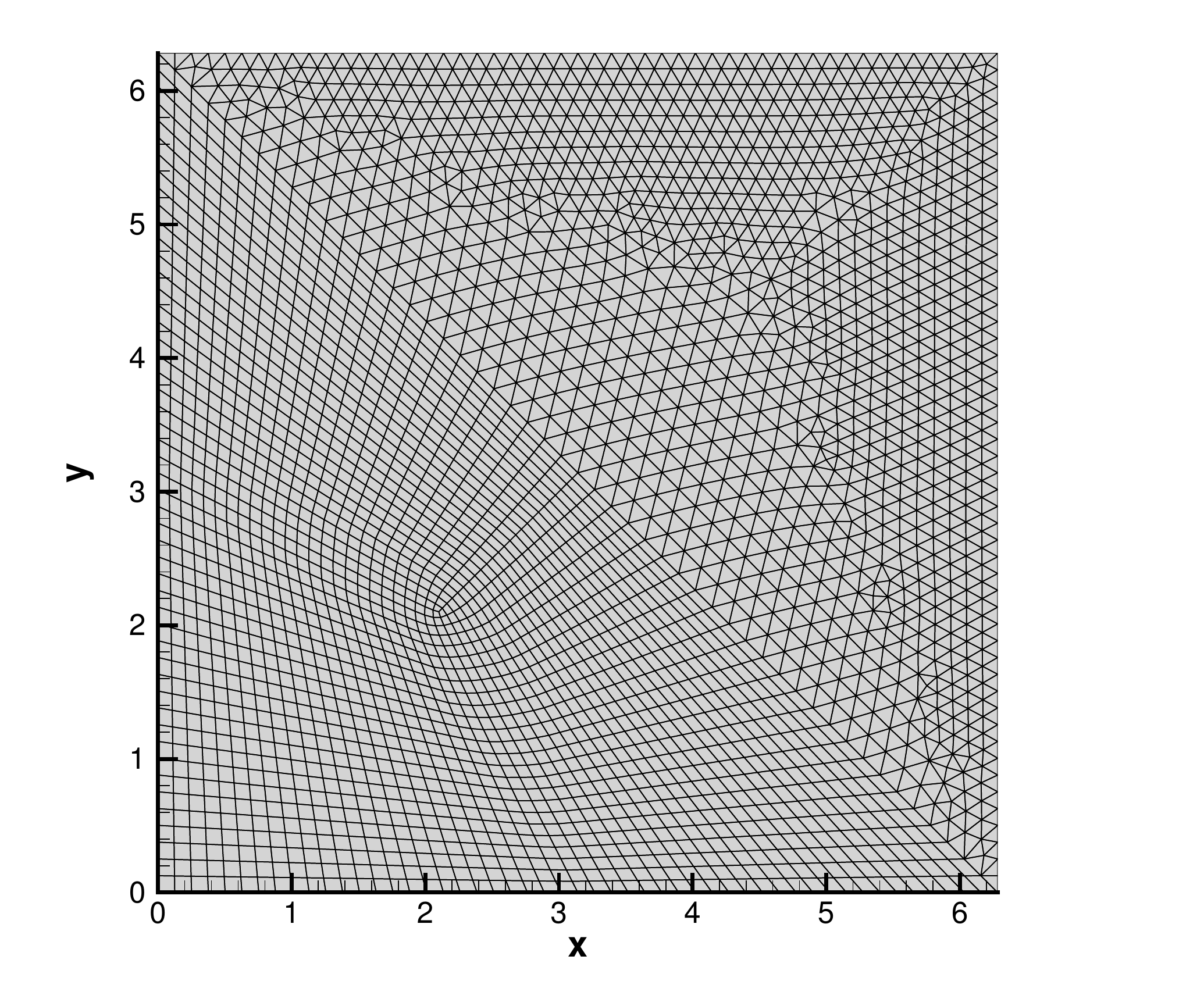}\hspace{-0.5cm}
	\includegraphics[trim=5 40 20 10,clip,width=0.55\textwidth]{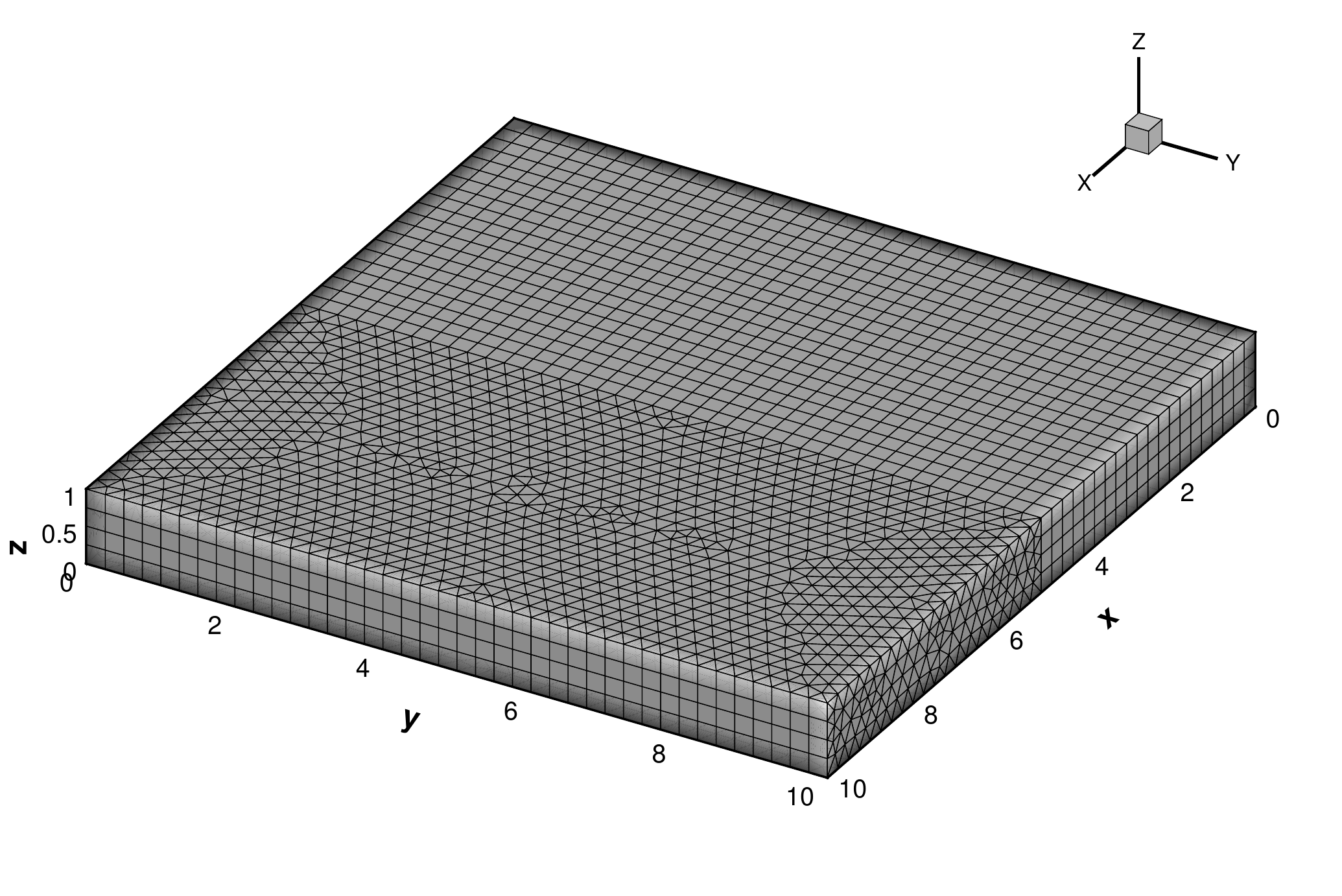}
	\end{center}
	\caption{Left: 2D unstructured mixed-element grid M$_{50}$, used in the domain $\Omega=[0,2\pi]^2$, with 4079 elements in total: 2204 triangles and 1875 quadrilateral elements. Right: 3D unstructured grid M$^{3\mathrm{D}}_{40}$, used in the domain $\Omega=[0,10]^2\times[0,1]$, with 21433 elements in total: 3200 Cartesian elements, 320 pyramids and 17913 tetrahedra.}
	\label{fig:mixmesh}
\end{figure}

First, we solve the incompressible Euler equations with $\mu=0$ and consider the Taylor-Green vortex benchmark in the two-dimensional domain $\Omega=\left[0,2\pi\right]^2$. The well-known exact solution for this test is 
\begin{equation}
	\bvel \left(\xx,t\right) = \left( \begin{array}{l} 
		\phantom{-}\sin(x)\cos(y) \\ 
		-\cos(x)\sin(y)\end{array} \right), \qquad 
	\press \left(\xx,t\right) = \frac{1}{4} \left(\cos(2x)+\cos(2y) \right).
	\label{eqn.TGV_ex}
\end{equation}
To run the simulations, the computational domain is discretized by using unstructured mixed-element primal grids formed by triangles and quadrilateral elements. The left plot of Figure~\ref{fig:mixmesh} shows, as an example, one of the used meshes (M$_{50}$), with a total amount of 4079 primal elements, 2204 of which are triangles and 1875 are quadrilateral elements. 
Table~\ref{tab:2D_mesh} contains some details about the grids used for the convergence analysis: the number of elements per edge in each coordinate direction, $N_x=N_y$, and the number of vertices and elements of the primal mesh, including the number of triangles and quadrangles of the grid. The subscript in mesh name states the number of elements in each direction ($N_x$). Periodic boundary conditions are set everywhere. Moreover, the time step is fixed for this test according to the mesh refinement. Table~\ref{tab:TGV_errors} contains the $L^2$ error norms for each mesh at time $t=0.1$ for the pressure, $\press$, at the vertex of the primal mesh and for both components of the velocity $\vel_1$ and $\vel_2$. These results were obtained using the scheme without well-balancing and considering the pressure correction approach. The differences between the errors computed by using either Rusanov or Ducros numerical fluxes are negligible, so the errors are shown only for the Rusanov flux. Since we are using a second order MUSCL-Hancock method in the explicit finite volume part of the scheme, the second order of accuracy for the pressure and the velocity fields is reached, as expected.  
\begin{table}[th]
	\caption{Main features of the computational grids and time steps considered to perform the convergence analysis in 2D of the hybrid FV/FE method.} 
	\label{tab:2D_mesh}	
	\renewcommand{\arraystretch}{1.2}
	\begin{center}		
		\begin{tabular}{ccccccc}
			\hline 
			Mesh & $N_x=N_y$ & Elements & Vertices & Triangles & Quadrangles & $\Dt$\\\hline
			M$_{20}$  & 20  & $666 $ & $524$ & $366$ & $300$ & $5.00\cdot 10^{-3} $\\ 
			M$_{40}$  & 40  & $2626 $ & $1994$ & $1426 $ & $1200 $ & $2.50\cdot 10^{-3} $\\ 
			M$_{50}$  & 50  & $4079 $ & $3078 $ & $2204 $ & $ 1875$ & $2.00\cdot 10^{-3} $\\ 
			M$_{60}$  & 60  & $5872 $ & $4407 $ & $3172 $ & $ 2700$ & $1.67\cdot 10^{-3} $\\ 
			M$_{80}$  & 80  & $10380 $ & $7751 $ & $5580 $ & $4800 $ & $1.25\cdot 10^{-3} $\\ 
			M$_{100}$ & 100 & $16174 $ & $12038 $ & $8674$ & $7500 $ & $1.00\cdot 10^{-3} $\\ 
			\hline 
		\end{tabular}
	\end{center}
\end{table}
\begin{table}[ht]
		\caption{Spatial $L^{2}$ error norms and convergence rates at time $t=0.1$ for the Taylor-Green vortex benchmark in 2D obtained using the hybrid FV/FE scheme with Rusanov flux function, without well-balancing and with pressure correction.}
	\label{tab:TGV_errors} 	
	\renewcommand{\arraystretch}{1.2}
	\begin{center}
		\begin{tabular}{ccccccc}
			\hline 
			Mesh 
			&$L^{2}_{\Omega}\left(\press \right)$ & $\mathcal{O}\left(\press \right)$                            
			&$L^{2}_{\Omega}\left(\vel_1\right)$ & $\mathcal{O}\left(\vel_1\right)$ 
			&$L^{2}_{\Omega}\left(\vel_2\right)$ & $\mathcal{O}\left(\vel_2\right)$
			\\ \hline
			M$_{20}$ & $1.0575\cdot 10^{-1} $ & $ $     &$3.4849 \cdot 10^{-2} $ & $ $  &$3.4809\cdot 10^{-2} $ & $ $ \\
			M$_{40}$ & $2.5991\cdot 10^{-2} $ & $2.02 $ &$8.4824\cdot 10^{-3} $ & $2.04 $  &$8.5197\cdot 10^{-3} $ & $2.03 $ \\
			M$_{50}$ & $1.6225\cdot 10^{-2} $ & $2.11 $ &$5.4300\cdot 10^{-3} $ & $2.00 $ &$5.4583\cdot 10^{-3} $ & $2.00 $ \\
			M$_{60}$ & $1.1166\cdot 10^{-2} $ & $2.05 $ &$3.7740\cdot 10^{-3} $ & $2.00 $ &$3.7926\cdot 10^{-3} $ & $2.00 $ \\
			M$_{80}$ & $6.2117\cdot 10^{-3} $ & $2.04 $ &$2.1337\cdot 10^{-3} $ & $1.98 $ &$2.1475\cdot 10^{-3} $ & $1.98 $ \\
			M$_{100}$ & $3.9129\cdot 10^{-3} $ & $2.07 $ &$1.3844\cdot 10^{-3} $ & $1.94 $ &$1.3873\cdot 10^{-3} $ & $1.96 $ \\
			\hline 
		\end{tabular}
	\end{center}
\end{table}

In order to illustrate that the well-balanced property of our new hybrid FV/FE scheme is satisfied and that the proposed method preserves stationary equilibrium solutions, some simulations of the 2D Taylor-Green vortex have been run using different machine precisions and imposing the exact solution \eqref{eqn.TGV_ex} as prescribed equilibrium solution. Table~\ref{tab:TGV_WB_Rusanov} reports the errors obtained when the simulation is run with single, double, and quadruple precision when considering the Rusanov and Ducros numerical flux functions. The errors are all of the order of machine precision, and therefore, clearly confirm the well-balanced property of the method. 

\begin{table}[ht!]
	\renewcommand{\arraystretch}{1.1}
	\caption{$L^{2}$ errors for the Taylor-Green vortex benchmark at time $t=0.1$ with different machine precisions for the well-balanced scheme.}
	\label{tab:TGV_WB_Rusanov}
	\begin{center}
		\begin{tabular}{cccc}
			\hline
			\multicolumn{4}{c}{Rusanov} \\
			\hline \hline 			
			Precision 
			&$L^{2}_{\Omega}\left(\press \right)$ 
			&$L^{2}_{\Omega}\left(\vel_1\right)$ 
			&$L^{2}_{\Omega}\left(\vel_2\right)$ 
			\\ \hline
			Single    & $9.4478 \cdot 10^{-6} $ & $1.0821 \cdot 10^{-6} $ & $1.0930 \cdot 10^{-6} $ \\
			Double    & $3.4167 \cdot 10^{-13}$ & $2.7140 \cdot 10^{-14}$ & $2.7030 \cdot 10^{-14}$ \\
			Quadruple & $3.8884 \cdot 10^{-27}$ & $2.6880 \cdot 10^{-28}$ & $2.6827 \cdot 10^{-28}$ \\
			\hline 

			\multicolumn{4}{c}{Ducros}\\ 
			\hline \hline 			
			Precision 
			&$L^{2}_{\Omega}\left(\press \right)$ 
			&$L^{2}_{\Omega}\left(u_1\right)$ 
			&$L^{2}_{\Omega}\left(u_2\right)$ 
			\\ \hline
			Single    & $1.1771 \cdot 10^{-6} $ & $1.0872 \cdot 10^{-6} $ & $1.0880 \cdot 10^{-6} $ \\
			Double    & $5.8481 \cdot 10^{-13}$ & $2.7159 \cdot 10^{-14}$ & $2.7014 \cdot 10^{-14}$ \\
			Quadruple & $2.5070 \cdot 10^{-27}$ & $2.6881 \cdot 10^{-28}$ & $2.6827 \cdot 10^{-28}$ \\ 
			\hline
		\end{tabular} 
	\end{center}
\end{table}

We now repeat the simulation, but this time prescribing as equilibrium solution a vortex defined by
\begin{align}
	\press = -\frac{1}{2} e^{-(r^2-1)}, \qquad
	\vel_1 	= \vel_0 - r e^{-\frac{1}{2}(r^2-1)}\sin(\phi), \qquad
	\vel_2 	= \vel_0 + r e^{-\frac{1}{2}(r^2-1)} \cos(\phi),
	\label{eqn.vortex}
\end{align}
where $\phi$ is the angle that satisfies $\tan \phi= (y - \vel_0 t - y_0)/(x - \vel_0 t - x_0)$, $r = \sqrt{(x-\vel_0 t-x_0)^2+(y-\vel_0 t-y_0)^2 }$ is the distance to the center of the computational domain $(x_0, y_0)$ in the $x-y$ plane and $\vel_0$ is a constant background velocity. In this case, we set $\vel_0=0$ and $x_0=y_0 = \pi$ and impose Dirichlet boundary conditions on all boundaries.   
The discrete initial condition and the discrete equilibrium solution are both initialized by evaluating the exact expressions given in \eqref{eqn.TGV_ex} and \eqref{eqn.vortex} and then performing five projection steps with time step $\Delta t = 10^{-7}$ in order to make sure that the discrete divergence errors of the velocity field in the initial condition are small.   
Table~\ref{tab:TGV_errors_EqVortex_Rusanov} contains the $L^2$ error norms for each mesh at time $t=0.1$ when we consider the initial condition \eqref{eqn.TGV_ex} together with the equilibrium \eqref{eqn.vortex}, which obviously do not coincide. As expected, when far from the prescribed equilibrium, the method converges to the exact solution of \eqref{eqn.TGV_ex} again with second order of accuracy.   

\begin{table}[ht]
		\caption{Spatial $L^{2}$ error norms and convergence rates at time $t=0.1$ obtained with the well-balanced hybrid FV/FE scheme with Rusanov flux for the Taylor-Green vortex in 2D for an equilibrium solution given by~\eqref{eqn.vortex} with $\vel_0=0$ and $x_0=y_0=\pi$.} 
	\label{tab:TGV_errors_EqVortex_Rusanov} 	
	\renewcommand{\arraystretch}{1.2}
	\begin{center}
		\begin{tabular}{ccccccc}
			\hline 
			Mesh 
			&$L^{2}_{\Omega}\left(\press\right)$ & $\mathcal{O}\left(\press \right)$                            
			&$L^{2}_{\Omega}\left(\vel_1\right)$ & $\mathcal{O}\left(\vel_1\right)$ 
			&$L^{2}_{\Omega}\left(\vel_2\right)$ & $\mathcal{O}\left(\vel_2\right)$
			\\ \hline
			M$_{20}$  & $8.5532\cdot 10^{-1} $ & $     $ &$2.2665\cdot 10^{-2} $ & $     $ &$2.3562\cdot 10^{-2} $ & $     $     \\
			M$_{40}$  & $2.2320\cdot 10^{-2} $ & $1.94 $ &$6.0195\cdot 10^{-3} $ & $1.91 $ &$6.2974\cdot 10^{-3} $ & $1.90 $  \\
			M$_{60}$  & $1.0051\cdot 10^{-2} $ & $1.97 $ &$2.8050\cdot 10^{-3} $ & $1.88 $ &$2.9340\cdot 10^{-3} $ & $1.88 $  \\
			M$_{80}$  & $5.7917\cdot 10^{-3} $ & $1.92 $ &$1.6042\cdot 10^{-3} $ & $1.94 $ &$1.7008\cdot 10^{-3} $ & $1.90 $  \\
			M$_{100}$ & $3.7256\cdot 10^{-3} $ & $1.98 $ &$1.0515\cdot 10^{-3} $ & $1.89 $ &$1.1271\cdot 10^{-3} $ & $1.84 $  \\
			\hline 
		\end{tabular}
	\end{center}
\end{table}

\subsubsection{3D Euler vortex}
Now, a stationary vortex in a three-dimensional domain is considered. We solve the incompressible Euler equations without magnetic field ($\magfield=0$, $\mu=0$) in the computational domain $\Omega=\left[0,10\right]^2\times[0,1]$, which is discretized by using unstructured primal grids formed by bricks, pyramids, and tetrahedra. Table~\ref{tab:3D_mesh} provides the main features of the three-dimensional grids used: the number of elements per edge in each coordinate direction, $N_x=N_y=N_z$, and the number of vertices and elements of the primal mesh, including the number of bricks, pyramids and tetrahedra. The right plot of Figure~\ref{fig:mixmesh} shows one of the used meshes M$^{3\mathrm{D}}_{40}$, with 21433 primal elements. The exact solution for this test is defined by \eqref{eqn.vortex} with $\vel_0=0$ and $x_0=y_0=5$. The simulation is performed with $\vel_0=0$ until $t=0.1$, and periodic boundary conditions are set everywhere. The simulations are run for the scheme without well-balancing and considering the pressure correction approach. The numerical convergence rates reported in Table \ref{tab:Vortex_errors}, show that second order of accuracy is reached for all variables, as expected. 

\begin{table}[th]
	\caption{Main features of the computational grids and time steps considered for the three-dimensional convergence tests.} 
	\label{tab:3D_mesh}	
	\renewcommand{\arraystretch}{1.5}
	\begin{center}		
		\begin{tabular}{cccccccc}
			\hline 
			Mesh & $N_x=N_y=N_z$ & Elements & Vertices & Bricks & Pyramids & Tetrahedra & $\Dt$\\\hline
			M$^{3\mathrm{D}}_{20}$  & 20  & $2659 $ & $1343$ & $400$ & $80$ & $2179$ & $5.00\cdot 10^{-3} $\\ 
			M$^{3\mathrm{D}}_{30}$  & 30  & $9445 $ & $3947$ & $1350 $ & $180$ & $7915 $ & $3.33\cdot 10^{-3} $\\ 
			M$^{3\mathrm{D}}_{40}$  & 40  & $21433 $ & $8413 $ & $3200 $ & $320$ & $ 17913$ & $2.50\cdot 10^{-3} $ \\ 
			M$^{3\mathrm{D}}_{50}$  & 50  & $41223 $ & $15454 $ & $6250 $ & $500$ & $ 34473$ & $2.00\cdot 10^{-3} $ \\ 
			M$^{3\mathrm{D}}_{60}$  & 60  & $68833 $ & $25308 $ & $10800 $ & $720$ & $57313 $ & $1.67\cdot 10^{-3} $ \\ 
			M$^{3\mathrm{D}}_{80}$ & 80 & $154580 $ & $55614 $ & $25600$ & $1280 $ & $127700 $ & $1.25\cdot 10^{-3} $\\ 
			\hline 
		\end{tabular}
	\end{center}
\end{table}
 
\begin{table}[ht]
	\caption{Spatial $L^{2}$ error norms and convergence rates at time $t=0.1$ obtained using the hybrid FV/FE scheme with Ducros flux for the stationary vortex in 3D without well-balancing and with pressure correction approach.} 
	\label{tab:Vortex_errors} 	
	\renewcommand{\arraystretch}{1.2}
	\begin{center}
		\begin{tabular}{ccccccc}
			\hline 
			Mesh 
			&$L^{2}_{\Omega}\left(\press\right)$ & $\mathcal{O}\left(\press\right)$                            
			&$L^{2}_{\Omega}\left(\vel_1\right)$ & $\mathcal{O}\left(\vel_1\right)$ 
			&$L^{2}_{\Omega}\left(\vel_2\right)$ & $\mathcal{O}\left(\vel_2\right)$
			\\ \hline
			M$^{3\mathrm{D}}_{20}$  & $1.3256\cdot 10^{-1} $ & $ $    &$5.8970\cdot 10^{-2} $ & $    $ &$5.9710\cdot 10^{-2} $ & $ $     \\
			M$^{3\mathrm{D}}_{30}$  & $4.7849\cdot 10^{-2} $ & $2.51 $ &$2.4526\cdot 10^{-2} $ & $2.16 $ &$2.6704\cdot 10^{-2} $ & $1.98 $  \\
			M$^{3\mathrm{D}}_{40}$  & $2.3492\cdot 10^{-2} $ & $2.47 $ &$1.4635\cdot 10^{-2} $ & $1.79 $ &$1.4955\cdot 10^{-2} $ & $2.02 $  \\
			M$^{3\mathrm{D}}_{50}$  & $1.4215\cdot 10^{-2} $ & $2.25 $ &$9.5686\cdot 10^{-3} $ & $1.90 $ &$9.7160\cdot 10^{-3} $ & $1.93 $  \\
			M$^{3\mathrm{D}}_{60}$ & $9.7048\cdot 10^{-3} $ & $2.09 $ &$6.7876\cdot 10^{-3} $ & $1.88 $ &$6.8443\cdot 10^{-3} $ & $1.92 $  \\
			M$^{3\mathrm{D}}_{80}$ & $5.2921\cdot 10^{-3} $ & $2.11 $ &$4.0781\cdot 10^{-3} $ & $1.77 $ &$4.1803\cdot 10^{-3} $ & $1.71 $  \\
			\hline 
		\end{tabular}
	\end{center}
\end{table}

\subsubsection{2D MHD vortex} 
The first MHD test considered is a stationary MHD vortex in the domain $\Omega = [0,10]^2$. The equilibrium solution is defined by 
\begin{equation}
	\begin{split}
	&
	\bvel \left(\xx,t\right) = e^{\frac{1}{2}(1-r^2)}\left( \begin{array}{l} 
		y_0 - y \\ 
		x - x_0\end{array} \right) + \left( \begin{array}{l} 
	v_0 \\ 
	v_0\end{array} \right) , \\
	&\press \left(\xx,t\right) = 1 + \frac{1}{2}e\left(1-r^2e^{-r^2}\right), \qquad \magfield\left(\xx,t\right) = e^{\frac{1}{2}(1-r^2)}\left( \begin{array}{l} 
		y_0 - y \\ 
		x - x_0\end{array} \right),
	\end{split}
	\label{eqn.MHD_ex}
\end{equation}
where  $r = \sqrt{(x-\vel_0 t-x_0)^2+(y-\vel_0 t-y_0)^2 }$ is the distance to the center of the computational domain $(x_0, y_0)$ in the $x-y$ plane, and $\vel_0$ is a constant background velocity. For this test, we set $\vel_0 = 0$ and $x_0 = y_0 = 5$. Periodic boundary conditions are imposed on all boundaries, and the time step is fixed according to the mesh refinement.

We start performing a numerical convergence analysis using the grids and the time steps already described in Table~\ref{tab:2D_mesh}. The results, reported in Table \ref{tab:MHD_errors}, show that the expected second order of accuracy is achieved for all variables. 

We now verify the well-balanced property by repeating the numerical experiment described in Section~\ref{sec.TGV2D} but for the MHD equations. We run some simulations of the MHD vortex with different machine precisions imposing \eqref{eqn.MHD_ex} as prescribed equilibrium solution. Table \ref{tab:MHD_WB_Rusanov_Ducros} reports the errors obtained with both the Rusanov and the Ducros fluxes with single, double, and quadruple precision. The errors are of the order of the machine precision, confirming the expected well-balanced property of the new algorithm proposed in this paper.

\begin{table}[ht]
	\caption{Spatial $L^{2}$ error norms and convergence rates at time $t=0.1$ obtained using the hybrid FV/FE scheme with Rusanov flux for the MHD vortex benchmark in 2D without well-balancing.}
	\label{tab:MHD_errors} 	
	\renewcommand{\arraystretch}{1.2}
	\begin{center}
		\begin{tabular}{ccccccc}
			\hline 
			Mesh 
			&$L^{2}_{\Omega}\left(\press \right)$ & $\mathcal{O}\left(\press \right)$                            
			&$L^{2}_{\Omega}\left(\vel_1\right)$ & $\mathcal{O}\left(\vel_1\right)$ 
			&$L^{2}_{\Omega}\left(\vel_2\right)$ & $\mathcal{O}\left(\vel_2\right)$ 
			\\ \hline
			M$_{20}$ & $1.6670\cdot 10^{-1}$&$$&$8.3987\cdot 10^{-2}$&$$&$8.4114\cdot 10^{-2}$&$$
			\\
			M$_{40}$ &  $4.0702\cdot 10^{-2}$&$2.0341$&$2.0105\cdot 10^{-2}$&$2.0626$&$2.0281\cdot 10^{-2}$&$2.0522$\\
			M$_{50}$ &  $2.5679\cdot 10^{-2}$&$2.0642$&$1.2612\cdot 10^{-2}$&$2.0897$&$1.2678\cdot 10^{-2}$&$2.1055$  \\
			M$_{60}$ & $1.7641\cdot 10^{-2}$&$2.0592$&$8.6049\cdot 10^{-3}$&$2.0971$&$8.6409\cdot 10^{-3}$&$2.1026$		 \\
			M$_{80}$ & $9.7836\cdot 10^{-3}$&$2.0492$&$4.7170\cdot 10^{-3}$&$2.0897$&$4.7666\cdot 10^{-3}$&$2.0678$\\
			M$_{100}$ & $6.2362\cdot 10^{-3}$&$2.0181$&$2.9862\cdot 10^{-3}$&$2.0487$&$3.0103\cdot 10^{-3}$&$2.0596$  \\
			\hline\hline 
			 	&$L^{2}_{\Omega}\left(B_1\right)$ & $\mathcal{O}\left(B_1\right)$ 
			 &$L^{2}_{\Omega}\left(B_2\right)$ & $\mathcal{O}\left(B_2\right)$  \\
			 \hline
			 M$_{20}$ &$1.3789\cdot 10^{-1}$&$$&$1.4064\cdot 10^{-1}$&$$\\
			 M$_{40}$ &$3.3165\cdot 10^{-2}$&$2.0558$&$3.3544\cdot 10^{-2}$&$2.0678$\\
			 M$_{50}$ &$2.0987\cdot 10^{-2}$&$2.0507$&$2.0930\cdot 10^{-2}$&$2.1137$\\
			 M$_{60}$ &$1.4441\cdot 10^{-2}$&$2.0504$&$1.4305\cdot 10^{-2}$&$2.0875$\\
			 M$_{80}$ &$7.9439\cdot 10^{-3}$&$2.0775$&$7.9804\cdot 10^{-3}$&$2.0286$ \\
			 M$_{100}$   &$5.0967\cdot 10^{-3}$&$1.9889$&$5.0781\cdot 10^{-3}$&$2.0258$\\
			\hline 
		\end{tabular}
	\end{center}
\end{table}

\begin{table}[ht!]
	\renewcommand{\arraystretch}{1.1}
	\caption{$L^{2}$ errors for the MHD vortex benchmark at time $t=0.1$ with different machine precisions when the simulation is carried out with the well-balanced scheme.}
	\label{tab:MHD_WB_Rusanov_Ducros}
	\begin{center}
		\begin{tabular}{cccccc}
			\hline
			\multicolumn{6}{c}{Rusanov} \\
			\hline \hline 			
			Precision 
			&$L^{2}_{\Omega}\left(\press \right)$ 
			&$L^{2}_{\Omega}\left(\vel_1\right)$ 
			&$L^{2}_{\Omega}\left(\vel_2\right)$ 
			&$L^{2}_{\Omega}\left(B_1\right)$ 
			&$L^{2}_{\Omega}\left(B_2\right)$ 
			\\ \hline
			Single    & $1.0089 \cdot 10^{-5} $ & $9.1549 \cdot 10^{-7} $ & $8.9882 \cdot 10^{-7} $  & $5.7389 \cdot 10^{-7}$ & $5.7724 \cdot 10^{-7}$\\
			Double    & $6.3177 \cdot 10^{-13}$ & $4.6600 \cdot 10^{-14}$ & $4.6360 \cdot 10^{-14}$  & $5.7500 \cdot 10^{-14}$  & $5.7894 \cdot 10^{-14}$ \\
		
			Quadruple & $1.1115 \cdot 10^{-27}$ & $4.6578 \cdot 10^{-28}$ & $4.6304 \cdot 10^{-28}$  & $5.7644 \cdot 10^{-28}$  & $5.7856 \cdot 10^{-28}$ \\
			\hline 
			
			\multicolumn{6}{c}{Ducros}\\ 
			\hline \hline 			
			Precision 
			&$L^{2}_{\Omega}\left(\press \right)$ 
			&$L^{2}_{\Omega}\left(\vel_1\right)$ 
			&$L^{2}_{\Omega}\left(\vel_2\right)$ 
			&$L^{2}_{\Omega}\left(B_1\right)$ 
			&$L^{2}_{\Omega}\left(B_2\right)$ 
			\\ \hline
			Single    & $9.4873 \cdot 10^{-6} $ & $8.7277 \cdot 10^{-7} $ & $8.9653 \cdot 10^{-7}$  & $5.7455 \cdot 10^{-7}$  & $5.7729 \cdot 10^{-7} $ \\
			Double    & $7.5255 \cdot 10^{-13}$ & $4.6622 \cdot 10^{-14}$ & $4.6357 \cdot 10^{-14}$  & $5.7497 \cdot 10^{-14}$  & $5.7889 \cdot 10^{-14}$ \\
			Quadruple & $4.6333 \cdot 10^{-27}$ & $4.6579 \cdot 10^{-28}$ & $4.6305 \cdot 10^{-28}$  & $5.7644 \cdot 10^{-28}$  & $ 5.7856 \cdot 10^{-28}$ \\ 
			\hline
		\end{tabular} 
	\end{center}
\end{table}

\subsubsection{3D MHD vortex}
We repeat the numerical convergence analysis for the same MHD vortex given in \eqref{eqn.MHD_ex}, but now inside a three-dimensional computational domain $\Omega = [0,10]^2 \times [0,1]$. We use the grids and the time steps listed in Table \ref{tab:3D_mesh}, starting from M$^{3\mathrm{D}}_{40}$. The exact solution of this problem is \eqref{eqn.MHD_ex} with $\vel_0 = 0$ and $x_0 = y_0 = 5$. The simulation is performed until a final time of $t=0.1$ and periodic boundary conditions are set everywhere. The differences between the results obtained with Rusanov and Ducros fluxes are negligible, therefore we report only the errors obtained from the former flux. In Table \ref{tab:MHD3dVortex_errors} we observe that second order of convergence is reached, as expected.  
\begin{table}[ht]
	\caption{Spatial $L^{2}$ error norms and convergence rates at time $t=0.1$ obtained with hybrid FV/FE scheme with Rusanov flux for the MHD stationary vortex in 3D.} 
	\label{tab:MHD3dVortex_errors} 	
	\renewcommand{\arraystretch}{1.2}
	\begin{center}
		\begin{tabular}{ccccccc}
			\hline 
			Mesh 
			&$L^{2}_{\Omega}\left(\press\right)$ & $\mathcal{O}\left(\press\right)$                            
			&$L^{2}_{\Omega}\left(\vel_1\right)$ & $\mathcal{O}\left(\vel_1\right)$ 
			&$L^{2}_{\Omega}\left(\vel_2\right)$ & $\mathcal{O}\left(\vel_2\right)$
			\\ \hline
				M$^{3\mathrm{D}}_{40}$  & $2.7716\cdot 10^{-2}$&$$&$1.4473\cdot 10^{-2}$&$$&$1.5287\cdot 10^{-2}$&$$\\
			M$^{3\mathrm{D}}_{50}$  &  $1.7843\cdot 10^{-2}$&$1.9737$&$9.5332\cdot 10^{-3}$&$1.8710$&$9.9961\cdot 10^{-3}$&$1.9038$\\
			M$^{3\mathrm{D}}_{60}$ &	$1.2930\cdot 10^{-2}$&$1.7666$&$6.8011\cdot 10^{-3}$&$1.8522$&$7.2149\cdot 10^{-3}$&$1.7883$\\
			M$^{3\mathrm{D}}_{80}$ &  $7.6015\cdot 10^{-3}$&$1.8464$&$4.2458\cdot 10^{-3}$&$1.6377$&$4.4381\cdot 10^{-3}$&$1.6891$\\
			\hline
				&$L^{2}_{\Omega}\left(B_1\right)$ & $\mathcal{O}\left(B_1\right)$ 
			&$L^{2}_{\Omega}\left(B_2\right)$ & $\mathcal{O}\left(B_2\right)$\\
			\hline
				M$^{3\mathrm{D}}_{40}$  &$1.9685\cdot 10^{-2}$&$$&$2.0359\cdot 10^{-2}$&$$\\
					M$^{3\mathrm{D}}_{50}$  &$1.2675\cdot 10^{-2}$&$1.9729$&$1.3154\cdot 10^{-2}$&$1.9574$\\
						M$^{3\mathrm{D}}_{60}$ &$9.0085\cdot 10^{-3}$&$1.8729$&$9.3240\cdot 10^{-3}$&$1.8877$\\
						M$^{3\mathrm{D}}_{80}$ &$5.2800\cdot 10^{-3}$&$1.8570$&$5.5552\cdot 10^{-3}$&$1.8001$\\
						
				\hline 
		\end{tabular}
	\end{center}
\end{table}

\subsection{Small perturbation of an equilibrium solution without and with magnetic field}

In order to show the clear benefit of the well-balanced scheme over the non well-balanced one, we now carry out a set of  numerical simulations where a small perturbation is added to a given equilibrium solution, one without magnetic field and one with magnetic field. The computational domain used for this test is $\Omega = [0,10]^2$ discretized with mesh M$_{40}$ and considering periodic boundary conditions everywhere. 

\begin{figure}[h!]
	\begin{center}
		\includegraphics[width=0.45\textwidth]{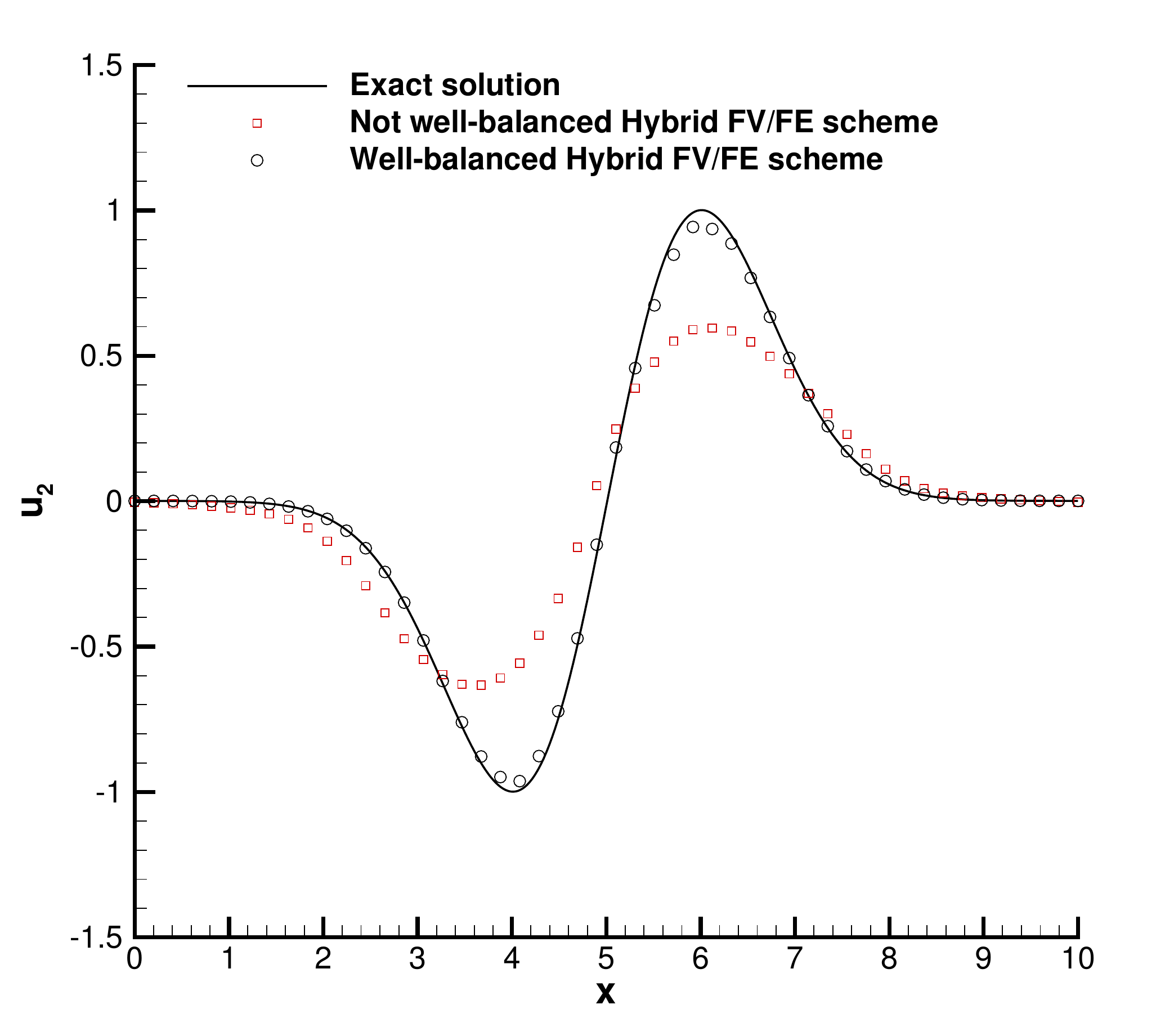}
		\includegraphics[width=0.45\textwidth]{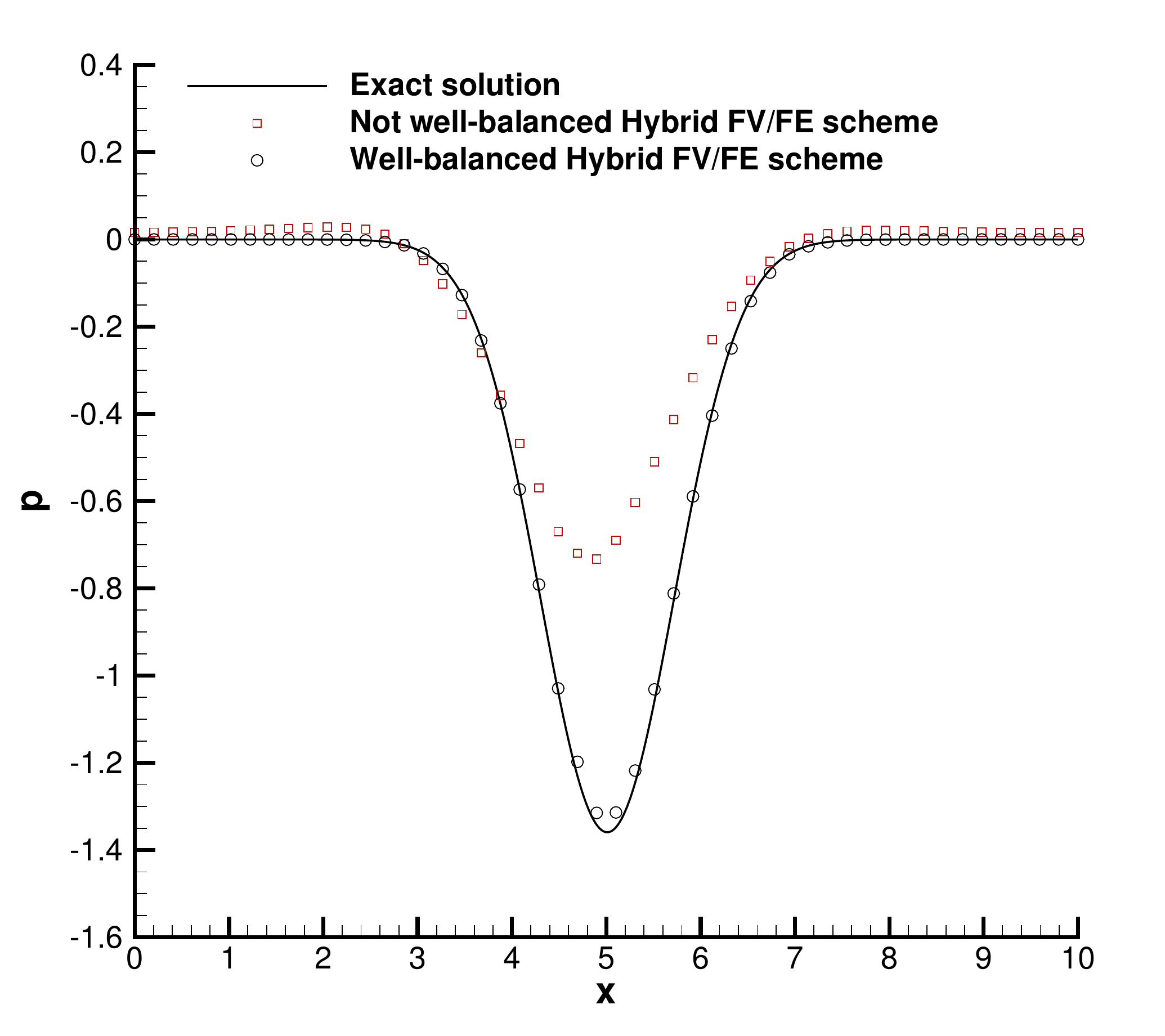}
	\end{center}
	\caption{Comparison of the exact solution with the non well-balanced and the well-balanced hybrid FV/FE scheme for a \textit{small perturbation} of the 2D vortex problem  \eqref{eqn.vortex} with $\vel_0=10^{-3}$ and $x_0=y_0=5$ at time $t=10$ using mesh M$_{40}$. Left: 1D cut along the $x-$axis through the velocity field $\vel_2$. Right: 1D cut along the $x-$axis through the pressure field $\press$.}
\label{fig:vortex_pert}
\end{figure}

\paragraph{2D vortex without magnetic field} The initial condition is now given by the vortex \eqref{eqn.vortex} with $x_0=y_0=5$, adding a small perturbation $\vel_0=10^{-3}$ to the velocity field. The equilibrium solution is given by the stationary vortex \eqref{eqn.vortex} without velocity  perturbation, i.e., with $\vel_0=0$. Simulations are run with the well-balanced hybrid FV/FE method and with the non well-balanced scheme with pressure correction until a final time of $t=10$, which corresponds to one full advection period through the domain. A 1D cut through the numerical solutions and the exact solution along the $x-$axis is presented in Figure \ref{fig:vortex_pert}, while the spatial $L^2$ error norms at time $t=10$ are listed in Table \ref{tab:euler_vortex_pert}. As expected, from the obtained numerical results, we clearly see the much better performance of the well-balanced scheme compared to the non well-balanced one when studying small perturbations around stationary equilibria. While the overall computational cost of both methods is comparable, the errors obtained with the well-balanced scheme are two orders of magnitude better than those of the simple non well-balanced method.     

\begin{table}[ht]
	\caption{Spatial $L^{2}$ error norms obtained at $t=10$ on mesh M$_{40}$ with the new well-balanced hybrid FV/FE scheme and with the classical non well-balanced scheme for the simulation of a small perturbation ($\vel_0=10^{-3}$) of the stationary 2D vortex without magnetic field \eqref{eqn.vortex}. }   
	\label{tab:euler_vortex_pert} 	
	\renewcommand{\arraystretch}{1.2}
	\begin{center}
		\begin{tabular}{cccc}
			\hline 
			numerical scheme  
			&$L^{2}_{\Omega}\left(\press\right)$ 
			&$L^{2}_{\Omega}\left(u_1\right)$ 
			&$L^{2}_{\Omega}\left(u_2\right)$
			\\ \hline
			well-balanced hybrid FV/FE scheme      & $4.4593\cdot 10^{-3} $ &$6.5385\cdot 10^{-3} $ &$5.3461\cdot 10^{-3}$   \\
			non well-balanced  hybrid FV/FE scheme & $7.3024\cdot 10^{-1} $ &$6.1672\cdot 10^{-1} $ &$7.7866\cdot 10^{-1}$   \\
			\hline 
		\end{tabular}
	\end{center}
\end{table}

\begin{figure}[h]
	\centering
	\includegraphics[width=0.42\linewidth]{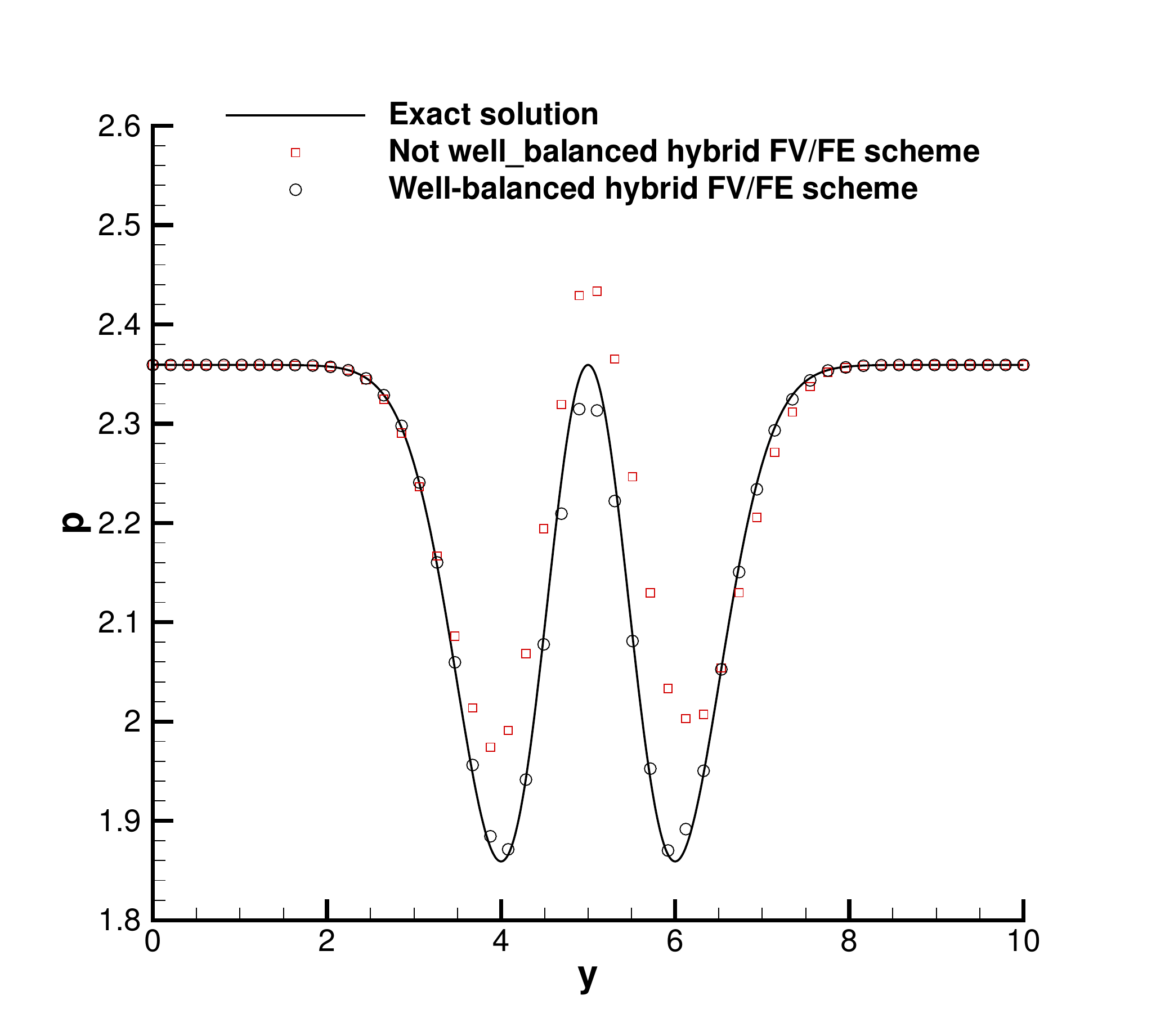}
	\includegraphics[width=0.42\linewidth]{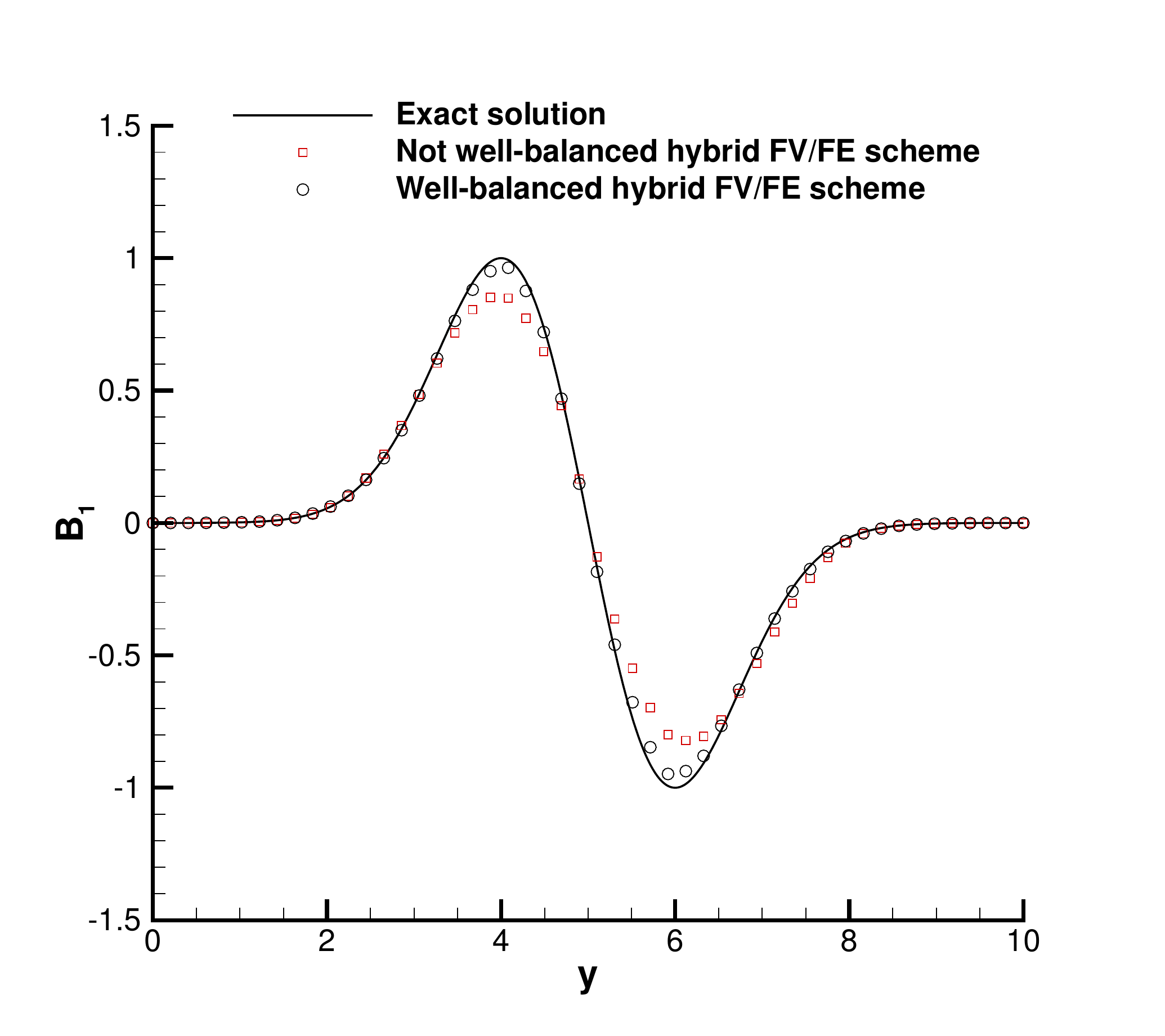}
	\caption{Comparison of the exact solution with the non well-balanced and the well-balanced hybrid FV/FE scheme for a \textit{small perturbation} of the 2D MHD vortex problem with $\vel_0=10^{-3}$ at time $t=10$ using mesh M$_{40}$. Left: 1D cut at $x=5$ through the pressure field $\press$. Right: 1D cut at $x=5$ through the first component of the magnetic field $B$.}
	\label{fig:MHD_vortex_perturbation}
\end{figure}

\paragraph{2D vortex with magnetic field}
We now repeat the same numerical experiment with the MHD vortex \eqref{eqn.MHD_ex}, setting again $\vel_0 = 10^{-3}$ for the initial condition and $\vel_0 = 0$ for the equilibrium solution. The simulations are run until a final time of $t = 10$ so that the vortex completes one entire advection period throughout the domain. A 1D cut of the numerical solution and the exact solution is shown at $x=5$ in Figure \ref{fig:MHD_vortex_perturbation}, while the spatial $L^2$ errors computed at the final time are listed in Table \ref{tab:MHD_vortex_pert}. 

\begin{table}[ht]
	\caption{Spatial $L^{2}$ error norms obtained at $t=10$ on mesh M$_{40}$ with the new well-balanced (WB) hybrid FV/FE scheme and with the non well-balanced (non WB) scheme for the simulation of a small perturbation ($\vel_0=10^{-3}$) of the stationary 2D MHD vortex. }   
	\label{tab:MHD_vortex_pert} 	
	\renewcommand{\arraystretch}{1.2}
	\begin{center}
		\begin{tabular}{cccccc}
			\hline 
			numerical scheme  
			&$L^{2}_{\Omega}\left(\press\right)$ 
			&$L^{2}_{\Omega}\left(\vel_1\right)$ 
			&$L^{2}_{\Omega}\left(\vel_2\right)$
			&$L^{2}_{\Omega}\left(B_1\right)$
			&$L^{2}_{\Omega}\left(B_2\right)$
			\\ \hline
			WB hybrid FV/FE scheme      & $2.9094\cdot 10^{-3} $ &$2.8213\cdot 10^{-3} $ &$3.4405\cdot 10^{-3}$ &  $2.8580\cdot 10^{-3}$ & $3.2472\cdot 10^{-3}$     \\
			not WB  hybrid FV/FE scheme & $3.3081\cdot 10^{-1} $ &$3.4881\cdot 10^{-1} $ &$3.5176\cdot 10^{-1}$ &$2.7129\cdot 10^{-1}$ &$2.6907\cdot 10^{-1}$  \\
			\hline 
		\end{tabular}
	\end{center}
\end{table}

Also for this second test problem, which adds the presence of a magnetic field compared to the first one, the error norms of the new well-balanced hybrid FV/FE scheme are \textit{two orders of magnitude} lower than those of the classical non well-balanced method. As expected, in the case of simulations of small perturbations around a stationary equilibrium, the obtained computational results  clearly show the substantially increased accuracy of the new well-balanced hybrid FV/FE scheme proposed in this paper compared to the same but non well-balanced scheme. 

\subsection{First problem of Stokes and current sheet test}  

\begin{figure}[h!]
	\begin{center}
	\includegraphics[trim=20 5 20 75,clip,width=1\textwidth]{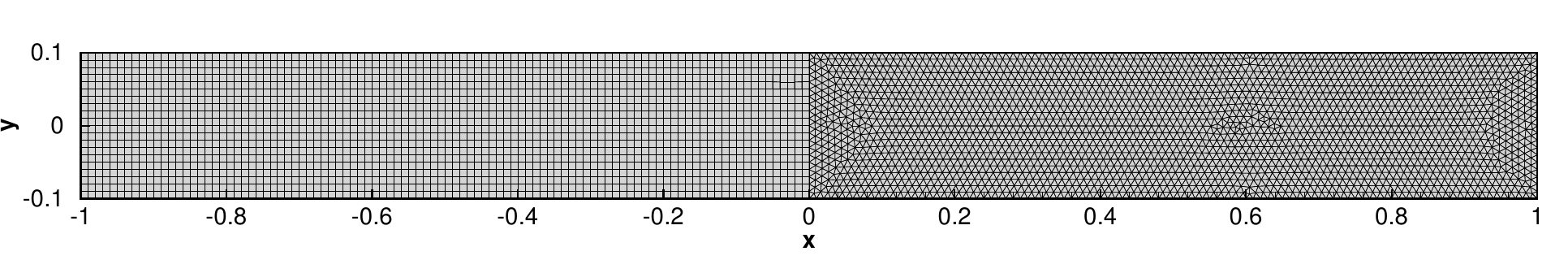}
	\end{center}
	\caption{Unstructured mix grid for the first problem of Stokes in 2D. }
	\label{fig:FPS2D_mixmesh}
\end{figure}

\begin{figure}[h]
	\centering
	\includegraphics[width=0.325\linewidth]{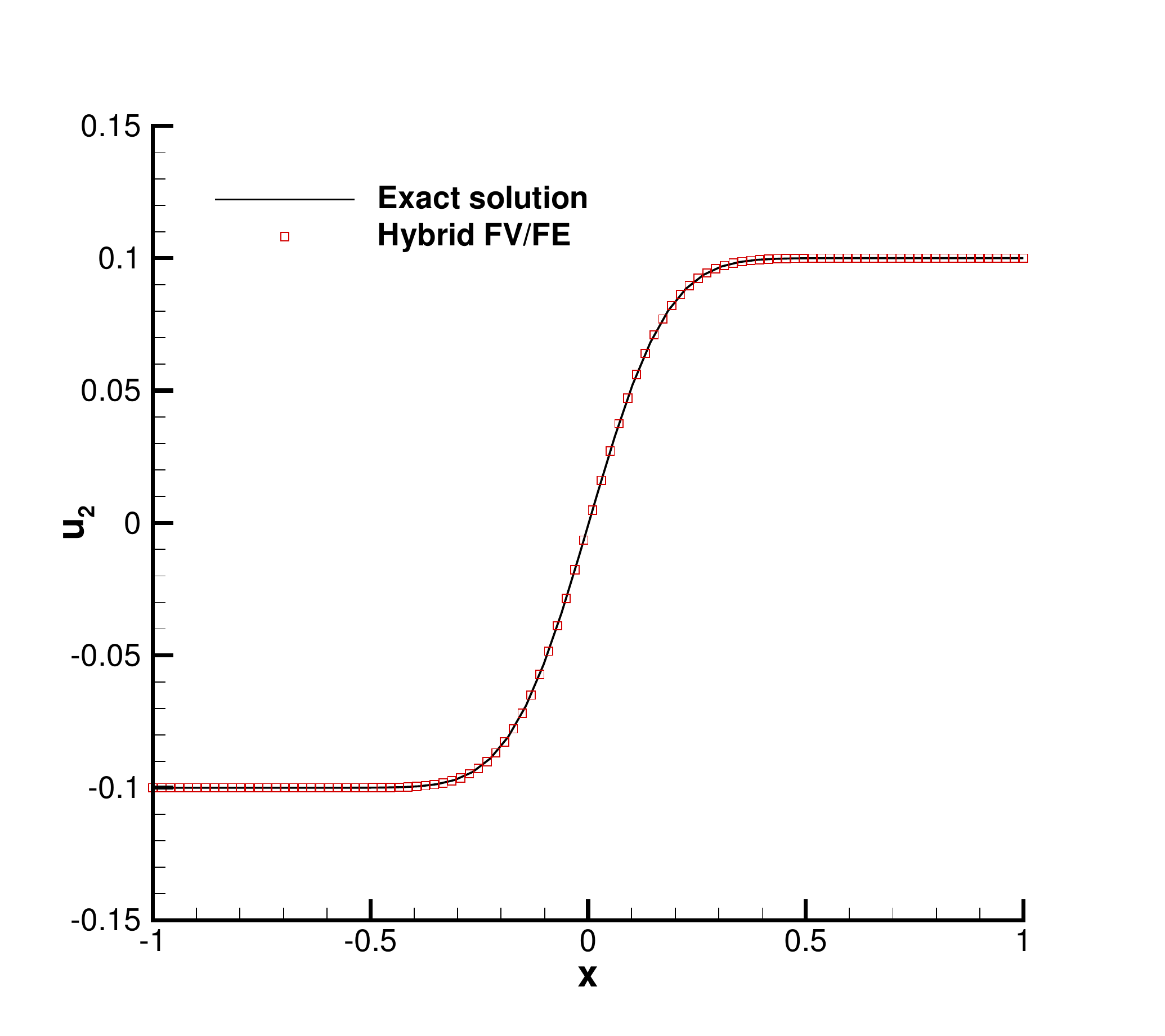}
	\includegraphics[width=0.325\linewidth]{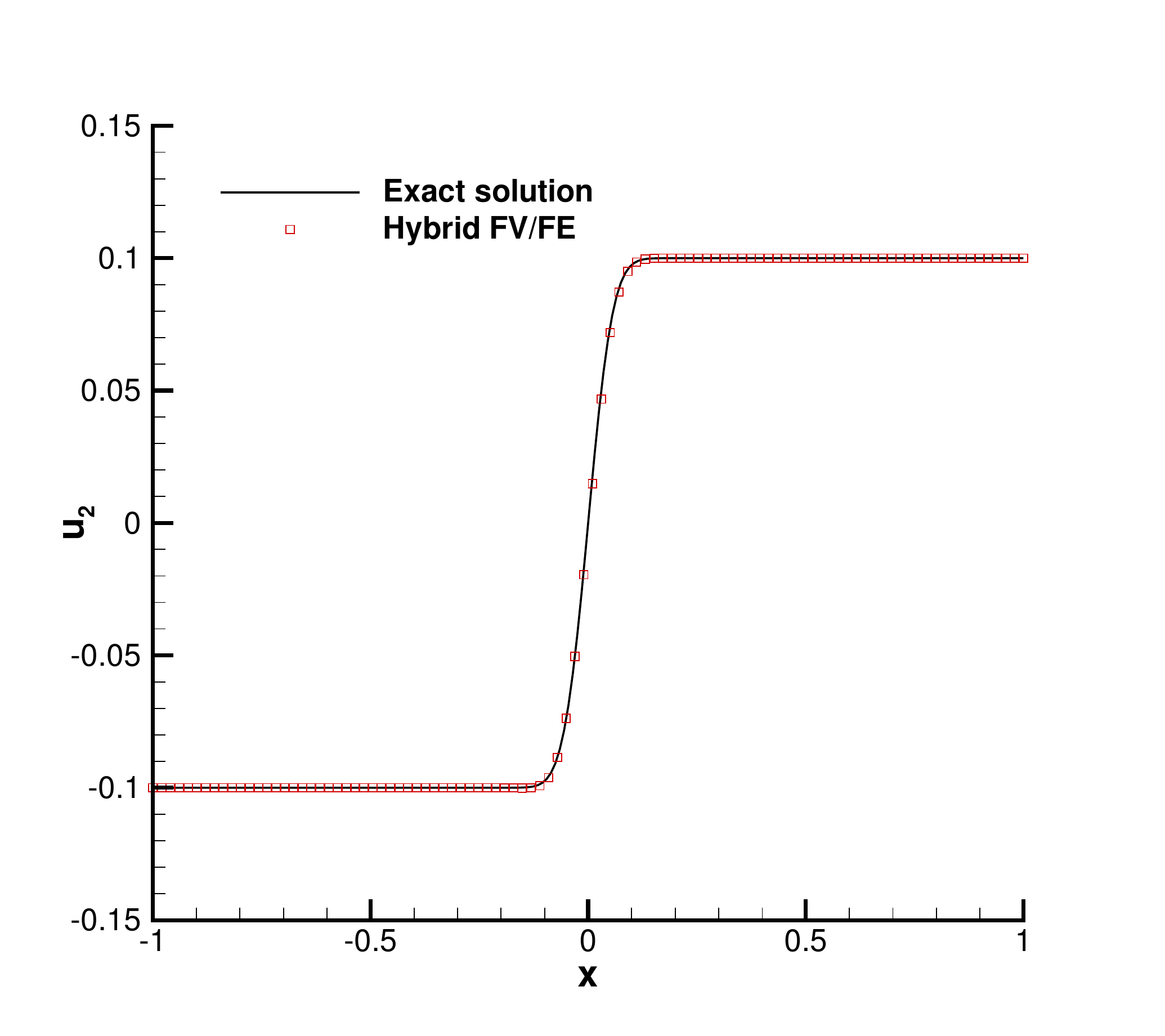}
	\includegraphics[width=0.325\linewidth]{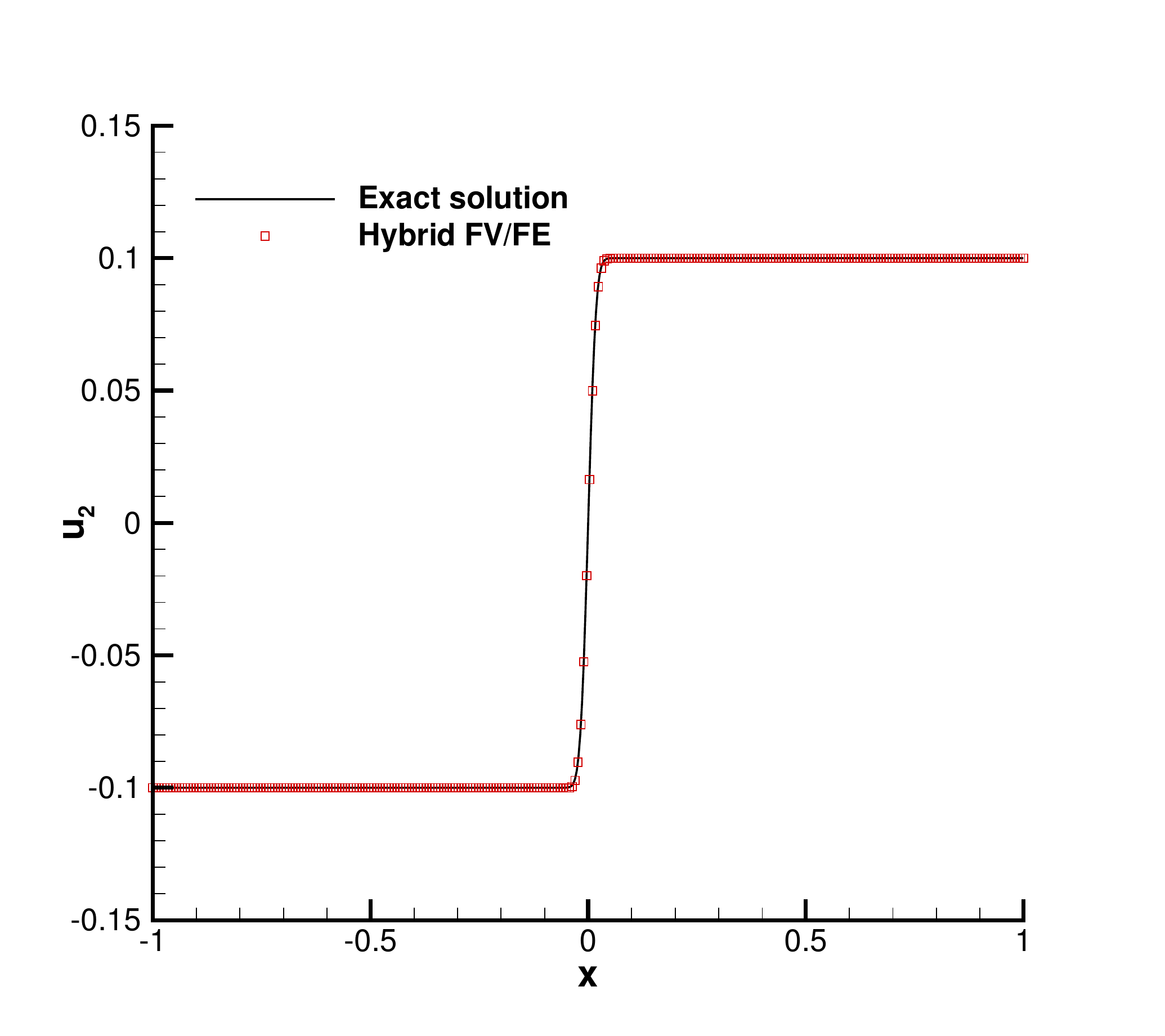}
	\caption{Velocity along the cut $y=0$ for the first problem of Stokes in 2D at time $t=1$. Comparison between the values computed with the non well-balanced scheme with pressure correction and the exact solution~\eqref{eqn.FPS_ex} for $\mu=10^{-2}$ (left), $\mu=10^{-3}$ (center), and $\mu=10^{-4}$ (right).}
	\label{fig:FPS2D}
\end{figure}

In this section, we solve the first problem of Stokes in the computational domain $\Omega=[-1,1]\times[-0.1,0.1]$. The initial condition is given by
\begin{equation*}
\press \left(\xx,0\right) = 1, \qquad
\vel_{1} \left(\xx,0\right) = 0, \qquad
\vel_{2} \left(\xx,0\right) = \left\lbrace \begin{array}{lc}
-0.1 & \mathrm{ if } \; x \le 0,\\
\phantom{-} 0.1 & \mathrm{ if } \; x > 0.
\end{array}\right.
\end{equation*}
The exact solution for the velocity in the $y$-direction reads, see \cite{BLTheory},  
\begin{equation}
\vel_{2} \left(\xx,t\right) = \frac{1}{10} \mathrm{erf}\left( \frac{x}{2\sqrt{\mu t}}\right).
\label{eqn.FPS_ex}
\end{equation}
Three different viscosities are considered $\mu\in\left\lbrace 10^{-2},10^{-3},10^{-4}\right\rbrace$. In the $y$-direction, periodic boundary conditions are considered, and Dirichlet conditions for the velocity and density are imposed on the left and right boundaries. In this case, we have employed a mixed-element mesh: for $x\in[-1,0]$ we have a Cartesian arrangement, and for $x\in[0,1]$, we have a triangular grid (see Figure~\ref{fig:FPS2D_mixmesh}). The simulations are carried out up to a final time of $t=1$, considering the non well-balanced scheme with pressure correction. Figure~\ref{fig:FPS2D} shows the numerical results for $\mu = 10^{-2}$, $\mu =10^{-3}$, and $\mu =10^{-4}$ in the left, center, and right plots, respectively. There, 1D cuts along $y=0$ are plotted and the results obtained with the proposed method are compared with the exact solution given by~\eqref{eqn.FPS_ex}. 

\begin{figure}[h]
	\centering
	\includegraphics[trim=15 5 20 10,clip,width=0.45\textwidth]{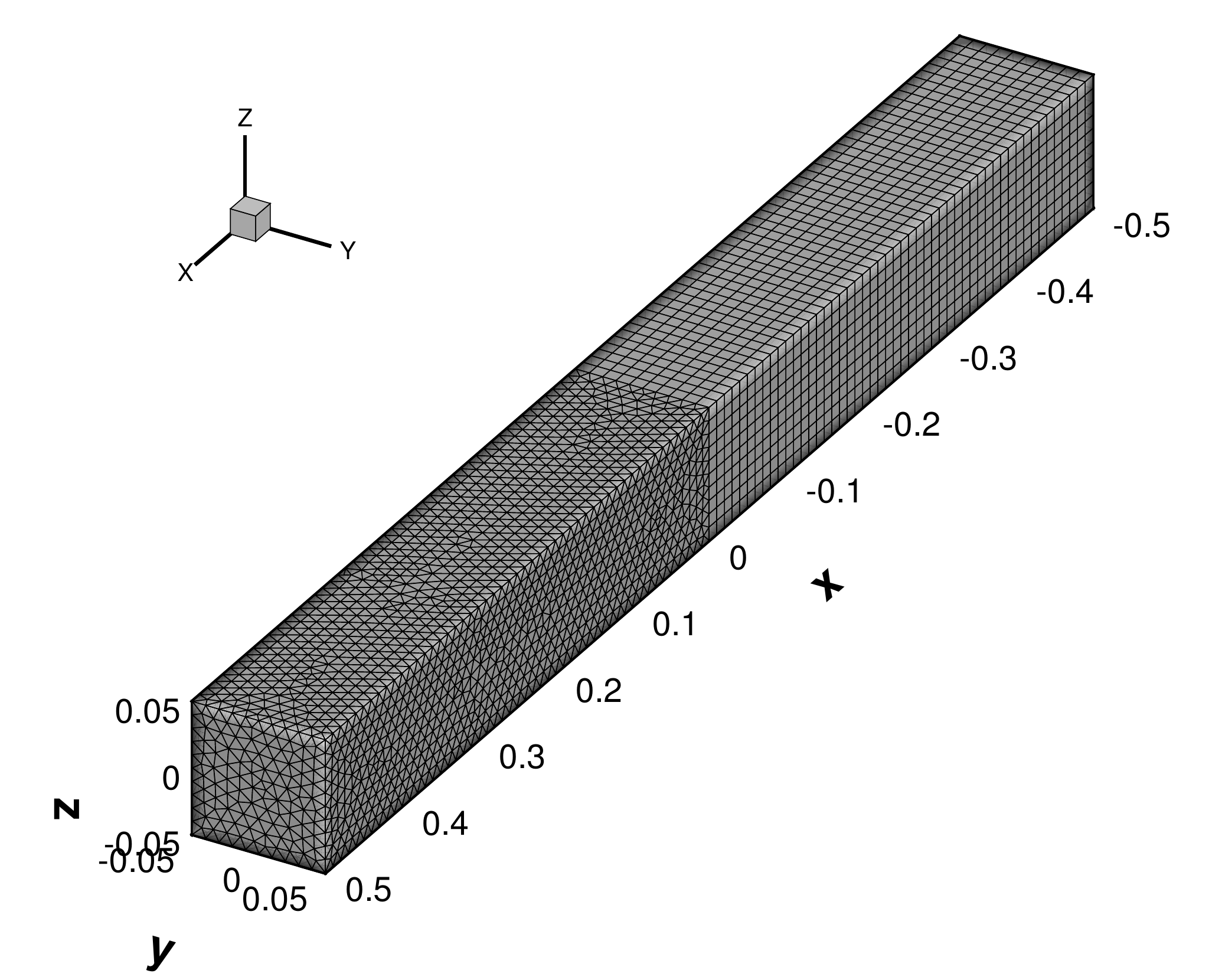}
	\includegraphics[width=0.45\linewidth]{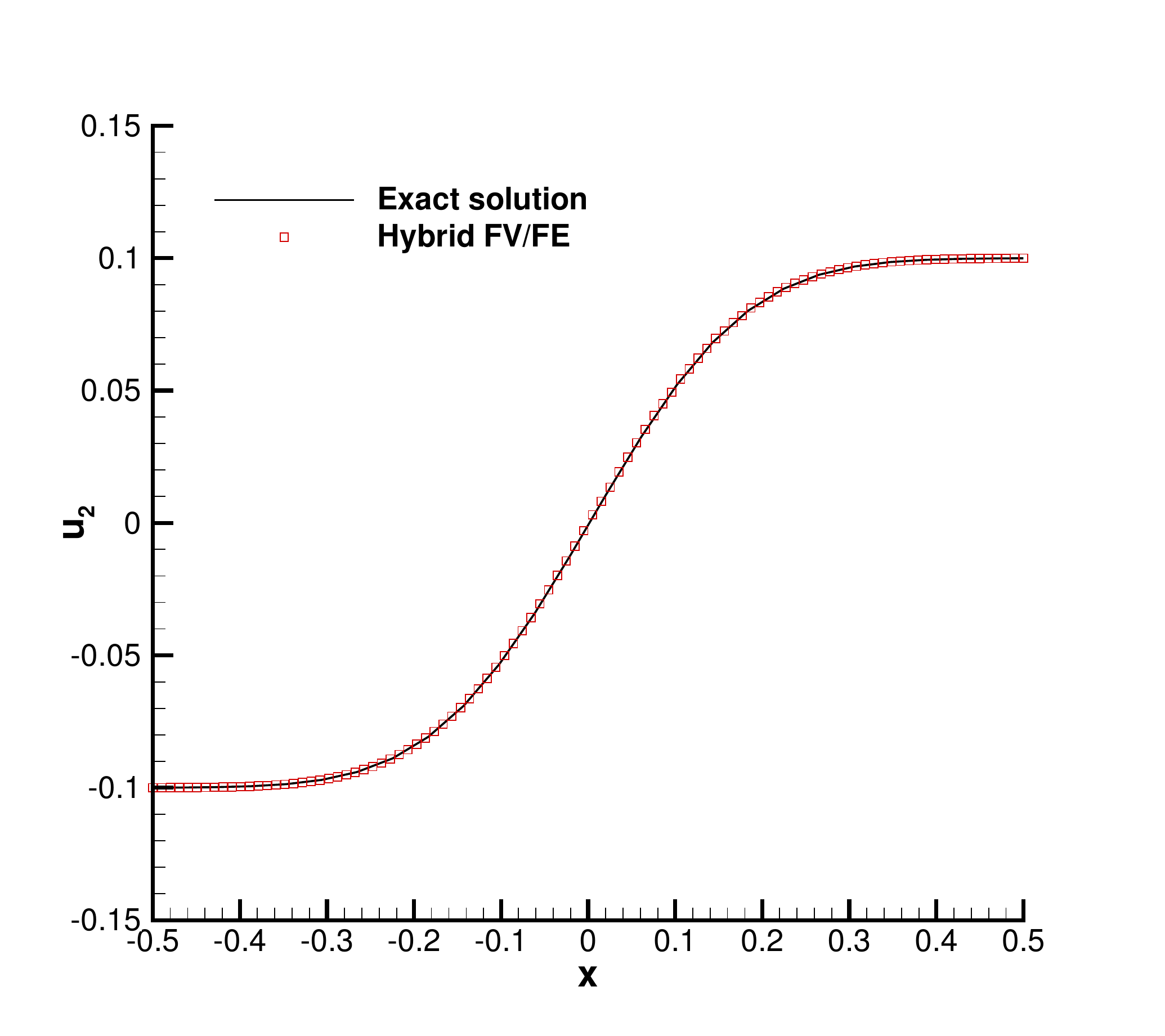}
	\caption{Left: Mesh considered for the first problem of Stokes in 3D. Right: Velocity along the cut $y=0$, $z=0$ for the first problem of Stokes in 3D at time $t=1$. Comparison between the values computed with the non well-balanced scheme, considering the pressure correction, and the exact solution~\eqref{eqn.FPS_ex} for $\mu=10^{-2}$.}
	\label{fig:FPS3D}
\end{figure}

We also consider the first problem of Stokes in 3D, where  $\Omega=[-0.5,0.5]\times[-0.05,0.05]\times[-0.05,0.05]$. Again, the simulations are performed with the non well-balanced scheme with the pressure correction approach. The left plot of Figure~\ref{fig:FPS3D} shows the grid considered, formed by hexahedra, tetrahedra, and pyramids, while the right plot reports the numerical results for $\mu = 10^{-2}$. More precisely, the 1D cut along $y=0$ and $z=0$ of the solution at time $t=1$ is shown in comparison with the exact solution given by~\eqref{eqn.FPS_ex}. 

\begin{figure}[h]
	\centering
	\includegraphics[width=0.325\linewidth]{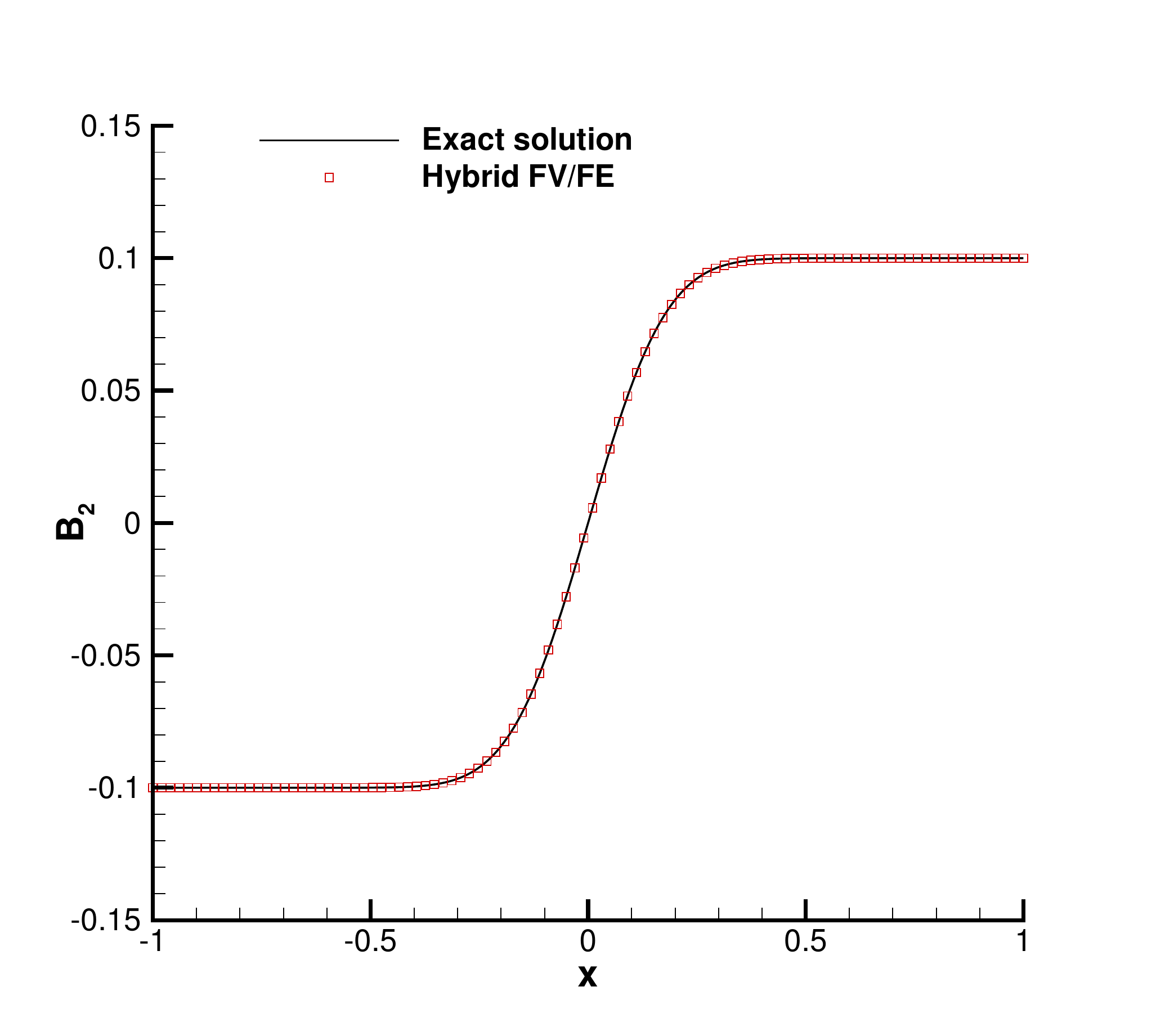}
	\includegraphics[width=0.325\linewidth]{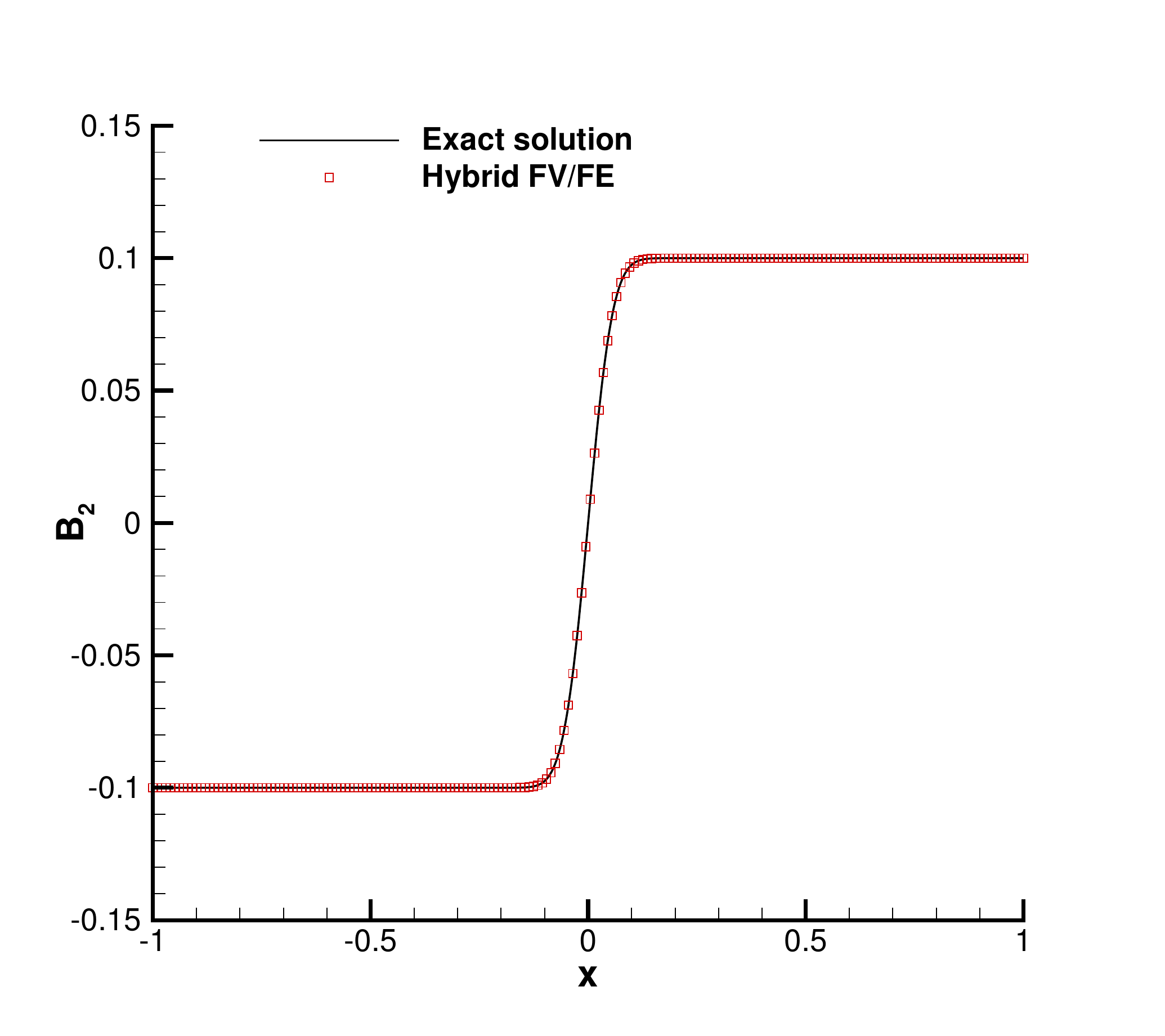}
	\includegraphics[width=0.325\linewidth]{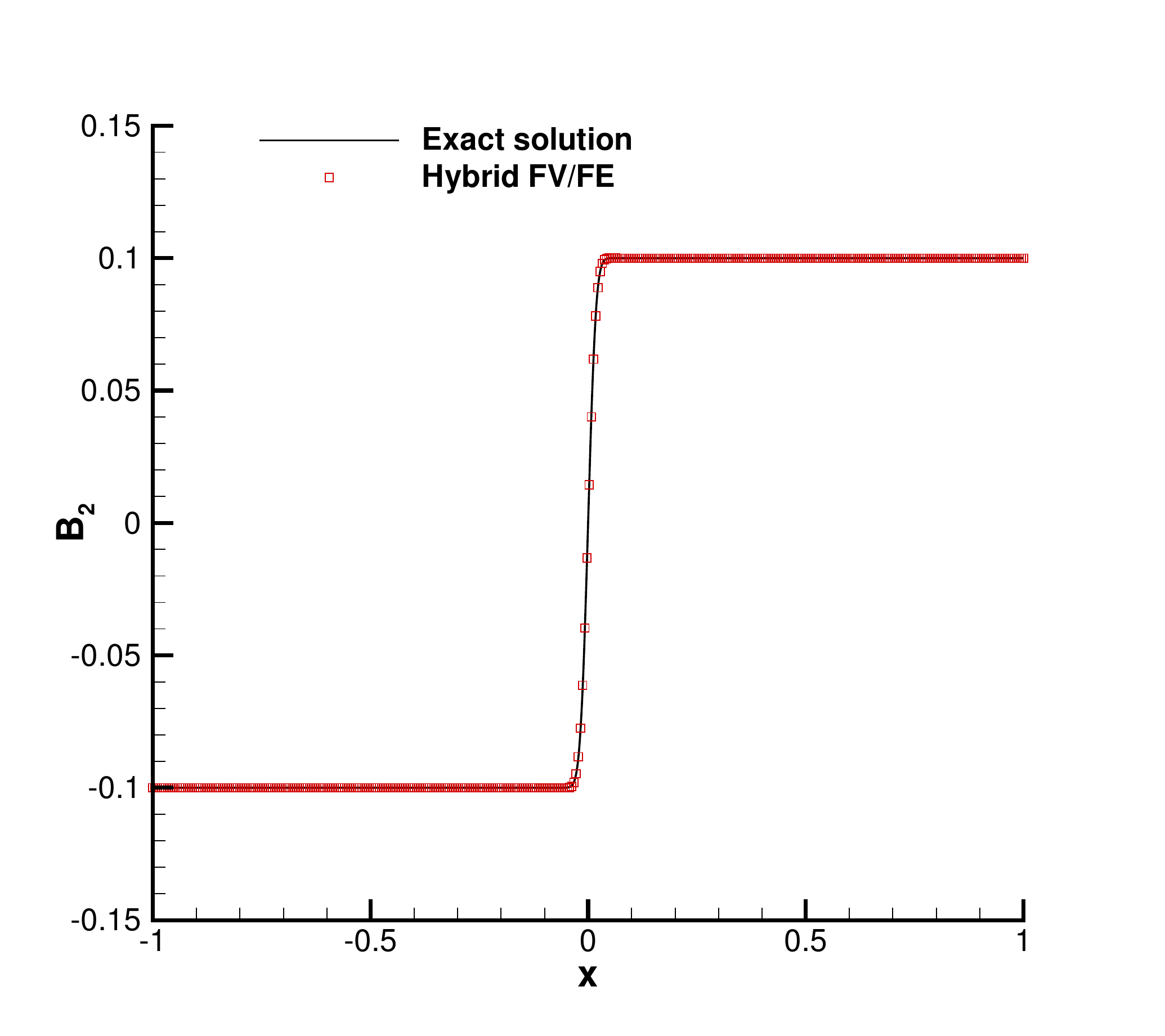}
	\caption{Magnetic field  along the cut $y=0$ for the current sheet problem in 2D at time $t=1$. Comparison between the values computed with the proposed methodology and the exact solution~\eqref{eqn.FPS_ex} with $\eta=10^{-2}$ (left), $\eta=10^{-3}$ (center), and $\eta=10^{-4}$ (right).}
	\label{fig:CS2D}
\end{figure}

As second test within this section, we consider a current sheet problem, similar to the first problem of Stokes, but for the magnetic field. The initial condition is given by
\begin{equation*}
	\press \left(\xx,0\right) = 1, \qquad \bvel\left(\xx,0\right) = {0}, \qquad 
	{B}_{1} \left(\xx,0\right) = 0, \qquad
	{B}_{2} \left(\xx,0\right) = \left\lbrace \begin{array}{lc}
		-0.1 & \mathrm{ if } \; x \le 0,\\
		\phantom{-} 0.1 & \mathrm{ if } \; x > 0.
	\end{array}\right.
\end{equation*}
The exact solution for the magnetic field in the $y$-direction is the right-hand side of \eqref{eqn.FPS_ex} with $\eta$ instead of $\mu$. We solve the problem in the computational domain $\Omega = [-1,1] \times [-0.1,0,1]$ for different magnetic resistivities: $\eta \in\{10^{-2}, 10^{-3}, 10^{-4}\}$. We consider the grid shown in Figure \ref{fig:FPS2D_mixmesh}, and impose periodic boundary conditions in the $y$-direction and Dirichlet boundary conditions in the $x$-direction. The results of this test are shown in Figure \ref{fig:CS2D}. 

\begin{figure}[h]
		\centering
		\includegraphics[width=0.5\linewidth]{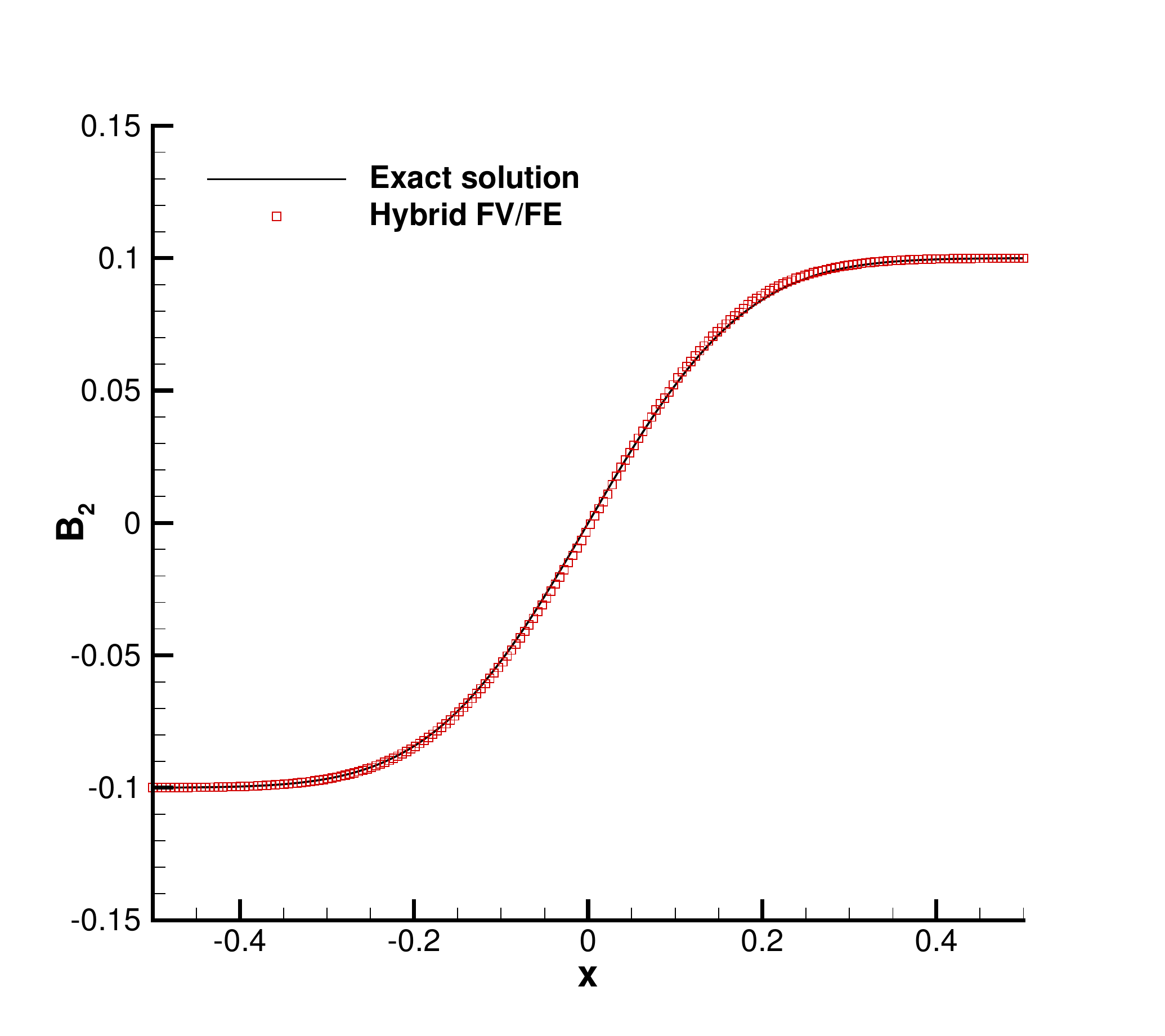}
		\caption{Magnetic field along the cut $y=0$, $z=0$ for the current sheet problem in 3D at time $t=1$. Comparison between the values computed with the proposed methodology and the exact solution~\eqref{eqn.FPS_ex} with $\eta=10^{-2}$.}
		\label{fig:CS3D}
\end{figure}

The test is repeated in the three-dimensional domain $\Omega = [-0.5,0.5]\times[-0.05,0.05]\times[-0.05,0.05]$ discretized with the mesh shown in Figure \ref{fig:FPS3D} on the left. For this test, we set $\eta = 10^{-2}$ and we impose Dirichlet boundary conditions in the $x$-direction and periodic boundary conditions in the $y$- and $z$-directions. The obtained results are shown in Figure \ref{fig:CS3D}. 
A very good agreement between the numerical and the exact solutions is observed for both the two and three-dimensional simulations.  

\subsection{Magnetic field loop advection in 2D and 3D}  

\begin{figure}[h]
 	\centering
 	\includegraphics[width=0.9\linewidth]{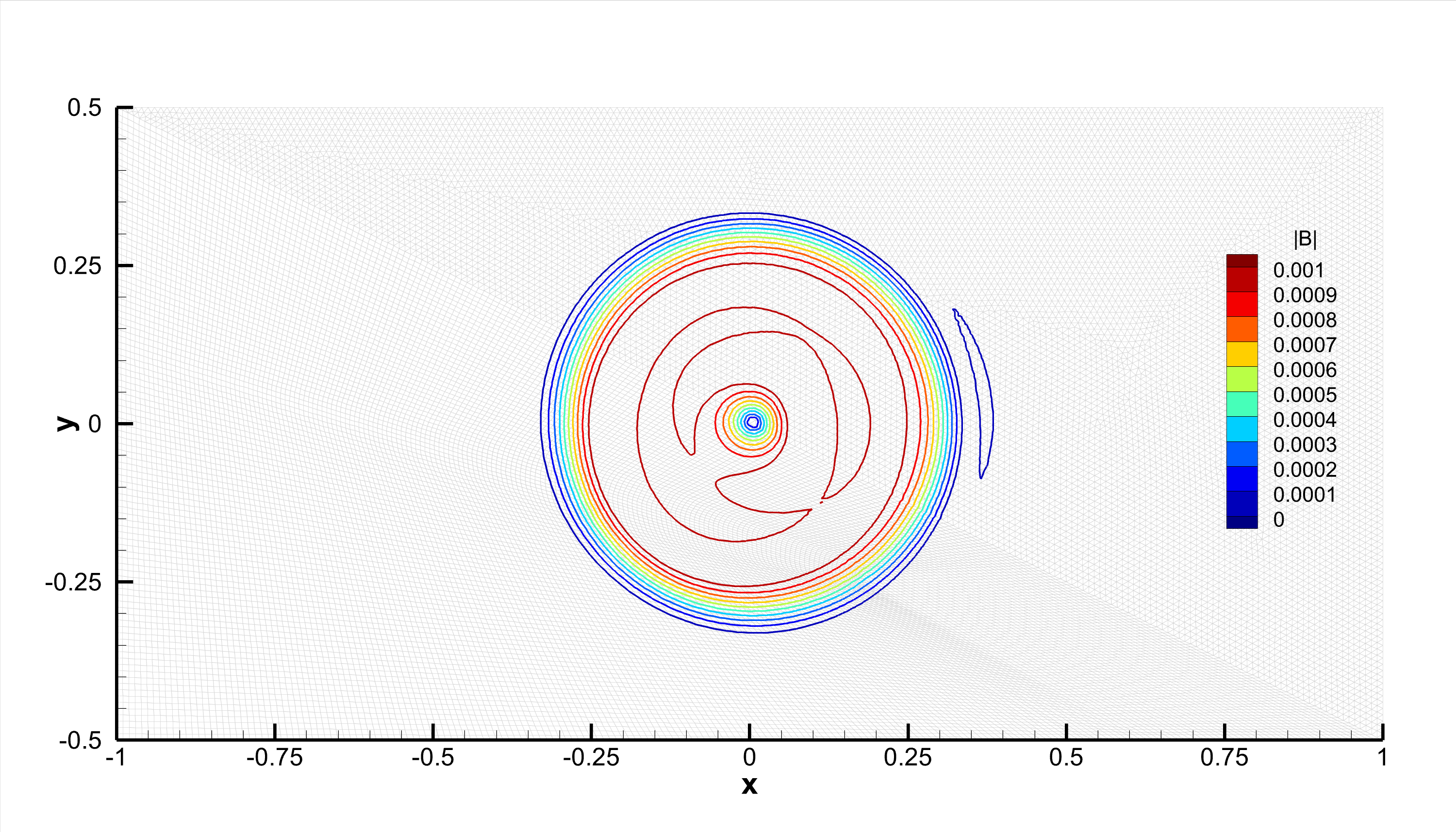}
 	\caption{Contours lines of the magnetic field strength at time $t = 1$ for the magnetic field loop advection in 2D.}
 	\label{fig:MHDloop2d}
 \end{figure}

We consider the classical magnetic field loop advection problem, see \cite{GardinerStone,divfree2015,SIMHD}. The initial condition is given by 
\begin{equation*}
	\begin{split}
	&
	\press \left(\xx,0\right) = 1, \qquad \vel_1\left(\xx,0\right) = 2, \qquad \vel_2\left(\xx,0\right) = 1, \\
	&A_{3}\left(\xx,0\right) = \left\lbrace \begin{array}{lc}
		 10^{-3}\cdot(0.3 - r\left(\xx\right))  & \mathrm{ if } \; r(\xx) \le 0.3,\\
		 0 & \mathrm{ if } \; r(\xx) > 0.3,
	\end{array}\right.
	\end{split}
\end{equation*}
where $r(\xx) = \sqrt{x^2 + y^2}$ and $A_{3}(\xx,0)$ is the $z$-component of the vector potential $\magpotential$, such that $\magfield = \nabla\times\mathbf{A}$. The domain $\Omega = [-1,1]\times [-0.5,0.5]$ is discretized with the unstructured mixed-element grid shown in Figure~\ref{fig:MHDloop2d}. Periodic boundary conditions are imposed everywhere. 
%

We run the simulation until $t=1$ in order to complete one advection period. The numerical solution at the end time, shown in Figure~\ref{fig:MHDloop2d}, is in good qualitative agreement with previous results reported in the literature, e.g. \cite{GardinerStone,divfree2015,SIMHD}. 
 
 \begin{figure}[h]
 	\centering
 	\includegraphics[width=0.53\linewidth]{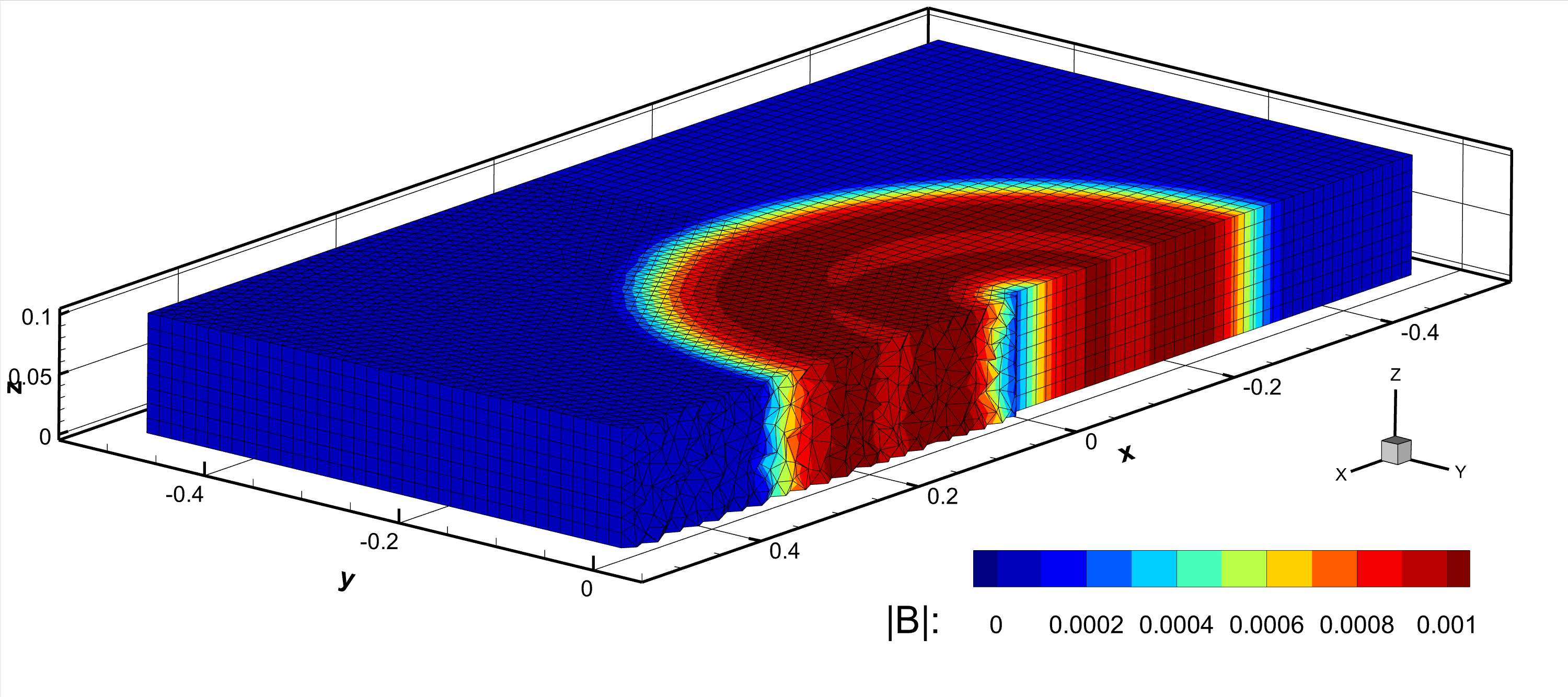} 
 	\includegraphics[width=0.43\linewidth]{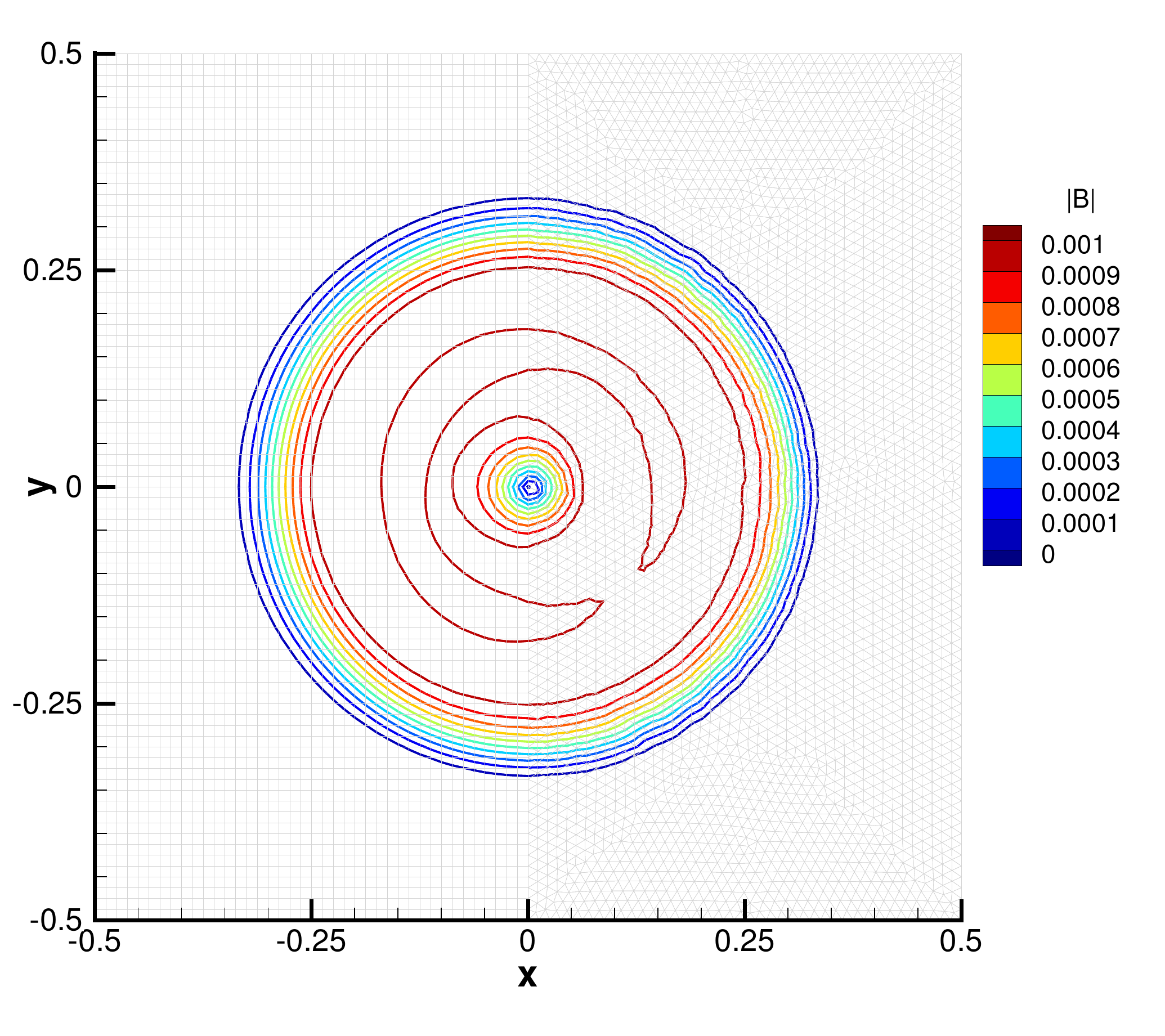}
 	\caption{Numerical results for the MHD advection loop in a 3D domain. Left: sketch of the mixed structured-unstructured mesh. Right: contour lines of the magnetic field strength at time $t = 1$.}
 	\label{fig:MHDloop3d}
 \end{figure}

\subsection{Lid-driven cavity} 
\label{sec:LDC}

We now solve the classical lid-driven cavity problem for the incompressible Navier-Stokes equations, i.e., without magnetic field, see~\cite{GGS82} for more details. This test is classically used to validate incompressible flow solvers and compressible ones in the low Mach number regime \cite{TD17,Hybrid1,Hybrid2}. 
We consider the two-dimensional computational domain $\Omega=[-0.5,0.5]\times[-0.5,0.5]$. The fluid is at rest at the initial time, $\bvel={0}$, and the lid velocity at the top boundary is chosen as $\bvel=(1,0)^{T}$. No-slip walls are set on the remaining boundaries, and the viscosity parameter is chosen as $\mu=10^{-2}$, which yields a Reynolds number equal to $100$ based on the size of the domain and the lid velocity.

\begin{figure}[!htbp]
	\begin{center}
		\begin{tabular}{cc} 
			\includegraphics[trim=20 20 20 20,clip,width=0.47\textwidth]{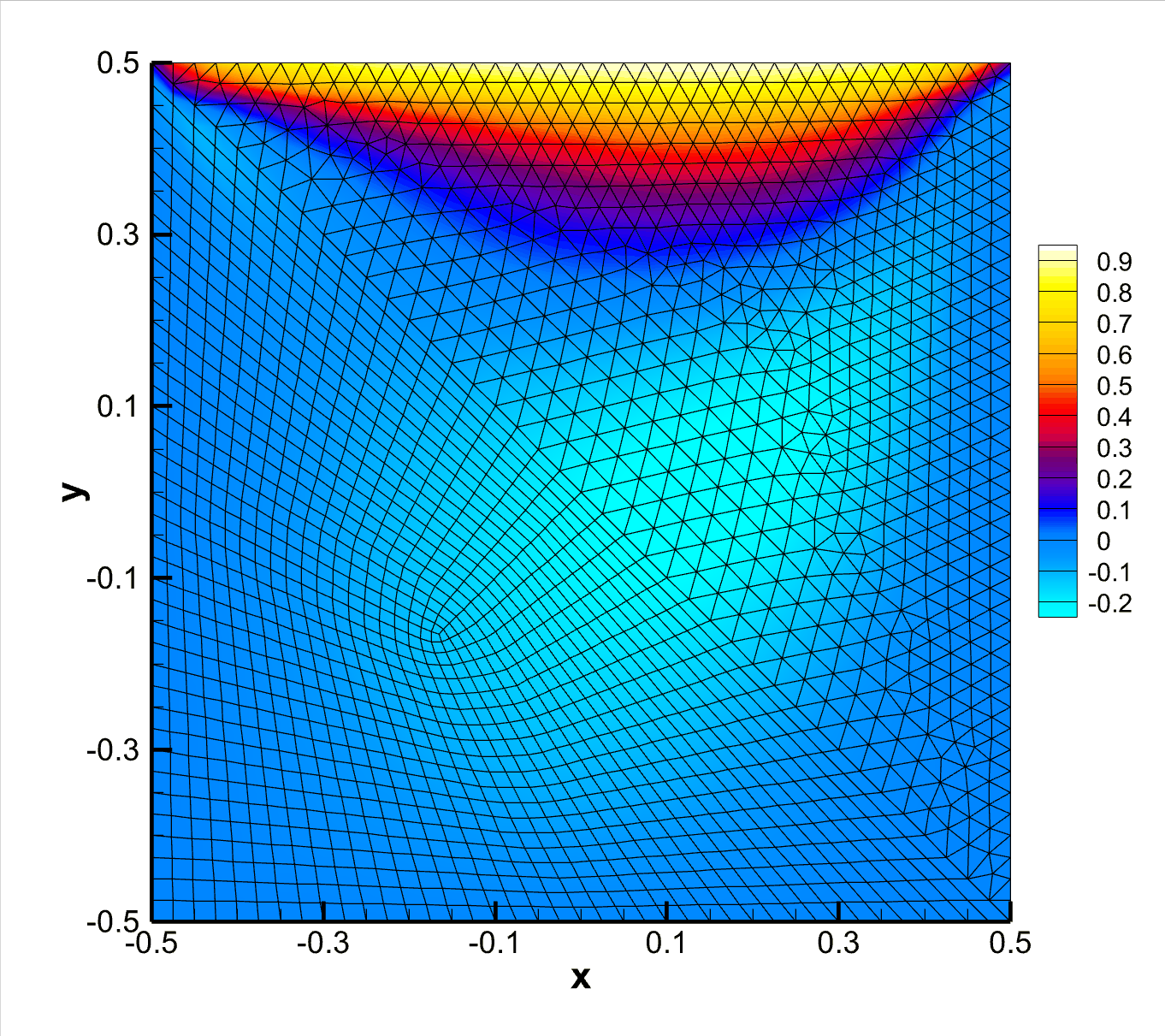}  & 
			\includegraphics[trim=20 20 20 
			20,clip,width=0.47\textwidth]{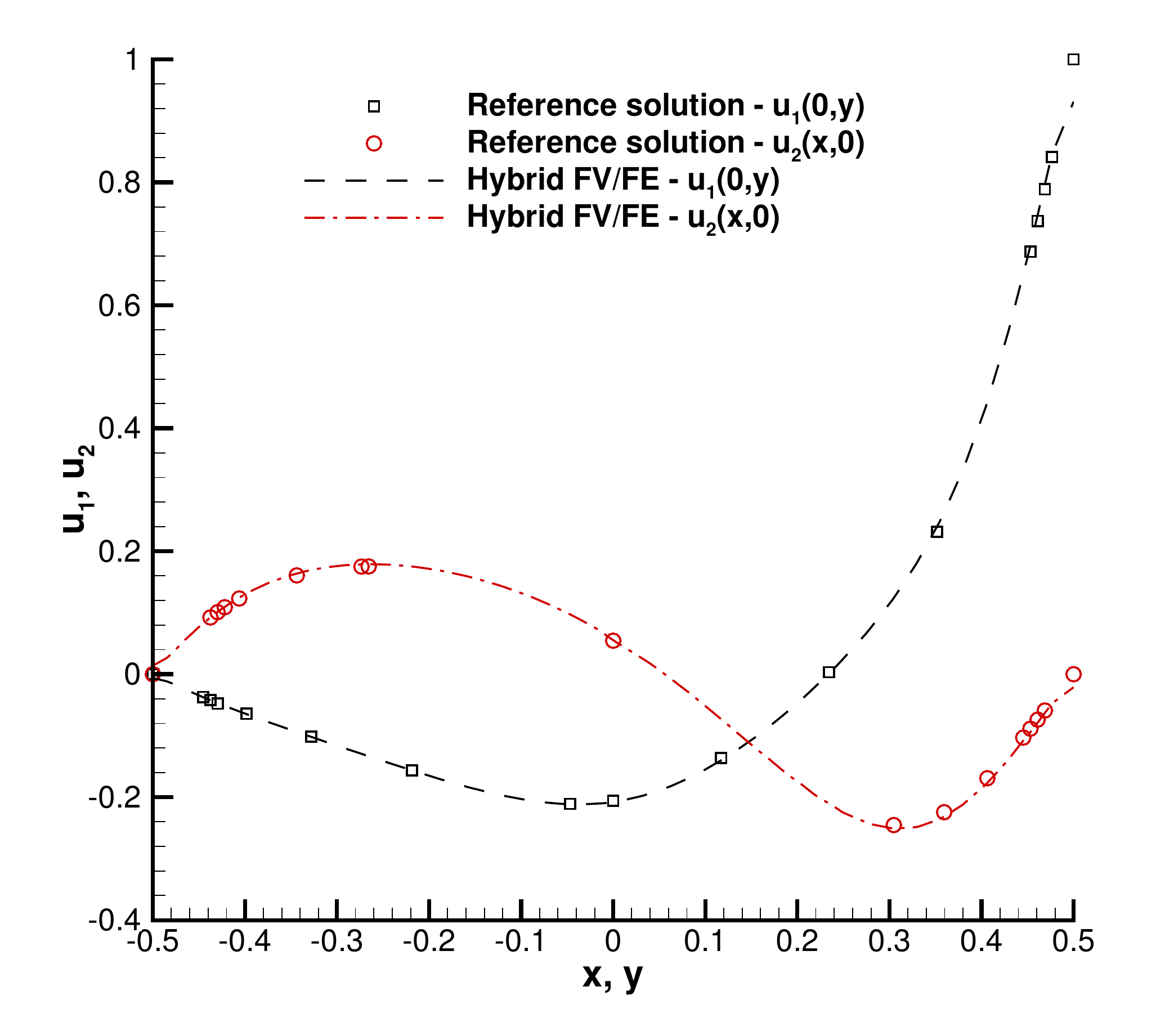}   
		\end{tabular} 
		\caption{Lid-driven cavity at Reynolds number $Re=100$ at time $t=10$. Contours of $\vel_1$ (left) and 1D cuts of the velocity along the $x$- and $y$-axis (right) in comparison with the reference solution~\cite{GGS82}.} 
		\label{fig:LDC}
	\end{center}
\end{figure}

The computational domain has been discretized with an unstructured mixed-element grid composed of triangular and quadrilateral elements, more specifically using M$_{40}$ (see the left plot of Figure~\ref{fig:mixmesh}). 
This simulation is run with the non well-balanced scheme using the pressure correction approach. The left graph in Figure~\ref{fig:LDC} reports the values of $u_1$ at time $t=10$, displayed over the mesh M$_{40}$. The right graph in Figure~\ref{fig:LDC} shows cuts for both velocity components along the $x$ and $y$ axes, respectively, compared to the reference solution provided in \cite{GGS82}. We can observe an excellent agreement between the reference solution and the numerical results obtained with the new hybrid FV/FE method proposed in this paper.

\subsection{Lid-driven cavity with magnetic field} \label{sec:MHD_LDC}
\begin{figure}
	\centering
	\includegraphics[width=0.45\linewidth]{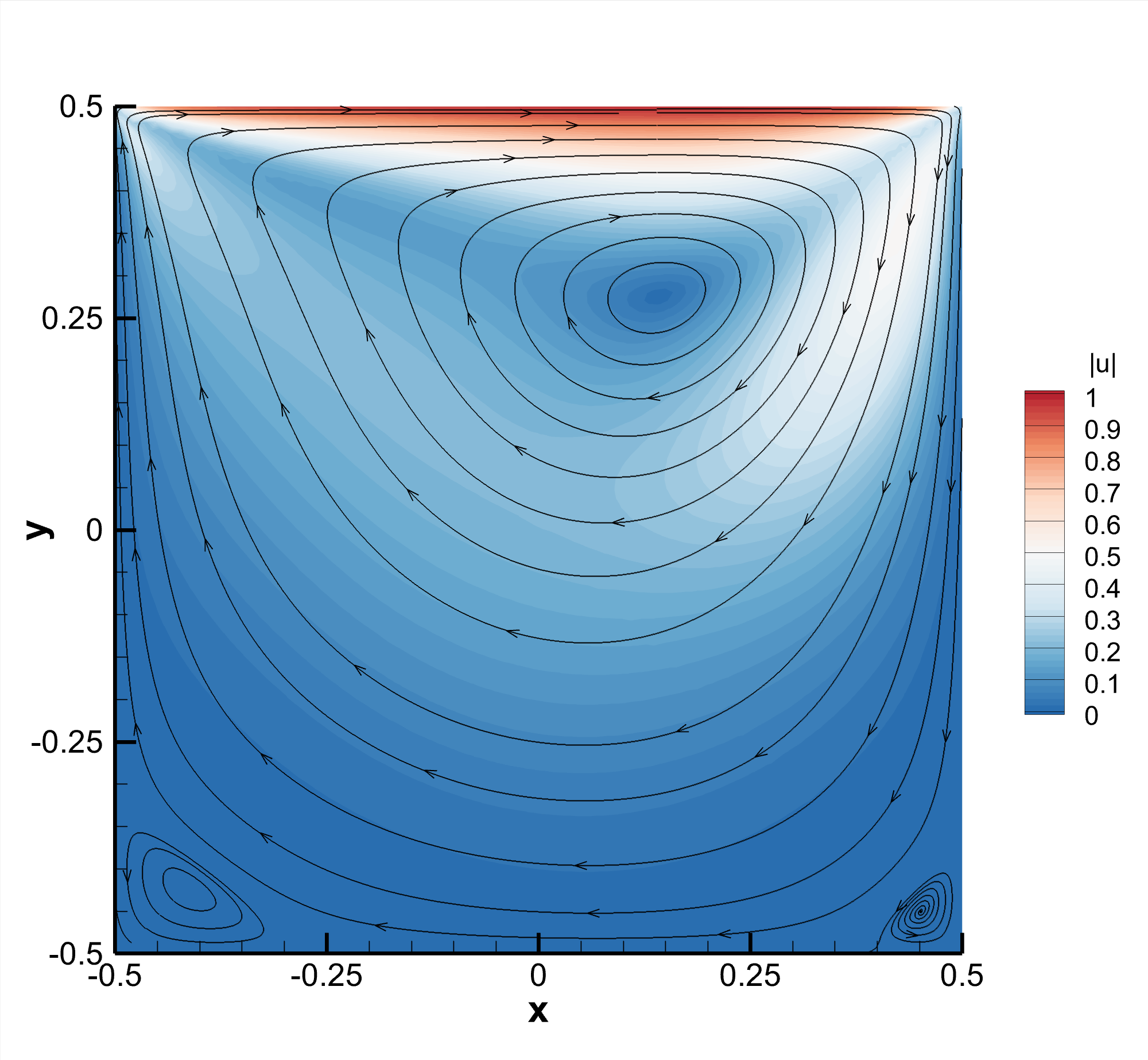} 
	\includegraphics[width=0.45\linewidth]{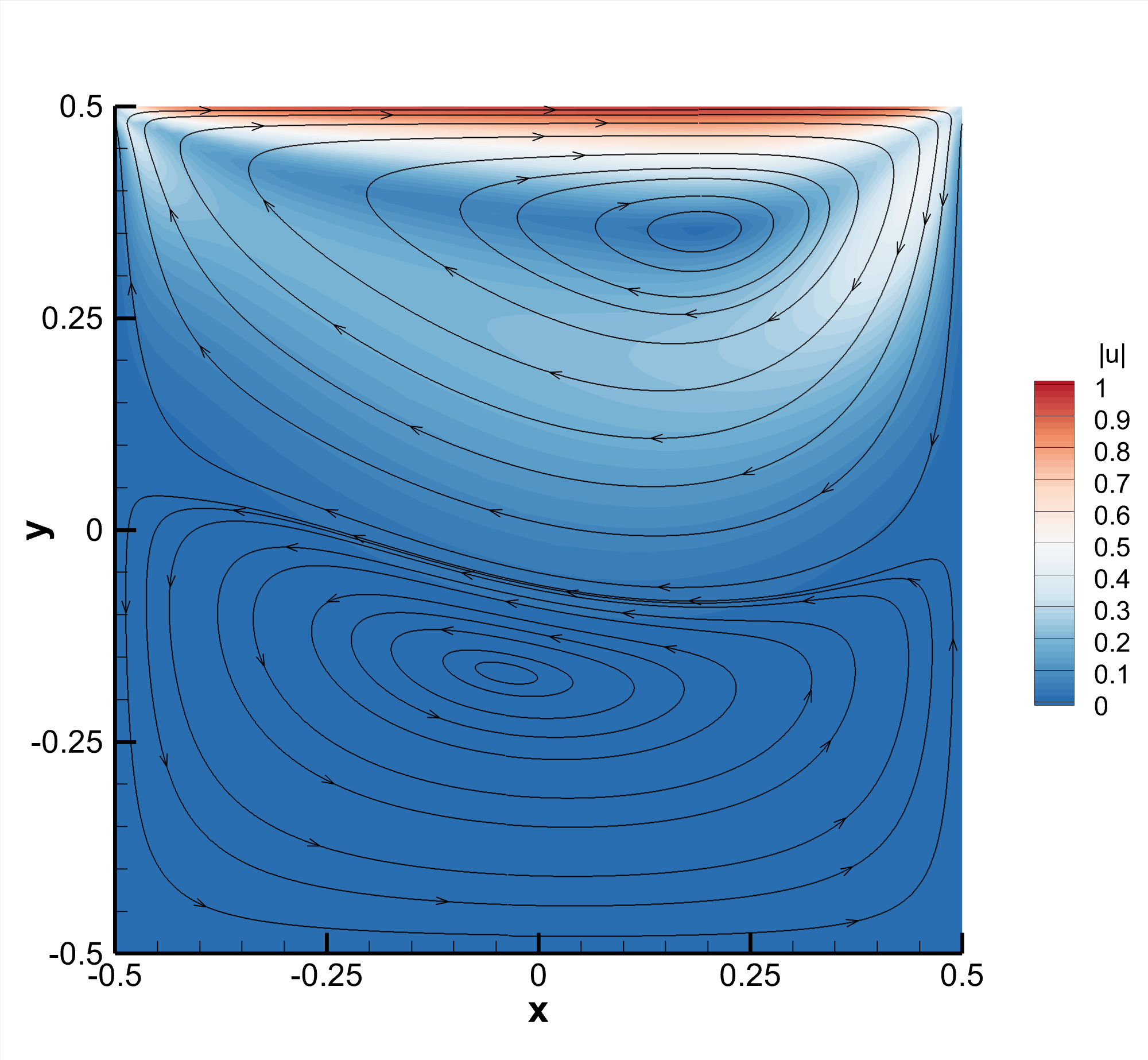} \\
	\includegraphics[width=0.45\linewidth]{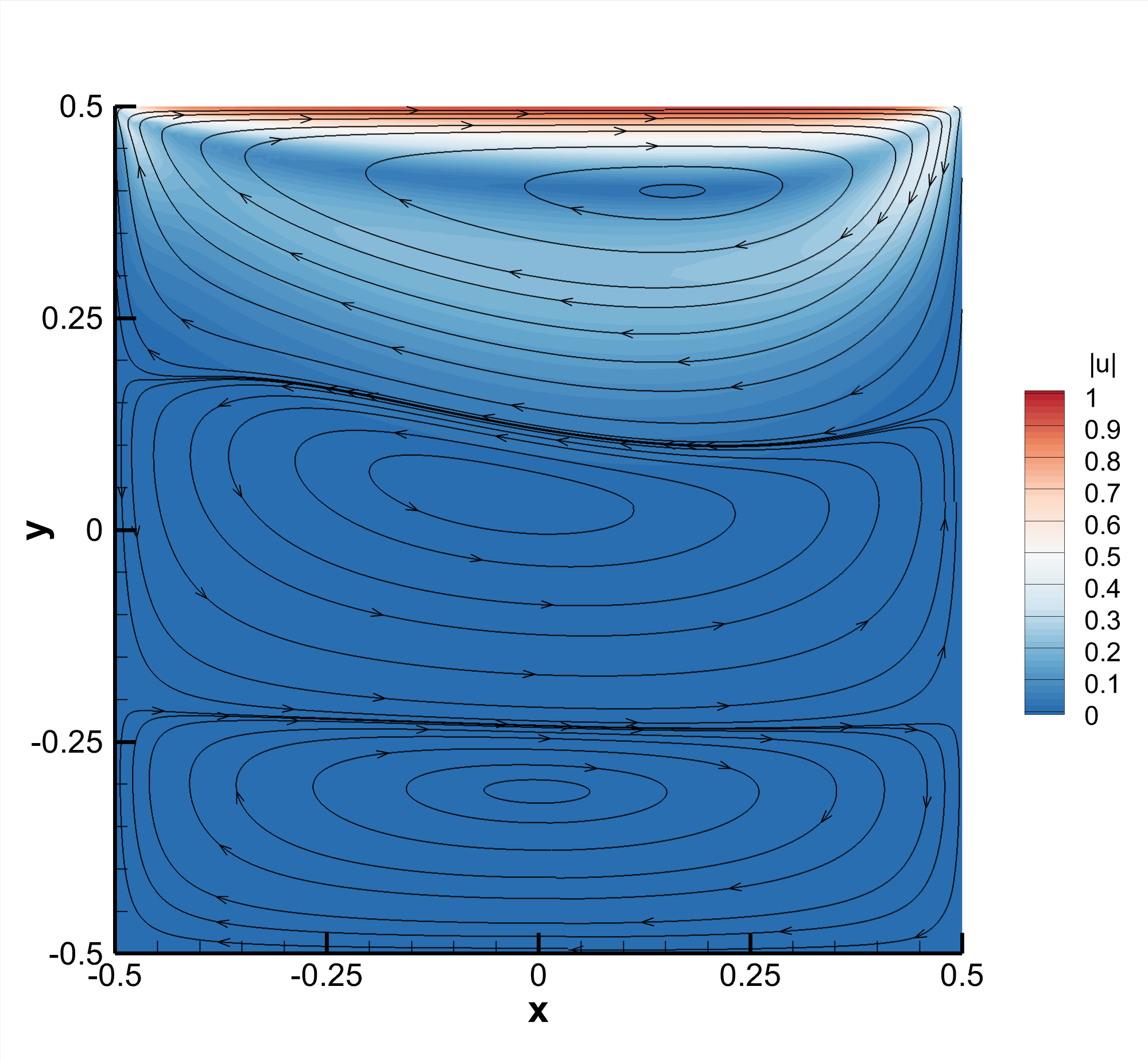} 
	\includegraphics[width=0.45\linewidth]{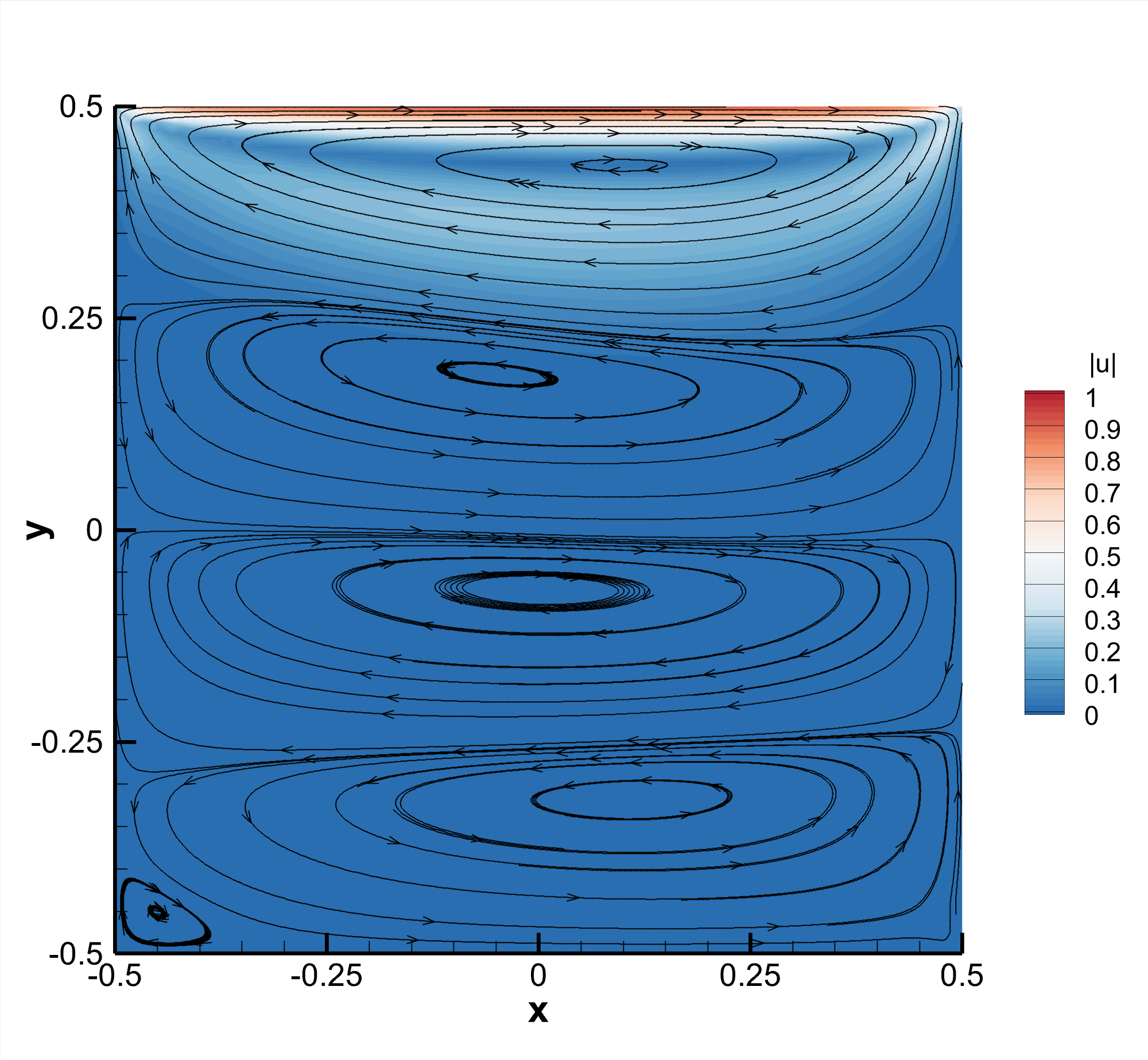} 
	\caption{Stream-lines of the velocity field at time $t=20$ for the MHD lid-driven cavity with horizontal magnetic field $\magfield_0 = (B_{0,x},0)^{T}$ for $B_{0,x} \in\left\lbrace 0.1, 0.25,0.5,1\right\rbrace$ (top-left, top-right, bottom-left, and bottom-right, respectively).}  
	\label{fig.MHDLDC_horizontal_modV}
\end{figure}

\begin{figure}
	\centering
	\includegraphics[width=0.45\linewidth]{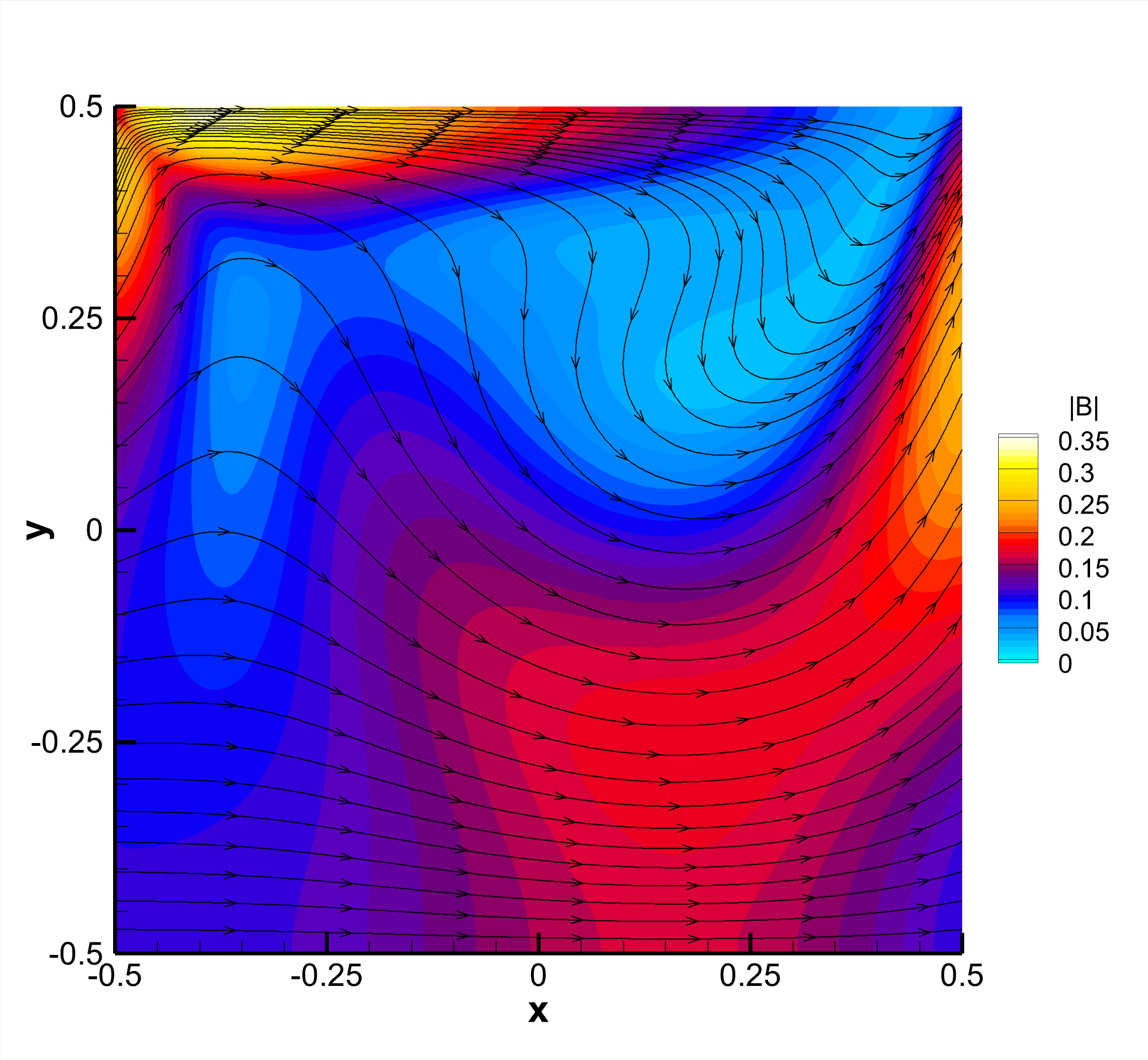} 
	\includegraphics[width=0.45\linewidth]{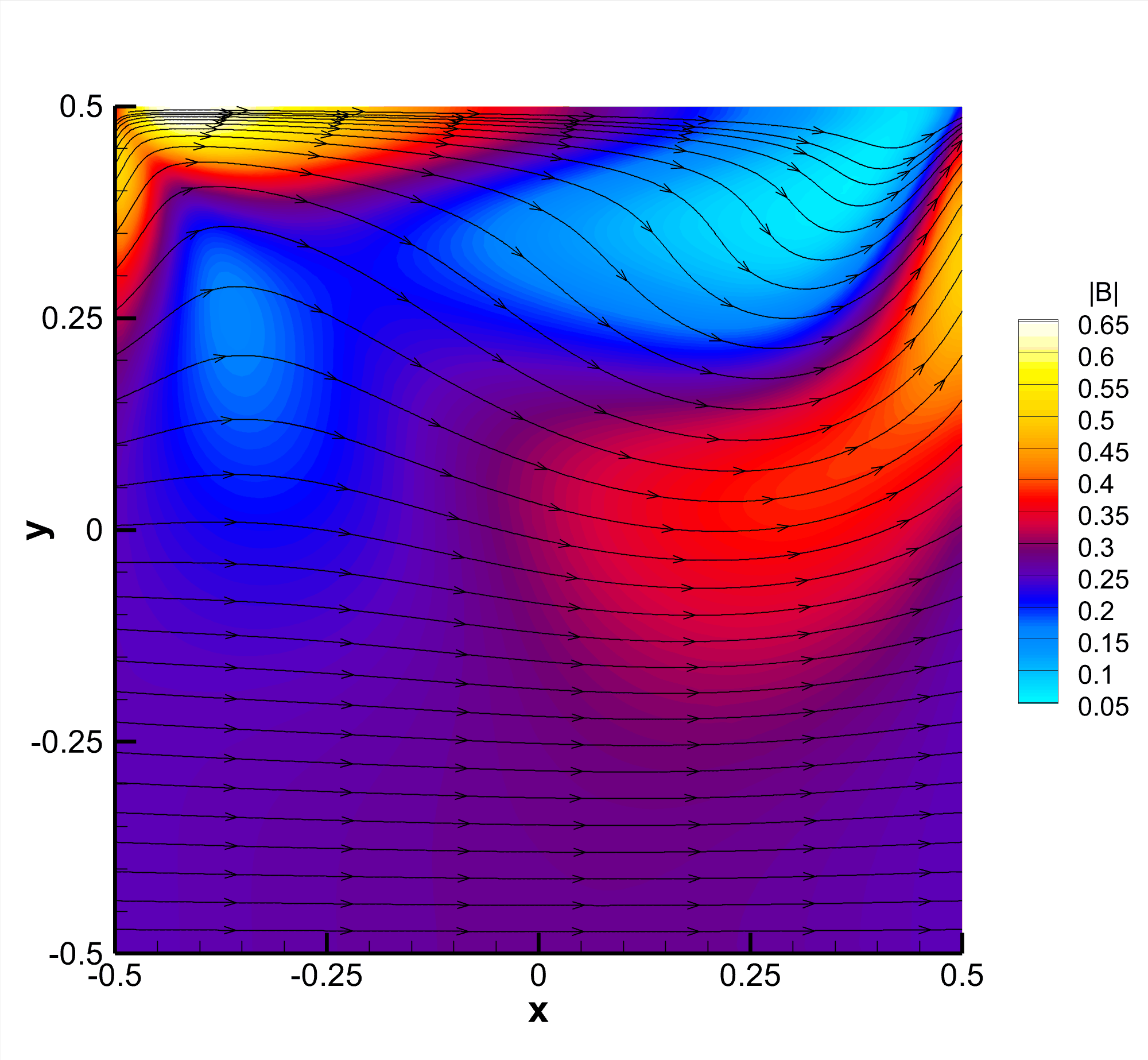} \\ 
	\includegraphics[width=0.45\linewidth]{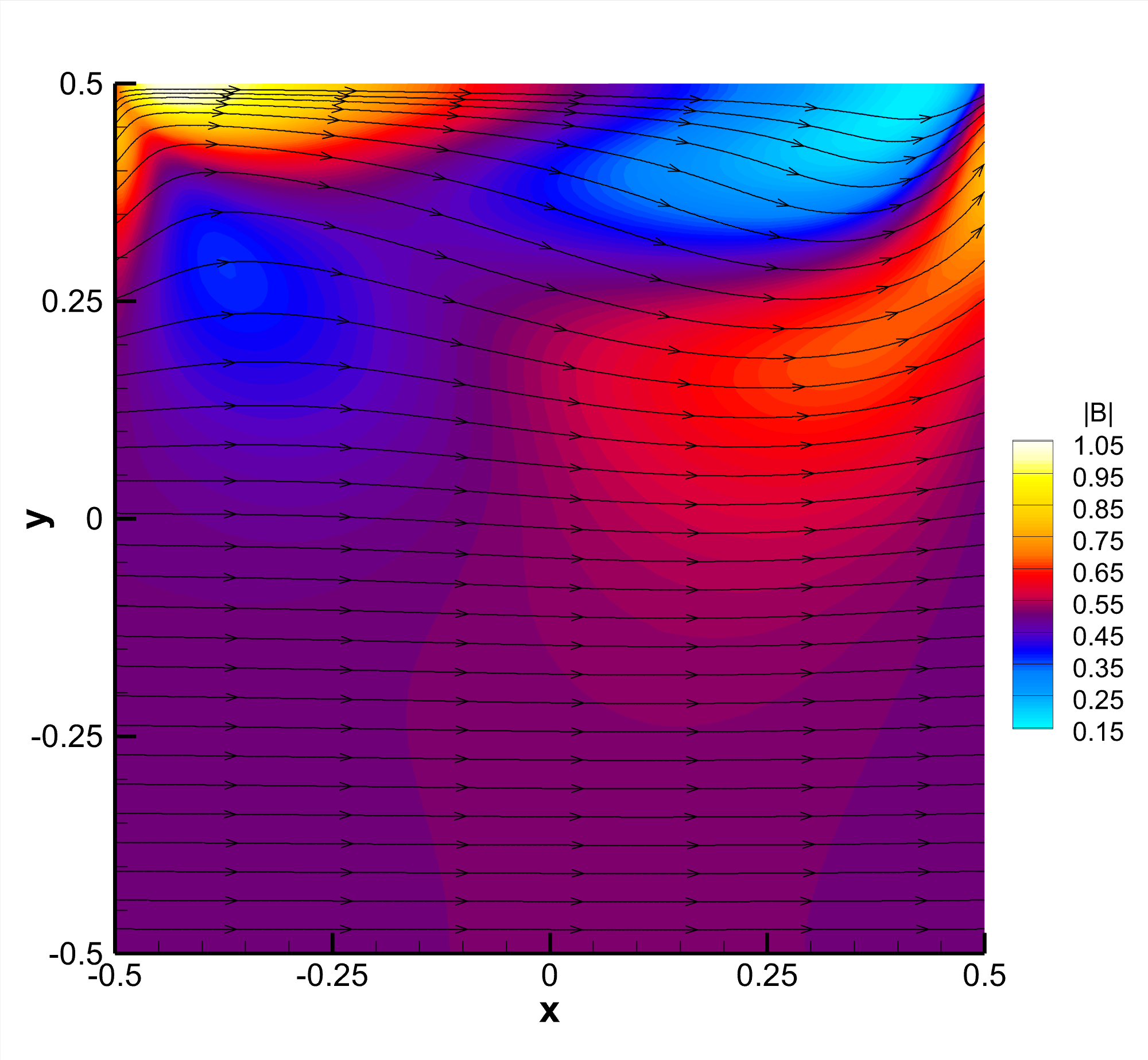} 
	\includegraphics[width=0.45\linewidth]{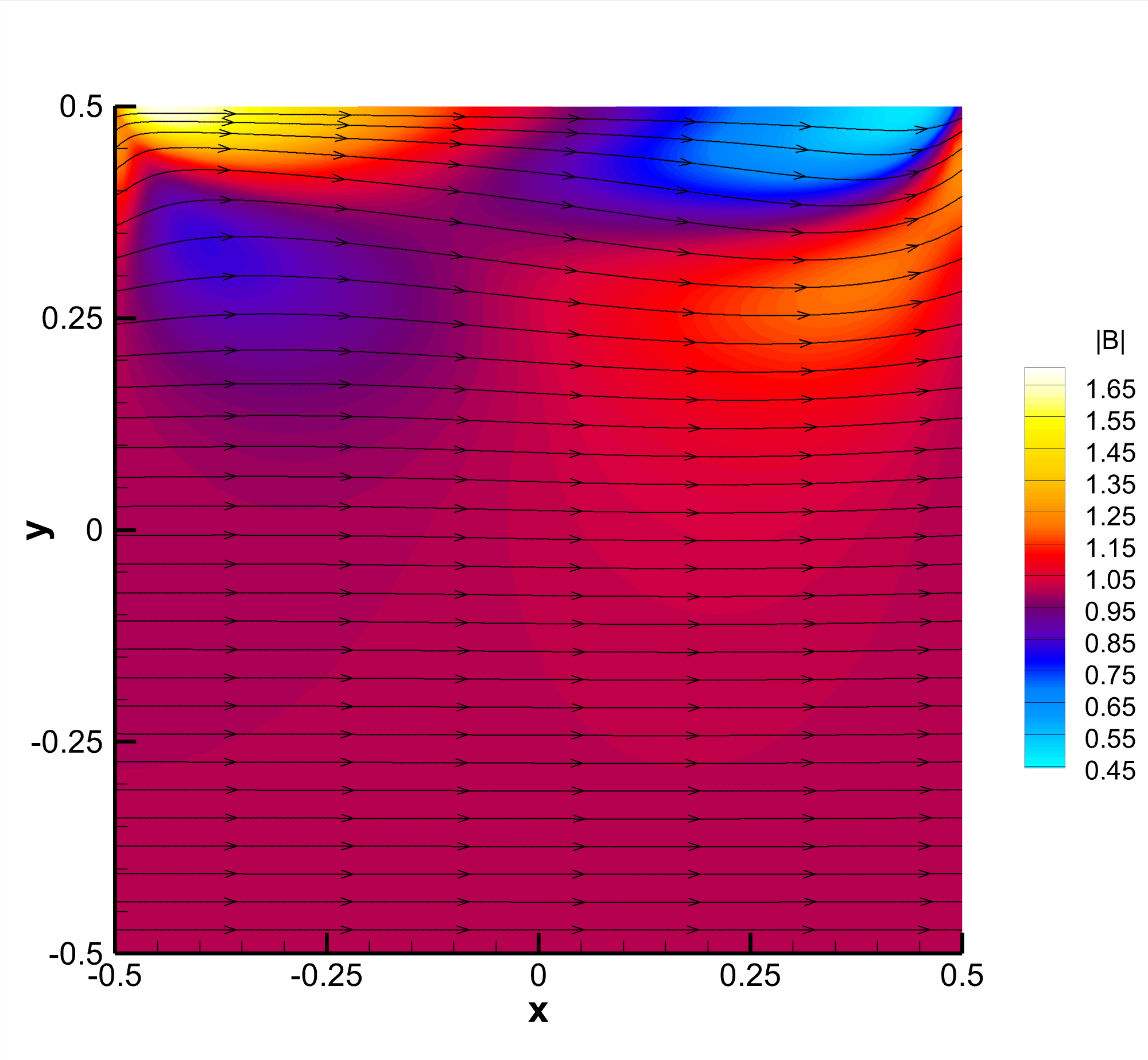} 
	\caption{Stream-lines of the magnetic field at time $t=20$ for the MHD lid-driven cavity with horizontal magnetic field $\magfield_0 = (B_{0,x},0)^{T}$ for $B_{0,x}  \in\left\lbrace 0.1, 0.25,0.5,1\right\rbrace$ (top-left, top-right, bottom-left, and bottom-right, respectively).}  
	\label{fig.MHDLDC_horizontal_modB}
\end{figure}

\begin{figure}
	\centering
	\includegraphics[width=0.45\linewidth]{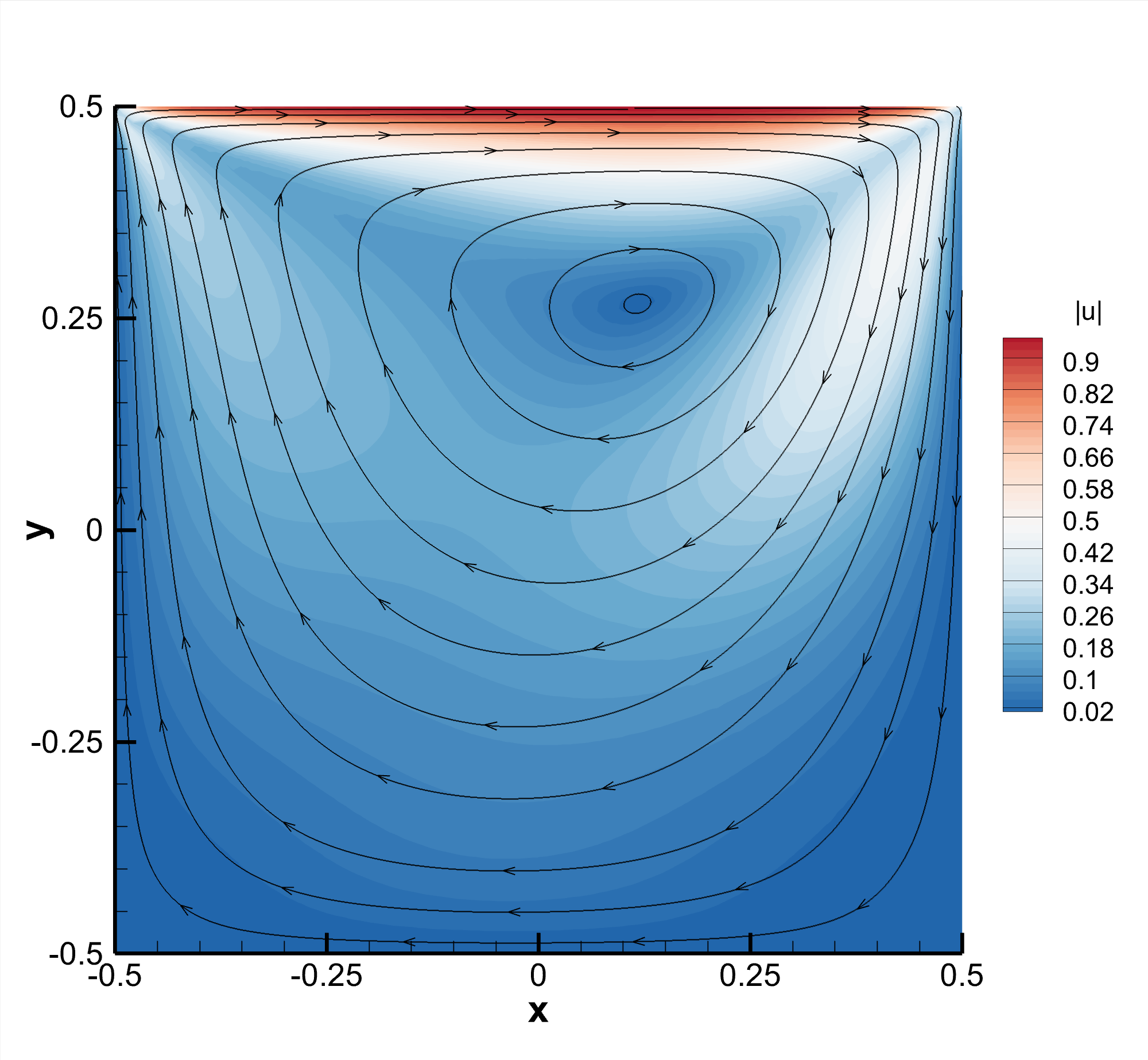}
	\includegraphics[width=0.45\linewidth]{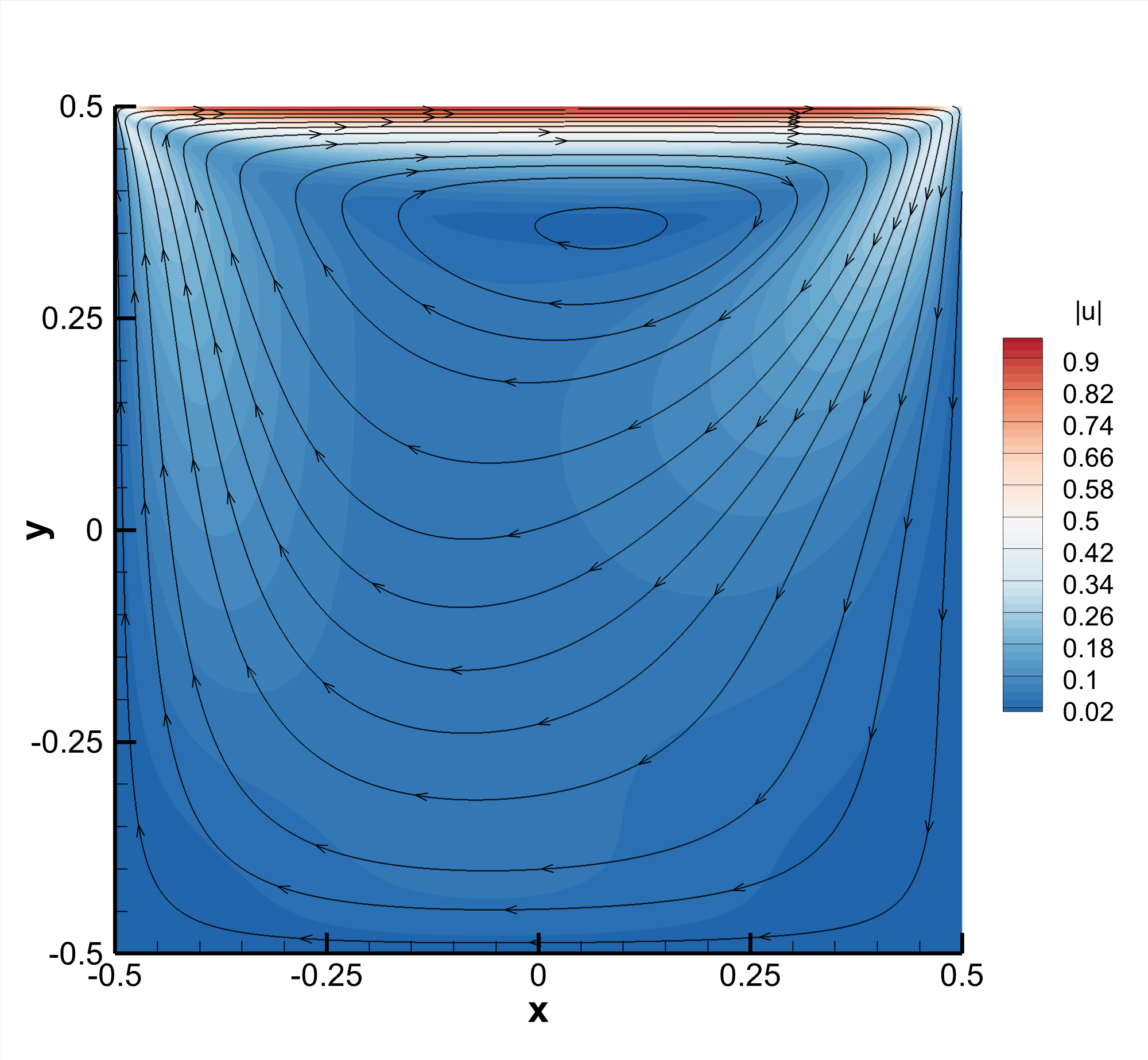} \\
	\includegraphics[width=0.45\linewidth]{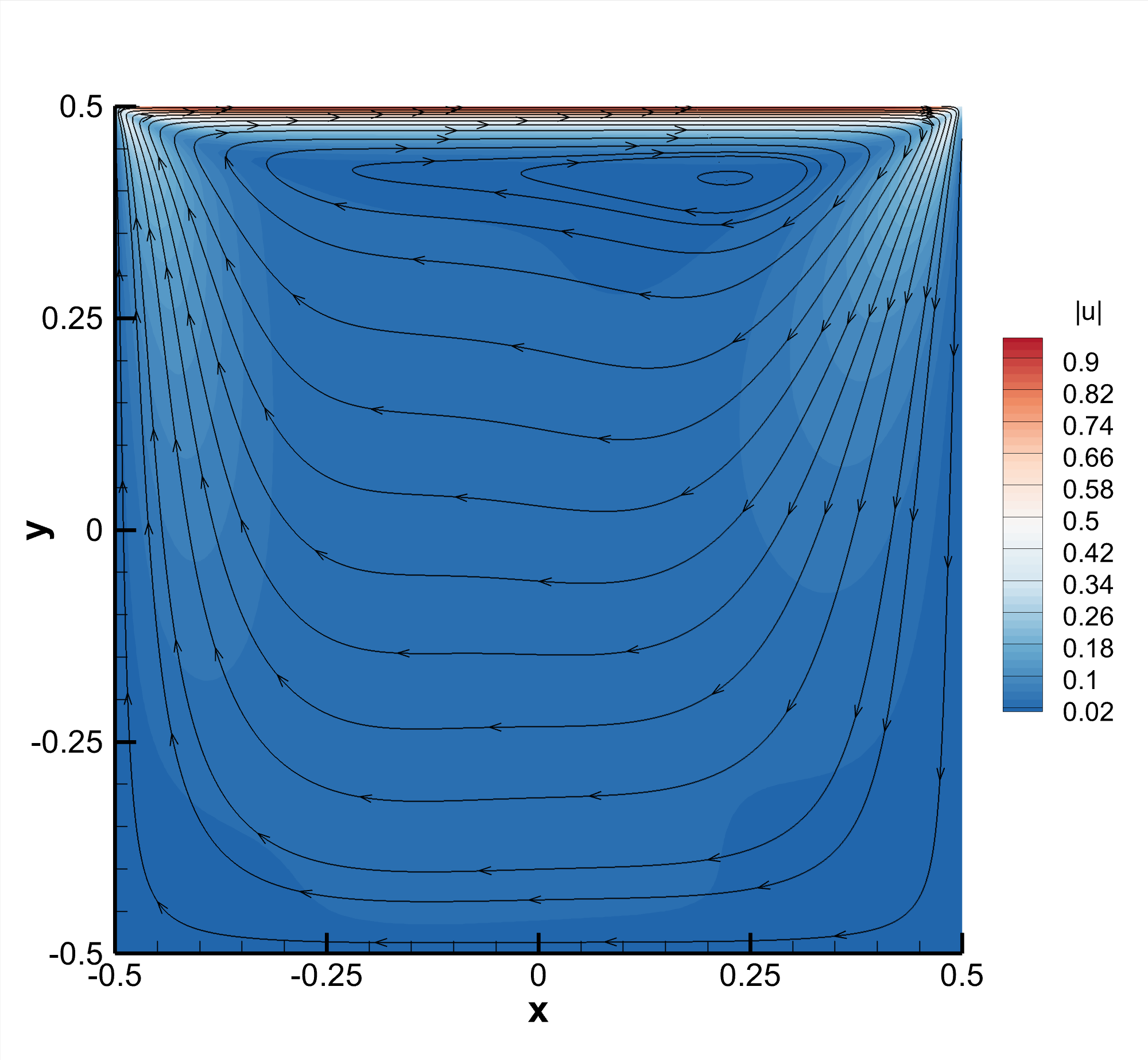}
	\includegraphics[width=0.45\linewidth]{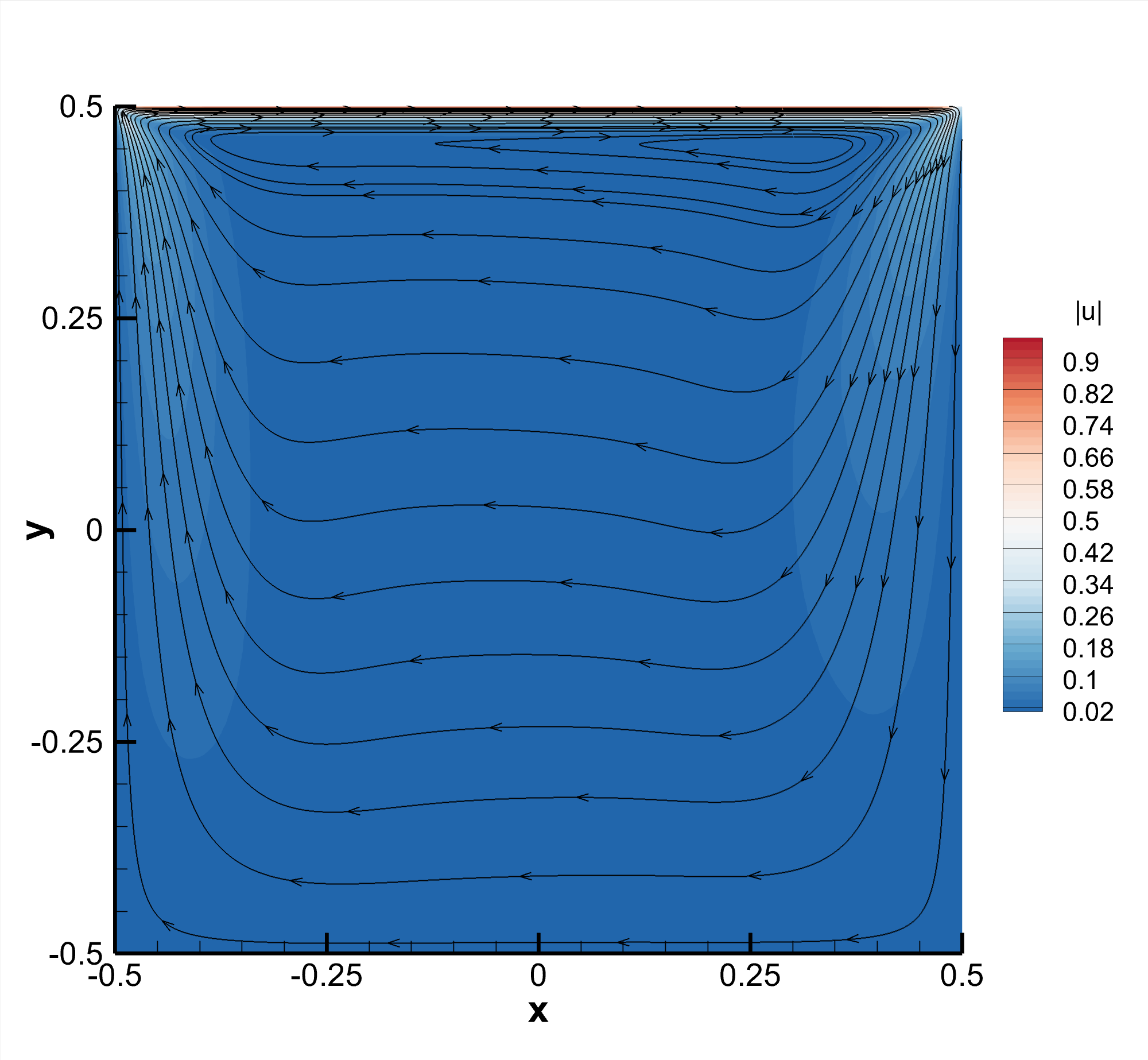} 
	\caption{Stream-lines of the velocity field at time $t=20$ for the MHD lid-driven cavity with vertical magnetic field $\magfield_0 = (0,B_{0,y})^{T}$ for $B_{0,y} \in\left\lbrace 0.1, 0.25,0.5,1\right\rbrace$ (top-left, top-right, bottom-left, and bottom-right, respectively).}  
	\label{fig.MHDLDC_vertical_modV}
\end{figure}

\begin{figure}
	\centering
	\includegraphics[width=0.45\linewidth]{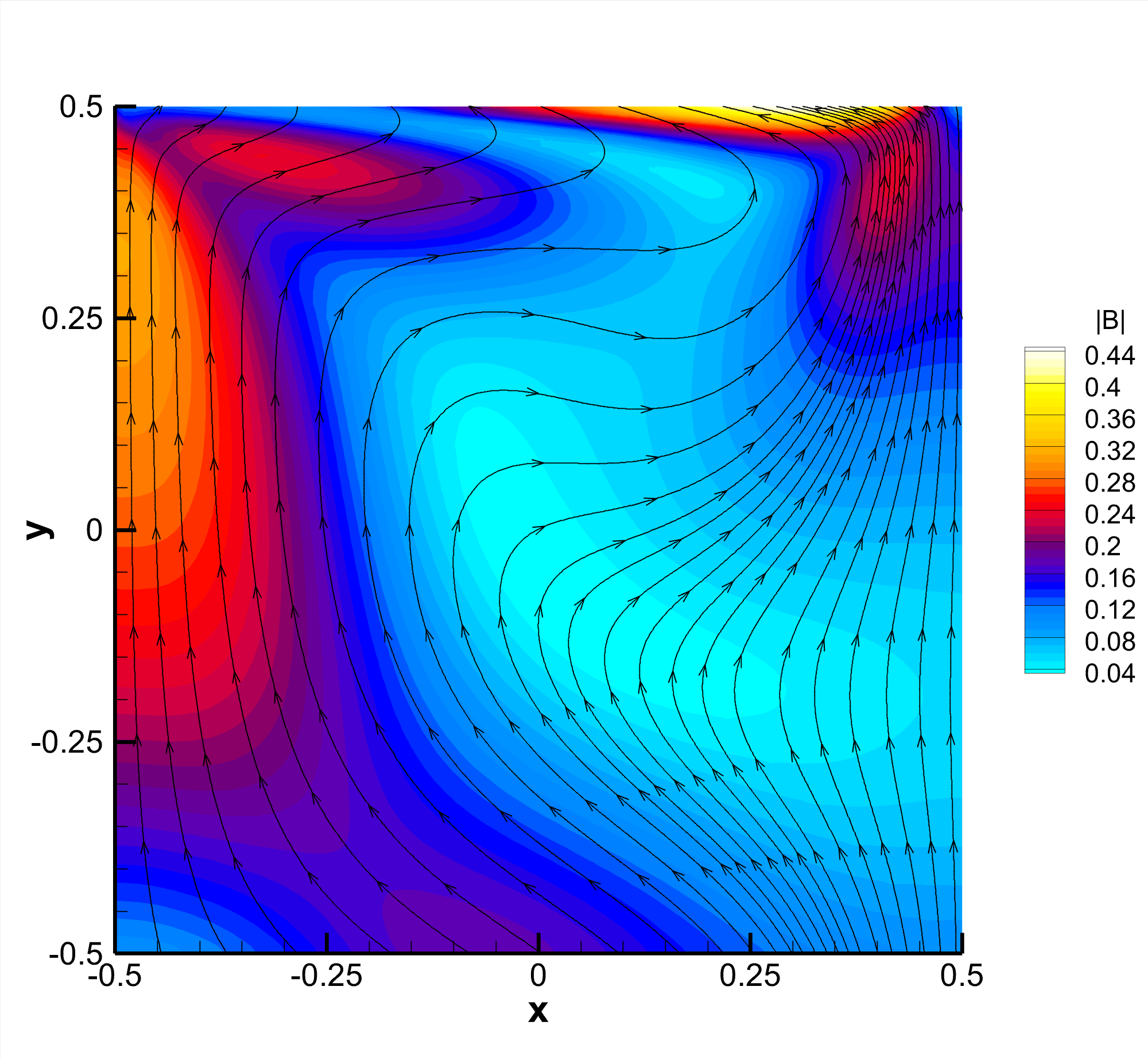} 
	\includegraphics[width=0.45\linewidth]{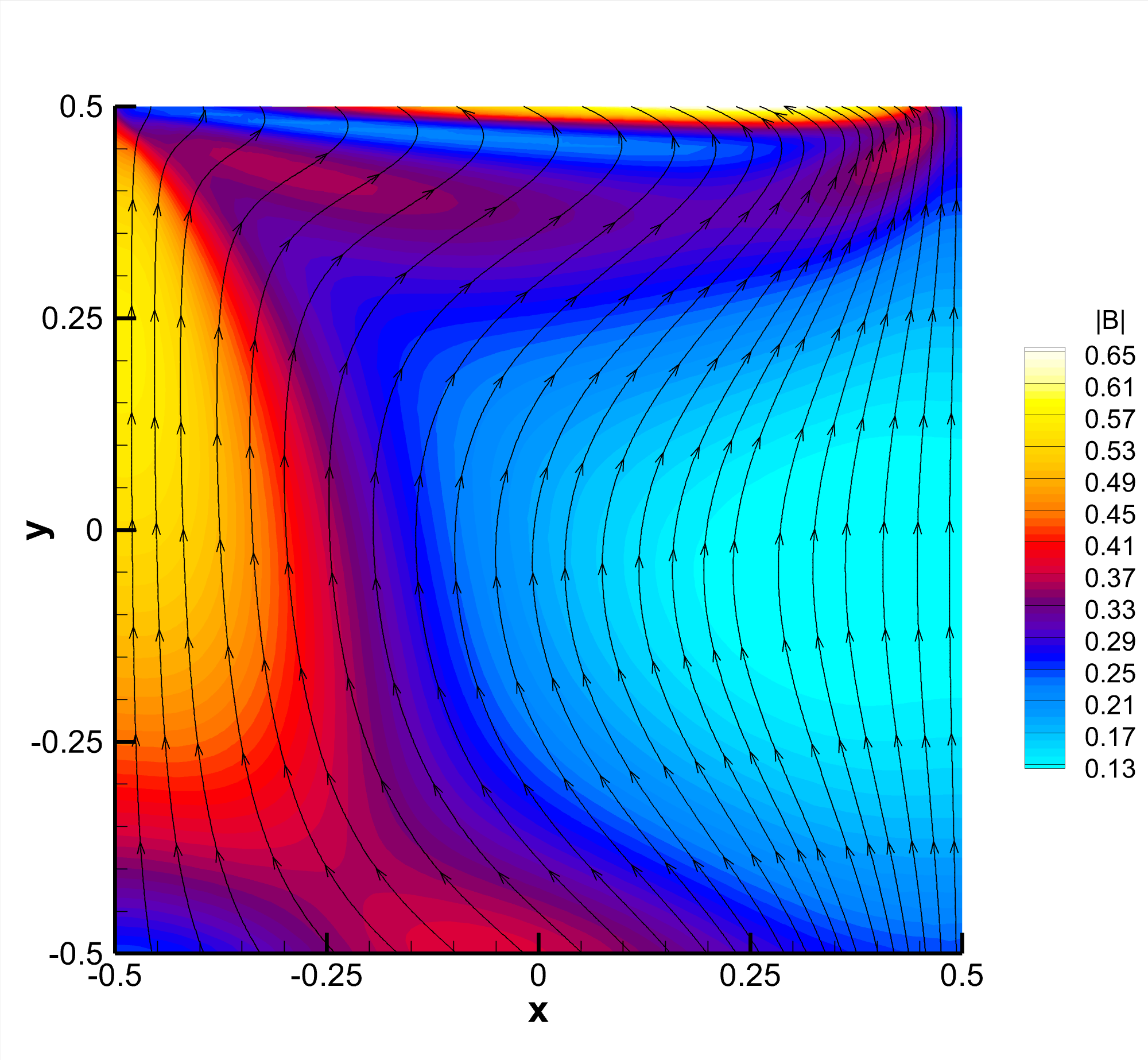} \\
	\includegraphics[width=0.45\linewidth]{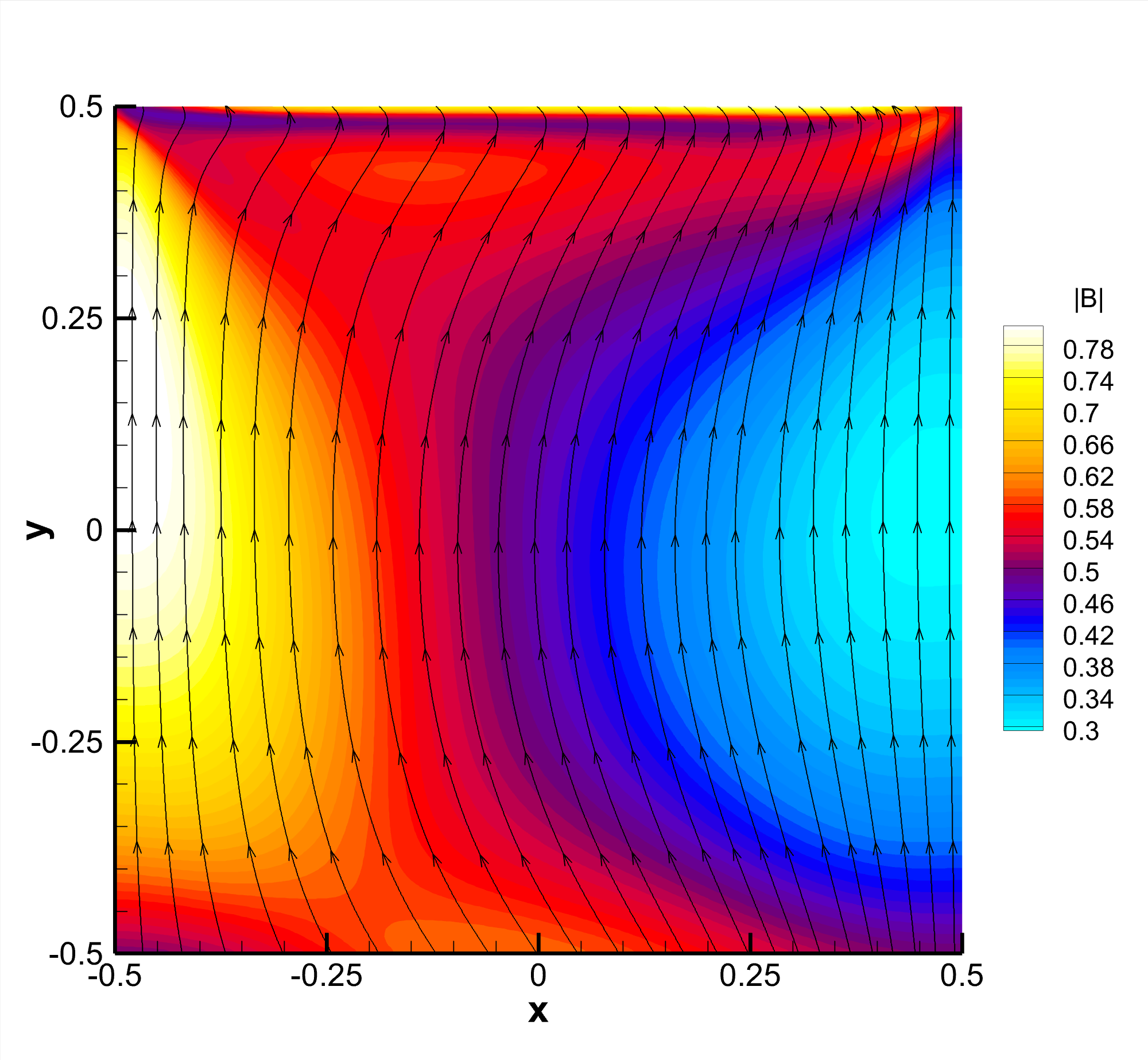} 
	\includegraphics[width=0.45\linewidth]{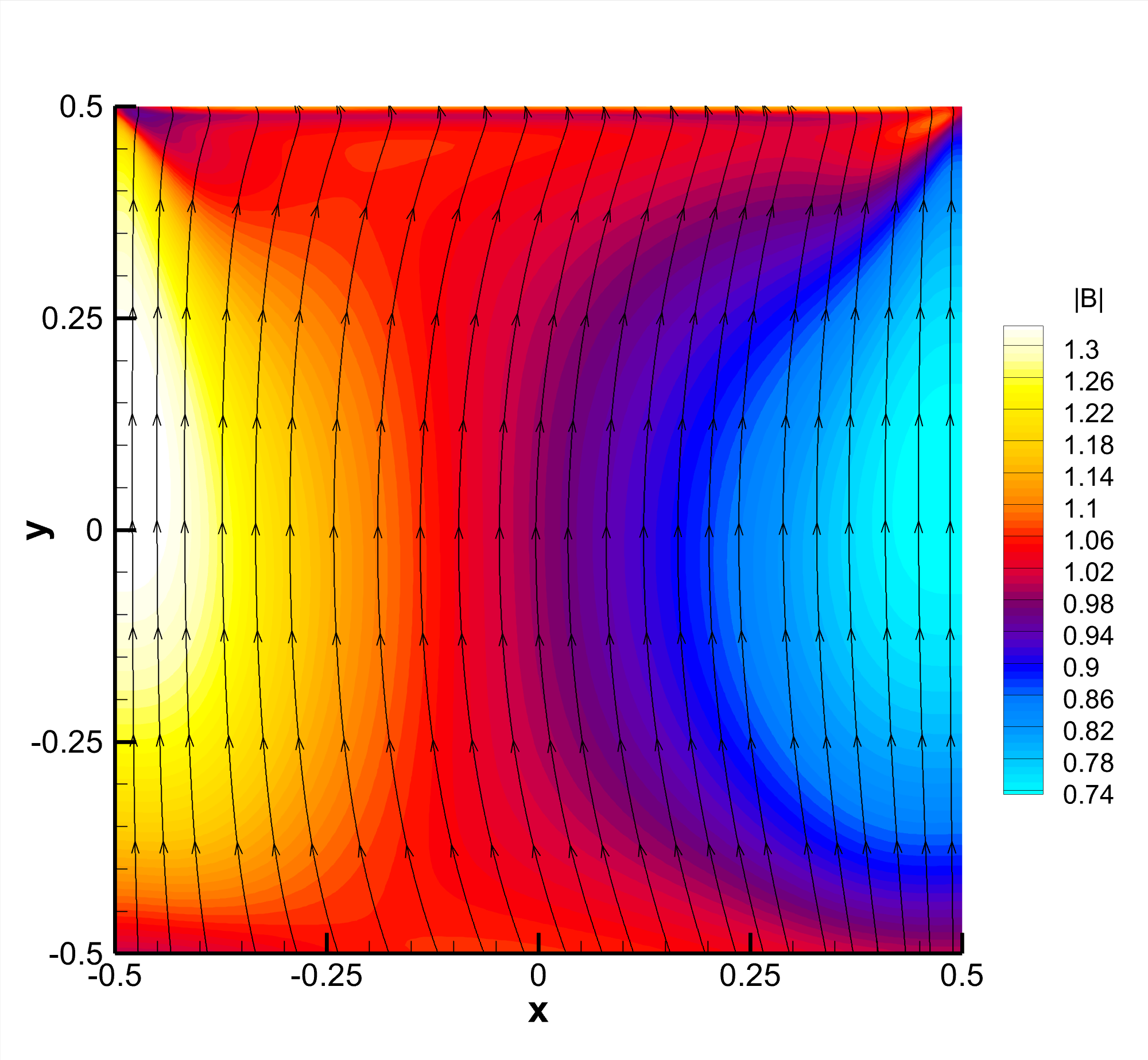} 
	\caption{Stream-lines of the magnetic field at time $t=20$ for the MHD lid-driven cavity with vertical magnetic field $\magfield_0 = (0,B_{0,y})^{T}$ for $B_{0,y} \in\left\lbrace 0.1, 0.25,0.5,1\right\rbrace$ (top-left, top-right, bottom-left, and bottom-right, respectively).}  
	\label{fig.MHDLDC_vertical_modB}
\end{figure}

This test is a variant of the traditional lid-driven cavity test in which a horizontal or a vertical magnetic field is imposed. The computational domain is again $\Omega=[-0.5,0.5]\times[-0.5,0.5]$, and the initial and boundary conditions for velocity, density, and pressure are the same of Section~\ref{sec:LDC}. The walls are assumed to be \emph{perfectly-conducting} and the following boundary conditions are imposed for the magnetic and electric field, respectively: 
$$\magfield\cdot\pmb{n} = \magfield_0\cdot\pmb{n}, \quad \quad \elfield\times\pmb{n} = \pmb{0}.$$
By using the Ohm law, we get 
$$ \left( - \bvel \times \bbvar + \eta \nabla \times \bbvar \right)\times\pmb{n} = \pmb{0}. $$
The initial magnetic field is either horizontal, i.e., $\bbvar_0 = (\bvar_{0,x},0)^{T}$, or vertical, i.e., $\bbvar_0 = (0,\bvar_{0,y})^{T}$. 
If we consider a vertical magnetic field, then at the lid the following Neumann boundary condition on the tangential component of the magnetic field holds 
\begin{equation*}
	\frac{\partial \bvar_1}{\partial y} = -\frac{v \bvar_{0,y}}{\eta}.
\end{equation*} 
Note that large values of $B_{0,y}$, or small values of $\eta$, will generate a very steep boundary layer at the wall. Because of this, the domain is discretized with the rather fine unstructured mixed-element mesh M$_{100}$ described in Section~\ref{sec:LDC}, except for $B_{0,y}\geq 0.5$, for which we use an even finer Cartesian $200\times200$ grid. The test is run for $\mu = \eta = 0.01$ and different values of $B_{0,x}$ and $B_{0,y}$ until a final time of $t = 20$ when the solution has been checked to become stationary. 
The stream-lines of velocity and magnetic field for the horizontal test are shown in Figures \ref{fig.MHDLDC_horizontal_modV} and \ref{fig.MHDLDC_horizontal_modB}, respectively, while those for the vertical tests are shown in Figures \ref{fig.MHDLDC_vertical_modV} and \ref{fig.MHDLDC_vertical_modB}. We can observe how the effect of the moving lid is damped as the strength of the magnetic field increases. 

\begin{figure}
	\centering
	\includegraphics[width=0.45\linewidth]{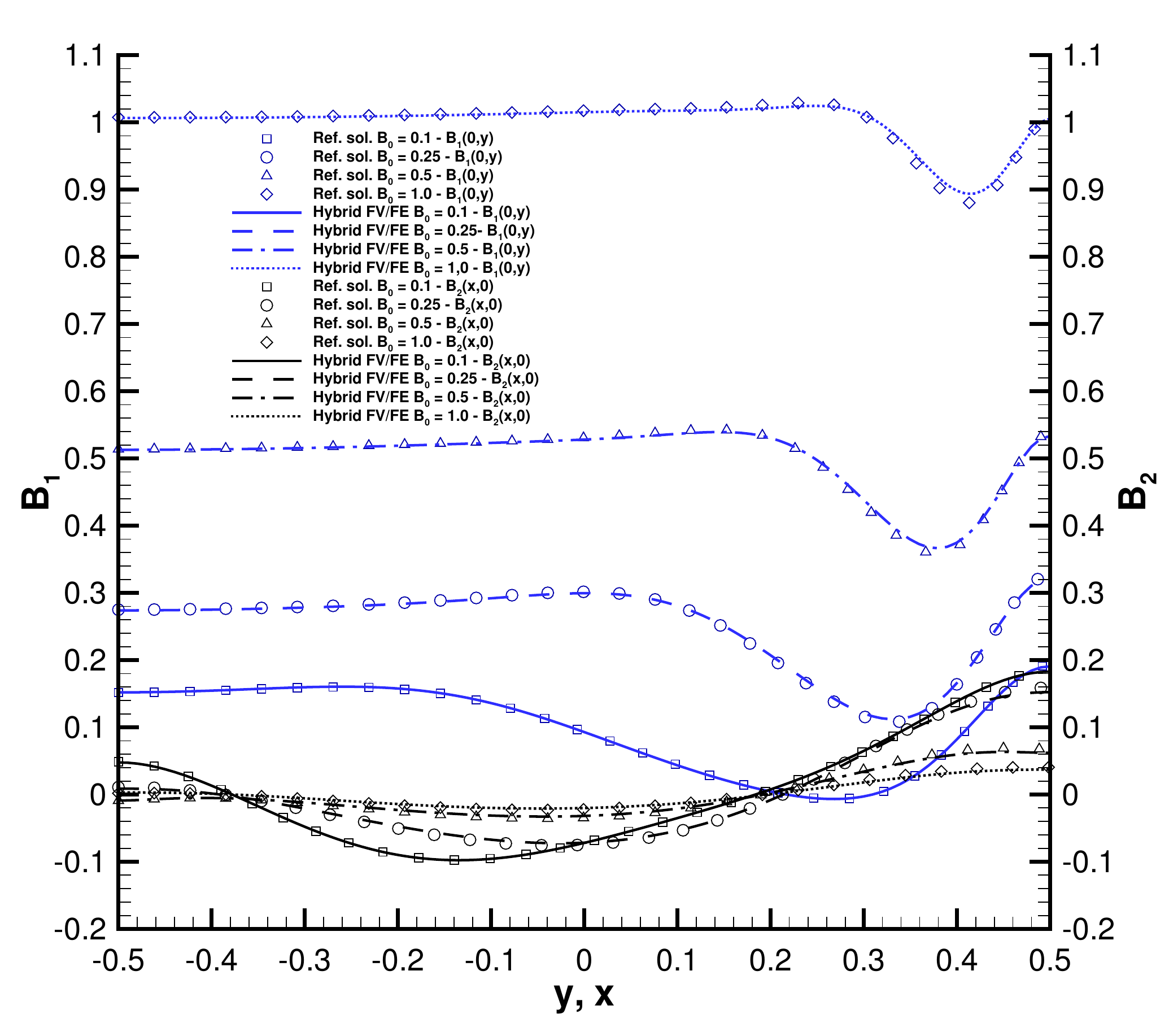}\hspace{0.02\linewidth}
	\includegraphics[width=0.45\linewidth]{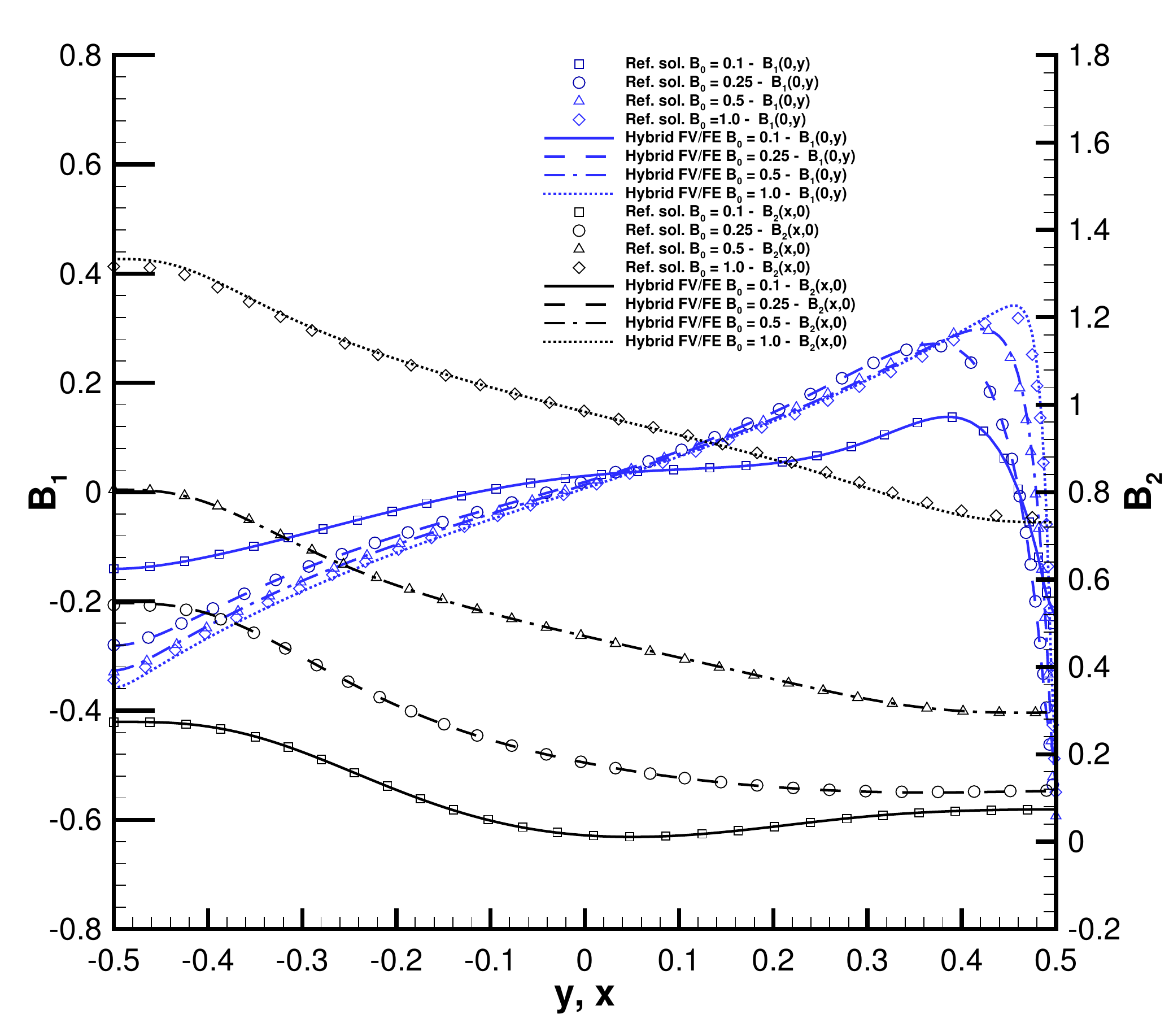} 
	\caption{1D cuts through the numerical solution at $x=0$ and $y=0$ for the MHD lid-driven cavity at time $t=20$ and comparison with the reference solution obtained with the numerical method introduced in \cite{SIMHD}. Left: initially horizontal magnetic field. Right: initially vertical magnetic field. In all cases, the lid velocity is $\vel=(1,0)$, and the viscosity and resistivity have been set to $\mu=\eta=10^{-2}$. }  
	\label{fig.MHDLDC_cuts}
\end{figure}

In Figure \ref{fig.MHDLDC_cuts}, 1D cuts of the solution along the $x$-axis and the $y$-axis are shown for the numerical solution obtained with the new hybrid FV/FE method, while the reference solution has been computed with the semi-implicit divergence-free method introduced in \cite{SIMHD}. A very good agreement between the numerical and reference solutions is obtained. In order to allow other research groups to reproduce our results and to compare with them, all results are available in electronic format\footnote[1]{The supplementary material is licensed under a Creative Commons ``Attribution-NonCommercial-ShareAlike 3.0 Unported'' license. \ccbyncsa} in the supplementary material provided with this paper. In~\ref{sec:appendix}, a detailed description of the available data files is given.

\subsection{Double shear layer} \label{sec:DSL}
\begin{figure}
	\centering
	\includegraphics[width=0.45\linewidth]{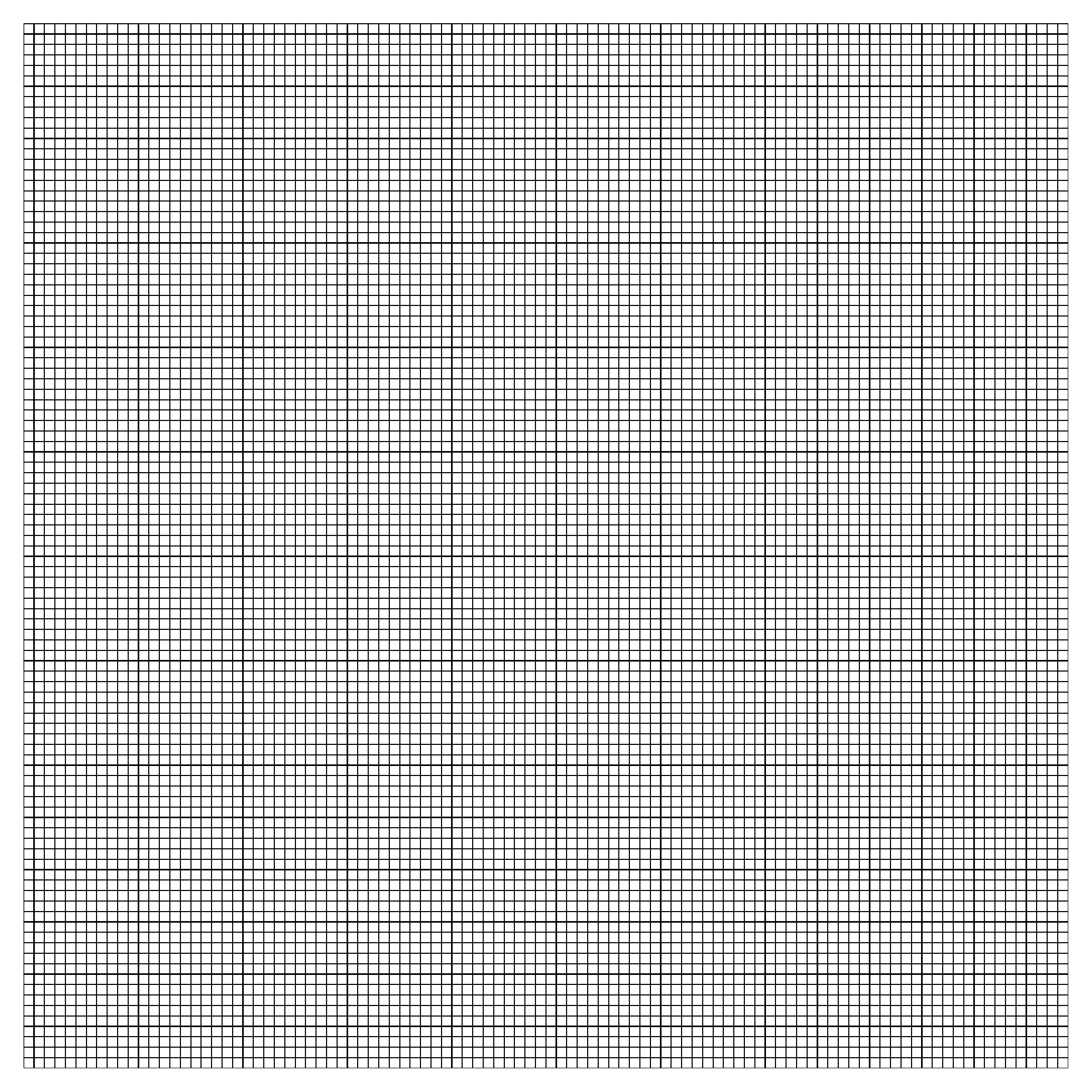}\hspace{0.05\linewidth}
	\includegraphics[width=0.45\linewidth]{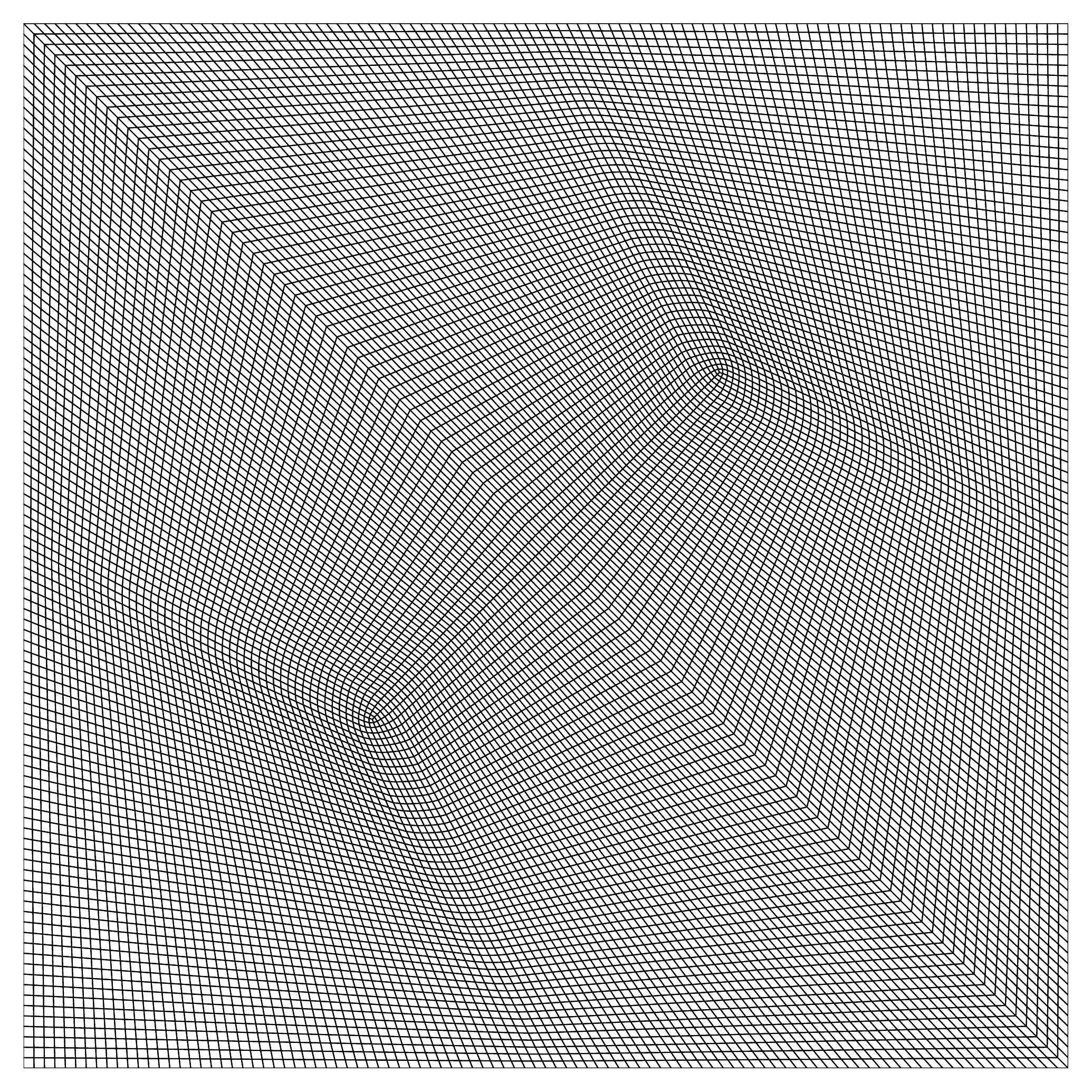}
	
	\vspace{0.05\linewidth}
	\includegraphics[width=0.45\linewidth]{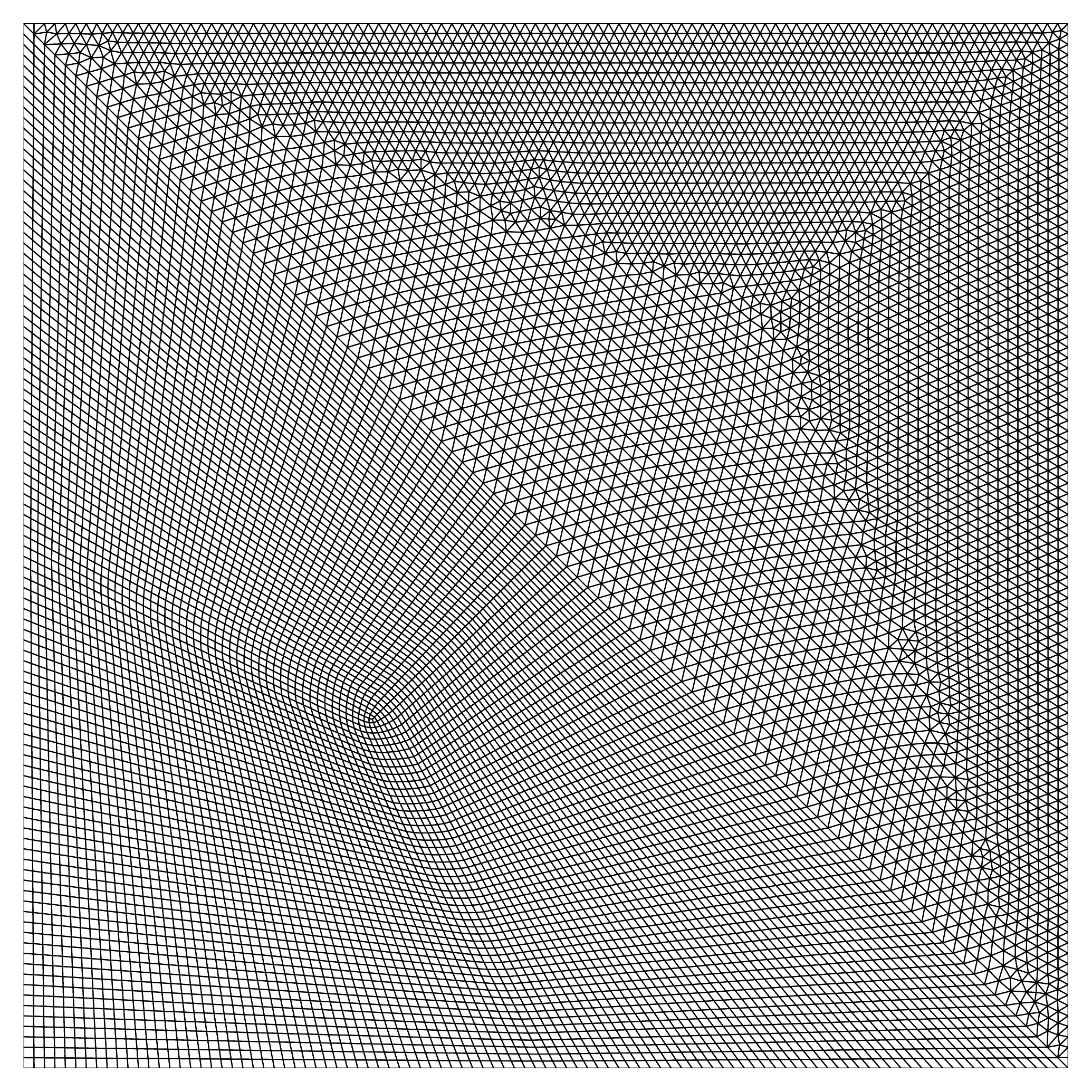}\hspace{0.05\linewidth}
	\includegraphics[width=0.45\linewidth]{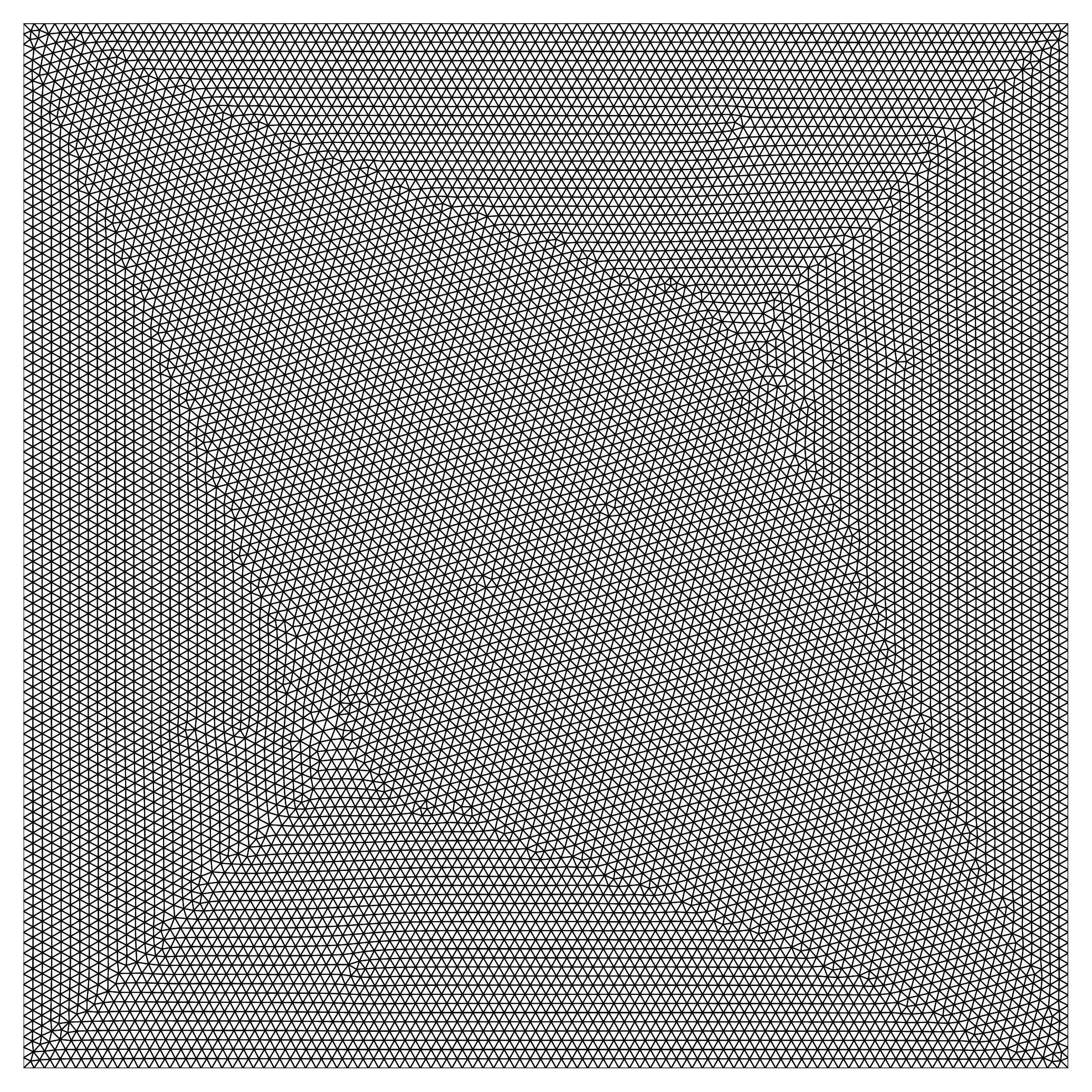}
	\caption{Sketch of the four different grids employed to run the double shear layer test case. Top left: Cartesian grid made of $10000$ primal elements, M$_{\mathrm{Cart}}$. Top right: grid made of $15000$ skewed quadrilaterals, M$_{\mathrm{Skew}}$. Bottom left: mixed mesh of $16174$ skewed quadrilaterals and triangles, M$_{\mathrm{Mix}}$. Bottom right: purely triangular grid made of $23056$ primal elements, M$_{\mathrm{Tria}}$.}
	\label{fig.DSL_grids}
\end{figure}

\begin{figure}
	\centering
	\includegraphics[width=0.43\linewidth]{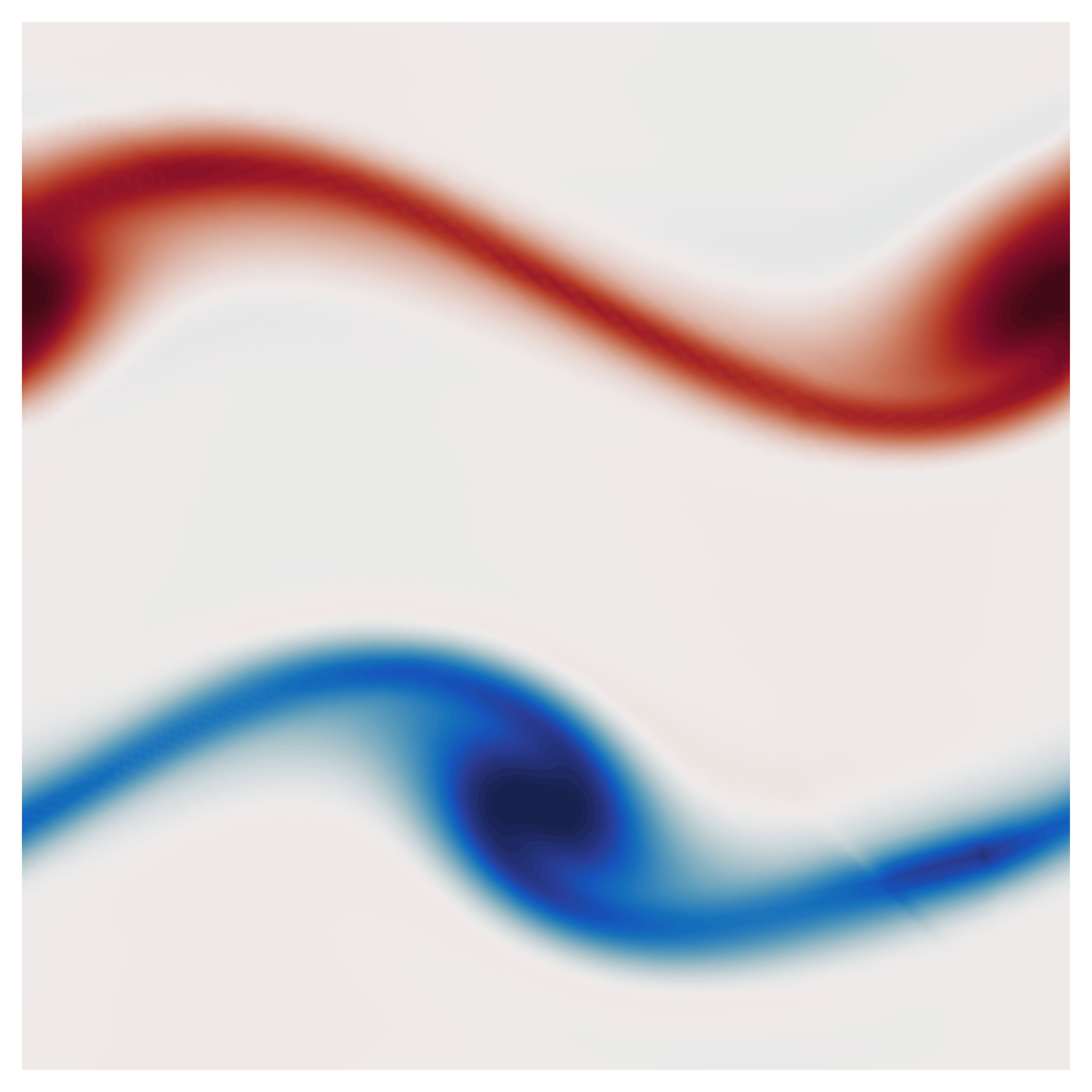} \hspace{0.01\linewidth}
	\includegraphics[width=0.43\linewidth]{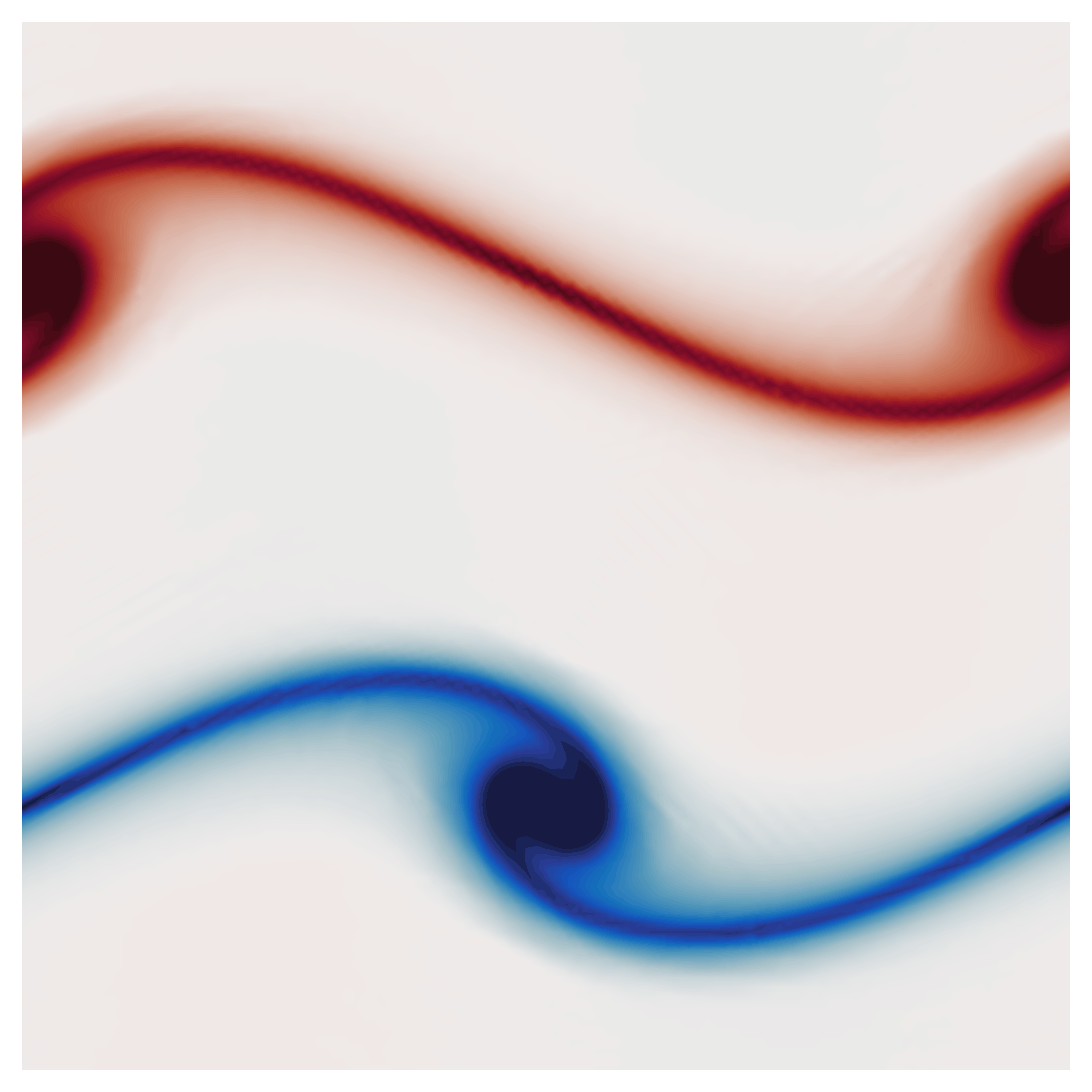}
	
	\includegraphics[width=0.43\linewidth]{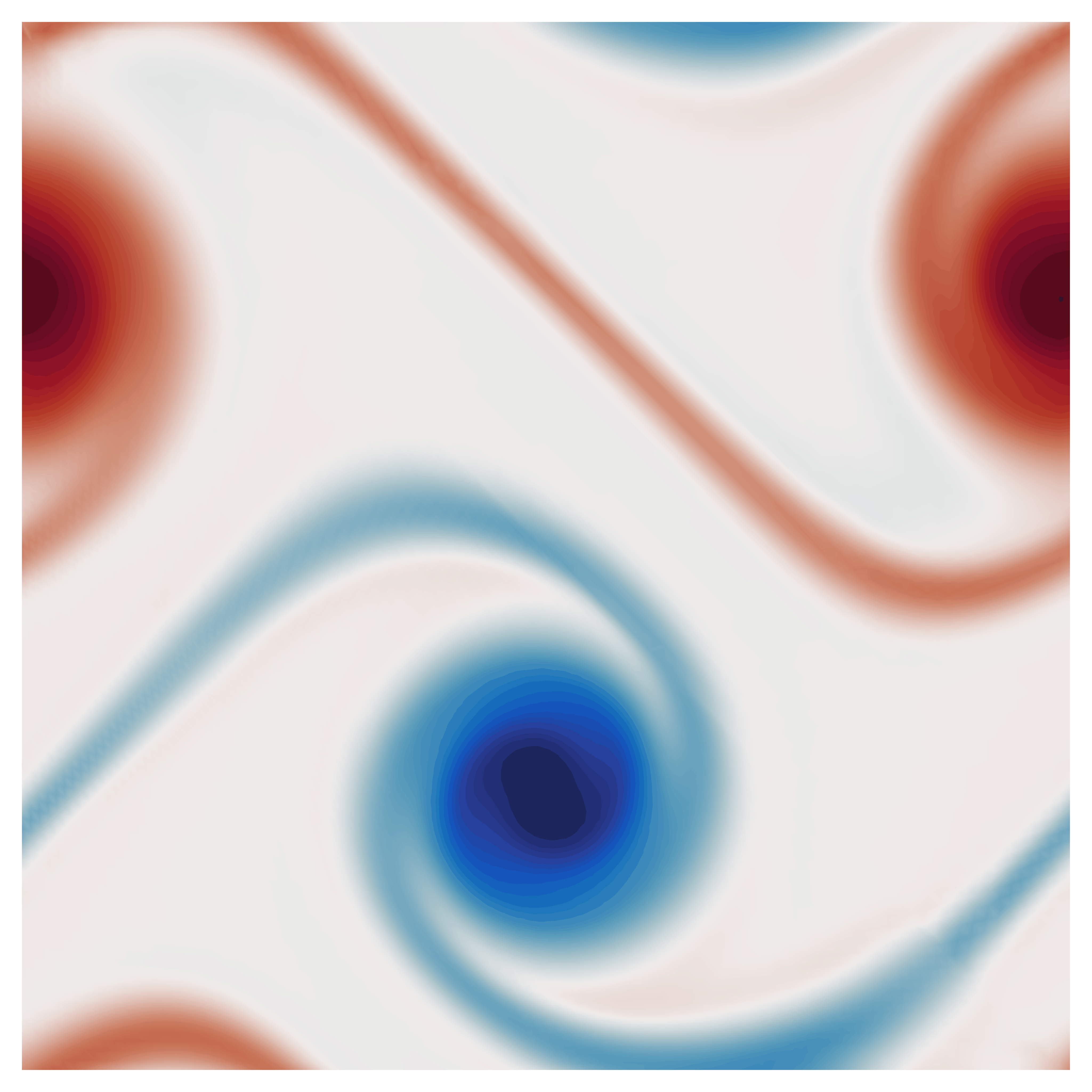}\hspace{0.01\linewidth}
	\includegraphics[width=0.43\linewidth]{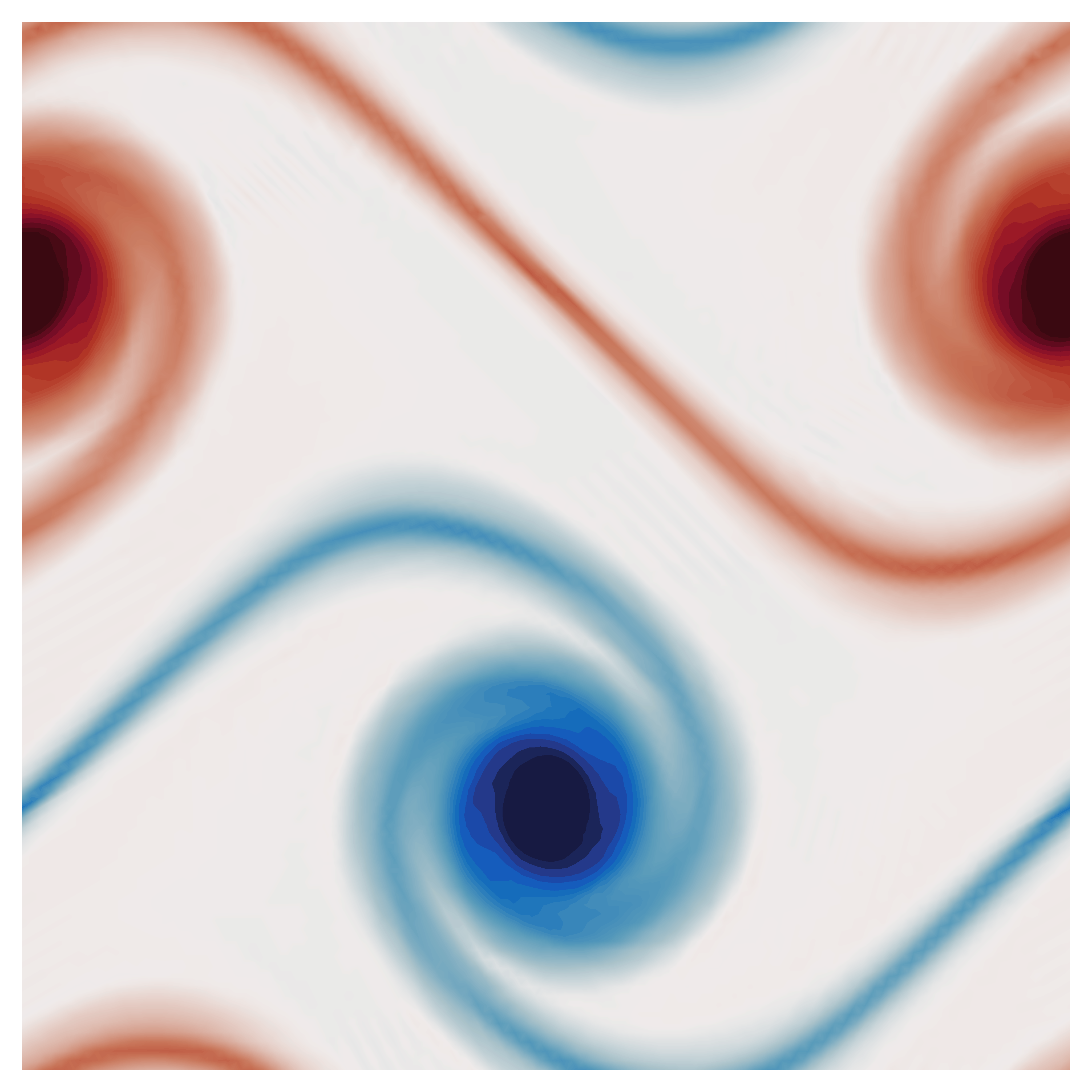}
	
	\includegraphics[width=0.43\linewidth]{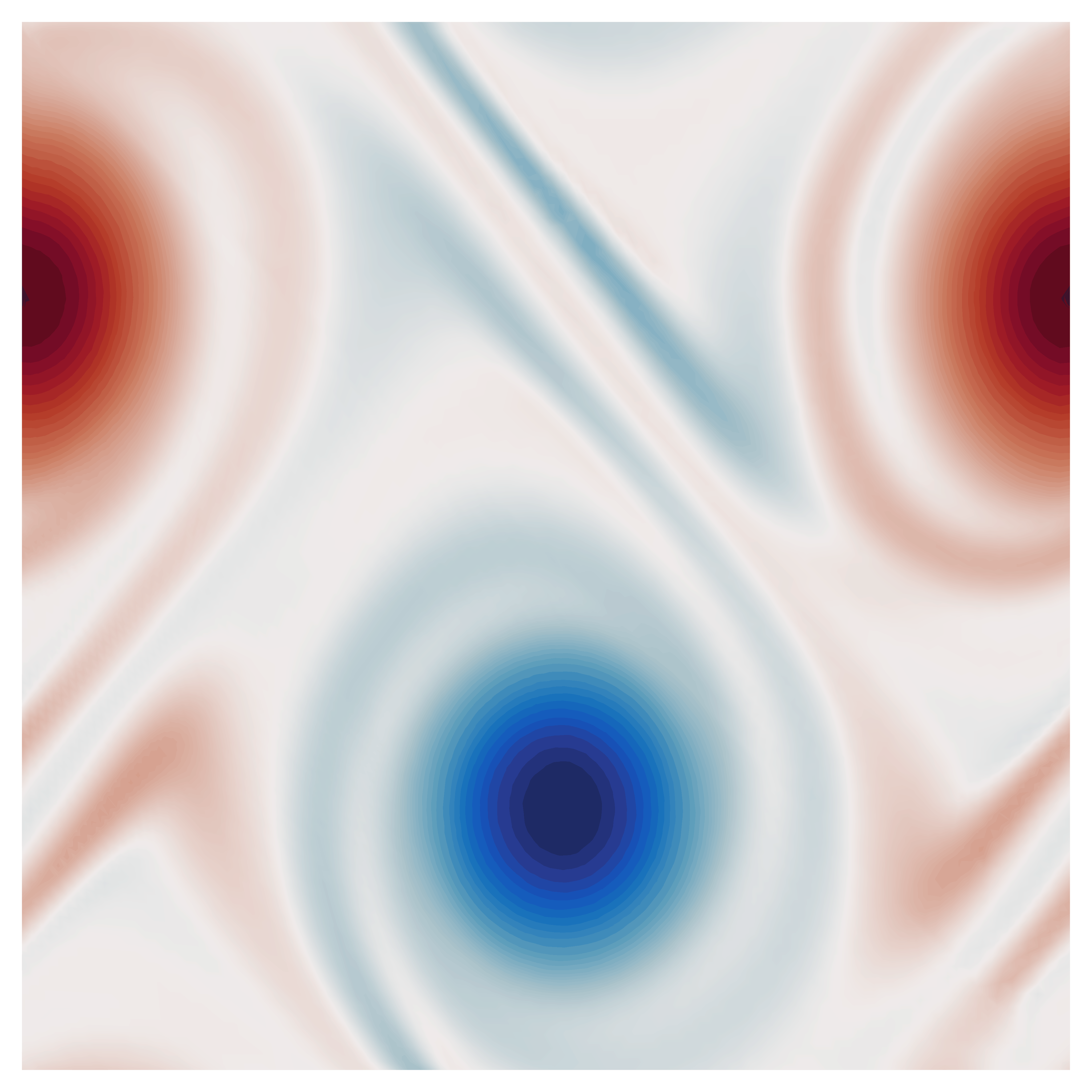}\hspace{0.01\linewidth}
	\includegraphics[width=0.43\linewidth]{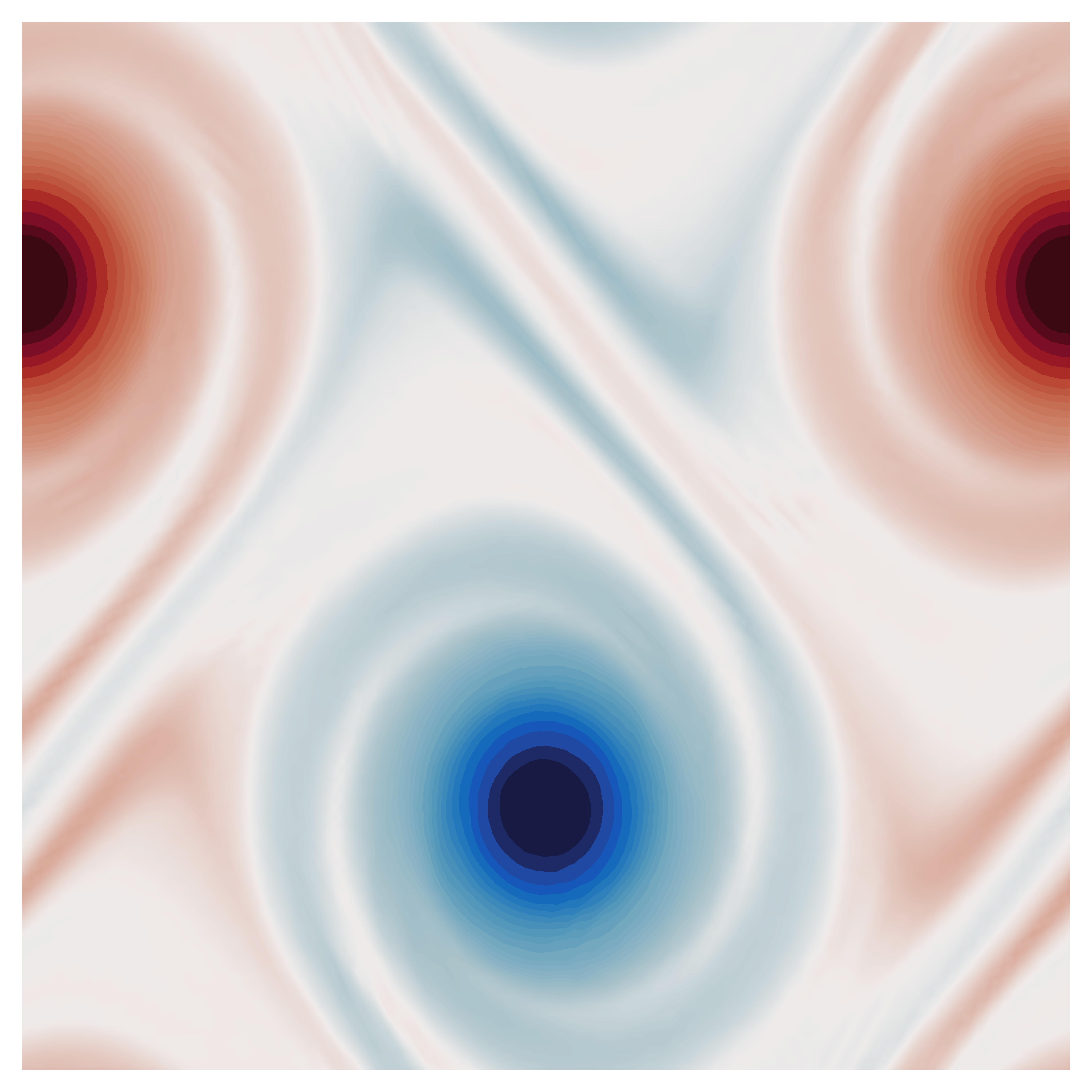}
	
	\includegraphics[width=0.4\linewidth]{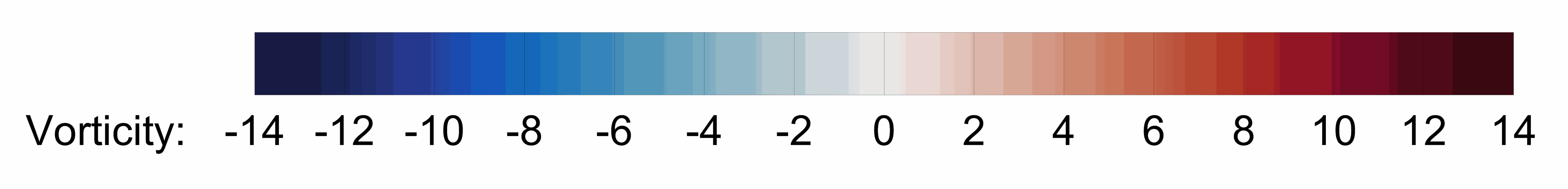}
	\caption{Vorticity contour plots of the double shear layer test case at times $t\in\left\lbrace 1.6,2.4,3.6 \right\rbrace$ (from top to bottom). Left: proposed hybrid FV/FE scheme on general grids, mesh M$_{\mathrm{Mix}}$. Right: reference solution \cite{HybridMPI}, mesh M$_{\mathrm{Tria}}$.}
	\label{fig.DSL_time}
\end{figure}

\begin{figure}
	\centering
	\includegraphics[width=0.45\linewidth]{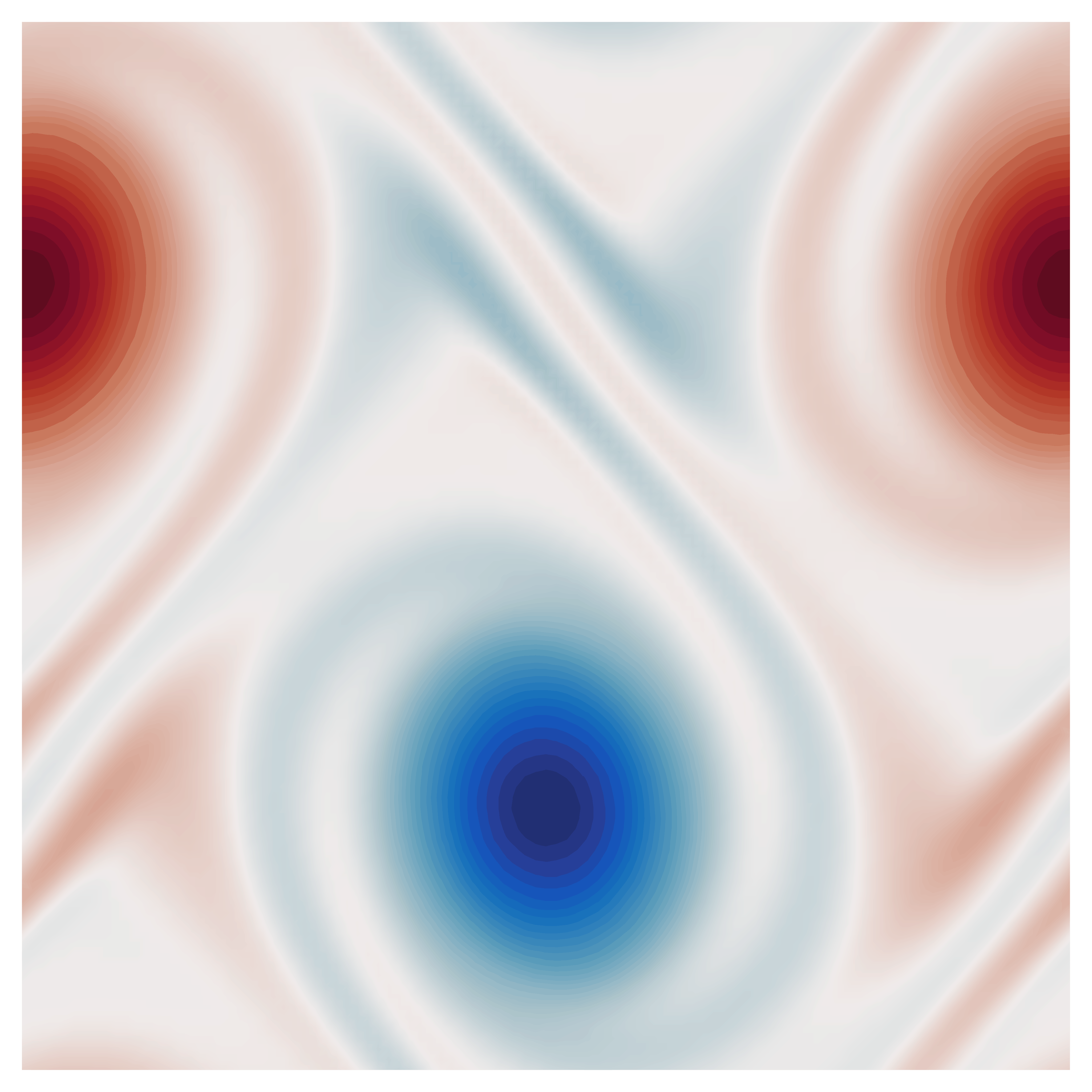}\hspace{0.02\linewidth}
	\includegraphics[width=0.45\linewidth]{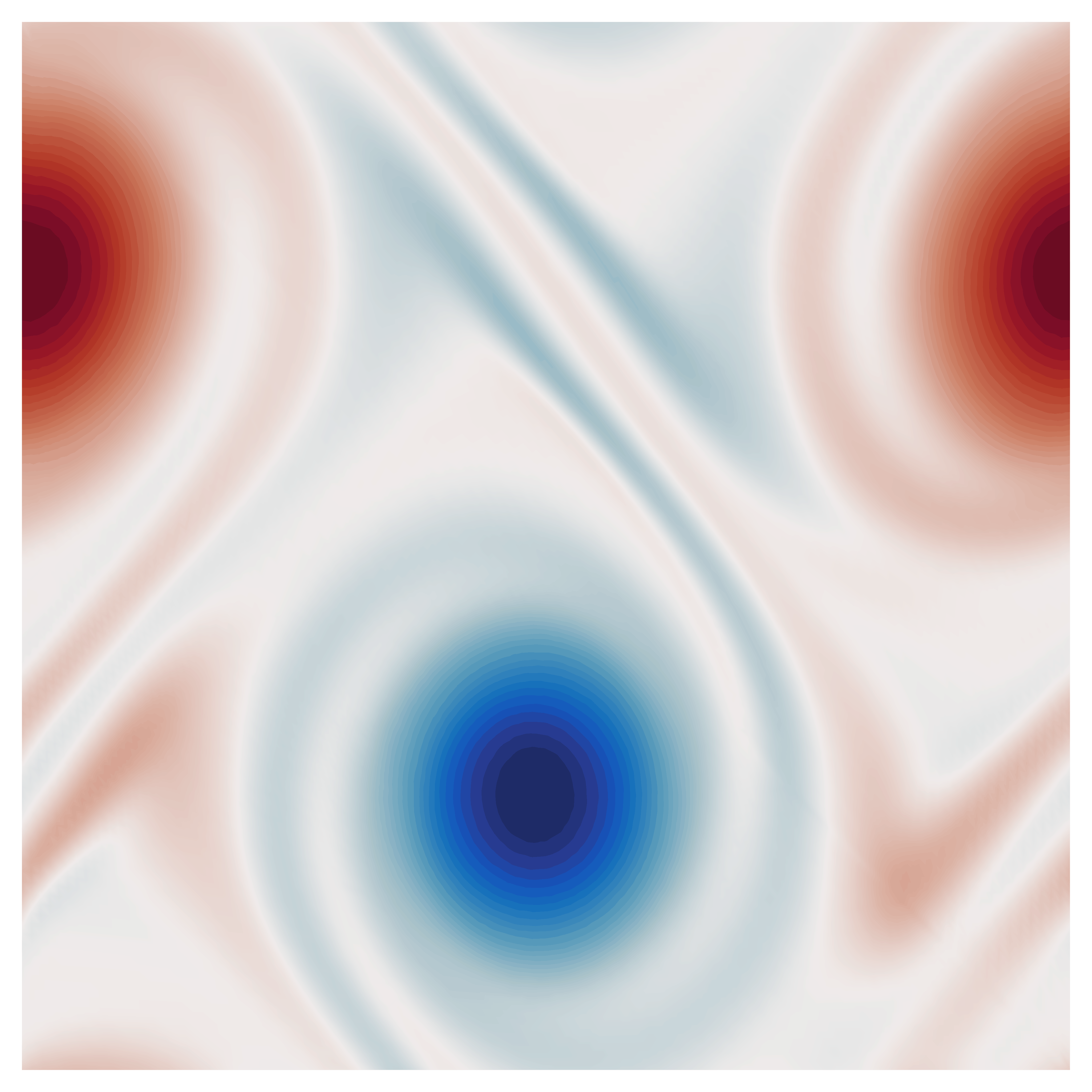}
	
	\vspace{0.02\linewidth}
	\includegraphics[width=0.45\linewidth]{DSL_vort36_MixedMesh}\hspace{0.02\linewidth}
	\includegraphics[width=0.45\linewidth]{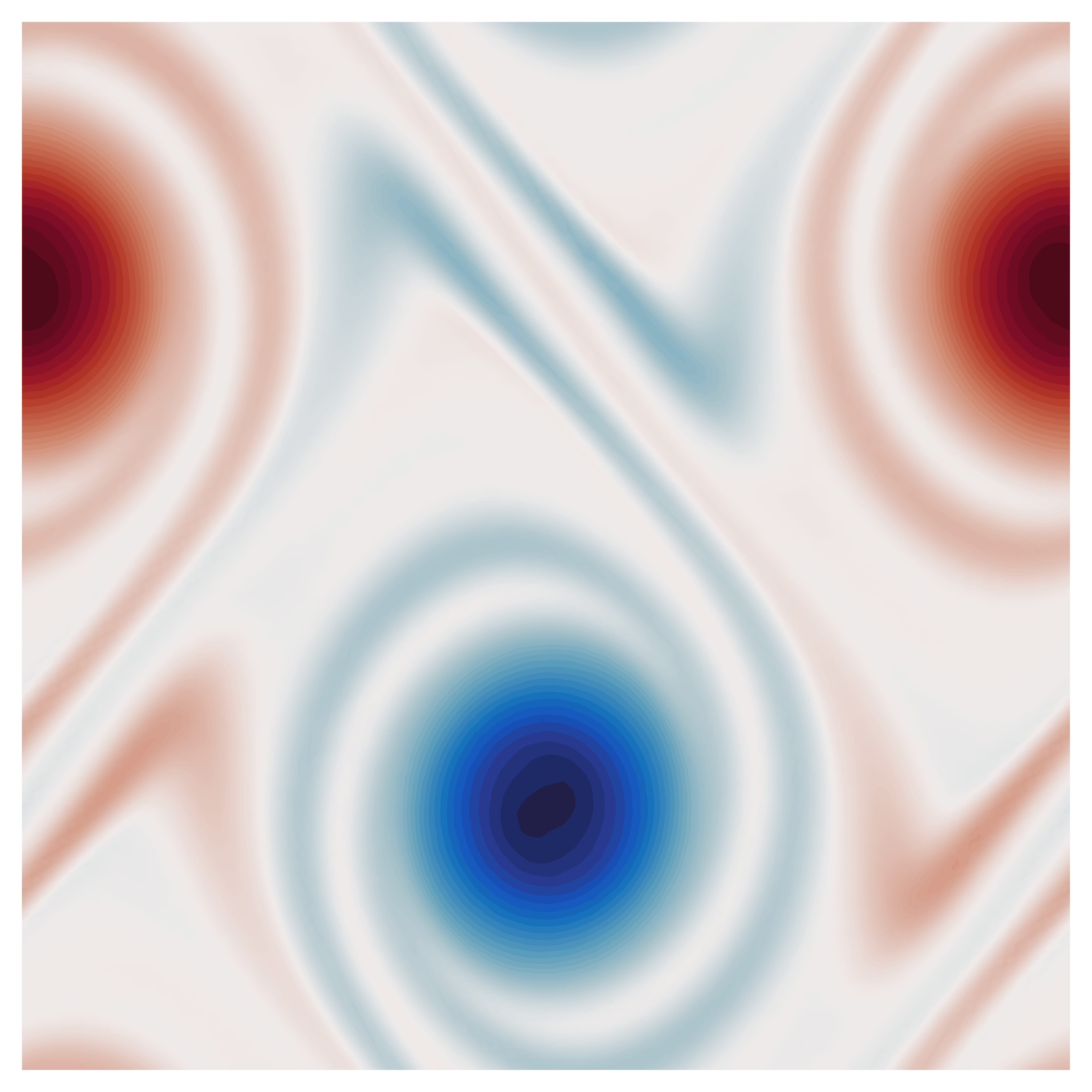}
	
	\includegraphics[width=0.4\linewidth]{DSL_vort_Legend}
	\caption{Vorticity contours of the double shear layer benchmark obtained at $t=3.6$ using the proposed hybrid FV/FE scheme on grids: M$_{\mathrm{Cart}}$, top left; M$_{\mathrm{Skew}}$, top right; M$_{\mathrm{Mix}}$, bottom left; M$_{\mathrm{Tria}}$, bottom right.}
	\label{fig.DSL_finaltimesol}
\end{figure}

We now consider the double shear layer test case, \cite{BCG89}, a classical benchmark used to analyse the behaviour of incompressible flow solvers. The initial flow is characterized by the presence of important velocity gradients, which lead to the formation of vortex-shaped patterns.
We consider an initial condition given by
\begin{gather*}
	\press \left(\xx,0\right) = \frac{10^{4}}{1.4}, \qquad
	\vel_{1} \left(\xx,0\right) = \left\lbrace \begin{array}{lc}
		\tanh \left[\hat{\rho}(\hat{y}-0.25)\right] & \mathrm{ if } \; \hat{y} \le 0.5,\\
		\tanh \left[\hat{\rho}(0.75-\hat{y})\right] & \mathrm{ if } \; \hat{y} > 0.5,
	\end{array}\right.
\\ 
	\vel_{2} \left(\xx,0\right) = \delta \sin \left(2\pi \hat{x}\right), \qquad
	\hat{x}= \frac{x+1}{2}, \qquad \hat{y}= \frac{y+1}{2},
\end{gather*}
in the computational domain $\Omega=[-1,1]^2$ with periodic boundary conditions everywhere. Moreover, we set the laminar viscosity to $\mu= 2\cdot 10^{-4}$ while $\hat{\rho} = 30$ and $\delta =0.05$ are the parameters that determine the slope of the shear layer and the amplitude of the initial perturbation.
The simulation is run up to time $t=3.6$ with $\mathrm{CFL}=0.5$ employing four different grids in order to check the versatility of the proposed scheme to deal with different shapes of primal elements. All meshes have been generated defining $100$ divisions along each boundary; the different grid configurations designed are depicted in Figure \ref{fig.DSL_grids}. The vorticity contours obtained for different time instants using an unstructured mixed-element grid are shown in Figure~\ref{fig.DSL_time}. A reference solution, computed using the original hybrid FV/FE scheme on triangular simplex grids based on a completely staggered approach where the momentum is effectively computed in the dual grid, see \cite{BFTVC17,Hybrid1,HybridMPI}, is also included for comparison. Furthermore, the solution obtained with all four grids at the final time, $t=3.6$, is reported in Figure \ref{fig.DSL_finaltimesol}. A good match between the numerical solution obtained with the new hybrid FV/FE scheme on four different unstructured mixed-element meshes and the reference solution is observed. 

%
%
%

\subsection{Viscous and resistive Orszag-Tang vortex}
\begin{figure}
	\centering
	\includegraphics[width=0.45\linewidth]{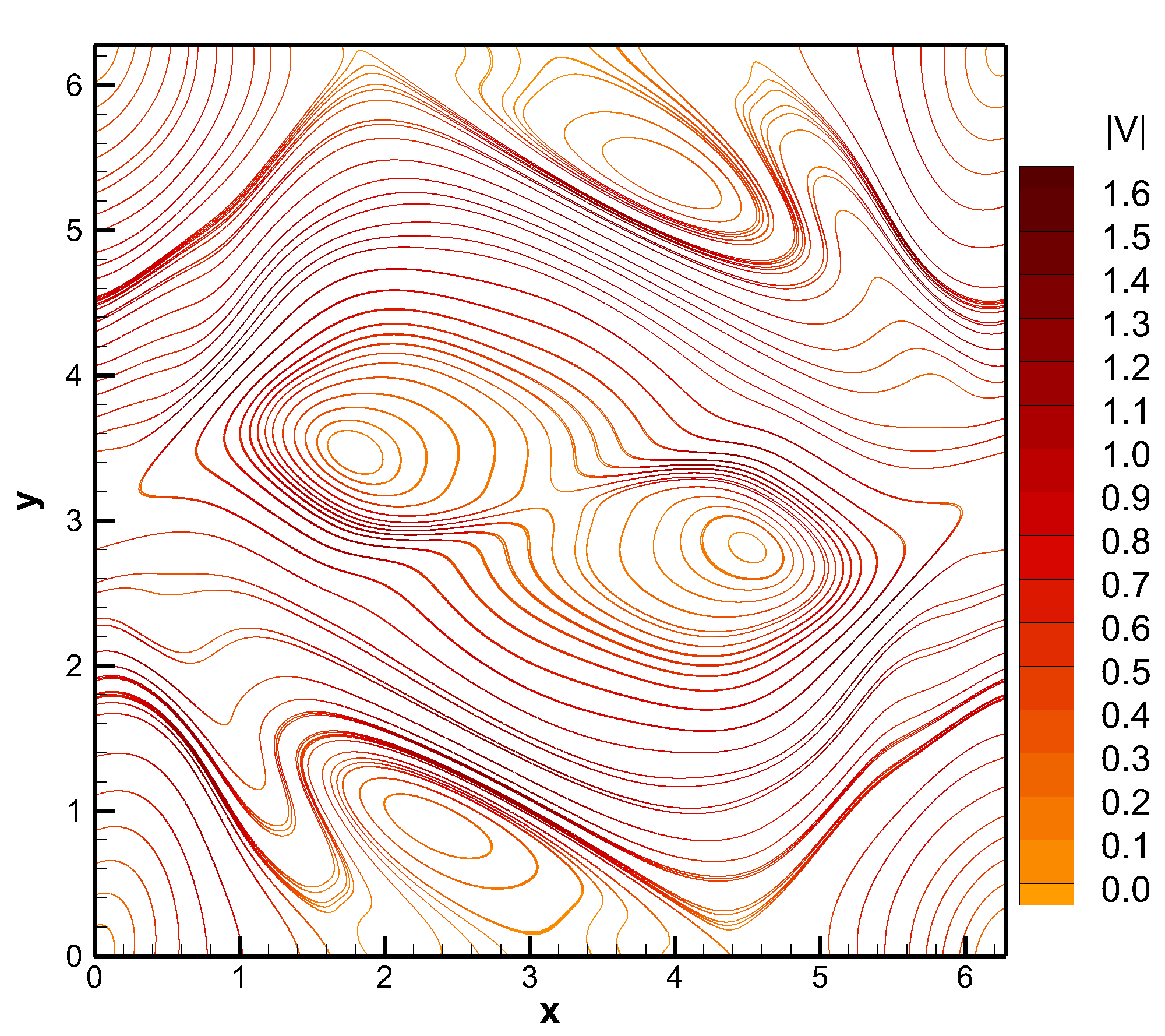}\hspace{0.02\linewidth}
	\includegraphics[width=0.45\linewidth]{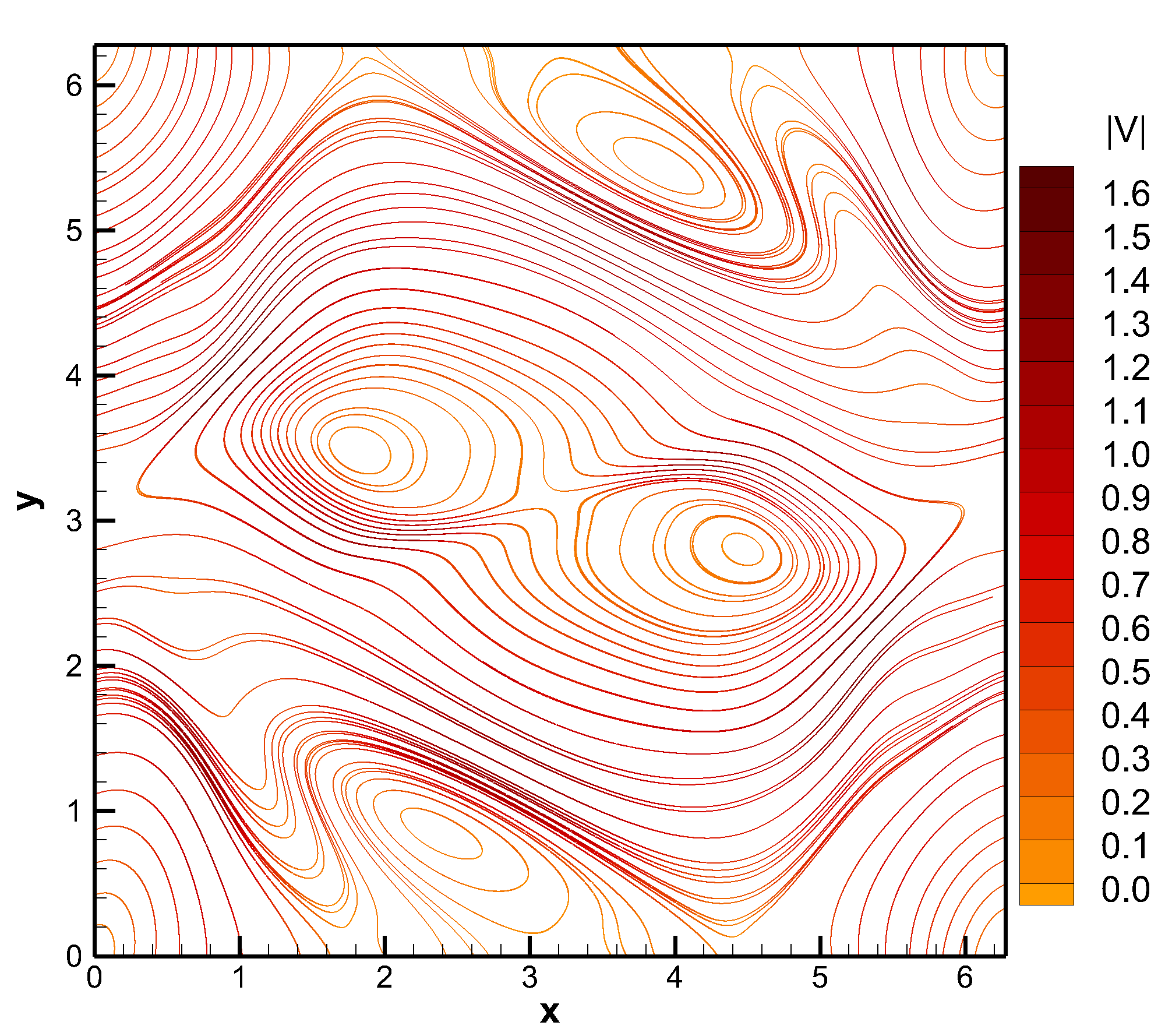} 
	\vspace{0.02\linewidth}
	\includegraphics[width=0.45\linewidth]{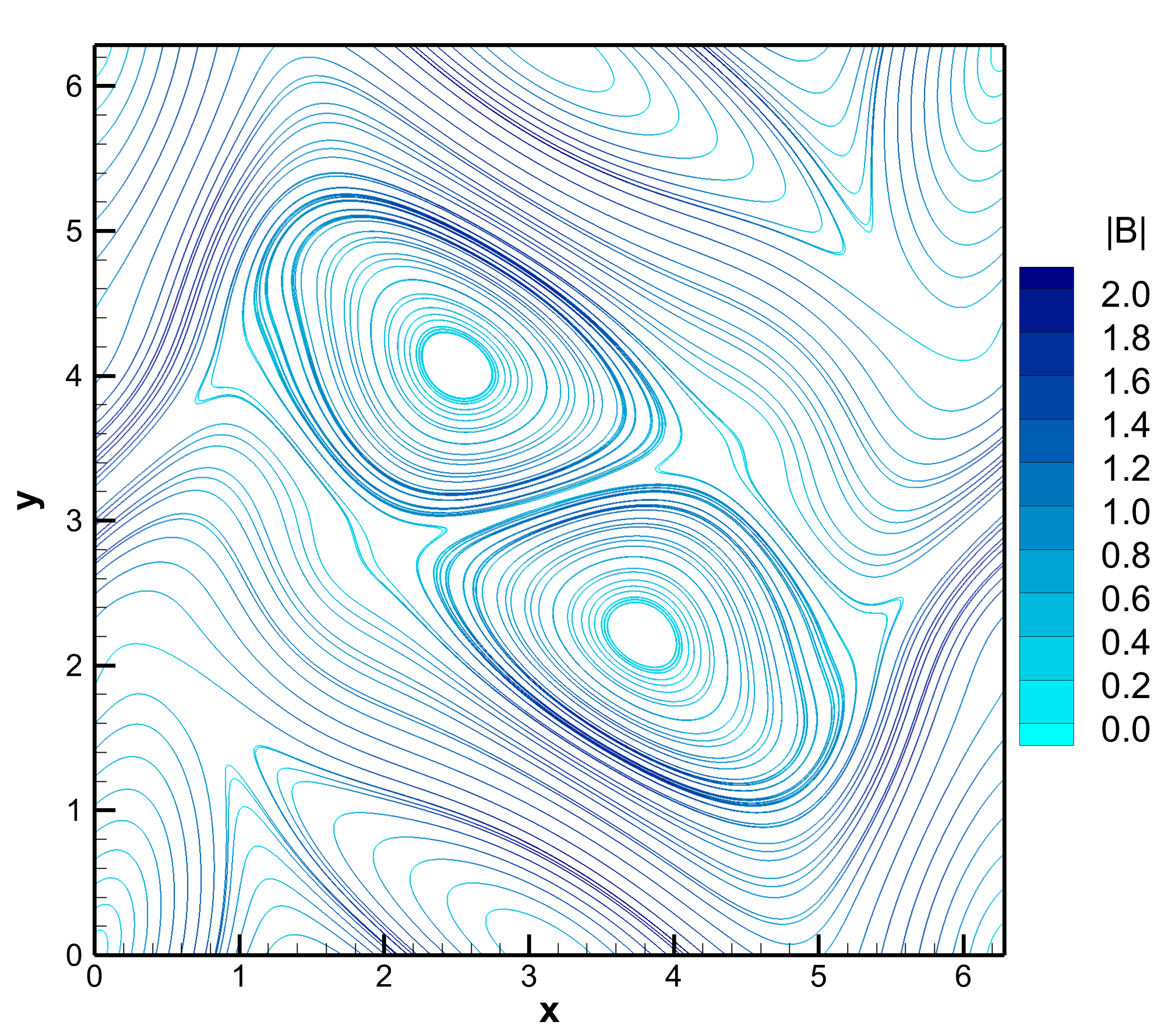}\hspace{0.02\linewidth}
	\includegraphics[width=0.45\linewidth]{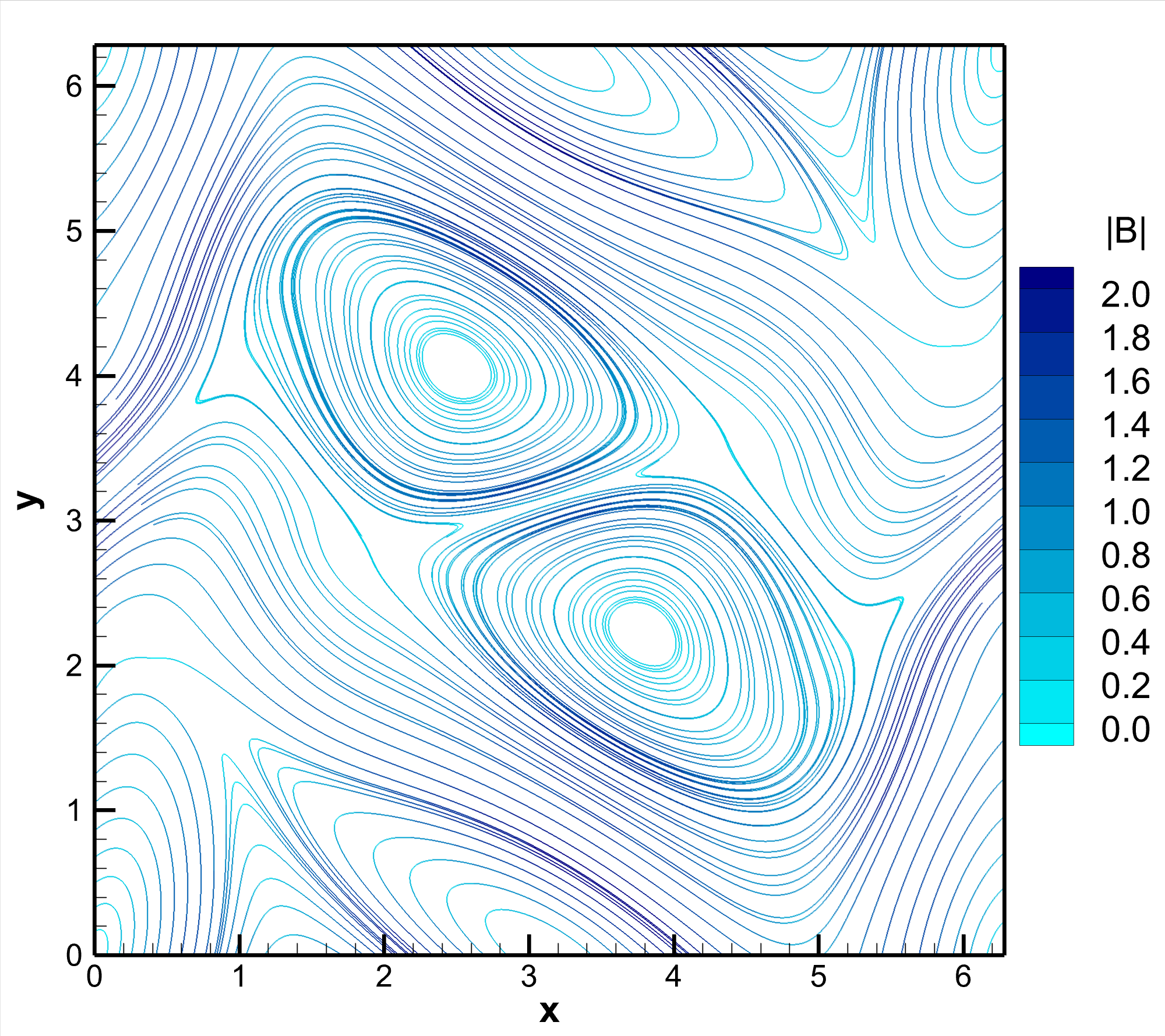} 
	\caption{Stream-lines of velocity (top row) and magnetic field (bottom row) at time $t=2$ for the viscous and resistive Orszag-Tang vortex problem. Numerical solution obtained with the proposed hybrid FV/FE scheme on a mixed skew-quadrilateral and triangular mesh (left column), compared to the reference solution computed with the divergence-free semi-implicit method outlined in \cite{SIMHD} on a $1000^2$ Cartesian grid (right column).}  
	\label{fig.OT_finaltimesol}
\end{figure}

The two-dimensional incompressible viscous and resistive Orszag-Tang vortex, see e.g. \cite{WarburtonVRMHD,SIMHD}, is a well-known configuration for testing both the robustness and the correct convergence of a numerical scheme for MHD. Indeed, in this test, the non-linear advection, viscous and resistive dissipations, incompressibility of the fluid and the divergence-free condition of $\magfield$ all play an important role.
The initial condition reads 
\begin{align*} 
	\bvel\left(\xx,0\right)  =  \sqrt{4 \pi} \left( - \sin\left(y\right), \sin \left(x \right), 0 \right)^{T}, \qquad
	\magfield\left(\xx,0\right)   =  \left( -  \sin\left(y\right), \sin \left(2x \right), 0\right)^{T},  \\
	\press\left(\xx,0\right) =\press_0 + \frac{15}{4} + \frac{1}{4} \cos(4x) + \frac{4}{5} \cos(2x) \cos(y) - \cos(x) \cos(y) + \frac{1}{4} \cos(2y),
 \end{align*} 
while the computational domain is the two-dimensional box $\Omega=[0,2\pi]^2$ that has been discretized with an unstructured mixed-element mesh composed of skewed quadrilaterals and triangles. The initial condition is taken from \cite{SIMHD,Fambri20}, but a very-high constant pressure $\press_0=10^5$ has been added in the background to reach the incompressible regime in the compressible flow solvers \cite{SIMHD,Fambri20}. The mesh is built by choosing $200$ points on every (periodic) boundary edge and about $200\sqrt{2}$ points on the diagonal that separate the skewed quads (south-east region) from the triangles (north-east region). A numerical reference solution has been computed by using the divergence-free semi-implicit method presented in \cite{SIMHD} on a uniform Cartesian grid composed of $1000^2$ elements. Here, the fluid viscosity and the magnetic resistivity are chosen to be $\nu=\mu=10^{-2}$.
The stream-lines of velocity and the magnetic field lines of the numerical solution obtained with the new hybrid FV/FE method are plotted in Figure \ref{fig.OT_finaltimesol}, achieving a very good agreement with the reference solution.
 
\subsection{Grad-Shafranov equilibrium in a torus} 

The aim of this last section is to show the applicability of the new well-balanced and divergence-free semi-implicit hybrid FV/FE scheme, presented in this paper, to more complex 3D problems related to continuum modeling of plasmas in 3D tokamak geometries and magnetic confinement fusion (MCF) applications. In particular, we show that the scheme is well-balanced and long-time stable for a non-trivial steady state solution on a 3D torus-shaped domain, including perfectly conducting wall boundary conditions. In the first part, we outline how the analytical solution is calculated and we then show the obtained numerical results. 

The so-called MHD equilibrium describes the steady  state of a plasma that is confined by a strong torus-shaped magnetic field, either in tokamak devices, being axisymmetric, or non-axisymmetric stellarator devices.  
It is derived from the ideal MHD equations ($\mu = \eta = 0$) with the assumption $\bvel=0$,  which subsequently reduce to the force balance equation 
$(\curl\bbvar)\times\bbvar=
\nabla\press$. Hence, the magnetic forces balance the pressure gradient inside the plasma, and the contours of constant pressure form a set of nested tori, on which the magnetic field is tangent ($\bbvar\cdot\nabla\press=0$).
Finding solutions of the force balance equation is not trivial, and is simplified by restricting to toroidal axisymmetric solutions in cylindrical coordinates $(R,Z,\phi)$, see Freidberg~\cite{FreidbergIdealMHD} for a derivation. The axisymmetric magnetic field is defined as
\begin{equation}
  \bbvar= \frac{F(\polflux)}{R} \mathrm{\mathbf{e}}_\phi + \frac{1}{R}\nabla \polflux \times \mathrm{\mathbf{e}}_\phi\,,
\end{equation}
with unit vector $\mathrm{\mathbf{e}}_\phi$, the net poloidal current profile $F(\polflux)=-\frac{I_p}{2\pi}$ and the poloidal magnetic flux $\polflux(R,Z)$, which is the solution of the Grad-Shafranov (GS) equation, see \cite{GradRubin,Shafranov},   
\begin{equation}
 R \partial_R\left(\frac{1}{R}\partial_R \polflux\right)+ \partial^2_Z \polflux = -\mu_0 R \frac{\partial p(\polflux)}{\partial\polflux}- \halb \frac{\partial F(\polflux)^2}{\partial \polflux}, \label{eq:GS}
\end{equation}
with the pressure profile $p(\polflux)$. In general, the GS equation is non-linear since the right-hand side depends on the solution. Under the assumption of Soloviev profiles \cite{solov}, which are linear $p\sim\polflux$, $F\sim\polflux$,
the GS equation becomes a linear inhomogeneous PDE. Under this assumption, Cerfon and Freidberg~\cite{CerfonFreidbergPOP2010} constructed a class of analytical equilibrium solutions of \eqref{eq:GS} as a series of polynomials. The coefficients of the series are defined by a profile parameter $C$, a global scaling factor $\Psi_0$, and a parametrized shape of the contour $\polflux=0$, which becomes the boundary of the plasma domain. The boundary is parametrized as
\begin{equation}
 R=R_0\left(1+\varepsilon\cos\left(\vartheta + \arcsin(\delta)\sin(\vartheta)\right)\right)\,,\quad Z=R_0\varepsilon\kappa\sin(\vartheta)\,,\quad \vartheta\in[0,2\pi], \label{eq:tok_bound}
\end{equation}
with major radius $R_0$, inverse aspect ratio $\varepsilon$, triangularity $\delta$, and elongation $\kappa$. All parameter symbols are chosen to match those in \cite{CerfonFreidbergPOP2010}.

For our numerical tests, we define two analytic equilibrium solutions with a circular and a D-shaped cross sections, which we will name the circular and D-shaped tokamak, respectively.
The parameters are summarized in Table~\ref{tab:tokamak_parameters}, with additional parameters for the profile definition. The profiles are chosen as
\begin{equation}
   F(\polflux):=F_0 \,,\quad p(\polflux):= p_0 + \frac{\polflux_0}{R_0^4}(\polflux_a-\polflux)\,,\quad \polflux_a=\polflux(R_a,Z_a), \label{eq:tok_profiles}
\end{equation}
with constants $p_0$, $F_0$ and the poloidal flux evaluated at the magnetic axis.
Furthermore, the Cartesian components of the vector potential can be expressed analytically as 
\begin{equation}
  \mathbf{A} = (\avar_1,\avar_2,\avar_3)=\left(\frac{-x_2\polflux}{R^2}\,,\,\frac{x_1\polflux}{R^2}\,,\,F_0\ln(R)\right)\,.
\end{equation}
The vector potential is needed for the initialization of a discretely divergence-free magnetic field, which must be computed by taking the discrete curl of the provided vector potential. 

\begin{figure}
 	\centering 
	 \includegraphics[width=0.45\linewidth]{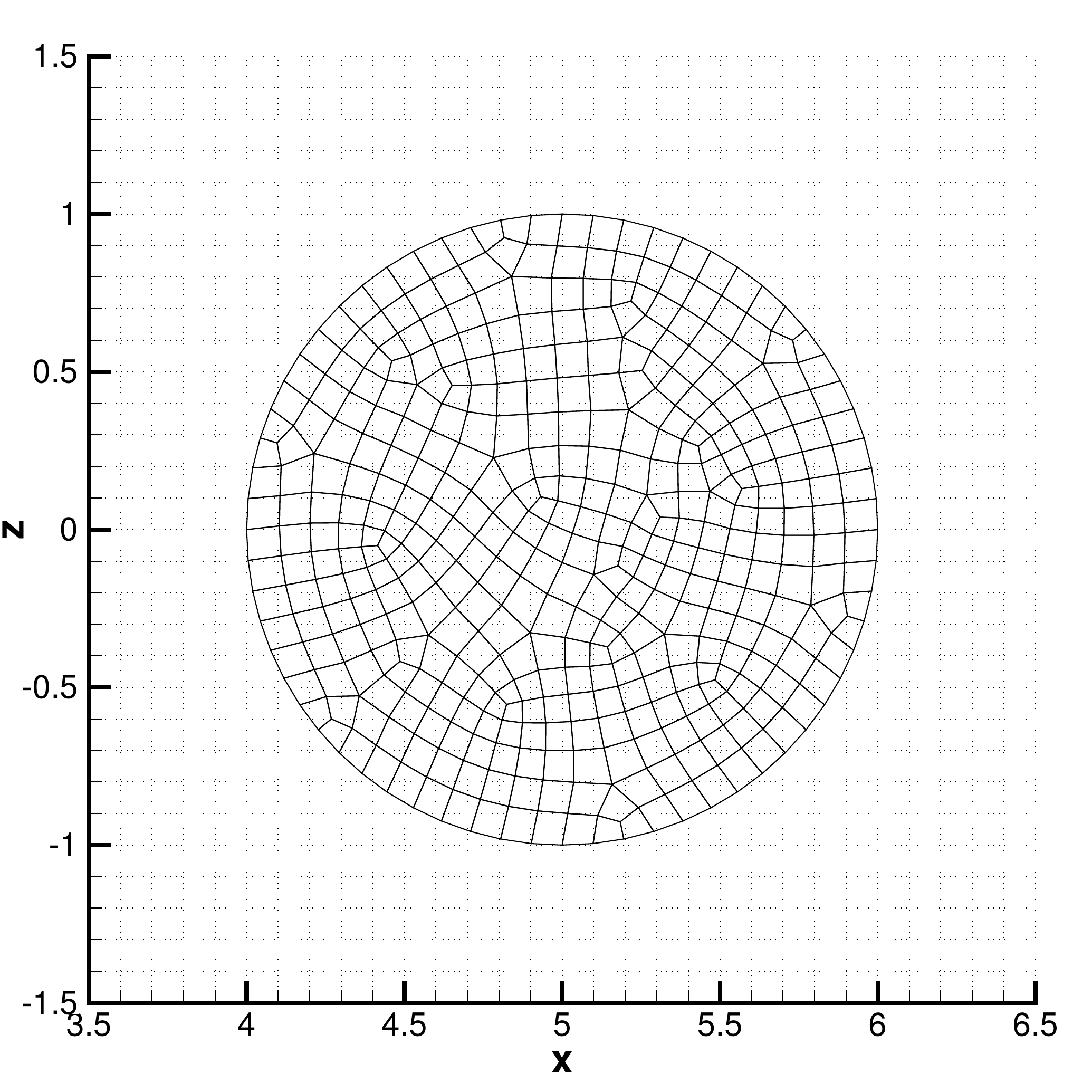}  
	\includegraphics[width=0.45\linewidth]{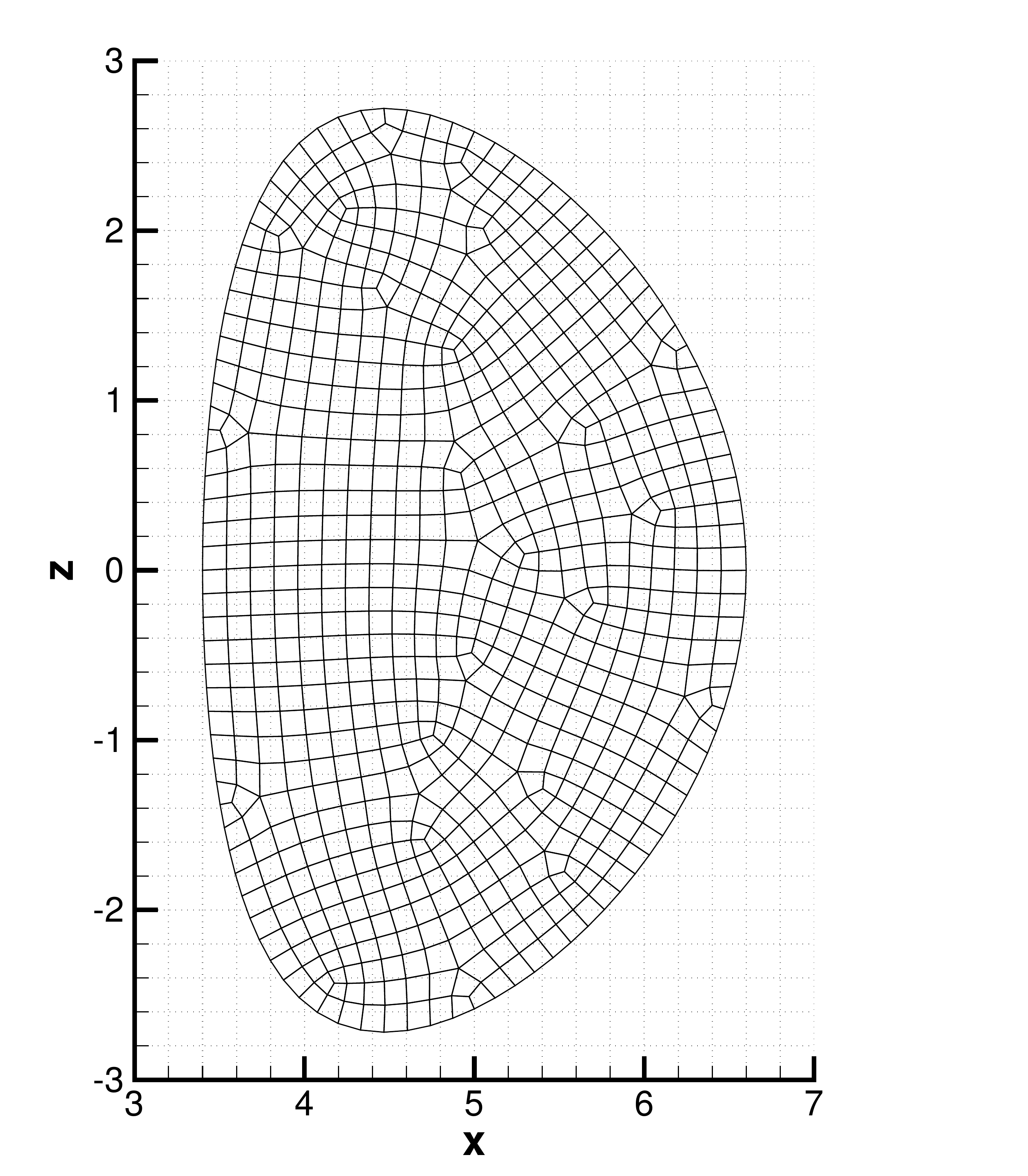} 
	\caption{ Cross section of the mesh of the circular (left) and D-shaped (right) tokamak, for which the analytical equilibrium solutions are defined. The boundary is also analytically defined, \eqref{eq:tok_bound}. }  
	\label{fig:tokamak_crossections}
\end{figure}

\begin{figure}
 	\centering 
	\includegraphics[clip,trim= 5 -500 5 5,width=0.35\linewidth]{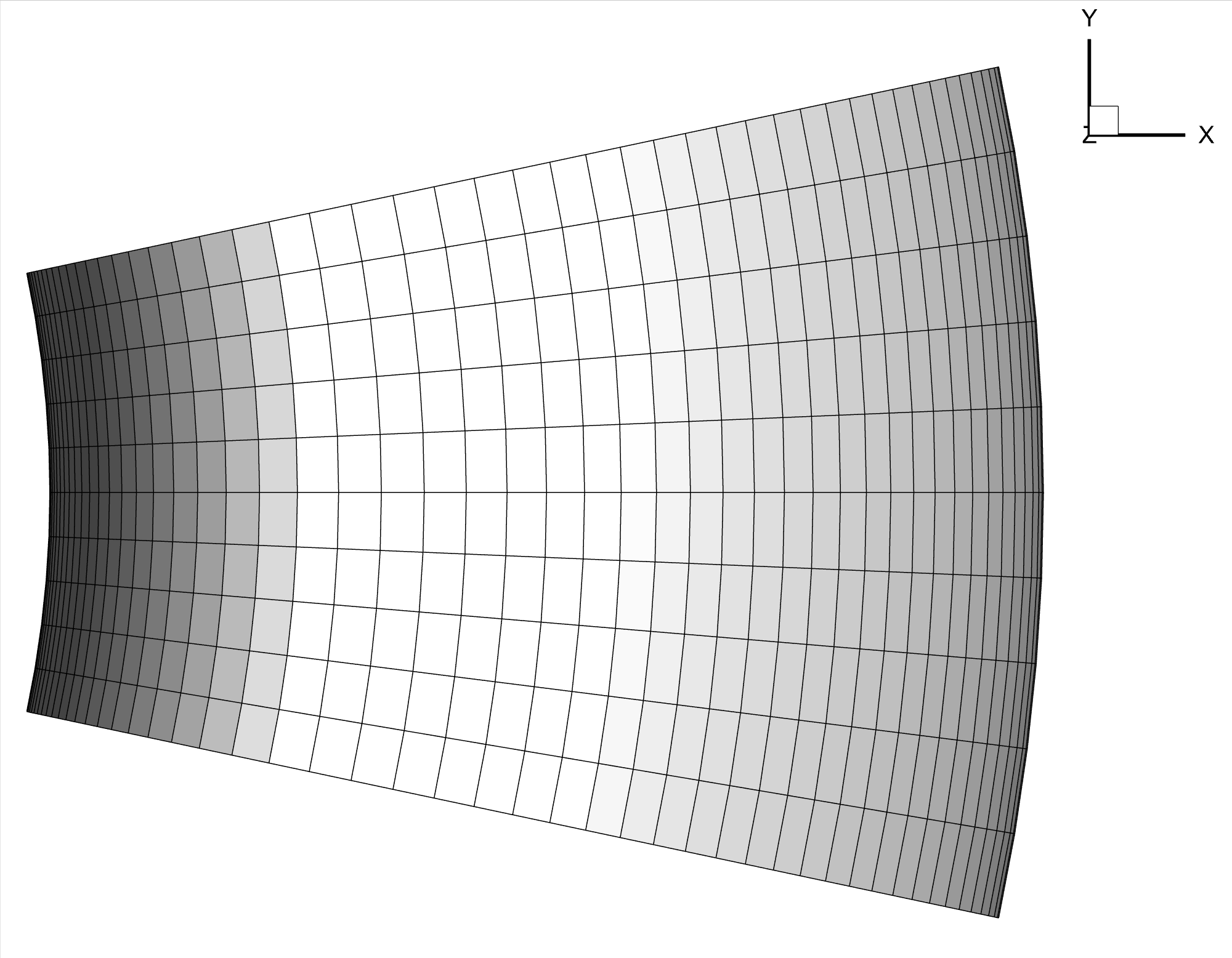}  \hspace{0.5cm}
	\includegraphics[clip,trim= 5 5 5 5,width=0.35\linewidth]{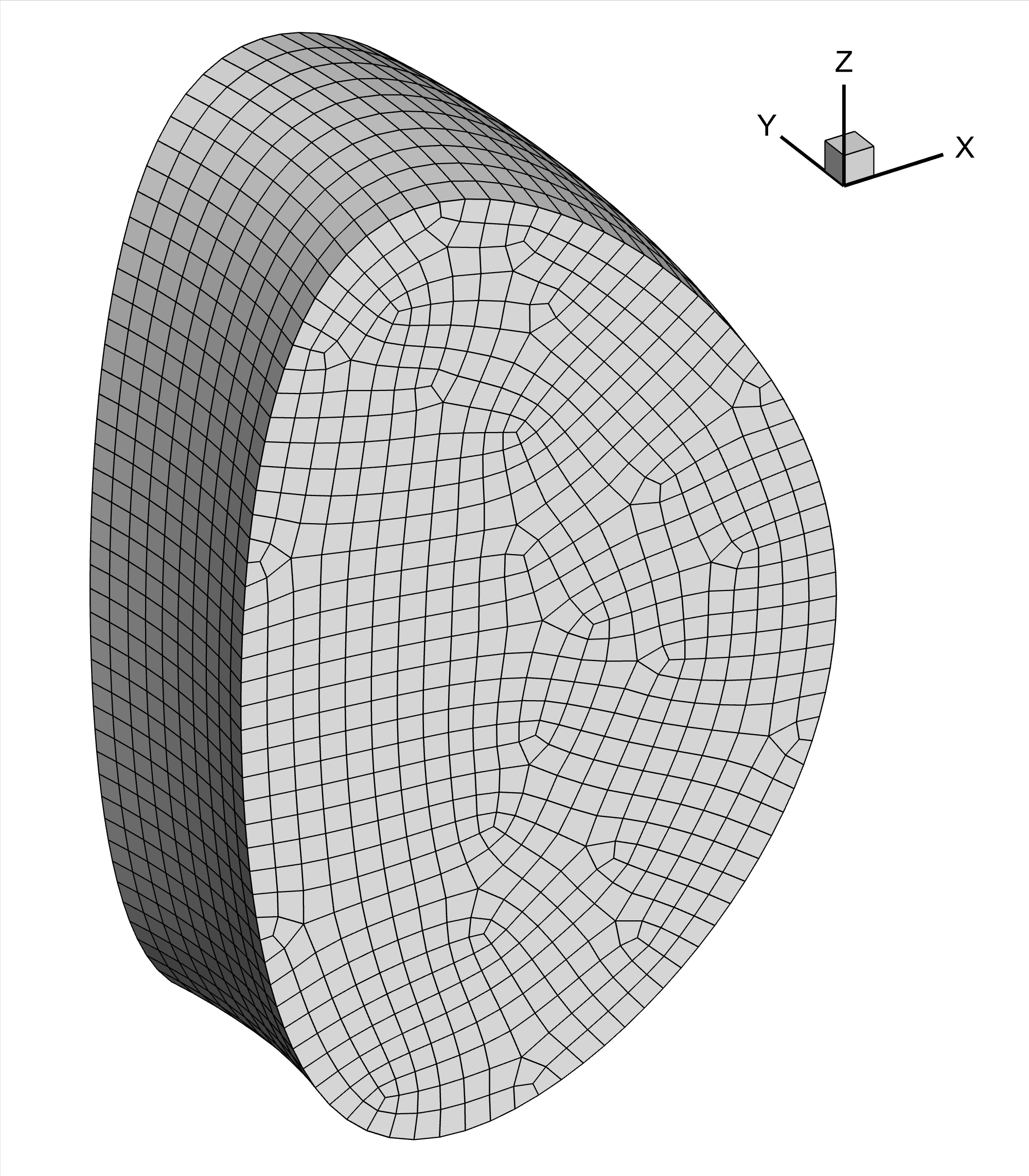} 
	\caption{Two different axonometries of the same 3D unstructured linear grid for the D-shaped tokamak, with $\phi \in [-\pi,\pi]/15$, a characteristic mesh size $h=0.1$, $N_\phi=10$.}  
	\label{fig:tokamak_mesh}
\end{figure}

The physical domain is a torus section with angle $\phi \in [\alpha, \beta]$ with (rotation-invariant) periodic boundary conditions, where $\beta-\alpha\in[0,2\pi]$. 
The computational domain is built starting from the discretization of the cross-section. A two-dimensional unstructured linear (non-curved) grid of quadrilaterals is built according to a prescribed characteristic mesh size $h$.
The discretized cross-sections shown in Figure~\ref{fig:tokamak_crossections}  correspond to a characteristic mesh-size $h=0.1$, for the circular and the D-shaped tokamak, respectively.
Then, a linear extrusion is performed in the toroidal direction by considering angles $\Delta\phi = (\beta-\alpha)/N_\phi$,   $N_\phi$ being the discretization number in the toroidal direction. The resulting 3D mesh is  a linear unstructured grid composed of hexahedra. In Figure \ref{fig:tokamak_mesh}, two different axonometries of the same 3D mesh for a section of a D-shaped tokamak are shown. In case the cross section is meshed via a mixed-element mesh of triangles and quadrilaterals, the resulting 3D mesh will be composed of hexahedra and triangular prisms.  

\begin{table}
\centering\begin{small}
 \begin{tabular}{|l|c|c|c|c|c|c|c||c|c|c|}
    \hline \rule[-1ex]{0pt}{3.5ex} 
    cross-section& $\mu_0$ & $R_0$ & $\varepsilon$ & $\kappa$ & $\delta$ & $C$ &$\polflux_0$ & $F_0$ & $p_0$ &  $(R_a,Z_a)$ \\\hline \rule[-1ex]{0pt}{3.5ex} 
    circular &  $1.$ & $5.$ & $0.20$ & $1.0$ & $0.0$ & $1.$ & $29.07886925004631$ & $5.$ & $0.02$ & $(5.12189635817554,0.)$ \\\hline\rule[-1ex]{0pt}{3.5ex} 
    D-shaped &  $1.$ & $5.$ & $0.32$ & $1.7$ & $0.33$ & $1.$ & $30.31531619122962$ & $5.$ & $0.08$ & $(5.22719821635024,0.)$ \\\hline
 \end{tabular}\end{small}
 \caption{Parameters that characterize the analytical equilibrium solutions from \cite{CerfonFreidbergPOP2010} and additional parameters for the profile definition, \eqref{eq:tok_profiles}.}  
 \label{tab:tokamak_parameters}
\end{table}

In the following numerical tests, long-time 3D simulations with and without well-balancing are compared with each other, for both, the circular and the D-shaped tokamaks. In all cases, perfectly-conducting wall boundary conditions are imposed. 

\begin{figure}
	\centering
	\includegraphics[width=0.7\linewidth]{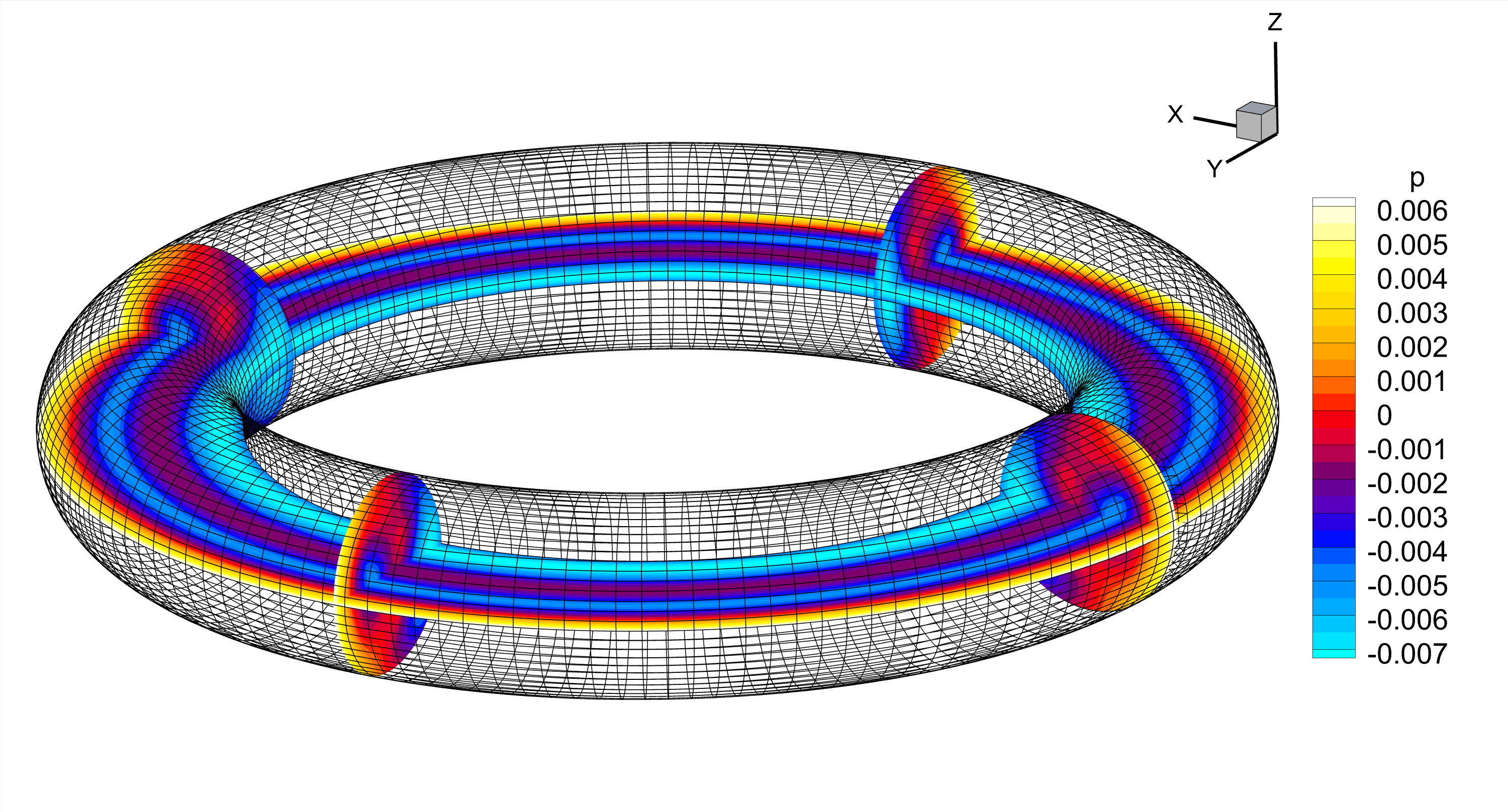} \\
	\includegraphics[width=0.7\linewidth]{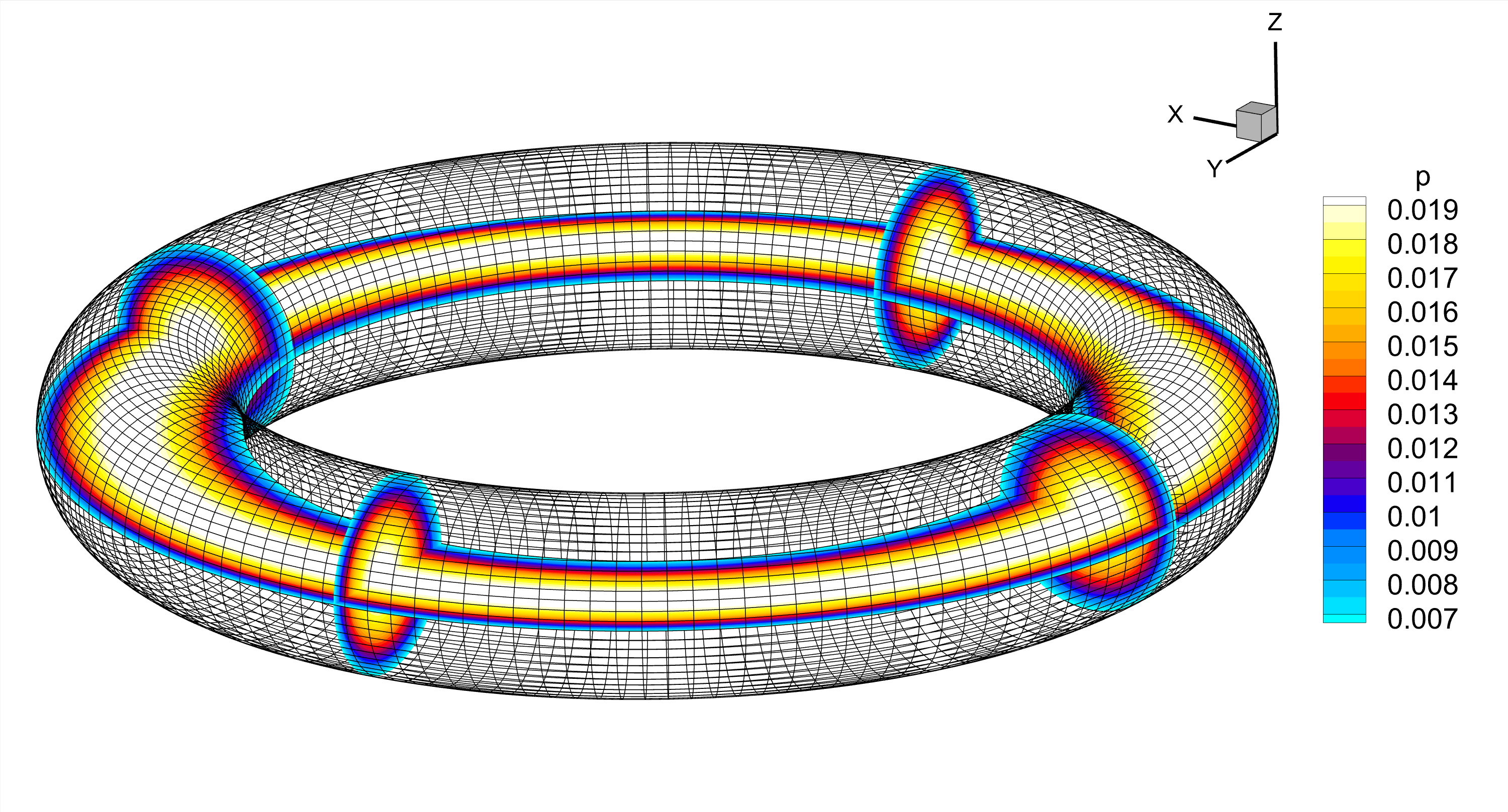}  
	\caption{Numerical results at $t=1000$ for the long time 3D simulation of the equilibrium solution in a circular tokamak. The numerical mesh is sketched together with the pressure contours: non well-balanced scheme (top) and the well-balanced scheme (bottom). For this test, the full geometry has been simulated.}  
	\label{fig:CSTokamak_3d}
\end{figure}

\begin{figure}
	\centering 
	\includegraphics[width=0.7\linewidth]{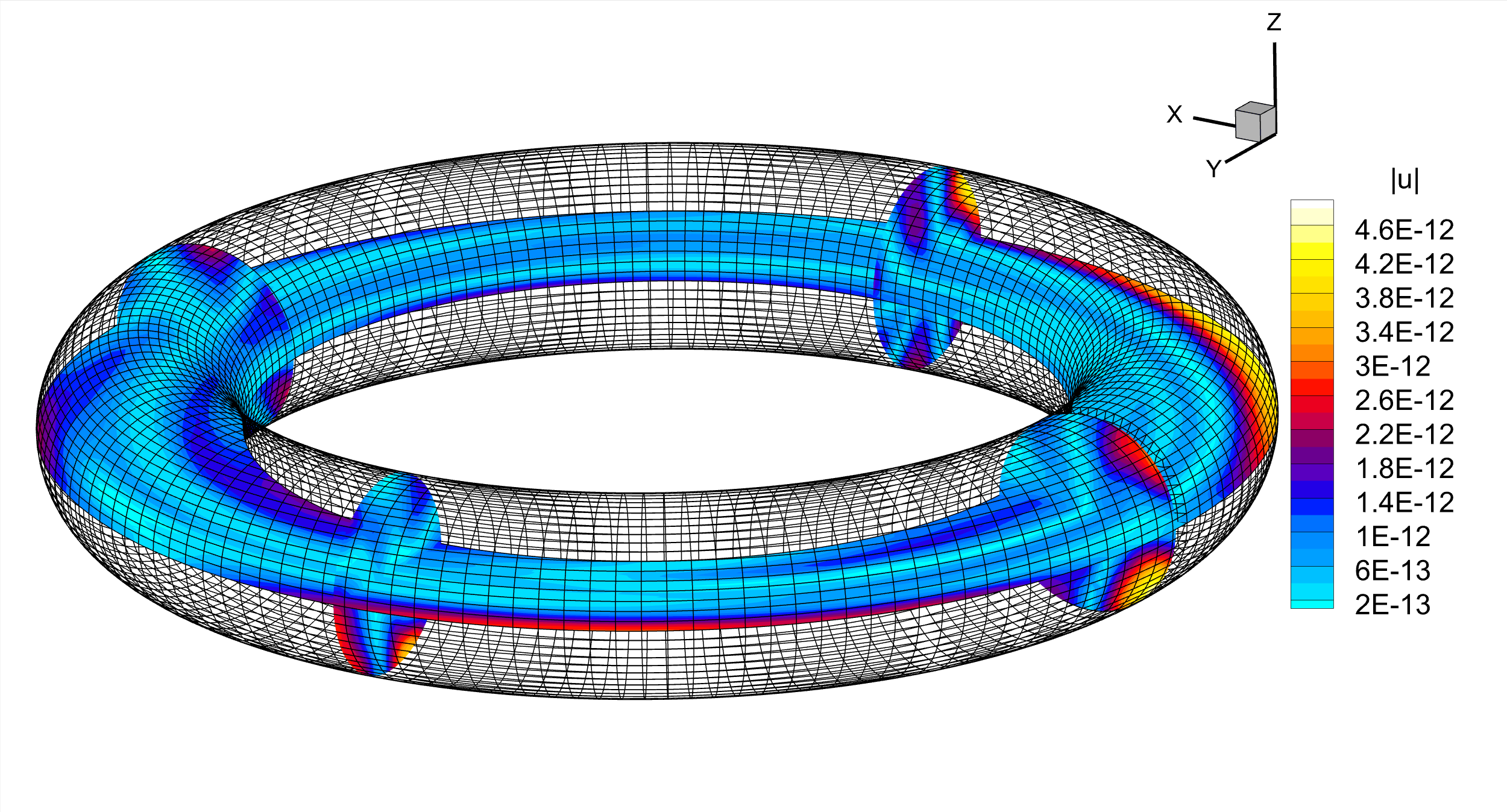} \\
	\includegraphics[width=0.7\linewidth]{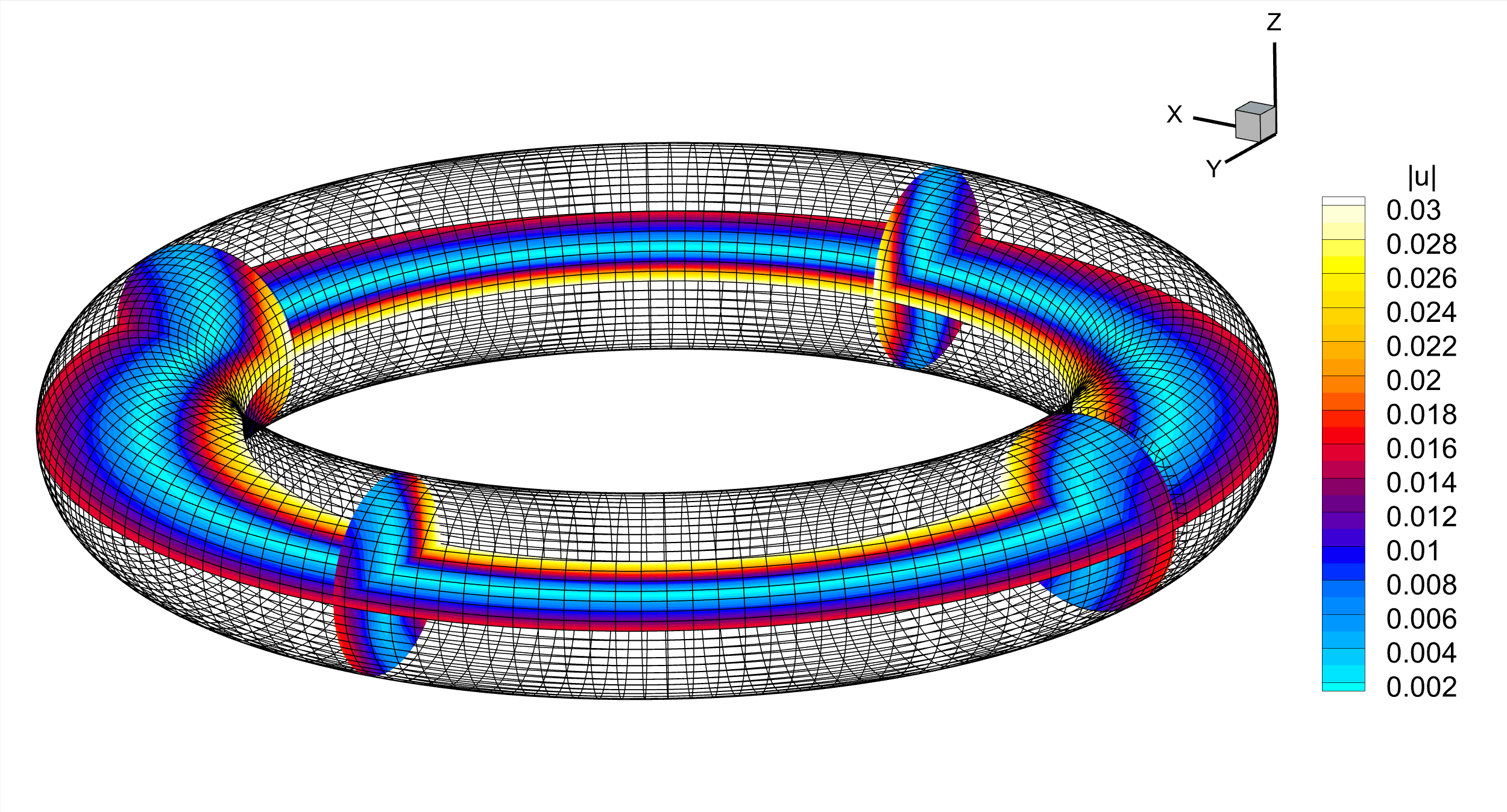} 
	\caption{Numerical results  at $t=1000$ for the long-time 3D simulation of the equilibrium solution in a circular tokamak. The numerical mesh is sketched together with the absolute value of the velocity: well-balanced scheme (top) and non well-balanced scheme (bottom). For this test, the full geometry has been simulated.}  
	\label{fig:CSTokamak_vel}
\end{figure}

\begin{figure}
	\centering
	\includegraphics[width=0.49\linewidth]{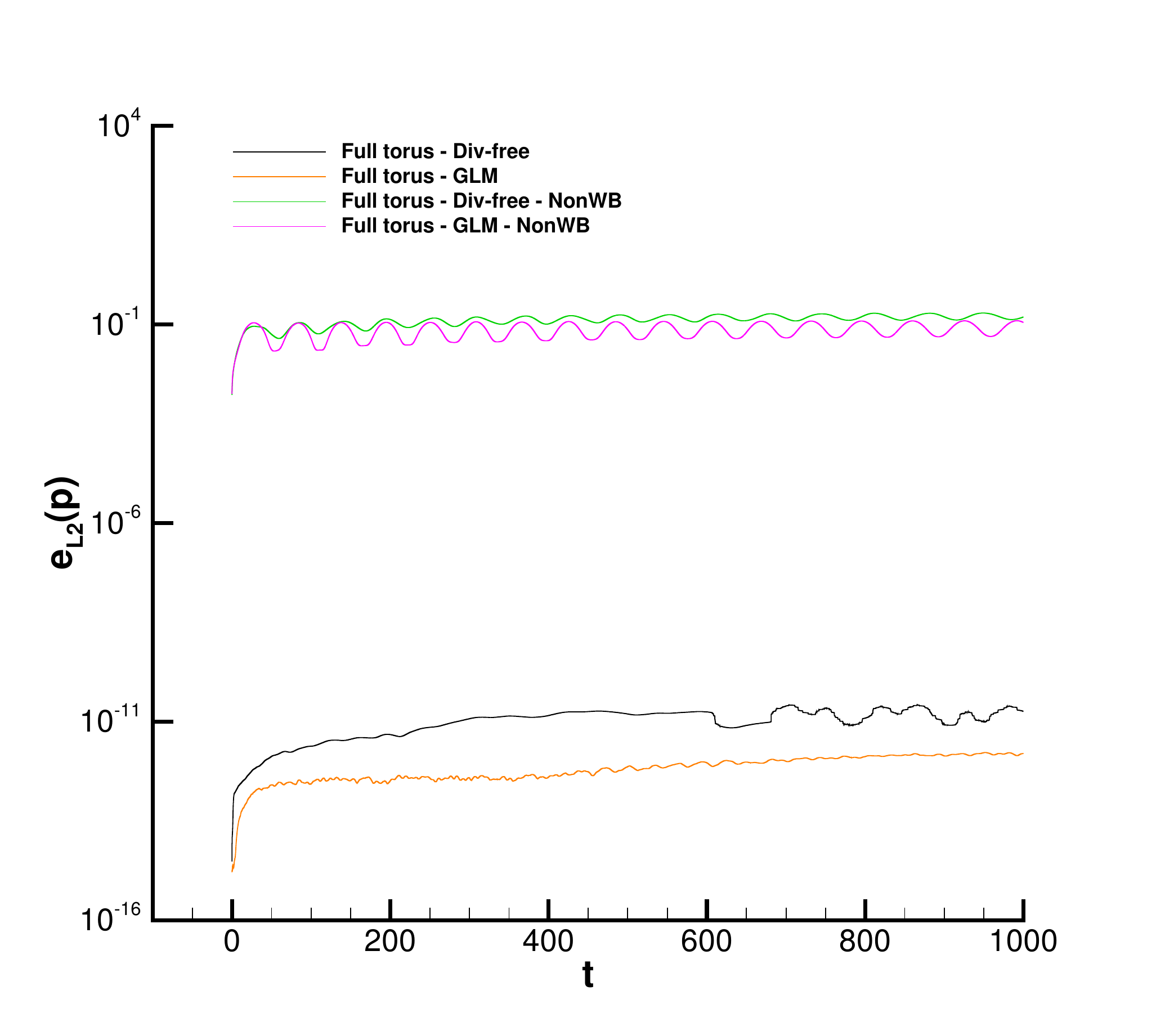} 
	\includegraphics[width=0.49\linewidth]{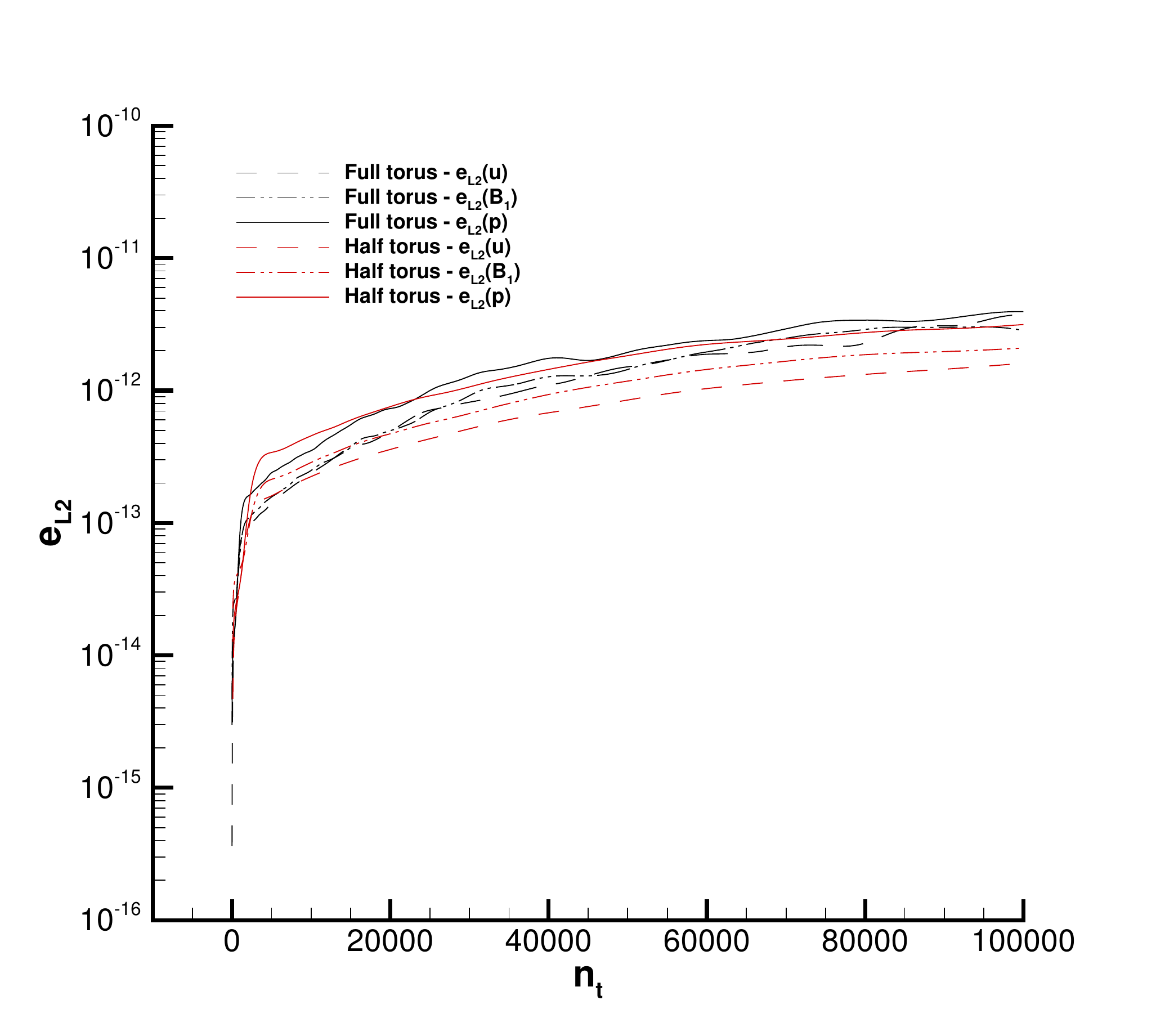}
	\caption{Numerical results for the long-time 3D simulation of the equilibrium solution in a circular tokamak. The time evolution for the $L^2$ norms of the error of some primitive variables is plotted. Left: $L^2$ norms of the error in the pressure field for the well-balanced versus the non-well-balanced scheme as a function of time up to $t=1000$, together with the corresponding solution obtained by using the hyperbolic GLM divergence cleaning approach. Right: $L^2$ norms of the error in $(u,p,B_1)$ for the well-balanced versus the non-well-balanced scheme as a function of the time-steps $n_t$ until $n_t=10^6$.}   
	\label{fig:CSTokamak_1d}
\end{figure}

In Figures \ref{fig:CSTokamak_3d}-\ref{fig:CSTokamak_vel} the contour plots of velocity and pressure along some poloidal and toroidal planes are shown for the numerical solution at the final time $t=1000$ of the prescribed equilibrium solution. 
 These solutions have been obtained after choosing a characteristic mesh size $h=0.1$ for the 2d discretization of the circular cross-section, while $N_\phi=150$ elements are used in the toroidal direction.
As expected, for long times the non well-balanced scheme significantly deviates from the equilibrium, while the well-balanced scheme, despite the use of a very coarse grid that is not necessarily aligned with the magnetic field lines, can properly preserve the stationary equilibrium solution up to machine precision. For a quantitative analysis, the time evolution of the $L^2$-error norms of some physical variables are plotted in Figure \ref{fig:CSTokamak_1d}. In the left subplot, we compare the results obtained for four different long-time simulations of the same equilibrium solution: the new exactly divergence-free scheme proposed in this paper versus the simpler non exactly divergence-free scheme based on the hyperbolic GLM-cleaning procedure \cite{MunzCleaning,Dedneretal} are tested and compared with each other, with and without well balancing. The stationary equilibrium solution is properly preserved for both well-balanced schemes, independently of the precise strategy that is used for the treatment of the divergence-free condition of the magnetic field, \eqref{eq:divB}, at the discrete level.
In the right subplot, the full torus simulation with $h=0.1$ is compared with a half-torus simulation that has a higher mesh resolution $h=0.05$, to show that the numerical error of the well-balanced scheme is stable in both cases and grows only linearly with the number of time-steps, due to the accumulation of round-off errors related to finite precision arithmetics. Indeed, for the finer mesh a more restrictive CFL condition applies, and a smaller time-step was used for the simulation. This means that the growth of the errors of the well-balanced simulation of an equilibrium is not improved by a mesh-refinement, but depends only on the number of operations performed within every single time-step.

\begin{figure}
	\centering
	\includegraphics[width=0.4\linewidth]{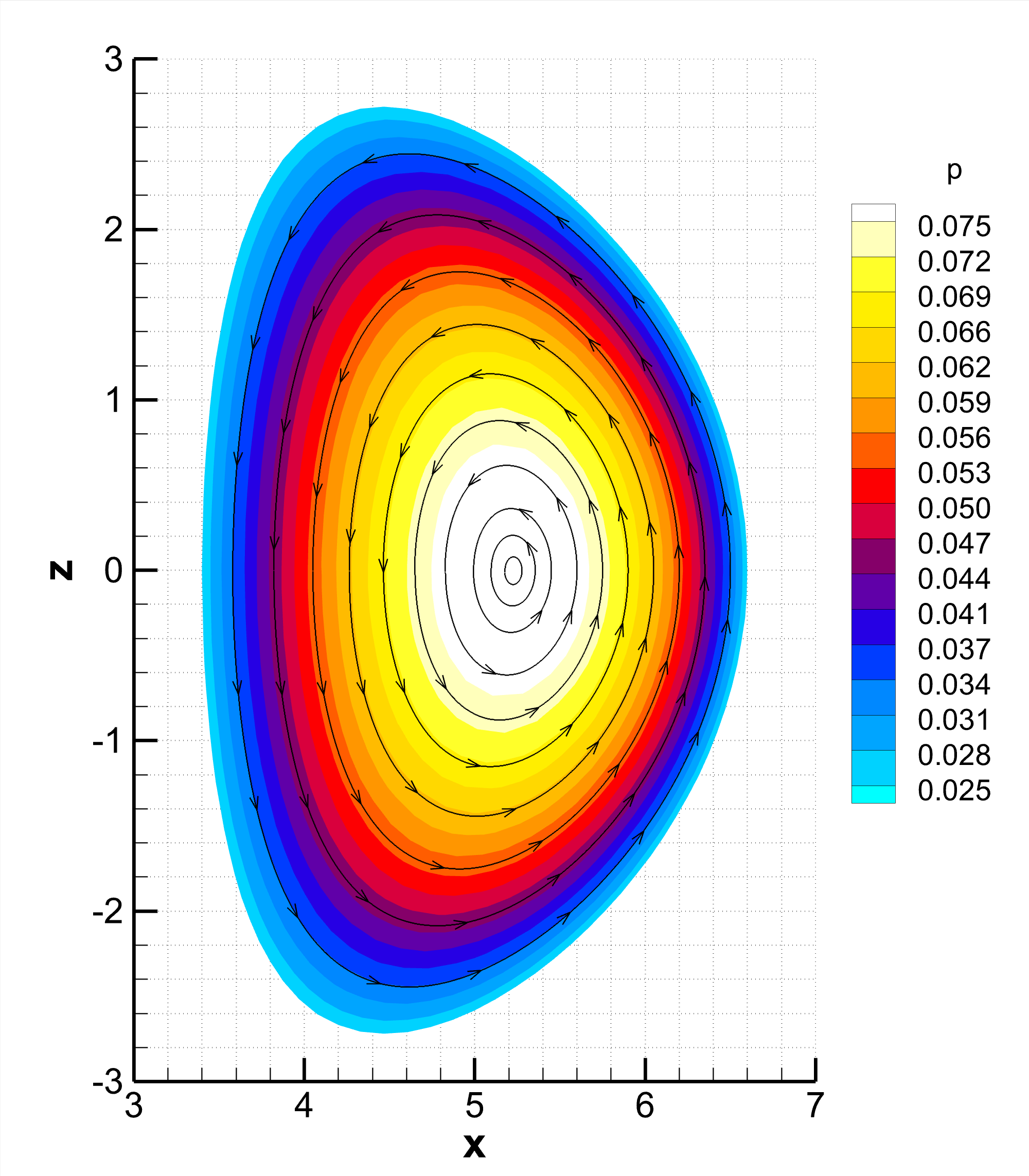}\hspace{0.02\linewidth}
	\includegraphics[width=0.4\linewidth]{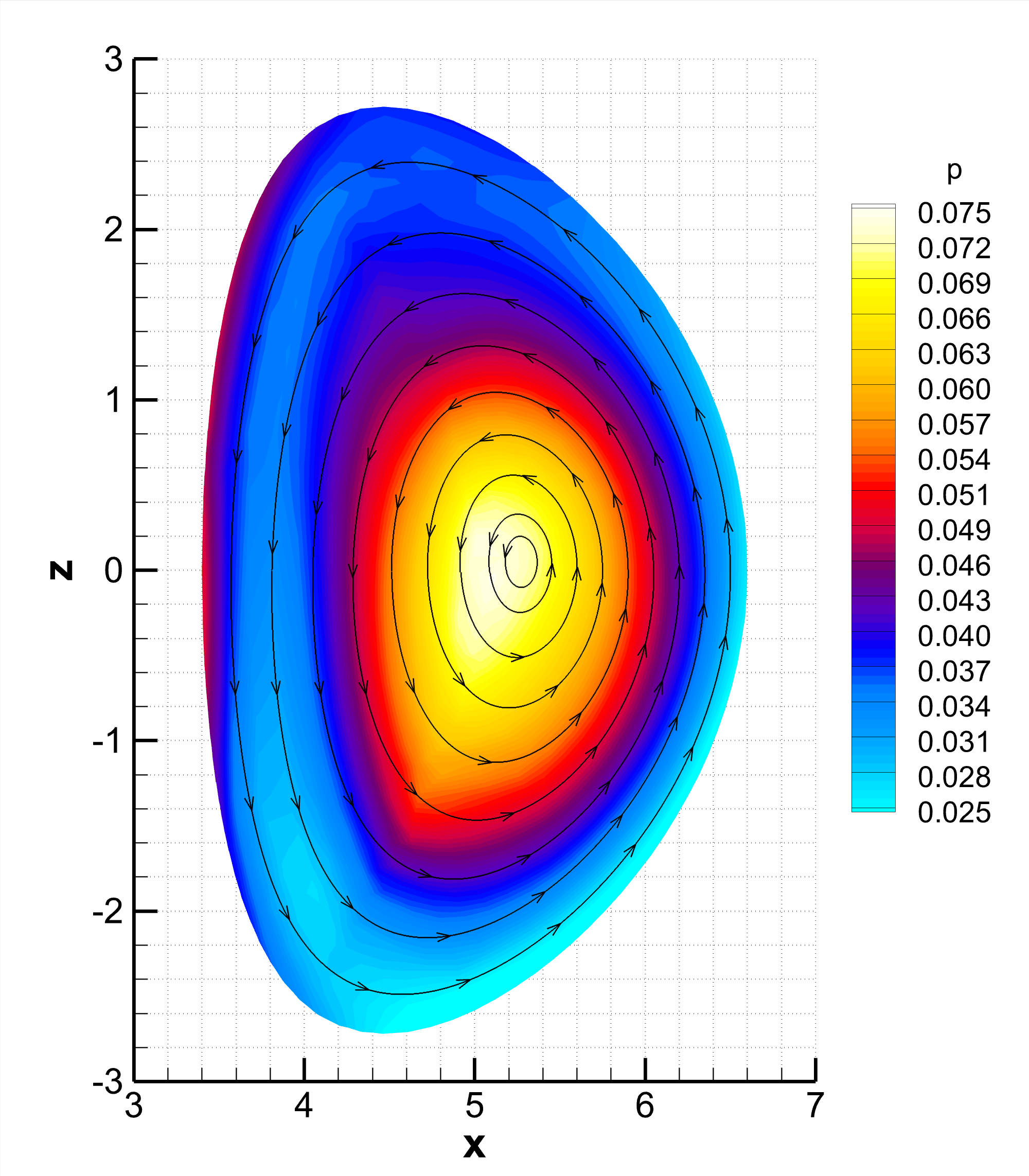} 
	\vspace{0.02\linewidth}\\
	\includegraphics[width=0.4\linewidth]{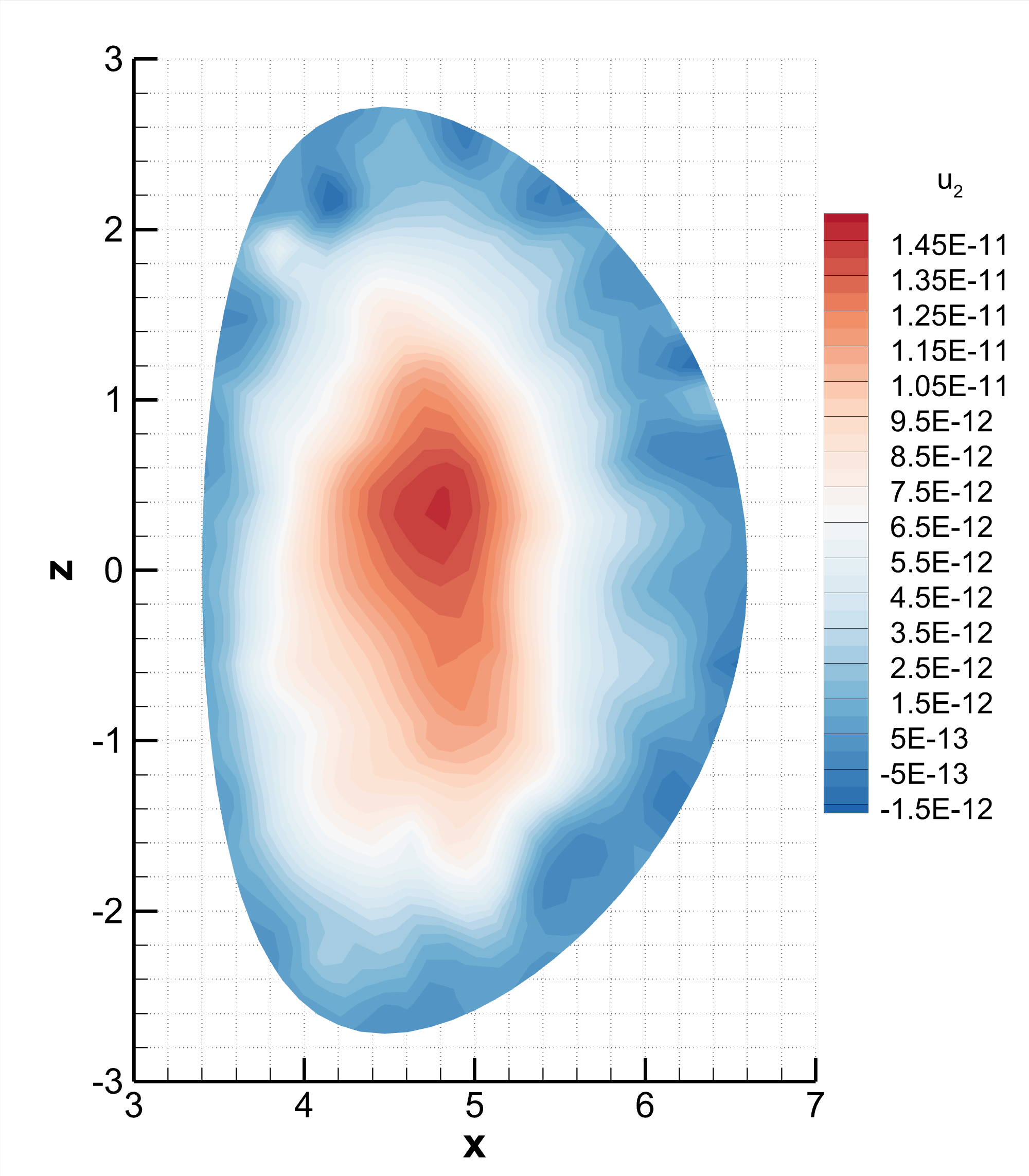}\hspace{0.02\linewidth}
	\includegraphics[width=0.4\linewidth]{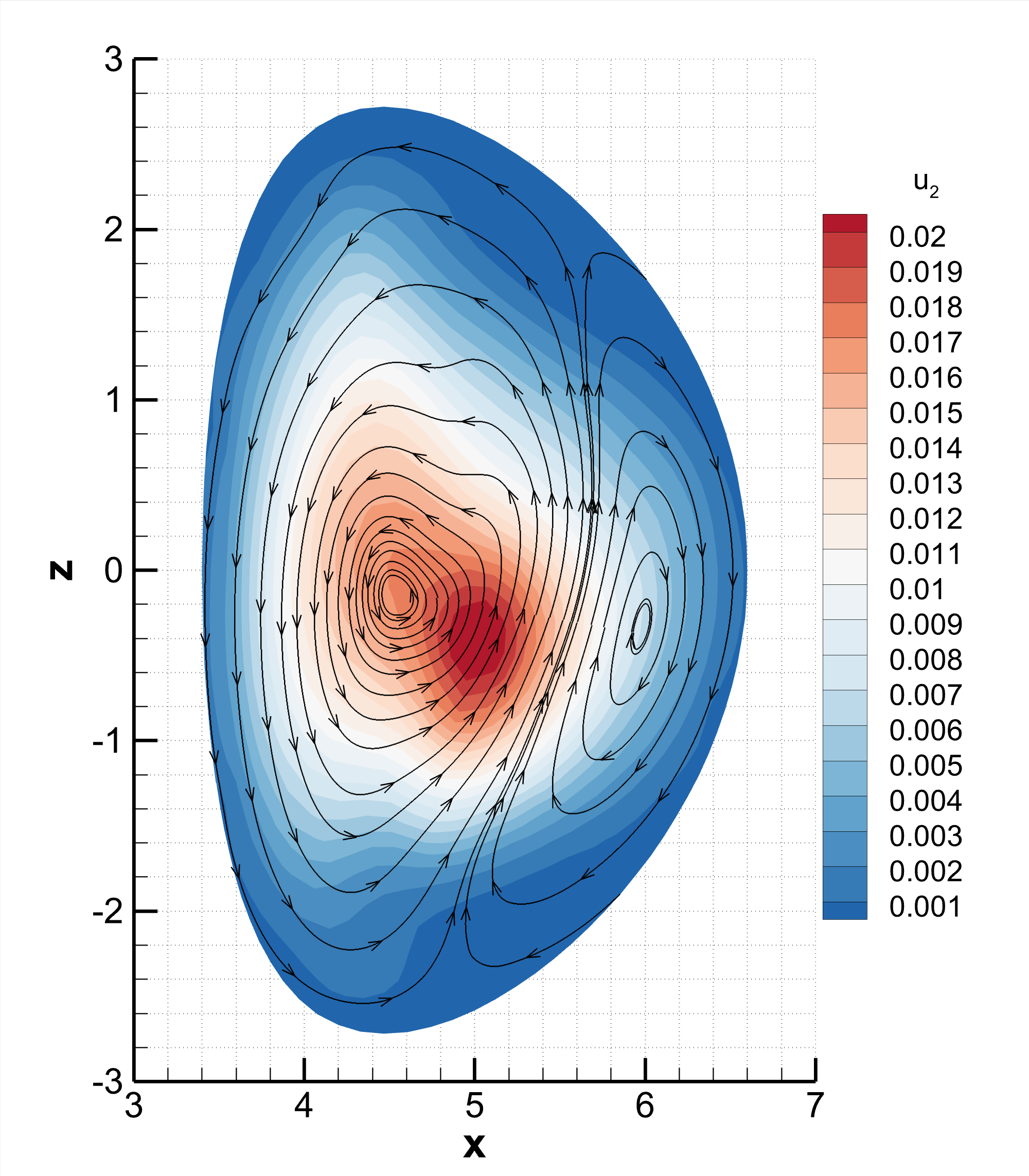} 
	\caption{Numerical results for the long-time 3D simulation of the equilibrium solution in a D-shaped tokamak at $t=1000$. Pressure contours together with the stream-lines of the poloidal magnetic field (top). Contours of the toroidal component of velocity ($u_2$) together with the stream-lines of the poloidal vector field $(u_1,u_3)$ (bottom). In both rows, the left plot corresponds to the well-balanced scheme, while the right image reports the solution for the non well-balanced scheme. For this test, only a small section, with periodic boundary conditions, is simulated. The plotted solution refers to a 2D cut at y=0. Velocity stream-lines are omitted for the well-balanced scheme because the velocity is zero up to round-off errors.}  
	\label{fig:DSTokamak_2d}
\end{figure}

\begin{figure}
	\centering
	\includegraphics[width=0.47\linewidth]{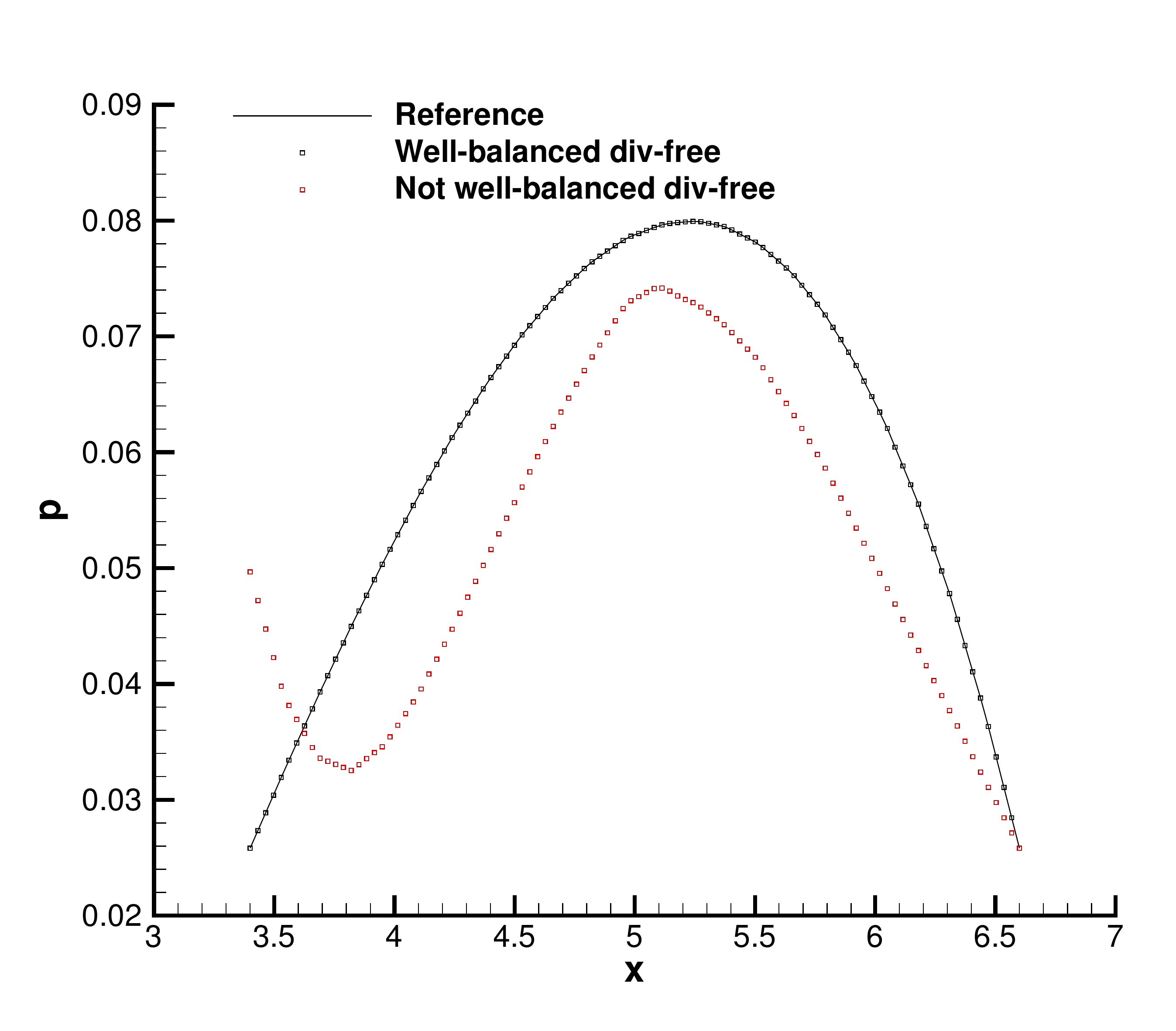}\hspace{0.02\linewidth}
	\includegraphics[width=0.47\linewidth]{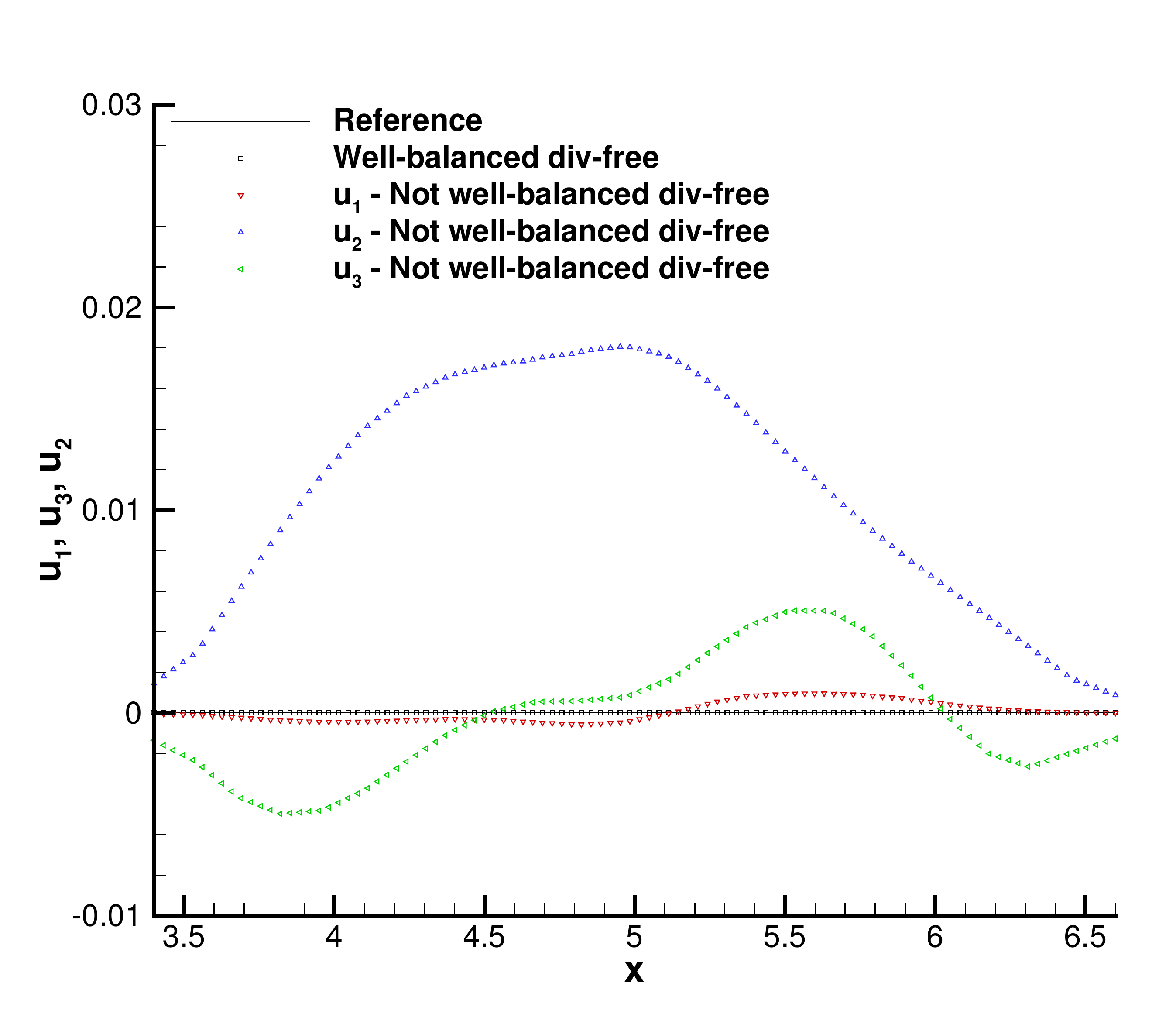} 
	\caption{Numerical results for the long-time 3D simulation of the equilibrium solution in a D-shaped tokamak at $t=1000$. The numerical solution has been interpolated over a uniform distribution of hundred points along the segment $x\in[3.4,6.6]$, $y=0$ and $z=0$. Pressure (left) and velocity (right) profiles are plotted to compare the well-balanced and the non-well-balanced scheme next to the reference solution. For the well-balanced scheme and the reference solution, only one curve is referred in the line-legend, because it coincides for the all three components up to round-off error.}  
	\label{fig:DSTokamak_1d_p}
\end{figure}

To complete this section, a long-time simulation of the equilibrium in a section of the D-shaped tokamak is depicted in Figure \ref{fig:tokamak_mesh}. Also in this case, the equilibrium is correctly preserved over long times when using the well-balanced scheme, with a stable linear growth of the $L^2$ error norms that are only due to round-off errors. In Figure \ref{fig:DSTokamak_2d}, the contour plots of pressure and velocity are shown. For the well-balanced scheme the magnetic field lines remain parallel to the contour lines of the pressure profile, while the velocity is a small perturbation of zero. On the contrary, for the non well-balanced scheme, the pressure profile and the contour lines of the magnetic-field are distorted, and spurious secondary-flows are generated. 
For a quantitative comparison, the solution is interpolated on a uniform distribution of hundred points along the segment in the poloidal plane with $x\in[3.4,6.6]$, $y=0$ and $z=0$. In Figure \ref{fig:DSTokamak_1d_p}, pressure and velocity profiles are plotted to compare the well-balanced and the non-well-balanced scheme next to the reference solution. The well-balanced solution cannot be distinguished from the prescribed reference equilibrium, as expected. Encouraged by these promising results, in the future we plan to apply our new scheme to non-ideal MHD in more complex 3D tokamak and stellarator geometries. 

\section{Conclusions}\label{sec:conclusions}

We have introduced a novel well-balanced and exactly divergence-free semi-implicit hybrid finite volume / finite element scheme for the solution of the incompressible viscous and resistive MHD equations on staggered unstructured mixed-element meshes in two and three space dimensions. The method is able to deal with mixed conforming triangular and quadrilateral meshes in 2D and with mixed conforming unstructured 3D meshes composed of tetrahedra, prisms, pyramids, and hexahedra. Since a staggered grid arrangement is quite common in both incompressible Navier-Stokes solvers, as well as in exactly divergence-free schemes for the Maxwell and MHD equations, the corresponding edge-based / face-based dual meshes are created by connecting each edge / face with the barycenters of the two elements sharing the same edge / face. In the paper, several examples have been shown on how to construct the appropriate dual meshes from a given mixed-element primary mesh. The approach presented here is a substantial generalization of the one introduced in \cite{BFSV14,BFTVC17,Hybrid1,Hybrid2}, which was so far limited to the use of simplex meshes. 

The overall governing PDE system is split into several subsystems, each of which is then discretized with a different numerical method that is appropriate for each subsystem. 
The non-linear convective and the viscous terms as well as the Maxwell stress tensor in the momentum equation have been discretized at the aid of a well-balanced explicit second order finite volume scheme that is run on the primal mesh. The resulting intermediate momentum, which does not yet take into account for the new pressure field, is then averaged onto the staggered dual mesh, where it is used in the computation of the right-hand side of the pressure Poisson equation. The latter is solved at the aid of classical conforming Lagrange finite elements, which are known to be well-suited for the solution of elliptic problems on general unstructured meshes. However, a peculiar feature of the pressure Poisson solver presented in this paper is that the scheme is exact for stationary equilibrium solutions. The magnetic field is evolved in an exactly divergence-free and well-balanced manner making use of a discrete form of the Stokes theorem inside the edges / faces of the primary mesh. The electric field, which is needed for this purpose, is computed in each vertex / edge by taking into account the data of all elements around a vertex / edge, adding the proper amount of numerical resistivity that is needed to stabilize the scheme and also adding the physical resistivity, present in the governing PDE system. To achieve high order of accuracy, piecewise linear polynomials of the magnetic field are reconstructed, which are then corrected via a constrained $L^2$ projection that guarantees that the obtained magnetic field is exactly divergence-free inside each element and that the integral average of its normal component is continuous across each edge / face, see \cite{divfree2015}.  
We stress again that each step of the algorithm presented in this paper is constructed in a well-balanced manner and is able to maintain any general but \textit{a priori} known stationary equilibrium solution exactly at the discrete level. 

We show numerous numerical test cases in two and three space dimensions in order to validate our new method carefully against known exact and numerical reference solutions. We also show numerical evidence that confirms the well-balancing of our new scheme and that the magnetic field remains exactly divergence-free. As a special feature, this paper includes a very thorough study of the lid-driven MHD cavity problem in the presence of different magnetic fields, and the obtained numerical solutions are provided as free supplementary electronic material to allow other research groups to reproduce our results and to compare with our data easily. 
In several numerical examples, we demonstrate that the well-balancing allows performing accurate and stable long-time simulations for simplified 3D tokamak geometries even on very coarse meshes and, in particular, on general unstructured meshes that do not need to be aligned with the magnetic field lines. In addition, we also show stable long-time simulations carried out by a simpler version of our scheme in which the magnetic field is not evolved in an exactly divergence-free manner but via the augmented generalized Lagrangian multiplier (GLM) divergence cleaning approach introduced in \cite{MunzCleaning,Dedneretal}. Our numerical results indicate that the well-balancing 
seems to be the key feature for accurate and stable long-time simulations on coarse meshes. 

Future work will concern the extension to higher order of accuracy in space and time \cite{TD16,BP2021,BDLTV2020} and to the weakly compressible case \cite{TD17,Hybrid1,Hybrid2}. We also plan to develop thermodynamically compatible and provably nonlinearly stable schemes following the framework introduced in \cite{HTCMHD,HTCAbgrall,ElenaHTC}, as well as a three-split semi-implicit extension that allows to get rid of the CFL time step restriction based on the Alfv\'en wave speed \cite{Fambri20}.    

%
\section*{Acknowledgements}
This work was financially supported by the Italian Ministry of Education, University 
and Research (MIUR) in the framework of the PRIN 2017 project \textit{Innovative numerical methods for evolutionary partial differential equations and  applications} and via the  Departments of Excellence  Initiative 2018--2027 attributed to DICAM of the University of Trento (grant L. 232/2016). S.B., M.D., L.R. and E.Z. also acknowledge the financial support from the Spanish Ministry of Science and Innovation, grant number PID2021-122625OB-I00. L.R. acknowledges funding from the Spanish Ministry of Universities and the European Union-Next GenerationEU under the project RSU.UDC.MS15. M.D., F.F. and L.R. are member of the GNCS-INdAM (Istituto Nazionale di Alta Matematica) group. M.D. was also co-funded by the European Union NextGenerationEU (PNRR, Spoke 7 CN HPC). Views and
opinions expressed are however those of the author(s) only and do not necessarily reflect those of the European Union or the European Research Council. Neither the European Union nor the granting authority can be
held responsible for them. 

\clearpage 

%
%
%
%
\bibliographystyle{plain}
\bibliography{../mibiblio} 

\begin{thebibliography}{100}

\bibitem{Abgrall1994}
R.~Abgrall.
\newblock On essentially non-oscillatory schemes on unstructured meshes:
  analysis and implementation.
\newblock {\em Journal of Computational Physics}, 114(1):45--58, 1994.

\bibitem{HTCAbgrall}
R.~Abgrall, S.~Busto, and M.~Dumbser.
\newblock {A simple and general framework for the construction of
  thermodynamically compatible schemes for computational fluid and solid
  mechanics}.
\newblock {\em Applied Mathematics and Computation}, 440:127629, 2023.

\bibitem{Abgrall99ENO}
R.~Abgrall, S.~Lant{\'e}ri, and T.~Sonar.
\newblock {ENO} approximations for compressible fluid dynamics.
\newblock {\em ZAMM-Journal of Applied Mathematics and Mechanics/Zeitschrift
  f{\"u}r Angewandte Mathematik und Mechanik: Applied Mathematics and
  Mechanics}, 79(1):3--28, 1999.

\bibitem{Amari}
T.~Amari, J.F. Luciani, and P.~Joly.
\newblock {Preconditioned semi--implicit method for magnetohydrodynamics
  equations}.
\newblock {\em SIAM Journal on Scientific Computing}, 21:970--986, 1999.

\bibitem{Balsara2004}
D.~Balsara.
\newblock Second-order accurate schemes for magnetohydrodynamics with
  divergence-free reconstruction.
\newblock {\em The Astrophysical Journal Supplement Series}, 151:149--184,
  2004.

\bibitem{BalsaraSpicer1999b}
D.~Balsara and D.~Spicer.
\newblock Maintaining pressure positivity in magneto\-hydrodynamic simulations.
\newblock {\em J. Comput. Phys.}, 148:133--148, 1999.

\bibitem{BalsaraSpicer1999}
D.~Balsara and D.~Spicer.
\newblock A staggered mesh algorithm using high order {G}odunov fluxes to
  ensure solenoidal magnetic fields in magnetohydrodynamic simulations.
\newblock {\em J. Comput. Phys.}, 149:270--292, 1999.

\bibitem{BalsaraAMR}
D.S. Balsara.
\newblock Divergence-free adaptive mesh refinement for magnetohydrodynamics.
\newblock {\em Journal of Computational Physics}, 174(2):614--648, 2001.

\bibitem{balsarahlle2d}
D.S. Balsara.
\newblock {Multidimensional HLLE Riemann solver: application to Euler and
  magnetohydrodynamic flows}.
\newblock {\em Journal of Computational Physics}, 229:1970--1993, 2010.

\bibitem{balsarahllc2d}
D.S. Balsara.
\newblock {A two-dimensional HLLC Riemann solver for conservation laws:
  application to Euler and magnetohydrodynamic flows}.
\newblock {\em Journal of Computational Physics}, 231:7476--7503, 2012.

\bibitem{MUSIC1}
D.S. Balsara.
\newblock {Multidimensional Riemann Problem with Self-Similar Internal
  Structure -- Part I -- Application to hyperbolic conservation laws on
  structured meshes}.
\newblock {\em Journal of Computational Physics}, 277:163--200, 2014.

\bibitem{balsarahlle3d}
D.S. Balsara.
\newblock {Three dimensional HLL Riemann solver for conservation laws on
  structured meshes: application to Euler and magnetohydrodynamic flows}.
\newblock {\em Journal of Computational Physics}, 295:1--23, 2015.

\bibitem{divfree2015}
D.S. Balsara and M.~Dumbser.
\newblock Divergence-free {MHD} on unstructured meshes using high order finite
  volume schemes based on multidimensional {R}iemann solvers.
\newblock {\em Journal of Computational Physics}, 299:687--715, 2015.

\bibitem{MUSIC2}
D.S. Balsara and M.~Dumbser.
\newblock {Multidimensional Riemann Problem with Self-Similar Internal
  Structure -- Part II -- Application to Hyperbolic Conservation Laws on
  Unstructured Meshes}.
\newblock {\em Journal of Computational Physics}, 287:269--292, 2015.

\bibitem{BalsaraMultiDRS}
D.S. Balsara, M.~Dumbser, and R.~Abgrall.
\newblock {Multidimensional HLLC Riemann solver for unstructured meshes - with
  application to Euler and MHD flows.}
\newblock {\em Journal of Computational Physics}, 261:172--208, 2014.

\bibitem{BJ89}
T.~Barth and D.~Jespersen.
\newblock The design and application of upwind schemes on unstructured meshes.
\newblock Technical report, 1989.

\bibitem{BarthFrederickson1990}
Timothy Barth and Paul Frederickson.
\newblock Higher order solution of the euler equations on unstructured grids
  using quadratic reconstruction.
\newblock In {\em 28th aerospace sciences meeting}, page~13, 1990.

\bibitem{BCG89}
J.~B. Bell, P.~Colella, and H.~M. Glaz.
\newblock A second-order projection method for the incompressible
  {N}avier--{S}tokes equations.
\newblock {\em Journal of Computational Physics}, 85(2):257--283, 1989.

\bibitem{BellColellaGlaz}
J.B. Bell, P.~Colella, and H.M. Glaz.
\newblock A second--order projection method for the incompressible
  {Navier--Stokes} equations.
\newblock {\em Journal of Computational Physics}, 85:257--283, 1989.

\bibitem{BCK21}
J.P. Berberich, P.~Chandrashekar, and C.~Klingenberg.
\newblock High order well-balanced finite volume methods for multi-dimensional
  systems of hyperbolic balance laws.
\newblock {\em Computers \& Fluids}, 219:104858, 2021.

\bibitem{Hybrid1}
A.~Berm\'udez, S.~Busto, M.~Dumbser, J.L. Ferr\'in, L.~Saavedra, and M.E.
  V\'azquez-Cend\'on.
\newblock A staggered semi-implicit hybrid {FV/FE} projection method for weakly
  compressible flows.
\newblock {\em J. Comput. Phys.}, 421:109743, 2020.

\bibitem{BFSV14}
A.~Berm\'udez, J.~L. Ferr\'in, L.~Saavedra, and M.~E. V\'azquez-Cend\'on.
\newblock A projection hybrid finite volume/element method for low-{M}ach
  number flows.
\newblock {\em Journal of Computational Physics}, 271:360--378, 2014.

\bibitem{Bermudez1994}
A.~Berm\'udez and M.~E. V\'azquez.
\newblock Upwind methods for hyperbolic conservation laws with source terms.
\newblock {\em Computers and Fluids}, 23:1049--1071, 1994.

\bibitem{BVC94}
A.~Berm{\'u}dez and M.~E. V{\'a}zquez-Cend{\'o}n.
\newblock Upwind methods for hyperbolic conservation laws with source terms.
\newblock {\em Comput. Fluids}, 23(8):1049--1071, 1994.

\bibitem{BDLTV2020}
W.~Boscheri, G.~Dimarco, R.~Loub{\`{e}}re, M.~Tavelli, and M.~H. Vignal.
\newblock A second order all mach number imex finite volume solver for the
  three dimensional euler equations.
\newblock {\em Journal of Computational Physics}, 415:109486, 2020.

\bibitem{BP2021}
W.~Boscheri and L.~Pareschi.
\newblock {High order pressure-based semi-implicit IMEX schemes for the 3D
  Navier--Stokes equations at all Mach numbers}.
\newblock {\em Journal of Computational Physics}, 434:110206, 2021.

\bibitem{BottaKlein}
N.~Botta, R.~Klein, S.~Langenberg, and S.~L\"utzenkirchen.
\newblock Well balanced finite volume methods for nearly hydrostatic flows.
\newblock {\em Journal of Computational Physics}, 196:539--565, 2004.

\bibitem{BrioWu}
M.~Brio and C.~C. Wu.
\newblock An upwind differencing scheme for the equations of ideal
  magnetohydrodynamics.
\newblock {\em Journal of Computational Physics}, 75:400--422, 1988.

\bibitem{HybridSWE}
S~Busto and M~Dumbser.
\newblock A staggered semi-implicit hybrid finite volume/finite element scheme
  for the shallow water equations at all {Froude} numbers.
\newblock {\em Applied Numerical Mathematics}, 175:108--132, 2022.

\bibitem{HTCMHD}
S.~Busto and M.~Dumbser.
\newblock {A new thermodynamically compatible finite volume scheme for
  magnetohydrodynamics}.
\newblock {\em SIAM Journal on Numerical Analysis}, 61:343--364, 2023.

\bibitem{HybridNNT}
S.~Busto, M.~Dumbser, and L.~R\'io-Mart\'in.
\newblock Staggered semi-implicit hybrid finite volume/finite element schemes
  for turbulent and non-{Newtonian} flows.
\newblock {\em Mathematics}, 9:2972, 2021.

\bibitem{HybridALE}
S~Busto, M~Dumbser, and L~R{\'\i}o-Mart{\'\i}n.
\newblock An {Arbitrary-Lagrangian-Eulerian} hybrid finite volume/finite
  element method on moving unstructured meshes for the {Navier--Stokes}
  equations.
\newblock {\em Applied Mathematics and Computation}, 437:127539, 2023.

\bibitem{BFTVC17}
S.~Busto, J.~L. Ferr\'in, E.~F. Toro, and M.~E. V\'azquez-Cend\'on.
\newblock A projection hybrid high order finite volume/finite element method
  for incompressible turbulent flows.
\newblock {\em Journal of Computational Physics}, 353:169--192, 2018.

\bibitem{Hybrid2}
S.~Busto, L.~R\'io-Mart\'in, M.E. V\'azquez-Cend\'on, and M.~Dumbser.
\newblock A semi-implicit hybrid finite volume / finite element scheme for all
  {Mach} number flows on staggered unstructured meshes.
\newblock {\em Appl. Math. Comput.}, 402:126117, 2021.

\bibitem{BTVC16}
S.~Busto, E.~F. Toro, and M.~E. V\'azquez-Cend\'on.
\newblock Design and analisis of {ADER}--type schemes for model
  advection--diffusion--reaction equations.
\newblock {\em Journal of Computational Physics}, 327:553--575, 2016.

\bibitem{CG23}
M.~Giuliano Carlino and E.~Gaburro.
\newblock Well balanced finite volume schemes for shallow water equations on
  manifolds.
\newblock {\em Applied Mathematics and Computation}, 441:127676, 2023.

\bibitem{Castro2008}
M.J. Castro, J.M. Gallardo, J.A. L\'opez, and C.~Par\'es.
\newblock Well-balanced high order extensions of godunov's method for
  semilinear balance laws.
\newblock {\em SIAM Journal of Numerical Analysis}, 46:1012--1039, 2008.

\bibitem{castro2007wellb}
M.J. Castro~D{\'\i}az, T.~Chac{\'o}n~Rebollo, E.D. Fern{\'a}ndez-Nieto, and
  C.~Pares.
\newblock On well-balanced finite volume methods for nonconservative
  nonhomogeneous hyperbolic systems.
\newblock {\em SIAM J. Sci. Comput.}, 29(3):1093--1126, 2007.

\bibitem{CerfonFreidbergPOP2010}
A.J. Cerfon and J.P. Freidberg.
\newblock ``{O}ne size fits all" analytic solutions to the {G}rad–{S}hafranov
  equation.
\newblock {\em Physics of Plasmas}, 17(3):032502, 2010.

\bibitem{Klingenberg2015}
P.~Chandrashekar and C.~Klingenberg.
\newblock {A second order well-balanced finite volume scheme for Euler
  Equations with gravity}.
\newblock {\em {Journal on Scientific Computing }}, 37:B382--B402, 2015.

\bibitem{chorin1}
A.J. Chorin.
\newblock A numerical method for solving incompressible viscous flow problems.
\newblock {\em Journal of Computational Physics}, 2:12--26, 1967.

\bibitem{chorin2}
A.J. Chorin.
\newblock Numerical solution of the {Navier--Stokes} equations.
\newblock {\em Mathematics of Computation}, 23:341--354, 1968.

\bibitem{Dedneretal}
A.~Dedner, F.~Kemm, D.~Kr\"oner, C.-D. Munz, T.~Schnitzer, and M.~Wesenberg.
\newblock Hyperbolic divergence cleaning for the {MHD} equations.
\newblock {\em Journal of Computational Physics}, 175:645--673, 2002.

\bibitem{GassnerEntropyGLM}
D.~Derigs, A.~R. Winters, G.~Gassner, S.~Walch, and M.~Bohm.
\newblock {Ideal GLM-MHD: About the entropy consistent nine-wave magnetic field
  divergence diminishing ideal magnetohydrodynamics equations}.
\newblock {\em J. Comput. Phys.}, 364:420--467, 2018.

\bibitem{Dingfelder2020}
B.~Dingfelder and F.~Hindenlang.
\newblock A locally field-aligned discontinuous {Galerkin} method for
  anisotropic wave equations.
\newblock {\em Journal of Computational Physics}, 408:109273, 2020.

\bibitem{DBTM08}
M.~Dumbser, D.~S. Balsara, E.~F. Toro, and C.-D. Munz.
\newblock A unified framework for the construction of one-step finite volume
  and discontinuous {G}alerkin schemes on unstructured meshes.
\newblock {\em Journal of Computational Physics}, 227(18):8209--8253, 2008.

\bibitem{SIMHD}
M.~Dumbser, D.S. Balsara, M.~Tavelli, and F.~Fambri.
\newblock A divergence-free semi-implicit finite volume scheme for ideal,
  viscous, and resistive magnetohydrodynamics.
\newblock {\em International Journal for Numerical Methods in Fluids},
  89(1-2):16--42, 2019.

\bibitem{cweno2017}
M.~Dumbser, W.~Boscheri, M.~Semplice, and G.~Russo.
\newblock Central weighted {ENO} schemes for hyperbolic conservation laws on
  fixed and moving unstructured meshes.
\newblock {\em SIAM Journal on Scientific Computing}, 39(6):A2564--A2591, 2017.

\bibitem{DET08}
M.~Dumbser, C.~Enaux, and E.~F. Toro.
\newblock Finite volume schemes of very high order of accuracy for stiff
  hyperbolic balance laws.
\newblock {\em Journal of Computational Physics}, 227(8):3971 -- 4001, 2008.

\bibitem{DumbserKaeser07}
M.~Dumbser, M.~K\"aser, V.~A. Titarev, and E.~F. Toro.
\newblock Quadrature-free non-oscillatory finite volume schemes on unstructured
  meshes for nonlinear hyperbolic systems.
\newblock {\em Journal of Computational Physics}, 226:204--243, 2007.

\bibitem{GPRmodelMHD}
M.~Dumbser, I.~Peshkov, E.~Romenski, and O.~Zanotti.
\newblock {H}igh order {ADER} schemes for a unified first order hyperbolic
  formulation of {N}ewtonian continuum mechanics coupled with
  electro--dynamics.
\newblock {\em J. Comput. Phys.}, 348:298--342, 2017.

\bibitem{WBGravity}
P.~V.~F. Edelmann, L.~Horst, J.~P. Berberich, R.~Andrassy, J.~Higl, G.~Leidi,
  C.~Klingenberg, and F.~K. Röpke.
\newblock Well-balanced treatment of gravity in astrophysical fluid dynamics
  simulations at low {Mach} numbers.
\newblock {\em Astronomy and Astrophysics}, 652, 2021.

\bibitem{fallemhd2}
S.~A. E.~G. Falle.
\newblock Rarefaction shocks, shock errors and low order of accuracy in {ZEUS}.
\newblock {\em The Astrophysical Journal}, 577:L123--L126, 2002.

\bibitem{fallemhd3}
S.A.E.G. Falle, S.S. Komissarov, and P.~Joarder.
\newblock {A multidimensional upwind scheme for magnetohydrodynamics}.
\newblock {\em Journal of Computational Physics}, 297:265--277, 1998.

\bibitem{Fambri20}
F.~Fambri.
\newblock A novel structure preserving semi-implicit finite volume method for
  viscous and resistive magnetohydrodynamics.
\newblock {\em International Journal for Numerical Methods in Fluids},
  93(12):3447--3489, 2021.

\bibitem{Finan}
C.H. Finan and J.~Killeen.
\newblock {Solution of the time--dependent, three--dimensional resistive
  magnetohydrodynamic equations}.
\newblock {\em Computer Physics Communications}, 24:441--463, 1981.

\bibitem{FreidbergIdealMHD}
Jeffrey~P. Freidberg.
\newblock {\em Ideal {MHD}}.
\newblock Cambridge University Press, 2014.

\bibitem{GCD18}
E.~Gaburro, M.~.J Castro, and M.~Dumbser.
\newblock Well-balanced {Arbitrary-Lagrangian-Eulerian} finite volume schemes
  on moving nonconforming meshes for the {E}uler equations of gas dynamics with
  gravity.
\newblock {\em Mon. Not. R. Astron. Soc.}, 477(2):2251--2275, 2018.

\bibitem{GCD21}
E.~Gaburro, M.J. Castro, and M.~Dumbser.
\newblock A well balanced finite volume scheme for general relativity.
\newblock {\em SIAM Journal on Scientific Computing}, 43(6):B1226--B1251, 2021.

\bibitem{GDC17}
E.~Gaburro, M.~Dumbser, and M.J. Castro.
\newblock Direct {Arbitrary-Lagrangian-Eulerian} finite volume schemes on
  moving nonconforming unstructured meshes.
\newblock {\em Computers and Fluids}, 159:254 – 275, 2017.

\bibitem{ElenaHTC}
E.~Gaburro, P.~\"Offner, M.~Ricchiuto, and D.~Torlo.
\newblock {High order entropy preserving ADER-DG schemes}.
\newblock {\em Applied Mathematics and Computation}, 440:127644, 2023.

\bibitem{mhdsi3}
H.~Gao and W.~Qiu.
\newblock A semi-implicit energy conserving finite element method for the
  dynamical incompressible magnetohydrodynamics equations.
\newblock {\em Computer Methods in Applied Mechanics and Engineering},
  346:982--1001, 2019.

\bibitem{GardinerStone}
T.A. Gardiner and J.M. Stone.
\newblock {An unsplit Godunov method for ideal MHD via constrained transport}.
\newblock {\em Journal of Computational Physics}, 205:509--539, 2005.

\bibitem{Gas07}
G.~Gassner, F.~Lorcher, and C.~D. Munz.
\newblock A contribution to the construction of diffusion fluxes for finite
  volume and discontinuous {G}alerkin schemes.
\newblock {\em Journal of Computational Physics}, 224(2):1049 -- 1063, 2007.

\bibitem{Gawlik2021}
E.S. Gawlik and F.~Gay-Balmaz.
\newblock A structure-preserving finite element method for compressible ideal
  and resistive magnetohydrodynamics.
\newblock {\em Journal of Plasma Physics}, 87:835870501, 2021.

\bibitem{Gawlik2022}
E.S. Gawlik and F.~Gay-Balmaz.
\newblock A finite element method for {MHD} that preserves energy,
  cross-helicity, magnetic helicity, incompressibility, and div {B}.
\newblock {\em Journal of Computational Physics}, 450, 2022.

\bibitem{GGS82}
U.~Ghia, K.N. Ghia, and C.T. Shin.
\newblock High-{R}e solutions for incompressible flow using the
  {N}avier--{S}tokes equations and a multigrid method.
\newblock {\em Journal of Computational Physics}, 48(3):387--411, 1982.

\bibitem{Giorgiani2020}
G.~Giorgiani, H.~Bufferand, F.~Schwander, E.~Serre, and P.~Tamain.
\newblock A high-order non field-aligned approach for the discretization of
  strongly anisotropic diffusion operators in magnetic fusion.
\newblock {\em Computer Physics Communications}, 254:107375, 2020.

\bibitem{God1961}
S.K. Godunov.
\newblock An interesting class of quasilinear systems.
\newblock {\em Dokl. Akad. Nauk SSSR}, 139(3):521--523, 1961.

\bibitem{God1972MHD}
S.K. Godunov.
\newblock Symmetric form of the magnetohydrodynamic equation.
\newblock {\em Numerical Methods for Mechanics of Continuum Medium},
  3(1):26--34, 1972.

\bibitem{gosse2000well}
L~Gosse.
\newblock A well-balanced flux-vector splitting scheme designed for hyperbolic
  systems of conservation laws with source terms.
\newblock {\em Computers \& Mathematics with Applications}, 39(9):135--159,
  2000.

\bibitem{GradRubin}
H.~Grad and H.~Rubin.
\newblock {Hydromagnetic Equilibria and Force-Free Fields}.
\newblock In {\em Proceedings of the 2nd UN Conf. on the Peaceful Uses of
  Atomic Energy, Geneva: IAEA}, volume~31, page 190, 1958.

\bibitem{greenbergleroux}
J.~M. Greenberg and A.~Y.~Le Roux.
\newblock A wellbalanced scheme for the numerical processing of source terms in
  hyperbolic equations.
\newblock {\em SIAM Journal on Numerical Analysis}, 33:1--16, 1996.

\bibitem{HLLCMHD}
K.~F. Gurski.
\newblock An {HLLC}-type approximate {Riemann} solver for ideal
  magnetohydrodynamics.
\newblock {\em SIAM Journal on Scientific Computing}, 25:2165--2187, 2004.

\bibitem{markerandcell}
F.H. Harlow and J.E. Welch.
\newblock Numerical calculation of time-dependent viscous incompressible flow
  of fluid with a free surface.
\newblock {\em Physics of Fluids}, 8:2182--2189, 1965.

\bibitem{Harned}
D.S. Harned and W.~Kerner.
\newblock {Semi--implicit method for three-dimensional resistive
  magnetohydrodynamic simulation of fusion plasmas}.
\newblock {\em Nuclear Science and Engineering}, 92:119--125, 1986.

\bibitem{HelzelRossmanith}
C.~Helzel, J.~A. Rossmanith, and B.~Taetz.
\newblock An unstaggered constrained transport method for the {3D} ideal
  magnetohydrodynamic equations.
\newblock {\em Journal of Computational Physics}, 230(10):3803--3829, 2011.

\bibitem{HiPa18}
R.~Hiptmair and C.~Pagliantini.
\newblock Splitting-based structure preserving discretization of
  magnetohydrodynamics.
\newblock {\em The SMAI Journal of Computational Mathematics}, 4:225--257,
  2018.

\bibitem{HuShuTri}
C.~Hu and {C.W.} Shu.
\newblock Weighted essentially non-oscillatory schemes on triangular meshes.
\newblock {\em Journal of Computational Physics}, 150:97--127, 1999.

\bibitem{Hu2021}
K.~Hu, Y.J. Lee, and J.~Xu.
\newblock {Helicity-conservative finite element discretization for
  incompressible MHD systems}.
\newblock {\em Journal of Computational Physics}, 436:110284, 2021.

\bibitem{HuMaXu2017}
K.~Hu, Y.~Ma, and J.~Xu.
\newblock Stable finite element methods preserving $\nabla\cdot {B} = 0$
  exactly for {MHD} models.
\newblock {\em Numerische Mathematik}, 135:371--396, 2017.

\bibitem{WBMHD}
F.~Kanbar, R.~Touma, and C.~Klingenberg.
\newblock Well-balanced central scheme for the system of {MHD} equations with
  gravitational source term.
\newblock {\em Communications in Computational Physics}, 32(3):878--898, 2022.

\bibitem{Kapelli2014}
R.~K{\"a}ppeli and S.~Mishra.
\newblock {Well-balanced schemes for the Euler equations with gravitation}.
\newblock {\em Journal of Computational Physics}, 259:199--219, 2014.

\bibitem{KlaMaj}
S.~Klainermann and A.~Majda.
\newblock Singular limits of quasilinear hyperbolic systems with large
  parameters and the incompressible limit of compressible fluid.
\newblock {\em Comm. Pure Appl. Math.}, 34:481--524, 1981.

\bibitem{KlaMaj82}
S.~Klainermann and A.~Majda.
\newblock Compressible and incompressible fluids.
\newblock {\em Communications on Pure and Applied Mathematics}, 35:629--651,
  1982.

\bibitem{Klein2001}
R.~Klein, N.~Botta, T.~Schneider, {C.D.} Munz, S.Roller, A.~Meister,
  L.~Hoffmann, and T.~Sonar.
\newblock Asymptotic adaptive methods for multi-scale problems in fluid
  mechanics.
\newblock {\em Journal of Engineering Mathematics}, 39:261--343, 2001.

\bibitem{KlingenbergPuppo}
C.~Klingenberg, G.~Puppo, and M.~Semplice.
\newblock {Arbitrary order finite volume well-balanced schemes for the Euler
  equations with gravity}.
\newblock {\em SIAM Journal on Scientific Computing}, 41:A695--A721, 2019.

\bibitem{mhdsi1}
M.~E. Kress and K.~S. Riedel.
\newblock Semi-implicit reduced magnetohydrodynamics.
\newblock {\em Journal of Computational Physics}, 83(1):237--239, 1989.

\bibitem{Leddy2017}
J.~Leddy, B.~Dudson, M.~Romanelli, B.~Shanahan, and N.~Walkden.
\newblock A novel flexible field-aligned coordinate system for tokamak edge
  plasma simulation.
\newblock {\em Computer Physics Communications}, 212:59--68, 2017.

\bibitem{Lerbinger}
K.~Lerbinger and J.F. Luciani.
\newblock {A new semi-implicit method for MHD computations}.
\newblock {\em Journal of Computational Physics}, 97:444--459, 1991.

\bibitem{leveque}
R.~J. LeVeque.
\newblock Balancing source terms and flux gradients in highresolution {Godunov}
  methods.
\newblock {\em Journal of Computational Physics}, 146:346--365, 1998.

\bibitem{divdg1}
F.~Li and C.W. Shu.
\newblock {Locally divergence-free discontinuous Galerkin methods for MHD
  equations}.
\newblock {\em Journal of Scientific Computing}, 22:413--442, 2005.

\bibitem{divdg2}
F.~Li and L.~Xu.
\newblock {Arbitrary order exactly divergence-free central discontinuous
  Galerkin methods for ideal MHD equations}.
\newblock {\em Journal of Computational Physics}, 231:2655--2675, 2012.

\bibitem{EntropyDGMHD}
Y.~Liu, C.W. Shu, and M.~Zhang.
\newblock {Entropy stable high order discontinuous Galerkin methods for ideal
  compressible MHD on structured meshes}.
\newblock {\em Journal of Computational Physics}, 354:163--178, 2018.

\bibitem{HybridImplicit}
A.~Lucca, S.~Busto, and M.~Dumbser.
\newblock An implicit staggered hybrid finite volume/finite element solver for
  the incompressible {Navier--Stokes} equations.
\newblock {\em East Asian Journal on Applied Mathematics}, 2023.

\bibitem{Ma2016}
Y.~Ma, K.~Hu, X.~Hu, and J.~Xu.
\newblock {Robust preconditioners for incompressible MHD models}.
\newblock {\em Journal of Computational Physics}, 316:721--746, 2016.

\bibitem{HLLD}
T.~Miyoshi and K.~Kusano.
\newblock A multi-state {HLL} approximate {Riemann} solver for ideal
  magnetohydrodynamics.
\newblock {\em Journal of Computational Physics}, 208:315--344, 2005.

\bibitem{MunzCleaning}
C.~D. Munz, P.~Omnes, R.~Schneider, E.~Sonnendr\"ucker, and U.~Voss.
\newblock Divergence correction techniques for {M}axwell solvers based on a
  hyperbolic model.
\newblock {\em Journal of Computational Physics}, 161:484--511, 07 2000.

\bibitem{munzMPV}
{C.D.} Munz, R.~Klein, S.~Roller, and {K.J.} Geratz.
\newblock The extension of incompressible flow solvers to the weakly
  compressible regime.
\newblock {\em Computers \& Fluids}, 32:173--196, 2003.

\bibitem{Noelle1}
S.~Noelle, N.~Pankratz, G.~Puppo, and J.R. Natvig.
\newblock Well-balanced finite volume schemes of arbitrary order of accuracy
  for shallow water flows.
\newblock {\em J. Comput. Phys.}, 213:474--499, 2006.

\bibitem{Noelle2}
S.~Noelle, Y.L. Xing, and C.W. Shu.
\newblock High-order well-balanced finite volume {WENO} schemes for shallow
  water equation with moving water.
\newblock {\em J. Comput. Phys.}, 226:29--58, 2007.

\bibitem{Pares2006}
C.~Par{\'e}s.
\newblock Numerical methods for nonconservative hyperbolic systems: a
  theoretical framework.
\newblock {\em SIAM Journal on Numerical Analysis}, 44(1):300--321, 2006.

\bibitem{pares2004well}
C.~Par{\'e}s and M.~Castro.
\newblock On the well-balance property of {Roe}'s method for nonconservative
  hyperbolic systems. applications to shallow-water systems.
\newblock {\em ESAIM Math. Model. Numer. Anal.}, 38(5):821--852, 2004.

\bibitem{patankar}
V.S. Patankar.
\newblock {\em Numerical {Heat} {Transfer} and {Fluid} {Flow}}.
\newblock Hemisphere Publishing Corporation, 1980.

\bibitem{patankarspalding}
V.S. Patankar and B.~Spalding.
\newblock A calculation procedure for heat, mass and momentum transfer in
  three-dimensional parabolic flows.
\newblock {\em International Journal of Heat and Mass Transfer}, 15:1787--1806,
  1972.

\bibitem{PowellMHD1}
K.G. Powell.
\newblock {An approximate Riemann solver for magnetohydrodynamics (that works
  in more than one dimension)}.
\newblock Technical Report ICASE-Report 94-24 (NASA CR-194902), NASA Langley
  Research Center, Hampton, VA, 1994.

\bibitem{PowellMHD2}
K.G. Powell, P.L. Roe, T.J. Linde, T.I. Gombosi, and D.L.~De Zeeuw.
\newblock A solution-adaptive upwind scheme for ideal magnetohydrodynamics.
\newblock {\em J. Comput. Phys.}, 154:284--309, 1999.

\bibitem{mhdsi4}
W.~Qiu and K.~Shi.
\newblock Analysis of a semi-implicit structure-preserving finite element
  method for the nonstationary incompressible magnetohydrodynamics equations.
\newblock {\em Computers and Mathematics with Applications}, 80(10):2150--2161,
  2020.

\bibitem{HybridMPI}
L.~R\'io-Mart\'in, S.~Busto, and M.~Dumbser.
\newblock A massively parallel hybrid finite volume/finite element scheme for
  computational fluid dynamics.
\newblock {\em Mathematics}, 9:2316, 2021.

\bibitem{Rom1998}
E.I. Romenski.
\newblock Hyperbolic systems of thermodynamically compatible conservation laws
  in continuum mechanics.
\newblock {\em Mathematical and computer modelling}, 28(10):115--130, 1998.

\bibitem{RyuJones}
D.~Ryu and T.~W. Jones.
\newblock Numerical magnetohydrodynamics in astrophysics: algorithm and tests
  for one-dimensional flow.
\newblock {\em Astrophysical Journal}, 442:228--258, 1995.

\bibitem{BLTheory}
H.~Schlichting and K.~Gersten.
\newblock {\em Boundary layer theory}.
\newblock Springer, 2016.

\bibitem{Shafranov}
V.D. Shafranov.
\newblock {Plasma equilibrium in a magnetic field}.
\newblock {\em Reviews of Plasma Physics}, 2:103, 1966.

\bibitem{EULAG}
P.K. Smolarkiewicz and P.~Charbonneau.
\newblock {EULAG, a computational model for multiscale flows: An MHD
  extension}.
\newblock {\em Journal of Computational Physics}, 236:608--623, 2013.

\bibitem{solov}
L.~S. Soloviev.
\newblock {T}he theory of hydromagnetic stability of toroidal plasma
  configurations.
\newblock {\em Sov. Phys. JETP}, 26:400--407, 1968.

\bibitem{Tamain2016}
P.~Tamain, H.~Bufferand, G.~Ciraolo, C.~Colin, D.~Galassi, Ph. Ghendrih,
  F.~Schwander, and E.~Serre.
\newblock {The TOKAM3X code for edge turbulence fluid simulations of tokamak
  plasmas in versatile magnetic geometries}.
\newblock {\em Journal of Computational Physics}, 321:606--623, 2016.

\bibitem{TD14}
M.~Tavelli and M.~Dumbser.
\newblock A staggered semi-implicit discontinuous {Galerkin} method for the two
  dimensional incompressible {N}avier--{S}tokes equations.
\newblock {\em Appl. Math. Comput.}, 248:70 -- 92, 2014.

\bibitem{TD16}
M.~Tavelli and M.~Dumbser.
\newblock A staggered space-time discontinuous {Galerkin} method for the
  three-dimensional incompressible {Navier--Stokes} equations on unstructured
  tetrahedral meshes.
\newblock {\em Journal of Computational Physics}, 319:294 -- 323, 2016.

\bibitem{TD17}
M.~Tavelli and M.~Dumbser.
\newblock A pressure-based semi-implicit space-time discontinuous {G}alerkin
  method on staggered unstructured meshes for the solution of the compressible
  {N}avier--{S}tokes equations at all {Mach} numbers.
\newblock {\em Journal of Computational Physics}, 341:341 -- 376, 2017.

\bibitem{ThomannWB1}
A.~Thomann, G.~Puppo, and C.~Klingenberg.
\newblock {An all speed second order well-balanced IMEX relaxation scheme for
  the Euler equations with gravity}.
\newblock {\em Journal of Computational Physics}, 420, 2020.

\bibitem{ThomannWB2}
A.~Thomann, M.~Zenk, and C.~Klingenberg.
\newblock {A second-order positivity-preserving well-balanced finite volume
  scheme for Euler equations with gravity for arbitrary hydrostatic
  equilibria}.
\newblock {\em International Journal for Numerical Methods in Fluids},
  89(11):465--482, 2019.

\bibitem{titarevtoro}
V.~A. Titarev and E.~F. Toro.
\newblock {ADER} schemes for three-dimensional nonlinear hyperbolic systems.
\newblock {\em Journal of Computational Physics}, 204:715--736, 2005.

\bibitem{MixedWENO2D}
V.A. Titarev, P.~Tsoutsanis, and D.~Drikakis.
\newblock {WENO schemes for mixed--element unstructured meshes}.
\newblock {\em Communications in Computational Physics}, 8:585--609, 2010.

\bibitem{Toro}
E.~F. Toro.
\newblock {\em Riemann solvers and numerical methods for fluid dynamics: A
  practical introduction}.
\newblock Springer, 2009.

\bibitem{TMN01}
E.~F. Toro, R.~C. Millington, and L.~A.~M. Nejad.
\newblock {\em Godunov methods}, chapter Towards very high order {G}odunov
  schemes.
\newblock Springer, 2001.

\bibitem{USFORCE}
E.F. Toro, A.~Hidalgo, and M.~Dumbser.
\newblock {FORCE} schemes on unstructured meshes {I}: Conservative hyperbolic
  systems.
\newblock {\em J. Comput. Phys.}, 228:3368--3389, 2009.

\bibitem{torrilhon}
M.~Torrilhon.
\newblock Non-uniform convergence of finite volume schemes for {Riemann}
  problems of ideal magnetohydrodynamics.
\newblock {\em Journal of Computational Physics}, 192:73--94, 2003.

\bibitem{torrilhonbalsara}
M.~Torrilhon and D.S. Balsara.
\newblock {High order WENO schemes: investigations on non-uniform convergence
  for MHD Riemann problems}.
\newblock {\em Journal of Computational Physics}, 201:586--600, 2004.

\bibitem{MixedWENO3D}
P.~Tsoutsanis, V.A. Titarev, and D.~Drikakis.
\newblock {WENO schemes on arbitrary mixed-element unstructured meshes in three
  space dimensions}.
\newblock {\em J. Comput. Phys.}, 230:1585--1601, 2011.

\bibitem{vanKan}
J.~van Kan.
\newblock {A second-order accurate pressure correction method for viscous
  incompressible flow}.
\newblock {\em SIAM Journal on Scientific and Statistical Computing},
  7:870--891, 1986.

\bibitem{VL97}
B.~van Leer.
\newblock Towards the ultimate conservative difference scheme.
\newblock {\em Journal of Computational Physics}, 135(2):229--248, 1997.

\bibitem{WarburtonVRMHD}
T.~Warburton and G.~Karniadakis.
\newblock A discontinuous {{G}alerkin} method for the viscous {MHD} equations.
\newblock {\em Journal of Computational Physics}, 152:608--641, 1999.

\bibitem{xing2013high}
Y.~Xing and C.W. Shu.
\newblock High order well-balanced {WENO} scheme for the gas dynamics equations
  under gravitational fields.
\newblock {\em J. Sci. Comput.}, 54(2-3):645--662, 2013.

\bibitem{divdg3}
Z.~Xu and Y.~Liu.
\newblock {New central and central discontinuous Galerkin schemes on
  overlapping cells of unstructured grids for solving ideal magnetohydrodynamic
  equations with globally divergence-free magnetic field}.
\newblock {\em Journal of Computational Physics}, 327:203--224, 2016.

\bibitem{Zanotti2015c}
O.~Zanotti, F.~Fambri, M.~Dumbser, and A.~Hidalgo.
\newblock Space-time adaptive {ADER} discontinuous {{G}alerkin} finite element
  schemes with a posteriori sub-cell finite volume limiting.
\newblock {\em Computers and Fluids}, 118:204 -- 224, 2015.

\bibitem{mhdsi2}
X.~Zhao, Y.~Yang, and C.~E. Seyler.
\newblock A positivity-preserving semi-implicit discontinuous {Galerkin} scheme
  for solving extended magnetohydrodynamics equations.
\newblock {\em Journal of Computational Physics}, 278:400--415, 2014.

\end{thebibliography}

\appendix
\section{Supplementary material}\label{sec:appendix}
This appendix contains a brief description of the supplementary material 
provided to the reader to ease comparison of the obtained results for the MHD lid driven cavity problem. The supplementary files can be downloaded from the online version of the article and are licensed under a Creative Commons ``Attribution-NonCommercial-ShareAlike 3.0 Unported'' license \ccbyncsa.  
\\The data files contain the solution of the lid-driven cavity with magnetic field, presented in Section~\ref{sec:MHD_LDC}. We recall that this test is a variant of the traditional lid-driven cavity in which an initial horizontal or a vertical magnetic field is imposed , i.e., $\bbvar_0 = (\bvar_{0,x},0)^{T}$ or  $\bbvar_0 = (0,\bvar_{0,y})^{T}$. We consider the lid velocity at the upper boundary $\bvel=(1,0)$, a viscosity $\mu = 0.01$ and a resistivity $\eta = 0.01$. A more detailed description of the contents of each data file, as well as the scripts made available to plot these data, in Tecplot$^{\textrm{\textregistered}}$, gnuplot, Matlab$^{\textrm{\textregistered}}$ and Python\texttrademark, is given below.

\renewcommand{\arraystretch}{1.2}
\noindent\begin{tabular}{ll}
	\toprule
	\multicolumn{2}{c}{\parbox{16.cm}{\centering\textbf{Available data files}}}\\
	\bottomrule
	\toprule
	\multicolumn{2}{c}{\parbox{16.cm}{\textbf{Description}: Data files containing the 1D cut of the magnetic field $\bvar_{2}$ (and the velocity $\vel_{2}$) along the $x$-axis for the MHD lid-driven cavity test presented in Section~\ref{sec:MHD_LDC}. These results have been obtained using the \textbf{hybrid FV/FE approach} introduced in this paper imposing a horizontal magnetic field. Different values of $\bvar_{x,0}$ are considered.}}\\
	\midrule
	{\textbf{File names:}}  & Initial horizontal field $\bvar_{x,0}$ \\
	HorizontalB\_HybridFVFE\_xBy\_B0\_010.dat &  $\bvar_{x,0}=0.1$ \\ 
	HorizontalB\_HybridFVFE\_xBy\_B0\_025.dat &  $\bvar_{x,0}=0.25$ \\
	HorizontalB\_HybridFVFE\_xBy\_B0\_050.dat &  $\bvar_{x,0}=0.5$ \\
	HorizontalB\_HybridFVFE\_xBy\_B0\_100.dat &  $\bvar_{x,0}=1.0$ \\[1pt]
	\bottomrule
	\toprule
	\multicolumn{2}{c}{\parbox{16.cm}{\textbf{Description}: Data files containing the 1D cut of the magnetic field $\bvar_{1}$ (and the velocity $\vel_{1}$) along the $y$-axis for the MHD lid-driven cavity test presented in Section~\ref{sec:MHD_LDC}. These results have been computed with the proposed \textbf{hybrid FV/FE approach}  for an initially imposed horizontal magnetic field. Different values of $\bvar_{x,0}$ are considered.}}\\
	\midrule
	{\textbf{File names:}}  & Initial horizontal field $\bvar_{x,0}$ \\
	HorizontalB\_HybridFVFE\_yBx\_B0\_010.dat &  $\bvar_{x,0}=0.1$ \\ 
	HorizontalB\_HybridFVFE\_yBx\_B0\_025.dat &  $\bvar_{x,0}=0.25$ \\
	HorizontalB\_HybridFVFE\_yBx\_B0\_050.dat &  $\bvar_{x,0}=0.5$ \\
	HorizontalB\_HybridFVFE\_yBx\_B0\_100.dat &  $\bvar_{x,0}=1.0$ \\[1pt]
	\bottomrule
	\toprule
	\multicolumn{2}{c}{\parbox{16.cm}{\textbf{Description}: Data files containing the 1D cut of the magnetic field $\bvar_{2}$ (and the velocity $\vel_{2}$) along the $x$-axis for the MHD lid-driven cavity test presented in Section~\ref{sec:MHD_LDC}. These results have been obtained using the \textbf{semi-implicit divergence-free method} introduced in~\cite{SIMHD} when a horizontal magnetic field is imposed. Different values of $\bvar_{x,0}$ are considered.}}\\
	\midrule
	{\textbf{File names:}}  & Initial horizontal field $\bvar_{x,0}$ \\
	HorizontalB\_Reference\_xBy\_B0\_010.dat &  $\bvar_{x,0}=0.1$ \\ 
	HorizontalB\_Reference\_xBy\_B0\_025.dat &  $\bvar_{x,0}=0.25$ \\
	HorizontalB\_Reference\_xBy\_B0\_050.dat &  $\bvar_{x,0}=0.5$ \\
	HorizontalB\_Reference\_xBy\_B0\_100.dat &  $\bvar_{x,0}=1.0$ \\[1pt]
	\bottomrule
	\toprule
	\multicolumn{2}{c}{\parbox{16.0cm}{\textbf{Description}: Data files containing the 1D cut of the magnetic field $\bvar_{1}$ (and the velocity $\vel_{1}$) along the $y$-axis for the MHD lid-driven cavity test presented in Section~\ref{sec:MHD_LDC}. These results have been obtained using the \textbf{semi-implicit divergence-free method}, introduced in~\cite{SIMHD}, for an initial horizontal magnetic field. Different values of $\bvar_{x,0}$ are considered.}}\\
	\midrule
	{\textbf{File names:}}  & Initial horizontal field $\bvar_{x,0}$ \\
	HorizontalB\_Reference\_yBx\_B0\_010.dat &  $\bvar_{x,0}=0.1$ \\ 
	HorizontalB\_Reference\_yBx\_B0\_025.dat &  $\bvar_{x,0}=0.25$ \\
	HorizontalB\_Reference\_yBx\_B0\_050.dat &  $\bvar_{x,0}=0.5$ \\
	HorizontalB\_Reference\_yBx\_B0\_100.dat &  $\bvar_{x,0}=1.0$ \\[1pt]
	\bottomrule
\end{tabular}

\noindent\begin{tabular}{ll}
	\toprule
	\multicolumn{2}{c}{\parbox{16.cm}{\centering\textbf{Available data files}}}\\
	\bottomrule
	\toprule
	\multicolumn{2}{c}{\parbox{16.cm}{\textbf{Description}: Data files containing the 1D cut of the magnetic field $\bvar_{2}$ (and the velocity $\vel_{2}$) along the $x$-axis for the MHD lid-driven cavity test presented in Section~\ref{sec:MHD_LDC}. These results have been obtained using the \textbf{hybrid FV/FE approach} introduced in this paper when a vertical magnetic field is imposed. Different values of $\bvar_{0,y}$ are considered.}}\\
	\midrule
	{\textbf{File names}}  & Initial vertical field $\bvar_{0,y}$ \\
	VerticalB\_HybridFVFE\_xBy\_B0\_010.dat &  $\bvar_{0,y}=0.1$ \\ 
	VerticalB\_HybridFVFE\_xBy\_B0\_025.dat &  $\bvar_{0,y}=0.25$ \\
	VerticalB\_HybridFVFE\_xBy\_B0\_050.dat &  $\bvar_{0,y}=0.5$ \\
	VerticalB\_HybridFVFE\_xBy\_B0\_100.dat &  $\bvar_{0,y}=1.0$ \\[1pt]
	\bottomrule
	\toprule
	\multicolumn{2}{c}{\parbox{16.cm}{\textbf{Description}: Data files containing the 1D cut of the magnetic field $\bvar_{1}$ (and the velocity $\vel_{1}$) along the $y$-axis for the MHD lid-driven cavity test presented in Section~\ref{sec:MHD_LDC}. These results have been computed with the \textbf{hybrid FV/FE approach} introduced in this paper if a vertical magnetic field is initially considered. Different values of  $\bvar_{0,y}$ are studied.}}\\
	\midrule
	{\textbf{File names:}}  & Initial vertical field $\bvar_{0,y}$ \\
	VerticalB\_HybridFVFE\_yBx\_B0\_010.dat &  $\bvar_{0,y}=0.1$ \\ 
	VerticalB\_HybridFVFE\_yBx\_B0\_025.dat &  $\bvar_{0,y}=0.25$ \\
	VerticalB\_HybridFVFE\_yBx\_B0\_050.dat &  $\bvar_{0,y}=0.5$ \\
	VerticalB\_HybridFVFE\_yBx\_B0\_100.dat &  $\bvar_{0,y}=1.0$ \\[1pt]
	\bottomrule
	\toprule
	\multicolumn{2}{c}{\parbox{16.cm}{\textbf{Description}: Data files containing the 1D cut of the magnetic field $\bvar_{2}$ (and the velocity $\vel_{2}$) along the $x$-axis for the MHD lid-driven cavity test presented in Section~\ref{sec:MHD_LDC}. These results have been obtained using the \textbf{semi-implicit divergence-free method} introduced in~\cite{SIMHD} when a vertical magnetic field is imposed for different values of $\bvar_{0,y}$.}}\\
	\midrule
	{\textbf{File names:}}  & Initial vertical field $\bvar_{0,y}$ \\
	VerticalB\_Reference\_xBy\_B0\_010.dat &  $\bvar_{0,y}=0.1$ \\ 
	VerticalB\_Reference\_xBy\_B0\_025.dat &  $\bvar_{0,y}=0.25$ \\
	VerticalB\_Reference\_xBy\_B0\_050.dat &  $\bvar_{0,y}=0.5$ \\
	VerticalB\_Reference\_xBy\_B0\_100.dat &  $\bvar_{0,y}=1.0$ \\[1pt]
	\bottomrule
	\toprule
	\multicolumn{2}{c}{\parbox{16.0cm}{\textbf{Description}: Data files containing the 1D cut of the magnetic field $\bvar_{1}$ (and the velocity $\vel_{1}$) along the $y$-axis for the MHD lid-driven cavity test presented in Section~\ref{sec:MHD_LDC}. These results have been obtained using the \textbf{semi-implicit divergence-free method} introduced in~\cite{SIMHD} when a vertical magnetic fields is imposed. Different values of $\bvar_{0,y}$ are considered.}}\\
	\midrule
	{\textbf{File names:}}  & Initial vertical field $\bvar_{0,y}$ \\
	VerticalB\_Reference\_yBx\_B0\_010.dat &  $\bvar_{0,y}=0.1$ \\ 
	VerticalB\_Reference\_yBx\_B0\_025.dat &  $\bvar_{0,y}=0.25$ \\
	VerticalB\_Reference\_yBx\_B0\_050.dat &  $\bvar_{0,y}=0.5$ \\
	VerticalB\_Reference\_yBx\_B0\_100.dat &  $\bvar_{0,y}=1.0$ \\[1pt]
	\bottomrule
\end{tabular}

\noindent\begin{tabular}{l}
	\toprule
	\parbox{16.cm}{\centering\textbf{Scripts for data visualization}}\\
	\bottomrule
	\toprule
	\textbf{File name:} Hybrid\_vs\_RefSol\_horizontalMagField.lay\\
	\parbox{16.cm}{\textbf{Description}: \textbf{Tecplot$^{\textrm{\textregistered}}$ layout} to plot the 1D cuts of the magnetic fields $\bvar_{2}$ and $\bvar_{1}$ (and the velocity fields $\vel_{2}$ and $\vel_{1}$) along the $x-$ and the $y-$axis, respectively, for the MHD lid-driven cavity test, presented in Section~\ref{sec:MHD_LDC}. The plots compare the solution computed using the hybrid FV/FE approach presented in this document and the semi-implicit divergence-free method introduced in~\cite{SIMHD} when a horizontal magnetic field is imposed, considering different values of $\bvar_{x,0}$.}\\
	\textbf{Requirements}:  Tecplot$^{\textrm{\textregistered}}$ 360.\\
	\bottomrule
	\toprule
	\textbf{File name:} Hybrid\_vs\_RefSol\_verticalMagField.lay\\
	\parbox{16.cm}{\textbf{Description}: \textbf{Tecplot$^{\textrm{\textregistered}}$ layout} to plot the 1D cuts of the magnetic fields $\bvar_{2}$ and $\bvar_{1}$ (and the velocity fields $\vel_{2}$ and $\vel_{1}$) along the $x-$ and the $y-$axis, respectively, for the MHD lid-driven cavity test, described in Section~\ref{sec:MHD_LDC}. The plots display both the solution computed using the hybrid FV/FE approach presented in this work and the semi-implicit divergence-free method introduced in~\cite{SIMHD} when initially a vertical magnetic field is considered, for different values of $\bvar_{0,y}$}\\
	\textbf{Requirements:} Tecplot$^{\textrm{\textregistered}}$ 360.\\
	\bottomrule
	\toprule
	\textbf{File name:} plot\_LDC\_MHD\_gnuplot.txt \\
	\parbox{16.cm}{\textbf{Description}: Script to plot, with the program \textbf{gnuplot}, the 1D cuts of the magnetic fields $\bvar_{2}$ and $\bvar_{1}$ (and the velocity fields $\vel_{2}$ and $\vel_{1}$) for the MHD lid-driven cavity test, described in Section~\ref{sec:MHD_LDC}, if horizontal (left plot) and vertical (right plot) magnetic fields are imposed initially. Running the script, we can see the comparison between the solution calculated using the hybrid FV/FE approach and the one obtained with a semi-implicit divergence-free method.}\\
	\textbf{Requirements:} Gnuplot.\\
	\textbf{Running from a command window:} \texttt{gnuplot plot\_LDC\_MHD\_gnuplot.txt}\\
	\bottomrule
	\toprule
	\textbf{File name:} plot\_LDC\_MHD\_MATLAB.m\\
	\parbox{16.cm}{\textbf{Description}: \textbf{MATLAB$^{\textrm{\textregistered}}$ script} to plot the 1D cuts of the magnetic fields $\bvar_{2}$ and $\bvar_{1}$ (and the velocity fields $\vel_{2}$ and $\vel_{1}$) for the MHD lid-driven cavity test, described in Section~\ref{sec:MHD_LDC}, if horizontal (left plot) and vertical (right plot) magnetic fields are imposed initially. The plots display the solution calculated using the hybrid FV/FE approach and that obtained with a semi-implicit divergence-free method.}\\
	\textbf{Requirements:} MATLAB$^{\textrm{\textregistered}}$ 2016 or higher.\\
	\textbf{Running from a command window:} \texttt{matlab plot\_LDC\_MHD\_MATLAB.m}\\
	\bottomrule
    \toprule
    \textbf{File name:} plot\_LDC\_MHD\_python.py\\
	\parbox{16.cm}{\textbf{Description}: \textbf{Python\texttrademark\, script} to plot the 1D cuts of the magnetic fields $\bvar_{2}$ and $\bvar_{1}$ (and the velocity fields $\vel_{2}$ and $\vel_{1}$) for the MHD lid-driven cavity test, described in Section~\ref{sec:MHD_LDC}. while the left plot shows the results when a horizontal magnetic field is initially imposed, the right one displays the results if a vertical magnetic field is initially considered. In both cases, figures show the comparison between the solution computed following the hybrid FV/FE approach and the one calculated using a semi-implicit divergence-free method.}\\
	\textbf{Requirements:} Python\texttrademark\, 3 and the libraries \texttt{numpy} and \texttt{matplotlib}.\\
	\textbf{Running from a command window:} \texttt{python plot\_LDC\_MHD\_python.py}\\
	\bottomrule
\end{tabular}

\vspace{0.5cm}



\end{document}